\documentclass{article}

\usepackage{arxiv}

\usepackage[utf8]{inputenc} 
\usepackage[T1]{fontenc}    
\usepackage{url}            
\usepackage{booktabs}       
\usepackage{amsfonts}       
\usepackage{nicefrac}       
\usepackage{microtype}      
\usepackage{lipsum}

\usepackage{amssymb}
\usepackage{amsmath}
\usepackage{lineno}
\usepackage{bbm}
\usepackage{amsthm}
\usepackage{longtable}
\usepackage{accents}
\usepackage{multirow}
\usepackage{subcaption}
\usepackage{rotating}
\usepackage{lscape}
\usepackage[nottoc,notlot,notlof]{tocbibind}

\usepackage{graphicx}

\usepackage{array}

\usepackage[shortlabels]{enumitem}
\allowdisplaybreaks

\newtheorem{dfn}{Definition}
\newtheorem{thm}{Theorem}
\newtheorem*{thm*}{Theorem}
\newtheorem{corr}{Corrollary}
\newtheorem*{corr*}{Corrollary}
\newtheorem{prop}{Proposition}
\newtheorem{lemma}{Lemma}
\newtheorem*{lemma*}{Lemma}
\newtheorem{rmk}{Remark}
\newtheorem*{rmk*}{Remark}

\newtheorem{innercustomgeneric}{\customgenericname}
\providecommand{\customgenericname}{}
\newcommand{\newcustomtheorem}[2]{%
  \newenvironment{#1}[1]
  {%
   \renewcommand\customgenericname{#2}%
   \renewcommand\theinnercustomgeneric{##1}%
   \innercustomgeneric
  }
  {\endinnercustomgeneric}
}

\newcustomtheorem{customthm}{Theorem}
\newcustomtheorem{customlemma}{Lemma}
\newcustomtheorem{customprop}{Proposition}
\newcustomtheorem{customcorollary}{Corollary}
\newcustomtheorem{customremark}{Remark}

\makeatletter
\newcommand*{\rom}[1]{\expandafter\@slowromancap\romannumeral #1@}
\makeatother

\newcolumntype{x}[1]{%
>{\centering\hspace{0pt}}p{#1}}%

\title{Testing for complete spatial randomness on three dimensional bounded convex shapes}

\author{
  Scott Ward \\
  Department of Mathematics\\
  Imperial College London\\
  London, SW7 2AZ\\
  \texttt{scott.ward12@imperial.ac.uk} \\
  \And
 Edward A.K. Cohen \\
  Department of Mathematics\\
  Imperial College London\\
  London, SW7 2AZ\\
  \texttt{e.cohen@imperial.ac.uk} \\
  \And
 Niall Adams \\
 Department of Mathematics\\
  Data Science Institute\\
  Imperial College London\\
  London, SW7 2AZ\\
  \texttt{n.adams@imperial.ac.uk} \\
}

\begin{document}
\maketitle

\begin{abstract}
There is currently a gap in theory for point patterns that lie on the surface of objects, with researchers focusing on patterns that lie in a Euclidean space, typically planar and spatial data. Methodology for planar and spatial data thus relies on Euclidean geometry and is therefore inappropriate for analysis of point patterns observed in non-Euclidean spaces.  Recently, there has been extensions to the analysis of point patterns on a sphere, however, many other shapes are left unexplored. This is in part due to the challenge of defining the notion of \emph{stationarity} for a point process existing on such a space due to the lack of rotational and translational isometries. Here, we construct functional summary statistics for Poisson processes defined on convex shapes in three dimensions. Using the Mapping Theorem, a Poisson process can be transformed from any convex shape to a Poisson process on the unit sphere which has rotational symmetries that allow for functional summary statistics to be constructed. We present the first and second order properties of such summary statistics and demonstrate how they can be used to test whether an observed pattern exhibits complete spatial randomness or spatial preference on the original convex space. A study of the Type \rom{1} and \rom{2} errors of our test statistics are explored through simulations on ellipsoids of varying dimensions. 
\end{abstract}


\section{Introduction}\label{Section:1}
Research in spatial statistics has predominantly concentrated on the development of theory and methodology for point processes on $\mathbb{R}^d$, with a significant focus on planar ($\mathbb{R}^2$) and spatial ($\mathbb{R}^3$) data. Point processes existing on non-Euclidean spaces, however, are still relatively under-explored. Recently, with the advent of spatial data on a global scale, and modelling Earth as a sphere, there have been important developments in the theory and analysis of point processes on the surface of $d-1$ dimensional unit spheres, $\mathbb{S}^{d-1}\subset \mathbb{R}^d$ \cite{Lawrence2016,Moller2016,Robeson2014}. Yet patterns can still arise for which these methodologies are inappropriate as they lie on other bounded metric spaces that deviate significantly from $\mathbb{S}^{d-1}$. For example, microbiologists are concerned with the spatial arrangement of lipids and proteins on the cellular membranes of microorganisms that are not adequately modelled by spheres. In the case of bacteria, ellipsoids or capsules are far more appropriate candidate surfaces. Recent advances in 3D super-resolution imaging techniques \cite[e.g.]{cabriel2019,gustavsson2018} output point patterns of this type, and there is a demand for the correct statistical procedures to analyse them.

Key to the statistical analysis of spatial data is the ability to form functional summary statistics from an observed pattern, primarily for performing exploratory data analysis and testing for complete spatial randomness (CSR). On $\mathbb{R}^d$ and $\mathbb{S}^{d-1}$, there exists an infinite number of isometries, allowing for the notions of stationarity and isotropy to be well defined, which in turn allows for well defined functional summary statistics. However, on the surface of an arbitrary convex shape $\mathbb{D}\subset\mathbb{R}^3$, the set of available isometries is finite, and thus defining stationarity, isotropy, and summary statistics directly on $\mathbb{D}$ is non-trivial. Building on the current literature for spherical point patterns, in particular the discussion of inhomogeneous point processes on a sphere by both \cite{Lawrence2016} and \cite{Moller2016}, we show that it is possible to construct functional summary statistics for point processes on the surface of arbitrary convex shapes in $\mathbb{R}^3$, with our primary interest being to test for CSR.

Our approach is to map the point pattern from an arbitrary convex shape $\mathbb{D}$ onto $\mathbb{S}^2$. For any Poisson process on $\mathbb{D}$, the Mapping Theorem \cite{Kingman1993} determines that the mapped process on $\mathbb{S}^2$ remains Poisson with the intensity function dependent on the mapping. By working on $\mathbb{S}^2$, we operate on a space that is more amenable to constructing summary statistics. These then allow us to test for CSR on $\mathbb{D}$. The functional summary statistics we develop are based on the inhomogeneous counterparts of typical functional summary statistics already established in the spatial statistics literature. In particular, we focus on the inhomogeneous $K$-function, first discussed by \cite{Baddeley2000} for $\mathbb{R}^d$ and later extended to $\mathbb{S}^{d-1}$ by \cite{Lawrence2016,Moller2016}. Furthermore, we also construct the empty-space function, $F$, spherical contact distribution, $H$, and $J$-function (the ratio of the $H$- and $F$-functions) for point processes on arbitrary convex shapes by extending the inhomogeneous definitions of \cite{vanLieshout2011} from $\mathbb{R}^d$ to $\mathbb{S}^2$. 

Section \ref{Section:2} introduces the notation used throughout this work and formally states the hypothesis for testing CSR on an arbitrary convex shape. Section \ref{Section:4} discusses functional summary statistics on $\mathbb{S}^2$, key for the construction of summary statistics on more general bounded subspaces of $\mathbb{R}^3$, and the impracticalities of attempting to define functional summary statistics directly on $\mathbb{D}$. Section \ref{Section:12} extends the inhomogeneous $F$-, $H$-, and $J$-functions from $\mathbb{R}^d$ \cite{vanLieshout2011} to $\mathbb{S}^2$. Section \ref{Section:5} describes the construction of functional summary statistics on bounded subspaces of $\mathbb{R}^3$ for Poisson processes, discussing their first and second order properties in the event that the intensity function is known. Section \ref{Section:11} provides two worked examples constructing functional summary statistics for realisations of a Poisson process observed on a cube and an ellipsoid. Section \ref{Section:8} discusses how regular and cluster processes can be detected based on the deviations of the empirical functional summary statistics. Section \ref{Section:6} describes estimation procedures for the functional summary statistics when the intensity function is unknown and we propose a test statistic for CSR. Finally in Section \ref{Section:9} we conduct empirical power tests using Monte Carlo simulations to explore the properties of our proposed test statistic.

\section{Preliminaries}\label{Section:2}
In this section we outline the necessary spatial theory and notation used throughout this work. We start by introducing the notion of a bounded convex space in $\mathbb{R}^3$ and then define what it means for a point process to lie on such a surface. We end with the statement of the problem that this work is primarily focused on.

\subsection{Notation}
Let $\mathbf{x}\in\mathbb{R}^3$ such that $\mathbf{x}=(x_1,x_2,x_3)^T$ and define $||\mathbf{x}||=(x_1^2+x_2^2+x_3^2)^{1/2}$ to be the Euclidean norm with the origin of $\mathbb{R}^3$ denoted as $\boldsymbol{0}=(0,0,0)^T$. Denote a subset of $\mathbb{R}^3$ as $\mathbb{D}=\{\mathbf{x}\in\mathbb{R}^3:g(\mathbf{x})=0\}$, where $g:\mathbb{R}^3\mapsto\mathbb{R}$. We also suppose that $\mathbb{D}$ is compact (i.e. closed and bounded) and call $g$ the level-set function of $\mathbb{D}$. Define the set $\mathbb{D}_{int}=\{\mathbf{x}\in\mathbb{R}^3:g(\mathbf{x})<0\}$, i.e. the boundary of $\mathbb{D}_{int}$ is $\mathbb{D}$ and we refer to $\mathbb{D}_{int}$ as the interior of $\mathbb{D}$. The set $\mathbb{D}_{int}$ is said to be convex if and only if for all $\mathbf{x},\mathbf{y}\in \mathbb{D}_{int}$ such that $\mathbf{x}\neq\mathbf{y}$ then $\{\mathbf{z}\in\mathbb{R}^3:\mathbf{z}=\mathbf{x}+\gamma(\mathbf{y}-\mathbf{x}),\; \gamma\in (0,1)\}\in \mathbb{D}_{int}$. We thus define $\mathbb{D}$ to be convex if its interior, $\mathbb{D}_{int}$, is also convex. Examples of bounded convex sets of $\mathbb{R}^3$ are spheres, ellipsoids, and cubes. Further for any bounded convex set $\mathbb{D}$ with level-set function $g$, we will also define $\tilde{g}$ which rearranges $g(\mathbf{x})=0$, such that $x_3=\tilde{g}(x_1,x_2)$, i.e. we write $x_3$ as a function of $x_1$ and $x_2$. It may not always be possible to find $\tilde{g}$ explicitly since, as defined previously, it may be the case that the resultant $\tilde{g}$ is not a proper function. This issue can be rectified by partitioning $\mathbb{D}$ appropriately. For example take the case of a sphere with radius 1, then $g(\mathbf{x})=x_1^2+x_2^2+x_3^2-1$, hence $\tilde{g}(x_1,x_2)=\pm (1-x_1^2-x_2^2)^{1/2}$, which is not a proper function. In this case we partition the region $\mathbb{D}$ into the regions $x_3\geq 0$ and $x_3<0$. Then for $x_3\geq 0$, we define $\tilde{g}(x_1,x_2)=+(1-x_1^2-x_2^2)^{1/2}$ and $x_3<0,$ $\tilde{g}(x_1,x_2)=-(1-x_1^2-x_2^2)^{1/2}$. For any bounded convex sets, $\mathbb{D}$, we also define its geodesic as the shortest path between two points $\mathbf{x},\mathbf{y}\in\mathbb{D}$ such that every point in the path is also an element of $\mathbb{D}$ and denote the geodesic distance by $d:\mathbb{D}\times\mathbb{D}\mapsto \mathbb{R}_+$, where $\mathbb{R}_+$ is the positive real line including $0$, thus $(\mathbb{D},d(\cdot,\cdot))$ is a metric space. Additionally, we will frequently need to evaluate integrals over $\mathbb{D}$, which can be done using its infinitesimal area element defined as,
\begin{equation*}
d\mathbb{D}=\sqrt{1+\left(\frac{\partial \tilde{g}}{\partial x_1}\right)^2+\left(\frac{\partial \tilde{g}}{\partial x_2}\right)^2}dx_1 dx_2.
\end{equation*} 
We assume that these convex subspaces of $\mathbb{R}^3$ are defined such that the origin is inside $\mathbb{D}$, that is $\mathbf{0}\in \mathbb{D}_{int}$, we then say the space $\mathbb{D}$ is centred. Our methodology can easily be adapted for non-centred spaces by making the appropriate translations to bring the origin inside $\mathbb{D}$.

Following the notation of \cite{Moller2004}, we define $\lambda_{\mathbb{D}}(\mathbf{x})$ as the Lebesgue measure restricted to the surface of the convex shape $\mathbb{D}$. Consider point processes which lie on some bounded convex metric space $(\mathbb{D},d(\cdot,\cdot))$. We define the notation $A_K=A\cap K$, $A,K\subseteq \mathbb{D}$. This nomenclature is often used when the set $A$ has finite cardinality and $K$ is any subset of $\mathbb{D}$. The cardinality of a set is denoted by $|\cdot|$. Further define $B_{\mathbb{D}}(\mathbf{x},r)=\{\mathbf{y}\in \mathbb{D},  d(\mathbf{x},\mathbf{y})\leq r\}$ and the set $N_{lf}=\{A\subset \mathbb{D}:|A_K|<\infty, K\subseteq \mathbb{D} \}$ where $K$ is any subset of $\mathbb{D}$. In other words, $N_{lf}$ is the set of subsets of $\mathbb{D}$ that have finite cardinality. To distinguish between points in a point process $X$ and any point in the space $\mathbb{D}$, we shall refer to elements of our point process $\mathbf{x}\in X$ as \emph{events} whilst retaining the term \emph{point} for any point in $\mathbb{D}$. We consider point processes which are \emph{locally finite} and \emph{simple}. A point process, $X$, lying on $\mathbb{D}$ is said to be locally finite, if for any bounded set $K\subseteq \mathbb{D}$, the number of events of $X$ in $K$ is finite almost surely, i.e. $X\in N_{lf}$ almost surely. A simple point process is one in which no coincident events exist almost surely, in other words if $\mathbf{x}_i,\mathbf{x}_j\in X$ such that $i\neq j$ then $\mathbf{x}_i\neq\mathbf{x}_j$. We also define the counting measure of $X$ as $N_X(K)=|X\cap K|=|X_K|$. We denote the reduced Palm distribution of the point process $X$ by $P_{X^!_{\mathbf{x}}}$ and define $X^!_{\mathbf{x}}$ to be the point process following this density, referring to this as the reduced Palm process \cite{Moller2004}. For a point process $X$ on $\mathbb{D}$ we define the \emph{intensity measure} as the expected number of events of $X$ for any $K\subseteq \mathbb{D}$, i.e. $\mu(K)=\mathbb{E}[N_{X}(K)]$, whilst the \emph{intensity function} $\rho(\mathbf{x})$ for all $\mathbf{x}\in \mathbb{D}$, if it exists, is given by,
\begin{equation*}
\mu(K)=\int_K \rho(\mathbf{x})\lambda_{\mathbb{D}}(d\mathbf{x}),
\end{equation*}
where $\rho(\mathbf{x})\lambda_{\mathbb{D}}(d\mathbf{x})$ can be interpreted heuristically as the probability of an event of $X$ being in the infinitesimal area $\lambda_{\mathbb{D}}(d\mathbf{x})$.

We also define, for $n\in\mathbb{N}$, $\alpha^{(n)}(K_1,\dots,K_n)=\mathbb{E}\sum^{\neq}_{\mathbf{x}_1,\dots,\mathbf{x}_n\in X}\mathbbm{1}[\mathbf{x_1}\in K_1,\dots,\mathbf{x_n}\in K_n]$ as the \emph{$n^{th}$-order factorial moment measures} where the summation is taken over pairwise distinct sets of $\{\mathbf{x}_1,\dots,\mathbf{x}_n\}$. Further we shall assume there exists $\rho^{(n)}(\mathbf{x}_1,\dots,\mathbf{x}_n)$ such that,
\begin{equation}\label{alpha:rho:relation}
\alpha^{(n)}(K_1,\dots,K_n)=\int_{K_1}\cdots\int_{K_n} \rho^{(n)}(\mathbf{x}_1,\dots,\mathbf{x}_n)\lambda_{\mathbb{D}}(d\mathbf{x}_1)\cdots \lambda_{\mathbb{D}}(d\mathbf{x}_n),
\end{equation}
where $K_1,\dots,K_n\subseteq \mathbb{D}$. We can interpret $\rho(\mathbf{x}_1,\dots,\mathbf{x}_n)^{(n)}\lambda_{\mathbb{D}}(d\mathbf{x}_1)\cdots \lambda_{\mathbb{D}}(d\mathbf{x}_n)$ as the probability that events of $X$ lie jointly in the infinitesimal areas $\lambda_{\mathbb{D}}(d\mathbf{x}_i),\; i=1,\dots,n$ and call $\rho^{(n)}$ the \emph{$n^{th}$-order factorial moment density}. Notice that $\alpha^{(1)}=\mu$, and $\rho^{(1)}=\rho$. The \emph{pair correlation function} is defined as,
\begin{equation*}
h(\mathbf{x},\mathbf{y})=\frac{\rho^{(2)}(\mathbf{x},\mathbf{y})}{\rho(\mathbf{x})\rho(\mathbf{y})}, 
\end{equation*}
where it is taken that division by $0$ results in the pair correlation function equalling $0$. 

A useful alternative to the $n^{th}$ order product intensities are the \emph{$n^{th}$ order correlation functions} \cite{vanLieshout2011}. They are recursively defined for $n\in\mathbb{N}$, based on product densities, with $\xi_1=1$ and
\begin{equation*}
\frac{\rho^{(n)}(\mathbf{x}_1,\dots,\mathbf{x}_n)}{\rho(\mathbf{x}_1)\cdots\rho(\mathbf{x}_n)}=\sum_{k=1}^n\sum_{D_1,\dots,D_k}\xi_{|D_1|}(\mathbf{x}_{D_1})\cdots\xi_{|D_k|}(\mathbf{x}_{D_k}),
\end{equation*}
where the final sum ranges over all partitions $\{D_1,\dots,D_k\}$ of $\{1,\dots,n\}$ in $k$ non-empty, disjoint sets, $\mathbf{x}_{D_j}=\{\mathbf{x}_i:i\in D_j\},\; j=1,\dots,k$, and $\mathbf{x}_i\in \mathbb{S}^2$ \cite{vanLieshout2011}. Further, we define the \emph{generating functional} \cite{Moller2004} of a point process $X$ as 
\begin{equation*}
G_X(u)=\mathbb{E}\prod_{\mathbf{x}\in X}u(\mathbf{x}),
\end{equation*}
for a function $u:\mathbb{D}\mapsto [0,1]$, which are also useful when discussing the $F$-, $H$-, and $J$-functions in Section \ref{Section:12}. 

We define a Poisson process on $\mathbb{D}$ identically to one on $\mathbb{R}^2$. Let $X$ be a point process on $\mathbb{D}$ such that $N_X(\mathbb{D})\sim \text{Poisson}(\mu(\mathbb{D}))$ where,
\begin{equation*}
\mu(B)=\int_B \rho(\mathbf{x}) \lambda_{\mathbb{D}}(d\mathbf{x}), \quad \text{for } B\subseteq\mathbb{D},
\end{equation*}
where $\mu(\mathbb{D})<\infty$. Then given $N_X(\mathbb{D})=n$, $\mathbf{x}_i\in X,\;i=1,\dots,n$ are independent and identically distributed across $\mathbb{D}$ with density proportional to $\rho(\mathbf{x})$. We say that $X$ is a Poisson process on $\mathbb{D}$ with intensity function $\rho:\mathbb{D}\mapsto\mathbb{R}_+$. When $\rho\in\mathbb{R}_+$ is constant we say the process is homogeneous Poisson or completely spatial random (CSR).

\subsection{Statement of the problem}

We are now in a position to formally state the hypothesis of CSR we are interested in and for which this work provides an approach to testing.

\begin{quote}
Let $X$ be a spatial point process, such that $g(X)=0$ where $g(X)$ is a notational convenience for $g(\mathbf{x})=0, \text{ for all } \mathbf{x}\in X$ and $g$ is the level set for the convex shape $\mathbb{D}$. From a realisation of $X$ we wish to conduct the following hypothesis test,
\begin{equation*}
H_0: X \text{ is CSR on } \mathbb{D} \quad\quad\text{vs.}\quad\quad H_1: X \text{ is not CSR on } \mathbb{D}.
\end{equation*}
\end{quote}

\section{Summary statistics on $\mathbb{S}^2$ and the impracticalities of defining functional summary statistics on general convex shapes}\label{Section:4}

To analyse point patterns on a convex shape $\mathbb{D}$, we will be required to map the point pattern onto the unit sphere $\mathbb{S}^2$. Hence it is necessary to discuss functional summary statistics in this space. We also explain why it is non trivial to construct analogous functional summary statistics directly on $\mathbb{D}$, motivating a need for new methodology in order to test patterns which arise on such surfaces.

\subsection{Summary statistics on $\mathbb{S}^2$}

Mentioned in passing by \cite{Ripley1977}, spherical point patterns only garnered interest over the past five years \cite{Lawrence2016,Moller2016,Robeson2014}. \cite{Robeson2014} shows that, on a sphere of radius $R$, the spherical $K$-function for a homogeneous Poisson process is $K(r)=2\pi R^2(1-\cos\left(r/R\right))$, where $r$ is the geodesic distance from an arbitrary point from the process. Building on this, both \cite{Lawrence2016} and \cite{Moller2016} define a range of typical functional summary statistics for isotropic (rotationally invariant) point processes, including the empty-space and spherical contact distributions. They also extend the $K$-function to the class of inhomogeneous point processes that have rotationally invariant pair correlation functions. This is analogous to the inhomogeneous extension given by \cite{Baddeley2000} for point processes on $\mathbb{R}^d$.

The geodesic on $\mathbb{S}^2$ is commonly referred to as the great circle distance and for two points $\mathbf{x},\mathbf{y}\in\mathbb{S}^2$ has analytic form $d(\mathbf{x},\mathbf{y})=\cos^{-1}(\mathbf{x}\cdot\mathbf{y}),$ where $\mathbf{x}\cdot\mathbf{y}$ is the dot product between vectors $\mathbf{x}, \mathbf{y}\in\mathbb{S}^2$. We define a point process $X$ on $\mathbb{S}^2$ to be isotropic if its distribution is invariant under rotations, i.e. $X\stackrel{d}{=}OX,$ for all $O\in\mathcal{O}(3),$ where $\mathcal{O}(3)$ is the set of orthogonal $3\times 3$ matrices \cite{Moller2016}. Here we use $\stackrel{d}{=}$ to denote \emph{equal in distribution}. Such an isotropic process has constant intensity function $\rho\in\mathbb{R}$.

Functional summary statistics are frequently employed for both exploratory data analysis and model fitting, playing a pivotal role in the early stages of any in depth investigation of an observed point pattern. In the homogeneous case, let $X$ be an isotropic spheroidal point process with constant intensity function $\rho\in\mathbb{R}$ then the $F$-, $H$-, $J$-, and $K$-functions are defined as,
\begin{align}
F(r)&=P(X_{B(\boldsymbol{o},r)}\neq\emptyset)\label{eq:F:stationary}\\
H(r)&=P(X^{!}_{\boldsymbol{o},B(\boldsymbol{o},r)}\neq\emptyset)\label{eq:H:stationary}\\
J(r)&=\frac{1-H(r)}{1-F(r)}\label{eq:J:stationary}\\
K(r)&=\frac{1}{\rho}\mathbb{E}\sum_{\mathbf{x}\in X^{!}_{\boldsymbol{o}}}\mathbbm{1}[d(\mathbf{x},\mathbf{o})\leq r],\label{eq:K:stationary}
\end{align}
where $r\in[0,\pi]$ and $\mathbf{o}=(0,0,1)^T$ is defined as the origin of $\mathbb{S}^2$. Estimators of these functional summary statistics can be used to determine whether the underlying process of an observed point pattern follows a specific distribution. In particular, they can be used to test whether a pattern arises from a CSR process or whether the underlying process exhibits regularity or clustering. A treatment of the standard isotropic functional summary statistics on $\mathbb{S}^2$ is given in \cite{Lawrence2016} and \cite{Moller2016}. 

It will be necessary for us to consider the inhomogeneous extensions of these functional summary statistics as they will form the foundation for our functional summary statistics for point process on convex bounded shapes in $\mathbb{R}^3$. We begin by reviewing the inhomogeneous $K$-function, originally attributed to \cite{Baddeley2000} for a class of inhomogeneous point processes in $\mathbb{R}^d$, and then extended to non-isotropic point processes on $\mathbb{S}^2$ by \cite{Lawrence2016} and \cite{Moller2016}. In Section \ref{Section:12}, we will construct the inhomogeneous $F$-, $H$-, and $J$-functions for non-isotropic point processes on $\mathbb{S}^2$. This builds on the formulation of \cite{vanLieshout2006} for non-stationary point processes in $\mathbb{R}^d$.

\subsection{Inhomogeneous $K$-function}

For the extension of the $K$-function to inhomogeneous processes on $\mathbb{R}^{2,3}$, \cite{Baddeley2000} introduce the notion of a point process being \emph{second order intensity reweighted stationary} (SOIRWS). Here we focus on the extension for $\mathbb{S}^2$ where \cite{Lawrence2016} and \cite{Moller2016} define the notion of \emph{second order intensity reweighted isotropic} (SOIRWI). A point process $X$ on $\mathbb{S}^2$ is said to be SOIRWI if its pair correlation function, $h$, is rotationally invariant, that is $h(\mathbf{x},\mathbf{y}) = h(d(\mathbf{x},\mathbf{y}))$, where $d$ is the great circle distance on $\mathbb{S}^2$. \cite{Lawrence2016} and \cite{Moller2016} define the inhomogeneous $K$-function for a SOIRWI process as,
\begin{equation*}
K_{\text{inhom}}(r)=\frac{1}{\lambda_{\mathbb{S}^2}(A)}\mathbb{E}\sum_{\mathbf{x},\mathbf{y}\in X}^{\neq}\frac{\mathbbm{1}[\mathbf{x}\in A, O_{\mathbf{x}}(\mathbf{y})\in B_{\mathbb{S}^2}(\mathbf{o},r)]}{\rho(\mathbf{x})\rho(\mathbf{y})},\quad 0\leq r \leq \pi,
\end{equation*}
where $O_{\mathbf{x}}:\mathbb{S}^2\mapsto\mathbb{S}^2$ is a rotation that takes $\mathbf{x}$ to $\mathbf{o}$. Here, $K_{\text{inhom}}$ is independent of the choice of $A\subseteq \mathbb{S}^2$ with $\lambda_{\mathbb{S}^2}(A)>0$, and by convention $a/0=0,$ for $a\geq 0$. For a Poisson process it is easy to show that $K_{\text{inhom}}(r)=2\pi(1-\cos(r))$. Therefore the $K$-function is the same for all Poisson processes regardless of whether the intensity function is constant or not \cite{Moller2016}.

Both \cite{Lawrence2016} and \cite{Moller2016} propose the following estimator for $K_{\text{inhom}}$ for a fully observed point pattern on $\mathbb{S}^2$,
\begin{equation}\label{k:inhom:sphere:estimator}
\hat{K}_{\text{inhom}}(r)=\frac{1}{4\pi}\sum_{\mathbf{x},\mathbf{y}\in X}^{\neq}\frac{\mathbbm{1}[d(\mathbf{x},\mathbf{y})\leq r]}{\rho(\mathbf{x})\rho(\mathbf{y})},
\end{equation} 
which is unbiased if $\rho$ is known. In the more likely event that $\rho$ is unknown, \cite{Lawrence2016} and \cite{Moller2016} suggest using a plugin estimator for $\rho(\mathbf{x})\rho(\mathbf{y})$.

\subsection{Impracticalities of defining functional summary statistics directly on $\mathbb{D}$}

We now explain the subtle reasoning as to why constructing functional summary statistics directly on $\mathbb{D}$ is not a trivial extension from $\mathbb{S}^2$. The definitions given by Equations (\ref{eq:F:stationary})-(\ref{eq:K:stationary}) are well defined when considering stationary or isotropic point processes on $\mathbb{R}^d$ or $\mathbb{S}^2$ respectively. This is because the symmetries of the space admit well defined notions of stationarity/isotropy based on translations and rotations. Since an arbitrary convex space $\mathbb{D}$ does not, in general, have isometries these notions of stationarity/isotropy cannot be well defined. Therefore, defining functional summary statistics analogous to (\ref{eq:F:stationary})-(\ref{eq:K:stationary}) is not possible. 

Further, we also argue that we cannot define a point process to be SOIRWI on $\mathbb{D}$. On $\mathbb{S}^2$ being SOIRWI is equivalent to having a rotationally invariant pair correlation function. We may be tempted to equivalently define a point process to be SOIRWI on $\mathbb{D}$ if it has an invariant form for its pair correlation function. In particular this would make sense for a Poisson process on $\mathbb{D}$ as it would have pair correlation function, $g(\mathbf{x})=1$ for all $\mathbf{x}\in\mathbb{D}$. Closer inspection though leads us to conclude that this is not an appropriate definition for SOIRWI on $\mathbb{D}$. Based on \cite[Definition 4.5, p. 32]{Moller2004} we can take a point process $X$ with intensity function $\rho:\mathbb{S}^2\mapsto\mathbb{R}_+$ as being SOIRWI on $\mathbb{S}^2$ if the measure,
\begin{equation}\label{eq:second:measure}
\mathcal{K}(B)=\frac{1}{\lambda_{\mathbb{S}^2}(A)}\mathbb{E}\sum_{\mathbf{x},\mathbf{y}\in X}^{\neq}\frac{\mathbbm{1}[\mathbf{x}\in A,O_{\mathbf{x}}(\mathbf{y})\in B]}{\rho(\mathbf{x})\rho(\mathbf{y})}, \quad B\subseteq \mathbb{S}^2,
\end{equation}
does not depend on the choice of $A\subseteq \mathbb{S}^2$ for $0<|A|<\infty$, where we take $a/0=0$. $\mathcal{K}$ is then called the \emph{second order reduced moment measure}. If the pair correlation function exists and is invariant under rotations, then by the Campbell-Mecke Theorem \cite{Moller2004} it follows that
\begin{equation*}
\mathcal{K}(B)=\int_B h(\mathbf{x})d\mathbf{x}, \quad B\subseteq \mathbb{S}^2.
\end{equation*}
Thus on $\mathbb{S}^2$ a point process is SOIRWI if $h$ is invariant under rotations. Equation (\ref{eq:second:measure}) implicitly depends on rotations $O_{\mathbf{x}}(\mathbf{y})$. If we now consider a point process on $\mathbb{D}$, we cannot construct the second order reduced moment measure as, in general, we do not have an analogous isometry. This in turn means that we cannot define SOIRWI directly on $\mathbb{D}$ based on an invariance of the pair correlation function. 

Moreover, for a point process on $\mathbb{S}^2$, consider the more specific case when $B=B_{\mathbb{S}^2}(\mathbf{o},r), r>0$ in (\ref{eq:second:measure}). This is identically the inhomogeneous $K$-function. The indicator function of (\ref{eq:second:measure}) is still well-defined in the case of $\mathbb{S}^2$ such that we are counting the events of $X\setminus\{\mathbf{x}\}$ that are at most a distance $r$ from $\mathbf{x}\in X$. This same intuition could not equivalently be applied to point processes on a convex shape as the ball of radius $r$ from a point $\mathbf{x}$ on $\mathbb{D}$ also depends on $\mathbf{x}$, i.e. $B_{\mathbb{D}}(\mathbf{x},r)\subset\mathbb{D}$ is different for each $\mathbf{x}\in\mathbb{D}$. Thus it is not possible to directly define an inhomogeneous $K$-function on $\mathbb{D}$.

\section{Extending the inhomogeneous $F$-, $H$-, and $J$-functions to $\mathbb{S}^2$}\label{Section:12}

On $\mathbb{R}^d$, in the stationary case, it can be shown that the $F$-, and $H$-functions have infinite series representations \cite{White1979}, and further work by \cite{vanLieshout2006} also gives an infinite series representation for the $J$-function based on the $n^{th}$-order correlation functions. Theorem \ref{thm:infinite:series:f:h} gives an infinite series representation when the $n^{th}$-order reduced factorial moment measure of all $n$ exists, similar to \cite{White1979} but where the underlying space is $\mathbb{S}^2$.

\begin{thm}\label{thm:infinite:series:f:h}
Let $X$ be an isotropic spheroidal point process with constant intensity function $\rho$. Further we assume the existence of all $n^{th}$-order factorial moment measures for both $X$ and its reduced Palm process, $X^!_{\mathbf{x}}$. Then the $F$- and $H$-functions have the following series representation,
\begin{align*}
F(r)&=-\sum_{n=1}^\infty\frac{(-1)^n}{n!}\alpha^{(n)}(B_{\mathbb{S}^2}(\mathbf{o},r),\dots,B_{\mathbb{S}^2}(\mathbf{o},r))\\
H(r)&=-\sum_{n=1}^\infty\frac{(-1)^n}{n!}\alpha^{!(n)}_{\mathbf{o}}(B_{\mathbb{S}^2}(\mathbf{o},r)\dots,B_{\mathbb{S}^2}(\mathbf{o},r))
\end{align*}
where $\alpha^{(n)}$ and $\alpha^{!(n)}_{\mathbf{x}}$ are the factorial moment measure for $X$ and $X^!_{\mathbf{x}}$ and $B_{\mathbb{S}^2}(\mathbf{o},r)$ is the spherical cap of radius $r$ at the origin $\mathbf{o}\in \mathbb{S}^2$. These representations hold provided the series is absolutely convergent, that is if $\lim_{n\rightarrow\infty}|a_{n+1}/a_n| <1$ or $\limsup_{n\rightarrow\infty} (|a_n|)^{1/n} < 1$, where $a_n=((-1)^n/n!)\alpha^{(n)}(B_{\mathbb{S}^2}(\mathbf{o},r),\dots,B_{\mathbb{S}^2}(\mathbf{o},r))$ for the $F$-function or $a_n=((-1)^n/n!)$ $\alpha^{!(n)}_{\mathbf{o}}(B_{\mathbb{S}^2}(\mathbf{o},r)\dots,B_{\mathbb{S}^2}(\mathbf{o},r))$ for the $H$-function.
\end{thm}
\begin{proof}
See Theorem 1 in Appendix A.
\end{proof}
The following corrollary reduces the representations for the $F$-, and $H$-function for when the $n^{th}$-order product density exist. These representations are those used by \cite{vanLieshout2011}.

\begin{corr}\label{corr:f:h:infinite:series}
\cite{White1979} Under the same assumptions as Theorem \ref{thm:infinite:series:f:h}, let $X$ be an isotropic spheroidal point process with constant intensity function $\rho$. Further we assume the existence of all $n^{th}$-order product intensities for both $X$ and its reduced Palm process, $X^!_{\mathbf{x}}$. Then the $F$- and $H$-functions have the following series representation,
\begin{align*}
F(r)&=-\sum_{n=1}^\infty\frac{(-1)^n}{n!}\int_{B_{\mathbb{S}^2}(\mathbf{o},r)}\cdots\int_{B_{\mathbb{S}^2}(\mathbf{o},r)}\rho^{(n)}(\mathbf{x}_1,\dots,\mathbf{x}_n)\lambda_{\mathbb{S}^2}(d\mathbf{x}_1)\cdots \lambda_{\mathbb{S}^2}(d\mathbf{x}_n)\\
H(r)&=-\sum_{n=1}^\infty\frac{(-1)^n}{n!}\int_{B_{\mathbb{S}^2}(\mathbf{o},r)}\cdots\int_{B_{\mathbb{S}^2}(\mathbf{o},r)}\frac{\rho^{(n+1)}(\mathbf{o},\mathbf{x}_1,\dots,\mathbf{x}_n)}{\rho}\lambda_{\mathbb{S}^2}(d\mathbf{x}_1)\cdots \lambda_{\mathbb{S}^2}(d\mathbf{x}_n)
\end{align*}
provided the series is absolutely convergent, where $B_{\mathbb{S}^2}(\mathbf{o},r)$ is the spherical cap of radius $r$ at the origin $\mathbf{o}\in \mathbb{S}^2$.
\end{corr}
\begin{proof}
See Corollary 1 in Appendix A.
\end{proof}

Adapting the work of \cite{vanLieshout2011}, the $J$-function for an isotropic spheroidal point process, based on the series for the $F$-, and $H$-function given by Theorem \ref{thm:infinite:series:f:h}, has the following infinite series representation
\begin{equation*}
J(r)=1+\sum_{n=1}^\infty\frac{(-\rho)^n}{n!}J_n(r), \quad 0\leq r \leq \pi,
\end{equation*}
where $J_n(r)=\int_{B_{\mathbb{S}^2}(\mathbf{o},r)} \cdots  \int_{B_{\mathbb{S}^2}(\mathbf{o},r)}\xi^{(n+1)}(\mathbf{o},\mathbf{x}_1,\dots,\mathbf{x}_n)\lambda_{\mathbb{S}^2}(d\mathbf{x}_1)\cdots \lambda_{\mathbb{S}^2}(d\mathbf{x}_n)$, and $B_{\mathbb{S}^2}(\mathbf{o},r)$ is the spherical cap at the origin $\mathbf{o}\in \mathbb{S}^2$.

In order to define the inhomogeneous $F$-, $H$-, and $J$-function we first define the notion of \emph{iterative reweighted moment isotropic} (IRWMI) for a class of non-isotropic spheroidal point processes, similar to the notion of \emph{iterative reweighted moment stationary} in $\mathbb{R}^{2,3}$ \cite{vanLieshout2011}.
\begin{dfn}
A spheroidal point process $X$ is said to be IRWMI if, $\text{for all } n\in\mathbb{N}$, the $n^{th}$-order correlation functions are rotationally invariant. That is $\xi_n(\mathbf{x}_1,\dots,\mathbf{x}_n)=\xi_n(O\mathbf{x}_1,\dots,O\mathbf{x}_n)$ for all $n\in\mathbb{N}$ and $O\in\mathcal{O}(3)$.
\end{dfn}
Identically to the inhomogeneous $J$-function in $\mathbb{R}^{2,3}$ \cite{vanLieshout2011}, we define the inhomogeneous $J$-function on $\mathbb{S}^2$.
\begin{dfn}
For an IRWMI point process $X$ with intensity function $\rho:\mathbb{S}^2\mapsto\mathbb{R}$ such that $\bar{\rho}\equiv\inf_{\mathbf{x}\in\mathbb{S}^2}\rho(\mathbf{x})>0$,
\begin{equation*}
J_{\text{inhom}}(r)=1+\sum_{n=1}^{\infty}\frac{(-\bar{\rho})^n}{n!}J_n(t), \quad 0\leq r \leq \pi
\end{equation*}
where $J_n(r)=\int_{B_{\mathbb{S}^2}(\mathbf{o},r)}\cdots\int_{B_{\mathbb{S}^2}(\mathbf{o},r)}\xi_{n+1}(\mathbf{o},\mathbf{x}_1,\dots,\mathbf{x}_n)\lambda_{\mathbb{S}^2}(d\mathbf{x}_1)\cdots \lambda_{\mathbb{S}^2}(d\mathbf{x}_n)$ and the series is absolutely convergent. 
\end{dfn}
Notice that since the point process is IRWMI then the $J$-function does not depend on the origin $\mathbf{o}$, and furthermore when the point process is isotropic $J_{\text{inhom}}$ collapses down to $J$ since the intensity function is constant. In the context of $\mathbb{R}^d$, \cite{vanLieshout2011} shows that the inhomogeneous $J$-function can be written as the ratio of generating functionals of the point process. Here we easily adapt the theorem for IRWMI point processes on $\mathbb{S}^2$.
\begin{thm*}
For all $r\in [0,\pi]$ and $\mathbf{y}\in\mathbb{S}^2$,
\begin{equation*}
u^{\mathbf{y}}_r(\mathbf{x})=\frac{\bar{\rho}\mathbbm{1}[O_{\mathbf{y}}(\mathbf{x})\in B_{\mathbb{S}^2}(\mathbf{o},r)]}{\rho(\mathbf{x})},\quad \mathbf{x}\in\mathbb{S}^2,
\end{equation*}
where $O_{\mathbf{y}}:\mathbb{S}^2\mapsto\mathbb{S}^2$ is a rotation that maps $\mathbf{y}$ to $\mathbf{o}$. Assuming that the series $\sum_{n=1}^{\infty}\frac{\bar{\rho}^n}{n!}\int_{B_{\mathbb{S}^2}(\mathbf{o},r)}\cdots\int_{B_{\mathbb{S}^2}(\mathbf{o},r)}\frac{\rho^{(n)}(\mathbf{x}_1,\dots,\mathbf{x}_n)}{\rho(\mathbf{x}_1)\cdots\rho(\mathbf{x}_n)}$ $\lambda_{\mathbb{S}^2}(d\mathbf{x}_1)\cdots \lambda_{\mathbb{S}^2}(d\mathbf{x}_n)$ is absolutely convergent. Then under the further assumptions associated with the inhomogeneous $J$-function and the existence of all $n^{th}-$order intensity function $\rho^{!(n)}_{\mathbf{y}}$ for the reduced Palm distribution $X^!_{\mathbf{y}}$, $\forall \mathbf{y}\in\mathbb{S}^2,$
\begin{equation*}
J_{\text{inhom}}(r)=\frac{G^!_{\mathbf{y}}(1-u^{\mathbf{y}}_r)}{G(1-u^{\mathbf{y}}_r)}, \quad 0\leq r \leq \pi,
\end{equation*}
for when $G(1-u^{\mathbf{y}}_r)>0$, where $G^!_{\mathbf{y}}$ and $G$ are the generating functionals for $X^!_\mathbf{y}$ and $X$ respectively.
\end{thm*}
\begin{proof}
See Theorem 1 of \cite{vanLieshout2011}.
\end{proof}
From the proof given by \cite{vanLieshout2011}, it can be shown that the numerator and denominator do not depend on the arbitrary point $\mathbf{y}$. Further, in the case of an isotropic point process the numerator can be shown to be $G^!_{\mathbf{y}}(1-u^{\mathbf{y}}_r)=1-H(r),$ whilst the denominator is $G(1-u^{\mathbf{y}}_r)=1-F(r)$, and so the $F$-, and $H$-functions can be extended to the inhomogeneous case,
\begin{align*}
F_{\text{inhom}}(r)&=1-G(1-u^{\mathbf{y}}_r)\\
H_{\text{inhom}}(r)&=1-G^!_{\mathbf{y}}(1-u^{\mathbf{y}}_r),
\end{align*}
where the functions do not depend on the arbitrary point $\mathbf{y}$ of the point process.

Similar to \cite{vanLieshout2011} in $\mathbb{R}^{2,3}$, we propose the following estimators for the inhomogeneous $F$-, and $H$-functions for spheroidal IRWMI point processes as,
\begin{align}
\hat{F}_{\text{inhom}}(r)&=1-\frac{\sum_{\mathbf{p}\in P}\prod_{\mathbf{x}\in X\cap B_{\mathbb{S}^2}(\mathbf{p},r)}\left(1-\frac{\bar{\rho}}{\rho(\mathbf{x})}\right)}{|P|}\label{F:inhom:sphere:estimator}\\
\hat{H}_{\text{inhom}}(r)&=1-\frac{\sum_{\mathbf{x}\in X}\prod_{\mathbf{y}\in (X\setminus \{\mathbf{x}\})\cap B_{\mathbb{S}^2}(\mathbf{x},r)}\left(1-\frac{\bar{\rho}}{\rho(\mathbf{y})}\right)}{N_X(\mathbb{S}^2)},\label{h:inhom:sphere:estimator}
\end{align}
where $P\subseteq \mathbb{S}^2$ is a finite grid of points. The properties of the $\hat{F}_{\text{inhom}}$-function are independent of the choice of $P$ \cite{vanLieshout2011}. In this work we choose $P$ such that the points on $\mathbb{S}^2$ are equidistant. \cite{vanLieshout2011} show that $\hat{F}_{\text{inhom}}(r)$ is unbiased whilst $\hat{H}_{\text{inhom}}(r)$ is ratio-unbiased. Then since $\hat{F}_{\text{inhom}}$ 
is unbiased and $\hat{H}_{\text{inhom}}$ is ratio-unbiased, constructing $\hat{J}_{\text{inhom}}$ as
\begin{equation}\label{j:inhom:sphere:estimator}
\hat{J}_{\text{inhom}}(r)=\frac{1-\hat{H}_{\text{inhom}}(r)}{1-\hat{F}_{\text{inhom}}(r)},
\end{equation}
gives a ratio-unbiased estimator for $J_{\text{inhom}}(r)$.

\section{Summary statistics for Poisson processes on convex shapes}\label{Section:5}

Here, we construct summary statistics for Poisson processes on general convex shapes. We show that a Poisson process on a general convex shape, $\mathbb{D}$, can be mapped to a Poisson process on a sphere, and then define functional summary statistics for such processes. We discuss properties of these functional summary statistics in the more general setting of inhomogeneous Poisson processes on $\mathbb{S}^2$.

\subsection{Mapping from $\mathbb{D}$ to $\mathbb{S}^2$}

To circumvent the geometrical restrictions of $\mathbb{D}$ we show, in this section, that we can map Poisson processes from $\mathbb{D}$ to $\mathbb{S}^2$ and construct functional summary statistics in this space. Theorem \ref{thm:mapping:general:text} shows that a Poisson process on $\mathbb{D}$ can be transformed to a Poisson process on a sphere where we can take advantage of the rotational symmetries. The invariance of Poisson processes between metric spaces is known as the Mapping Theorem \cite{Kingman1993}. We use the function $f(\mathbf{x})=\mathbf{x}/||\mathbf{x}||$ to map point patterns from $\mathbb{D}$ to $\mathbb{S}^2$. Lemma \ref{lemma:f:bijective} shows that this function is bijective and hence measurable.

\begin{lemma}\label{lemma:f:bijective}
Let $\mathbb{D}$ be a convex subspace of $\mathbb{R}^3$ such that the origin in $\mathbb{R}^3$ is in the interior of $\mathbb{D}$, i.e. $\mathbf{o}\in\mathbb{D}_{int}$. Then the function $f(\mathbf{x})=\mathbf{x}/||\mathbf{x}||, f:\mathbb{D}\mapsto\mathbb{S}^2$ is bijective.
\end{lemma}  
\begin{proof}
See Lemma 1 in Appendix B. 
\end{proof}

Rather than using the Mapping Theorem \cite{Kingman1993}, we utilise Proposition 3.1 of \cite{Moller2004} to show that mapping a Poisson process from $\mathbb{D}$ to $\mathbb{S}^2$ results in a new Poisson process on $\mathbb{S}^2$ and also derive the intensity function of the mapped process on $\mathbb{S}^2$. 

\begin{thm}\label{thm:mapping:general:text}
Let $X$ be a Poisson process on an arbitrary bounded convex shape $\mathbb{D}\subset\mathbb{R}^3$ with intensity function $\rho:\mathbb{D}\mapsto \mathbb{R}$. We assume that $\mathbb{D}=\{\mathbf{x}\in\mathbb{R}^3: g(\mathbf{x})=0\}$ where $g(\mathbf{x})=0$ is the level-set function and is defined as,
\begin{equation*}
g(\mathbf{x})=
\left\{\begin{alignedat}{2}
    g_1(\mathbf{x})=0,&\quad \mathbf{x}\in\mathbb{D}_1\\
    &\vdots \\
    g_n(\mathbf{x})=0,&\quad \mathbf{x}\in\mathbb{D}_n
  \end{alignedat}\right.
\end{equation*}
such that $\cup_{i=1}^n\mathbb{D}_i=\mathbb{D}$ and $\mathbb{D}_i\cap\mathbb{D}_j=\emptyset,\; \forall i\neq j$. Let $Y=f(X)$, where $f(\mathbf{x})=\mathbf{x}/||\mathbf{x}||$ with $f(X)=\{\mathbf{y}\in\mathbb{S}^2: \mathbf{y}=\mathbf{x}/||\mathbf{x}||, \mathbf{x}\in X\}$. Then $Y$ is a Poisson process on $\mathbb{S}^2$, with intensity function,
\begin{equation}\label{rho:mapping}
\rho^*(\mathbf{x})=
\left\{\begin{alignedat}{2}
\rho(f^{-1}(\mathbf{x}))l_1(f^{-1}(\mathbf{x})) &J_{(1,f^*)}(\mathbf{x})\sqrt{1-x_1^2-x_2^2}, \quad \mathbf{x}\in f(\mathbb{D}_1)\\
&\vdots\\
\rho(f^{-1}(\mathbf{x}))l_n(f^{-1}(\mathbf{x})) &J_{(n,f^*)}(\mathbf{x})\sqrt{1-x_1^2-x_2^2}, \quad \mathbf{x}\in f(\mathbb{D}_n)
  \end{alignedat}\right.
\end{equation}
where,
\begin{align*}
x_3 &=\tilde{g}_i(x_1,x_2)\\
l_i(\mathbf{x}) &=\left[1+\left(\frac{\partial \tilde{g}_i}{\partial x_1}\right)^2+\left(\frac{\partial \tilde{g}_i}{\partial x_2}\right)^2\right]^{\frac{1}{2}}\\
J_{(i,f^{*-1})}(\mathbf{x})&=\frac{1}{(x_1^2+x_2^2+\tilde{g}_i^2(x_1,x_2))^3}\\
&\det\left[
\begin{pmatrix}
x_2^2+\tilde{g}_i^2(x_1,x_2)-x_1\tilde{g}_i(x_1,x_2)\frac{\partial \tilde{g}_i}{\partial x_1} & -x_1\left(x_2+\tilde{g}_i(x_1,x_2)\frac{\partial \tilde{g}_i}{\partial x_2}\right) \\
-x_2\left(x_1+\tilde{g}_i(x_1,x_2)\frac{\partial \tilde{g}_i}{\partial x_1}\right) & x_1^2+\tilde{g}_i^2(x_1,x_2)-x_2\tilde{g}_i(x_1,x_2)\frac{\partial \tilde{g}_i}{\partial x_2}
\end{pmatrix}\right]\\
J_{(i,f^*)}(\boldsymbol{\mathbf{x}}) &= \frac{1}{J_{(i,f^{*-1})}(f^{-1}(\boldsymbol{\mathbf{x}}))},
\end{align*}
where $f^{-1}$ is the inverse of $f$, $\det(\cdot)$ is the determinant operator, and $f^*:\mathbb{R}^2\mapsto\mathbb{R}^2$ is the function which maps $x_1\mapsto x_1/||\mathbf{x}||$ and $x_2\mapsto x_2/||\mathbf{x}||$.
\end{thm}
\begin{proof}
See Theorem 2 in Appendix B.
\end{proof}

\begin{rmk}\label{remark:bijectivity}
A notion of bijectivety arises from this theorem. Consider the set of all Poisson processes on $\mathbb{D}$ such that their intensity functions exist, label this set $T_{\mathbb{D}}$. Also define $T_{\mathbb{S}^2}$ as all the Poisson processes on $\mathbb{S}^2$ such that their intensity functions exits. Then for any $X\in T_{\mathbb{D}}$ implies that $f(X)\in T_{\mathbb{S}^2}$. Similarly by considering the inverse operation $f^{-1}$, which exists by Lemma \ref{lemma:f:bijective}, for all $Y\in T_{\mathbb{S}^2}$ implies that $f(Y)\in T_{\mathbb{D}}$. Hence the mapping $f:T_{\mathbb{D}}\mapsto T_{\mathbb{S}^2}$ is surjective. By Theorem \ref{thm:mapping:general:text} if $X,Y$ are Poisson processes on $\mathbb{D}$ with intensity function $\rho_X$ and $\rho_Y$ respectively then $f(X)$ and $f(Y)$ are the same Poisson process if and only if $\rho_X=\rho_Y$ and so the mapping is also injective, and hence bijective. This means that analysis of a Poisson process, $X$, on $\mathbb{D}$ is equivalent to the analysis of $f(X)$ on $\mathbb{S}^2$. 
\end{rmk} 

\begin{rmk}\label{remark:facet}
Further, another useful result which follows directly from Theorem \ref{thm:mapping:general:text} is the construction of approximate Poisson processes on $\mathbb{S}^2$. More precisely consider a convex surface $\mathbb{D}$ for which instead of having a level-set function, $g$, we have an approximation to the space, for example consider we have a finite piecewise planar approximation to $\mathbb{D}$. Then $\mathbb{D}$ can be approximated by $\cup_{i=1}^n \mathbb{D}_i,$ where each $\mathbb{D}_i$ is a planar piece and there are $n\in\mathbb{N}$ pieces to the approximation. For each $\mathbb{D}_i$ the level-set function is $g_i(\mathbf{x})=a_i x_1+b_i x_2+c_i x_3+d_i=0$, and we can then use this approximation of $\mathbb{D}$ to map a Poisson process on $\mathbb{D}$ to $\mathbb{S}^2$. 
\end{rmk}

\subsection{Construction of functional summary statistics}
We are now in a position to construct functional summary statistics for a Poisson process which lies on some bounded convex space $\mathbb{D}$. Since all Poisson processes on $\mathbb{S}^2$ are SOIRWI \cite{Moller2016} and IRWMI \cite{vanLieshout2011}, the estimators for $F_{\text{inhom}}, H_{\text{inhom}}, J_{\text{inhom}},$ and $K_{\text{inhom}}$ (see Equations \ref{F:inhom:sphere:estimator}-\ref{j:inhom:sphere:estimator} and \ref{k:inhom:sphere:estimator} respectively) \cite{vanLieshout2011,Moller2016} can be combined with the mapped intensity function from Theorem \ref{thm:mapping:general:text} to construct estimators as follows,
\begin{align}
\hat{F}_{\text{inhom},\mathbb{D}}(r)&=1-\frac{\sum_{\mathbf{p}\in P}\prod_{\mathbf{x}\in Y\cap B_{\mathbb{S}^2}(\mathbf{p},r)}\left(1-\frac{\bar{\rho^*}}{\rho^*(\mathbf{x})}\right)}{|P|}\label{F:estimator}\\
\hat{H}_{\text{inhom},\mathbb{D}}(r)&=1-\frac{\sum_{\mathbf{x}\in Y}\prod_{\mathbf{y}\in (Y\setminus \{\mathbf{x}\})\cap B_{\mathbb{S}^2}(\mathbf{x},r)}\left(1-\frac{\bar{\rho^*}}{\rho^*(\mathbf{y})}\right)}{N_Y(\mathbb{S}^2)}\label{H:estimator}\\
\hat{J}_{\text{inhom},\mathbb{D}}(r)&=\frac{1-\hat{H}_{\text{inhom},\mathbb{D}}(r)}{1-\hat{F}_{\text{inhom},\mathbb{D}}(r)}\label{J:estimator}\\
\hat{K}_{\text{inhom},\mathbb{D}}(r)&=\frac{1}{4\pi}\sum_{\mathbf{x},\mathbf{y}\in Y}^{\neq}\frac{\mathbbm{1}[d(\mathbf{x},\mathbf{y})\leq r]}{\rho^*(\mathbf{x})\rho^*(\mathbf{y})},\label{K:estimator}
\end{align}
where $X$ is a Poisson process on $\mathbb{D}$ with intensity function $\rho$, $Y=f(X)$ is the mapped Poisson process onto $\mathbb{S}^2$, $\rho^*$ is given by (\ref{rho:mapping}) and $\bar{\rho^*}=\inf_{\mathbf{x}\in\mathbb{S}^2}\rho^*(\mathbf{x})$. In the event that $\rho:\mathbb{D}\mapsto\mathbb{R}_+$ is unknown and therefore $\rho^*$ is unknown, nonparametric plug-in estimates of $\rho^*$ can be constructed on $\mathbb{S}^2$ \cite{Lawrence2016,Moller2016}. 

\subsection{Properties of functional summary statistics}
Consider the general case of all Poisson processes on $\mathbb{S}^2$. Theorem \ref{mean:functional:summary:statistics} gives the expectations of $\hat{F}_{\text{inhom}}(r)$, $\hat H_{\text{inhom}}(r)$, and $\hat K_{\text{inhom}}(r)$. We restate the mean of $\hat{K}_{\text{inhom}}$ \cite{Lawrence2016,Moller2016} and adapt the proof to Proposition 1 in \cite{vanLieshout2011} for $\mathbb{R}^d$, to show that $\hat{F}_{\text{inhom}}$ is unbiased and $\hat{H}_{\text{inhom}}$ is ratio unbiased for $\mathbb{S}^2$. In addition we also provide the expectation of $\hat{H}_{\text{inhom}}(r)$.

\begin{thm}\label{mean:functional:summary:statistics}
Let $X$ be a spherical Poisson process on $\mathbb{S}^2$ with known intensity function $\rho:\mathbb{S}^2\mapsto \mathbb{R}_+$, such that $\bar{\rho}=\inf_{\mathbf{x}\in\mathbb{S}^2}\rho(\mathbf{x})>0$. Then the estimators for $\hat{F}_{\text{inhom}}(r)$, and $\hat K_{\text{inhom}}(r)$ are unbiased whilst $\hat{H}_{\text{inhom}}(r)$ is ratio-unbiased. More precisely, 
\begin{align*}
\mathbb{E}[\hat{F}_{\text{inhom}}(r)]&=1-\exp(-\bar{\rho}2\pi(1-\cos r))\\
\mathbb{E}[\hat{H}_{\text{inhom}}(r)]&=1-\frac{\mathrm{exp}(-\bar{\rho}2\pi(1-\cos r))-\mathrm{exp}(-\mu(\mathbb{S}^2))}{1-\frac{\bar{\rho}2\pi(1-\cos r)}{\mu(\mathbb{S}^2)}}\\
\mathbb{E}[\hat{K}_{\text{inhom}}(r)]&=2\pi (1-\cos r),
\end{align*}
where $r\in [0,\pi]$, and $\bar{\rho}=\inf_{\mathbf{x}\in\mathbb{S}^2}\rho(\mathbf{x})>0$. Further by unbiasedness and ratio-unbiasedness of $\hat{F}_{\text{inhom}}(r)$ and $\hat{H}_{\text{inhom}}(r)$, respectively, we immediately have ratio-unbiasedness of $\hat{J}_{\text{inhom}}(r)$. 
\end{thm}
\begin{proof}
See  \cite{Lawrence2016} for treatment of $\hat{K}_{\text{inhom}}(r)$. Results for $\hat{F}_{\text{inhom}}(r)$ and $\hat{H}_{\text{inhom}}(r)$ follow from a trivial adaptation of the proof for Proposition 1 in \cite{vanLieshout2011}. For the expectation of $\hat{H}_{\text{inhom}}(r)$ see Theorem 3 in Appendix C.
\end{proof}

Theorem \ref{mean:functional:summary:statistics} shows that $\hat{H}_{\text{inhom}}(r)$ is a biased estimator for $H_{\text{inhom}}(r)$. Although biased it can be bounded.

\begin{corr}\label{corr:bound:H:estimate}
With the same assumptions as Theorem \ref{mean:functional:summary:statistics}, let $X$ be a spherical Poisson process on $\mathbb{S}^2$ with intensity function $\rho:\mathbb{S}^2\mapsto \mathbb{R}_+$. Defining $\bar{\rho}=\inf_{\mathbf{x}\in\mathbb{S}^2}\rho(\mathbf{x})$, the bias of the estimator $\hat{H}_{\text{inhom}}(r)$ is bounded by
\begin{equation*}
|\text{Bias}(\hat{H}_{\text{inhom}}(r))|\leq \exp(-\mu(\mathbb{S}^2))\leq \exp(-4\pi\bar{\rho}),
\end{equation*}
for all $r\in[0,\pi]$.
\end{corr}
\begin{proof}
See Corollary 2 in Appendix C.
\end{proof}
Corollary \ref{corr:bound:H:estimate} shows that, depending on the intensity function and hence $\bar{\rho}=\inf_{\mathbf{x}\in\mathbb{S}^2\rho(\mathbf{x})}$, the bias can be considered negligible. In the examples to come we set the expected number of points of the process to be large enough for the bias to be considered negligible. Next we provide the variance of the estimators of the functional summary statistics.

\begin{thm}\label{thm:variance:K:inhom:text}
Let $X$ be a spherical Poisson process on $\mathbb{S}^2$ with known intensity function $\rho:\mathbb{S}^2\mapsto \mathbb{R}_+$, such that $\bar{\rho}=\inf_{\mathbf{x}\in\mathbb{S}^2}\rho(\mathbf{x})>0$. Then the estimators $\hat{K}_{\rm{inhom}}(r)$, $\hat{F}_{\rm{inhom}}(r)$, and $\hat{H}_{\rm{inhom}}(r)$ have variance,
\begin{align*}
&\emph{Var}(\hat{K}_{\rm{inhom}}(r)) = \frac{1}{8\pi^2}\int_{\mathbb{S}^2}\int_{\mathbb{S}^2}\frac{\mathbbm{1}[d(\mathbf{x},\mathbf{y})\leq r]}{\rho(\mathbf{x})\rho(\mathbf{y})}\lambda_{\mathbb{S}^2}(d\mathbf{x})\lambda_{\mathbb{S}^2}(d\mathbf{y}) +(1-\cos r)^2\int_{\mathbb{S}^2}\frac{1}{\rho(\mathbf{x})}\lambda_{\mathbb{S}^2}(d\mathbf{x}),\\
&\emph{Var}(\hat{F}_{\rm{inhom}}(r)) = \frac{\exp\left(-2\bar{\rho}\lambda_{\mathbb{S}^2}(B_{\mathbb{S}^2}(\mathbf{o},r))\right)}{|P|^2}\\
&\phantom{AAAAAAAA}\sum_{\mathbf{p}\in P}\sum_{\mathbf{p}'\in P} \exp\left(\int_{B_{\mathbb{S}^2}(\mathbf{p},r)\cap B_{\mathbb{S}^2}(\mathbf{p}',r)}\frac{\bar{\rho}^2}{\rho(\mathbf{x})}\lambda_{\mathbb{S}^2}(d\mathbf{x})\right)-\exp\left(-2\bar{\rho}\lambda_{\mathbb{S}^2}(B_{\mathbb{S}^2}(\mathbf{o},r))\right),\\
&\emph{Var}(\hat{H}_{\rm{inhom}}(r))\\
&=\frac{1}{\mu^2(\mathbb{S}^2)}\int_{\mathbb{S}^2}\int_{\mathbb{S}^2}\left(\rho(\mathbf{x})-\bar{\rho}\mathbbm{1}[\mathbf{x}\in B_{\mathbb{S}^2}(\mathbf{y},r)]\right)\left(\rho(\mathbf{y})-\bar{\rho}\mathbbm{1}[\mathbf{y}\in B_{\mathbb{S}^2}(\mathbf{x},r)]\right)\\
&\phantom{=+}\frac{\mathrm{e}^{-\mu(\mathbb{S}^2)}}{A_1^{2}(\mathbf{x},\mathbf{y})}\left(\mathrm{e}^{\mu(\mathbb{S}^2) A_1(\mathbf{x},\mathbf{y})}-1 -\rm{Ei}(\mu(\mathbb{S}^2) A_1(\mathbf{x},\mathbf{y}))+\gamma +\log(\mu(\mathbb{S}^2) A_1(\mathbf{x},\mathbf{y}))\right) \lambda_{\mathbb{S}^2}(d\mathbf{x}) \lambda_{\mathbb{S}^2}(d\mathbf{y})\\
&\phantom{=}+\frac{1}{\mu(\mathbb{S}^2)}\int_{\mathbb{S}^2}\frac{\mathrm{e}^{-\mu(\mathbb{S}^2)}}{A_2(\mathbf{x})}\left(\gamma +\log(\mu(\mathbb{S}^2) A_2(\mathbf{x}))-\rm{Ei}(\mu(\mathbb{S}^2) A_2(\mathbf{x}))\right)\rho(\mathbf{y})\lambda_{\mathbb{S}^2}(d\mathbf{y})\\
&\phantom{=}-\frac{\mathrm{e}^{-2\mu(\mathbb{S}^2)}}{\left(1-\frac{\bar{\rho}}{\mu(\mathbb{S}^2)}2\pi(1-\cos r)\right)^2}\left(\mathrm{e}^{\mu(\mathbb{S}^2) \left(1-\frac{\bar{\rho}}{\mu(\mathbb{S}^2)}2\pi(1-\cos r)\right)}-1\right)^2
\end{align*}
where,
\begin{align*}
A_1(\mathbf{x},\mathbf{y})&= 1-\frac{2\bar{\rho}}{\mu(\mathbb{S}^2)}2\pi(1-\cos r)+\frac{\bar{\rho}^2}{\mu(\mathbb{S}^2)}\int_{B_{\mathbb{S}^2}(\mathbf{x},r)\cap B_{\mathbb{S}^2}(\mathbf{y},r)}\frac{1}{\rho(\mathbf{z})}d\mathbf{z}\\
A_2(\mathbf{x})&= 1-\frac{2\bar{\rho}}{\mu(\mathbb{S}^2)}2\pi(1-\cos r)+\frac{\bar{\rho}^2}{\mu(\mathbb{S}^2)}\int_{B_{\mathbb{S}^2}(\mathbf{x},r)}\frac{1}{\rho(\mathbf{y})}\lambda_{\mathbb{S}^2}(d\mathbf{y})\\
\rm{Ei}(x) &= -\int_{-x}^{\infty}\frac{\mathrm{e}^{-t}}{t} dt
\end{align*}
and $\rm{Ei}(x)$ is the exponential integral and $r\in [0,\pi]$.
\end{thm}
\begin{proof}
See Theorem 4 in Appendix D.
\end{proof}

Due to the complexity of the estimator for the $\hat{J}_{\text{inhom}}$-function, its mean and variance are extremely complex and although can be derived in terms of integrals over $\mathbb{S}^2$, we instead give an approximation based on the Taylor series expansion of the function $f(x,y)=x/y$ around the means of the numerator and denominator. We first provide conditions for which the first two moments of $\hat{J}_{\text{inhom}}(r)$ exist and then proceed to show how it can be approximated.
\begin{thm}
Let $X$ be a spheroidal Poisson process with intensity function $\rho:\mathbb{S}^2\mapsto\mathbb{R}_+$ such that $\bar{\rho}\equiv\inf_{\mathbf{x}\in\mathbb{S}^2}\rho(\mathbf{x})>0$. Let $P$ be any finite grid on $\mathbb{S}^2$ and define $r_{\max}=\sup\{r\in [0,\pi]:\text{ there exists } \mathbf{p}\in P \text{ such that } \rho(\mathbf{x})\neq \bar{\rho} \text{ for all } \mathbf{x} \in B_{\mathbb{S}^2}(\mathbf{p},r)\}$. Then for any given $r\in [0,r_{\max}]$ both $\mathbb{E}[\hat{J}_{\rm{inhom}}(r)]$ and ${\rm{Var}}(\hat{J}_{\rm{inhom}}(r))$ exist.
\end{thm}
\begin{proof}
See Theorem 5 in Appendix E.
\end{proof}
\begin{prop}\label{prop:J:moment:approx}
  Let $X$ be a spheroidal Poisson process with known intensity function $\rho:\mathbb{S}^2\mapsto\mathbb{R}$. Then the covariance between $1-\hat{H}_{\text{inhom}}(r)$ and $1-\hat{F}_{\text{inhom}}(r)$ for $r\in [0,\pi]$ is,
  \begin{equation*}
  \begin{split}
  &\text{Cov}(1-\hat{H}_{\text{inhom}}(r),1-\hat{F}_{\text{inhom}}(r))\\
  &=\frac{1}{|P|}\sum_{\mathbf{p}\in P}\int_{\mathbb{S}^2}\left(1-\frac{\bar{\rho}\mathbbm{1}[\mathbf{x}\in B_{\mathbb{S}^2}(\mathbf{p},r)]}{\rho(\mathbf{x})}\right)\\
  &\phantom{=-}\frac{\mathrm{exp}\left\{-2\bar{\rho}2\pi(1-\cos r)-\int_{B_{\mathbb{S}^2}(\mathbf{x},r)\cap B_{\mathbb{S}^2}(\mathbf{p},r)}\frac{\bar{\rho}^2}{\rho(\mathbf{y})}\lambda_{\mathbb{S}^2}(d\mathbf{y})\right\}}{A(\mathbf{x},\mathbf{p})}\frac{\rho(\mathbf{x})}{\mu(\mathbb{S}^2)}\lambda_{\mathbb{S}^2}(d\mathbf{x})\\
  &\phantom{=}-\exp(-2\pi(1-\cos r)\bar{\rho})\left(\exp(-2\pi(1-\cos r)\bar{\rho}\right)-\exp(-\mu(\mathbb{S}^2))\frac{\mu(\mathbb{S}^2)}{\mu(\mathbb{S}^2)-2\pi(1-\cos r)\bar{\rho}},
  \end{split}
  \end{equation*}
  where $P$ is a finite grid of points on $\mathbb{S}^2$ and,
  \begin{equation*}
  A(\mathbf{x},\mathbf{p})=1-\frac{2\bar{\rho}}{\mu(\mathbb{S}^2)}2\pi(1-\cos r) +\frac{1}{\mu(\mathbb{S}^2)}\int_{ B_{\mathbb{S}^2}(\mathbf{x},r)\cap B_{\mathbb{S}^2}(\mathbf{p},r)}\frac{\bar{\rho}^2}{\rho(\mathbf{y})}\lambda_{\mathbb{S}^2}(d\mathbf{y}).
  \end{equation*}
\end{prop}
\begin{proof}
  See Proposition 1 in Appendix E.
\end{proof}
Using a Taylor series expansion (see Section S5.2 of the Supplementary Material), we can approximate the expectation and variance of $\hat{J}_{\text{inhom}}(r)$ as
\begin{align}
\mathbb{E}\left[\frac{X}{Y}\right] &\approx\frac{\mu_X}{\mu_Y}-\frac{\text{Cov}(X,Y)}{\mu_Y^2}+\frac{\text{Var}(Y)\mu_X}{\mu_Y^3}\label{mean:Taylor}\\
\text{Var}\left(\frac{X}{Y}\right) &\approx \frac{\mu_X}{\mu_Y}\left[\frac{\text{Var}(X)}{\mu_X^2}-2\frac{\text{Cov}(X,Y)}{\mu_X\mu_Y}+\frac{\text{Var}(Y)}{\mu_Y^2}\right],\label{variance:Taylor}
\end{align}
where $X=1-\hat{H}_{\text{inhom}}(r)$ and $Y=1-\hat{F}_{\text{inhom}}(r)$. The terms in Equations (\ref{mean:Taylor}) and (\ref{variance:Taylor}) are given in Theorems \ref{mean:functional:summary:statistics} and \ref{thm:variance:K:inhom:text}, and Proposition \ref{prop:J:moment:approx}.

\section{Examples}\label{Section:11}
We now look at two examples where we simulate homogeneous Poisson processes on their surfaces and construct the previously described functional summary statistics. 
\subsection{Cube}
\begin{figure}[t]
\begin{subfigure}{0.32\textwidth}
\centering
\includegraphics[trim={5em 5em 5em 5em},clip,width=\textwidth,height=!]{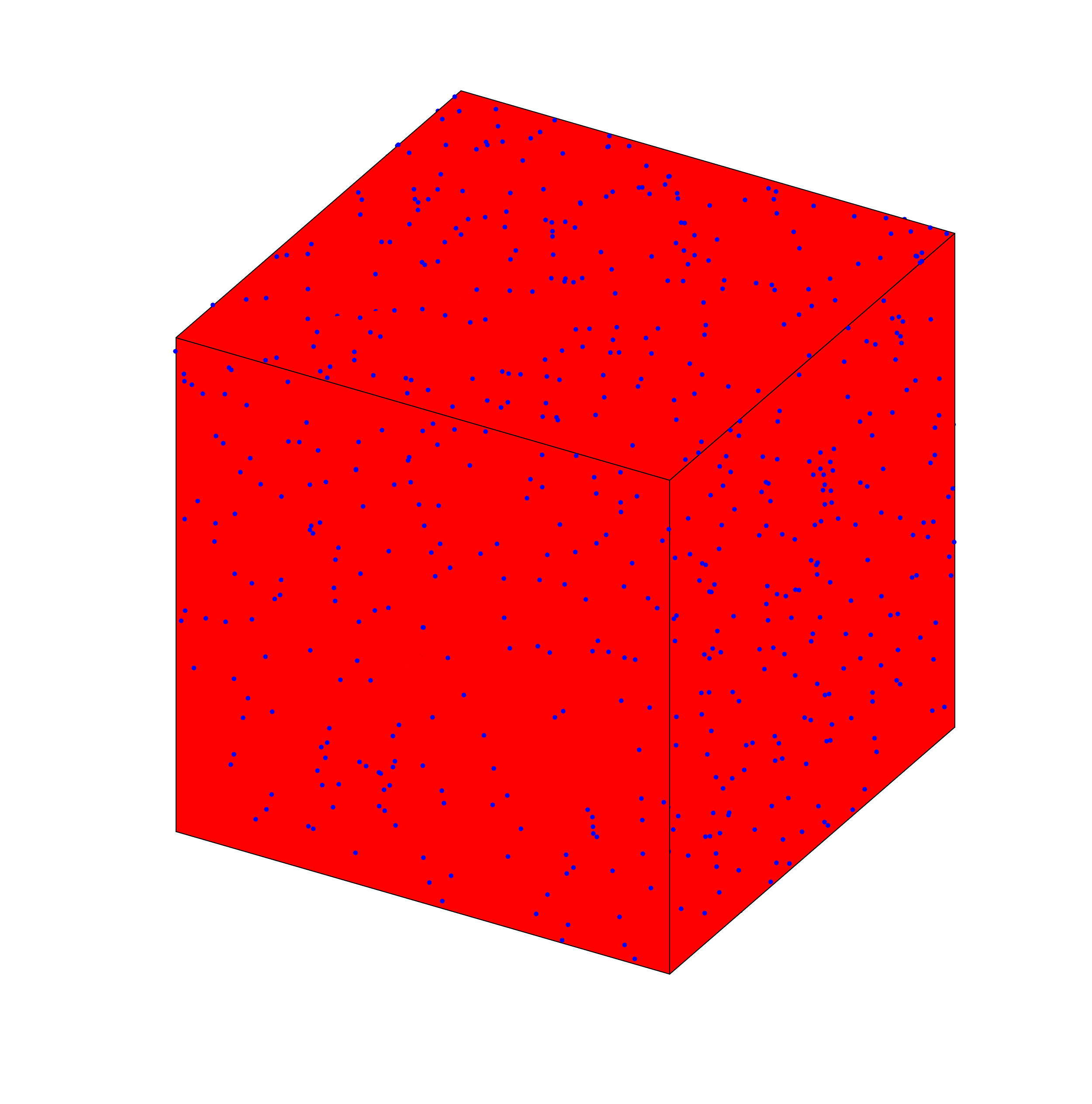}
\end{subfigure}
\begin{subfigure}{0.32\textwidth}
\centering
\includegraphics[trim={5em 5em 5em 4em},clip,width=\textwidth,height=!]{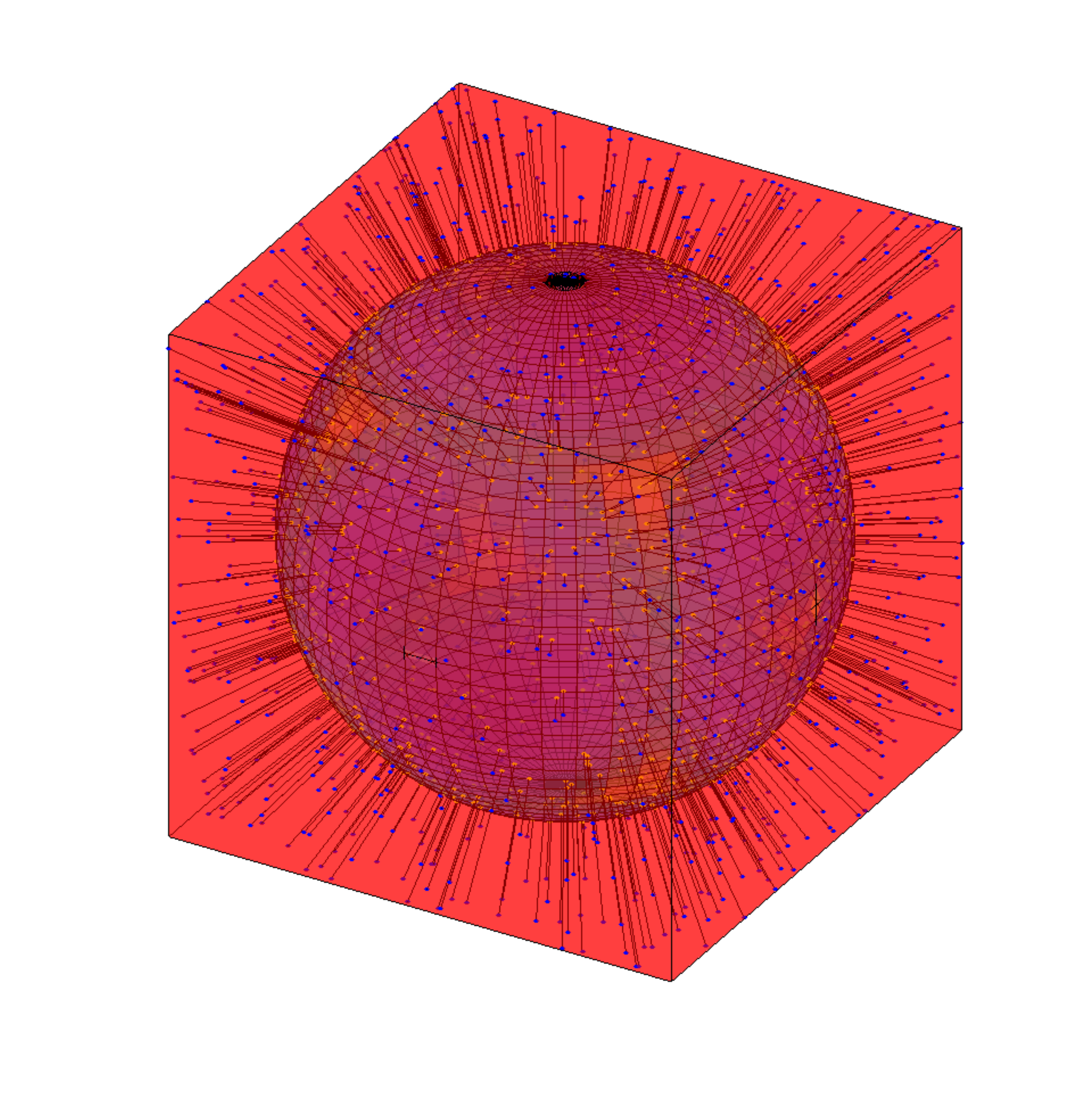}
\end{subfigure}
\begin{subfigure}{0.32\textwidth}
\centering
\includegraphics[trim={5em 5em 5em 4em},clip,width=\textwidth,height=!]{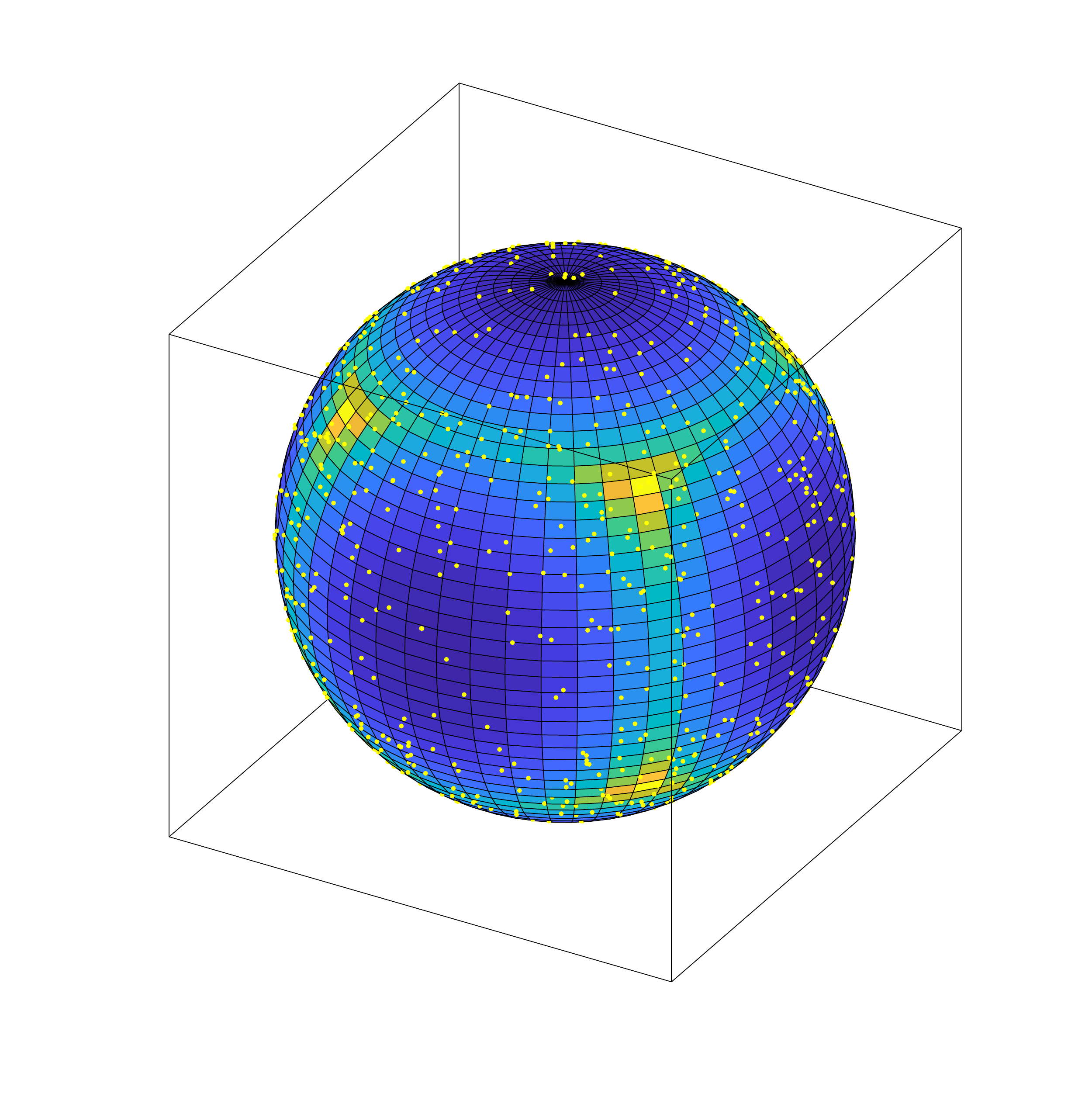}
\end{subfigure}
\caption{Example of simulating and mapping a CSR process on the cube to the sphere. \emph{Left:} example of a CSR process on a cube with $l=1$ and constant intensity function 50. \emph{Middle:} mapping of points from cube to the sphere by the function $f(\mathbf{x})=\mathbf{x}/||\mathbf{x}||$. \emph{Right:} mapped point pattern on the sphere with the new intensity function indicated by the colour on the sphere. High intensity is indicated in yellow whilst low intensity is indicated in blue.}
\label{fig:cube:map}
\hrulefill
\end{figure}

We define a centred cube over each of the six faces with a side length $2l$, where $l=1$. The level-set function for a cube is,
\begin{equation*}
g(\mathbf{x}) =
\begin{cases}
x_3 - l, \quad \text{for } -l \leq x_1,x_2 \leq l\\
x_3 + l, \quad \text{for } -l \leq x_1,x_2 \leq l\\
x_2 - l, \quad \text{for } -l \leq x_1,x_3 \leq l\\
x_2 + l, \quad \text{for } -l \leq x_1,x_3 \leq l\\
x_1 - l, \quad \text{for } -l \leq x_2,x_3 \leq l\\
x_1 + l, \quad \text{for } -l \leq x_2,x_3 \leq l.
\end{cases}
\end{equation*}
Using Theorem \ref{thm:mapping:general:text} we can derive the intensity function for the point process that is mapped to the sphere. By symmetry we need only consider one of the faces of the cube and by rotation we will be able to derive the intensity function on the sphere. Consider the bottom face, i.e. $z=-l$, and in the notation of Theorem \ref{thm:mapping:general:text} label this $\mathbb{D}_1$. Then,
\begin{align*}
l_1(\mathbf{x})&=1\\
J_{(1,f^*)}(\mathbf{x})&=(1+x_1+x_2)^2,
\end{align*}
and so the intensity function over $f(\mathbb{D}_1)$ is,
\begin{equation*}
\rho^*_1(\mathbf{x})=\rho(1+(f_1^{*-1}(x_1))^2+(f_2^{*-1}(x_2))^2)^2(1-x_1^2-x_1^2)^{\frac{1}{2}},
\end{equation*}
thus by the appropriate rotations the intensity function over the entire sphere is,
\begin{equation*}
\rho^*(\mathbf{x}) = 
\begin{cases}
\rho(1+(f_1^{*-1}(x_1))^2+(f_2^{*-1}(x_2))^2)^2(1-x_1^2-x_1^2)^{\frac{1}{2}}, \quad \mathbf{x}\in f(\mathbb{D}_1)\cup f(\mathbb{D}_2)\\
\rho(1+(f_1^{*-1}(x_1))^2+(f_3^{*-1}(x_3))^2)^2(1-x_1^2-x_3^2)^{\frac{1}{2}}, \quad \mathbf{x}\in f(\mathbb{D}_3)\cup f(\mathbb{D}_4)\\
\rho(1+(f_2^{*-1}(x_2))^2+(f_3^{*-1}(x_3))^2)^2(1-x_2^2-x_3^2)^{\frac{1}{2}}, \quad \mathbf{x}\in f(\mathbb{D}_5)\cup f(\mathbb{D}_6),
\end{cases}
\end{equation*}
where $\mathbb{D}_1,\mathbb{D}_2,\mathbb{D}_3,\mathbb{D}_4,\mathbb{D}_5,$ and $\mathbb{D}_6$ are the faces such that $z=-1,z=1,y=-1,y=1,x=-1,$ and $x=1$ respectively. Figure \ref{fig:cube:map} demonstrates mapping from a cube with $l=1$ and $\rho=50$ to the unit sphere where the colour over the sphere indicates areas of low (blue) and high (yellow) intensity. The figure also shows an example of a CSR pattern over the cube and how this pattern changes under the mapping. 

\begin{figure}[p]
\begin{subfigure}{0.5\textwidth}
\centering
\includegraphics[trim={0em 10em 0 10em},clip,width=0.95\textwidth,height=!]{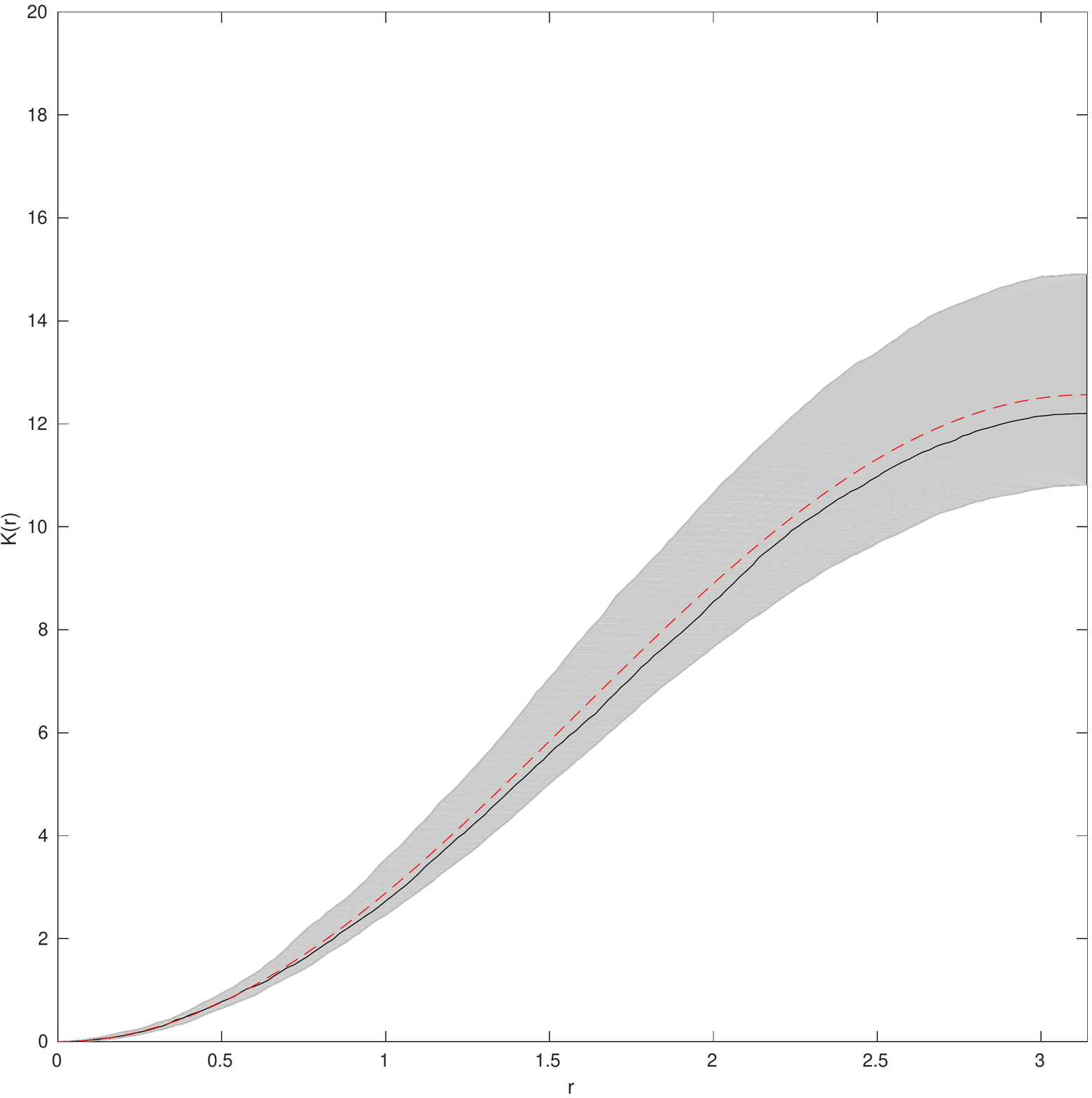}
\end{subfigure}
\begin{subfigure}{0.5\textwidth}
\centering
\includegraphics[trim={0 10em 0 10em},clip,width=0.95\textwidth,height=!]{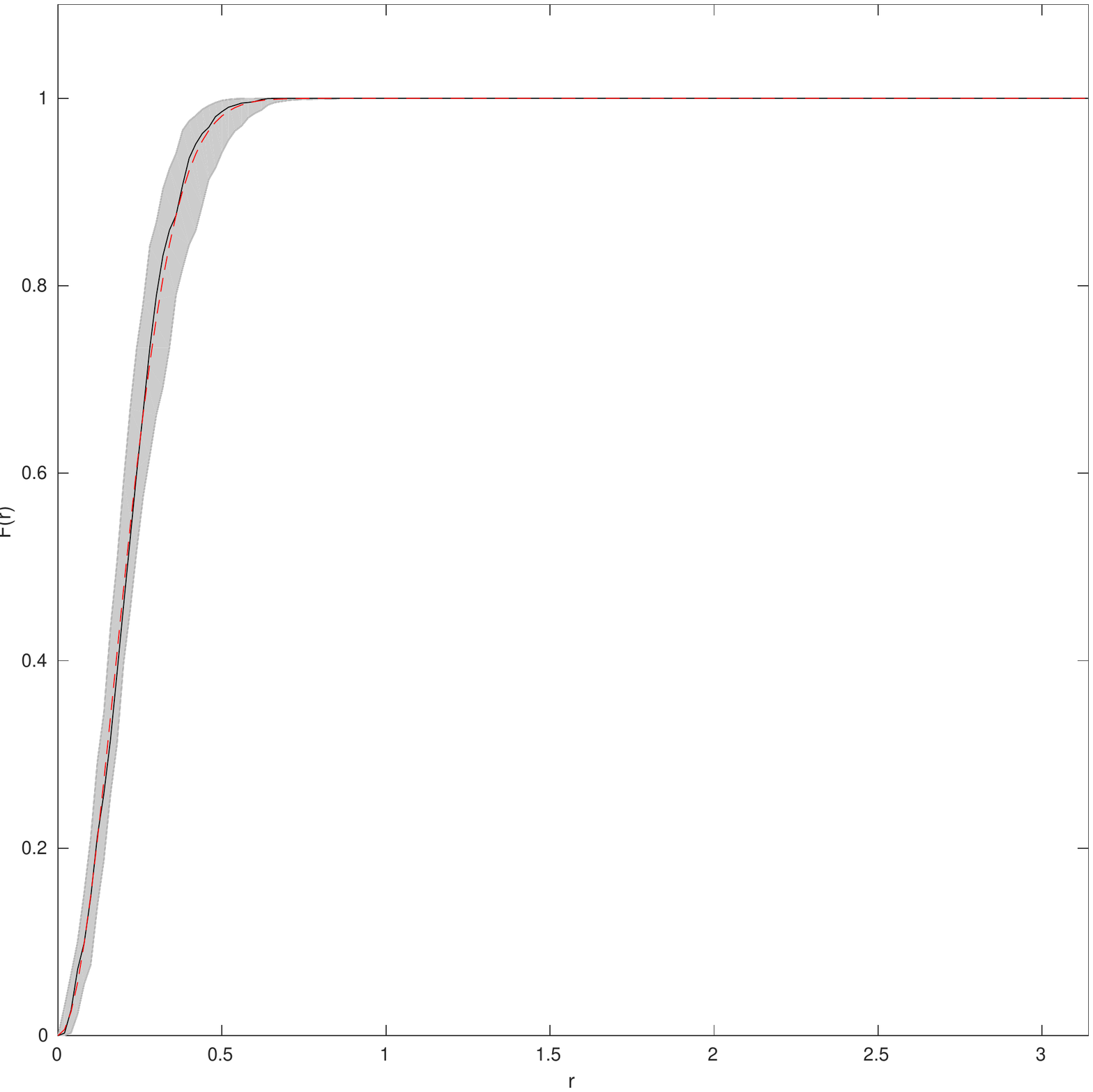}
\end{subfigure}
\begin{subfigure}{0.5\textwidth}
\centering
\includegraphics[trim={0 10em 0 10em},clip,width=0.95\textwidth,height=!]{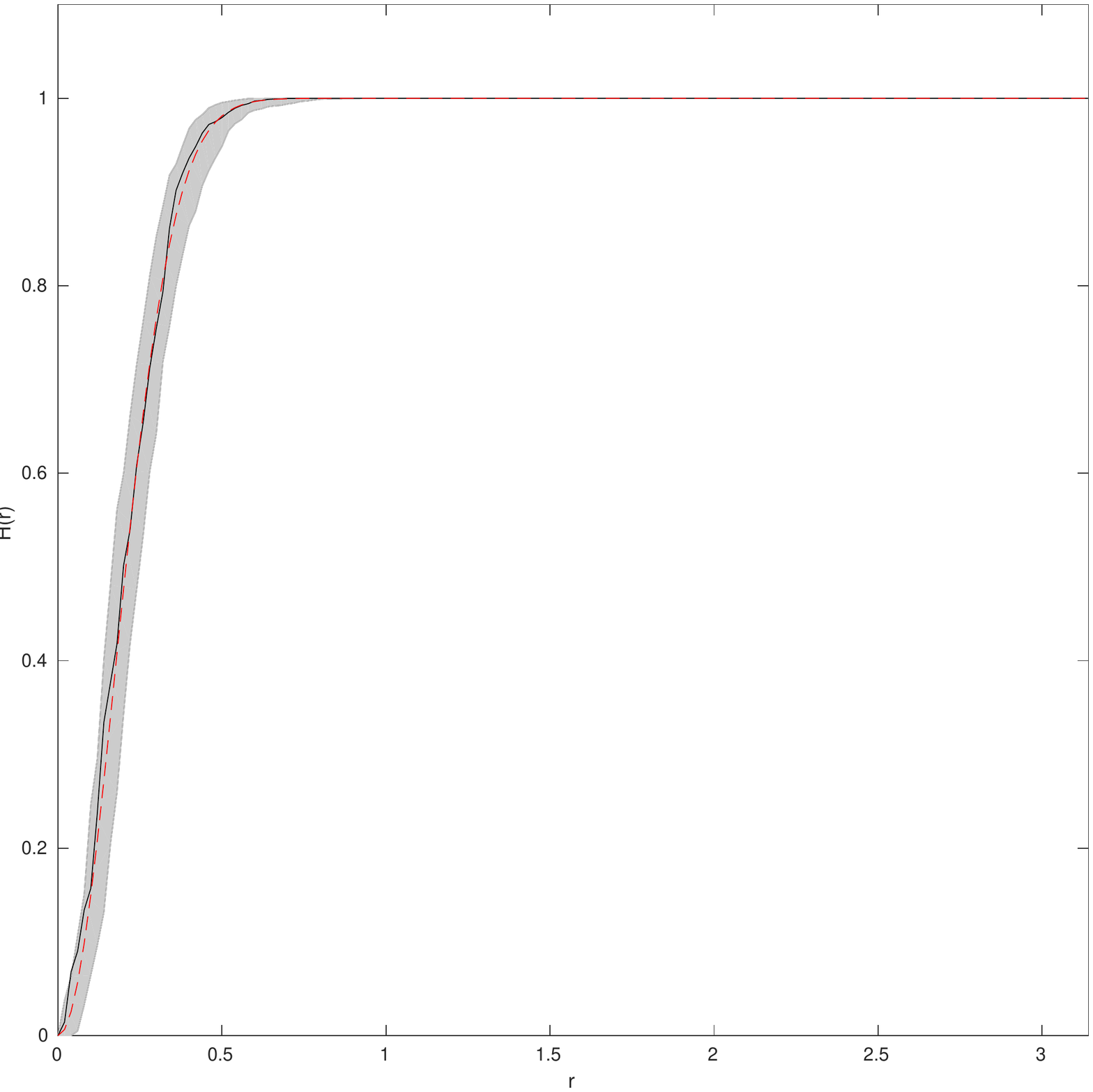}
\end{subfigure}
\begin{subfigure}{0.5\textwidth}
\centering
\includegraphics[trim={0 10em 0 10em},clip,width=0.95\textwidth,height=!]{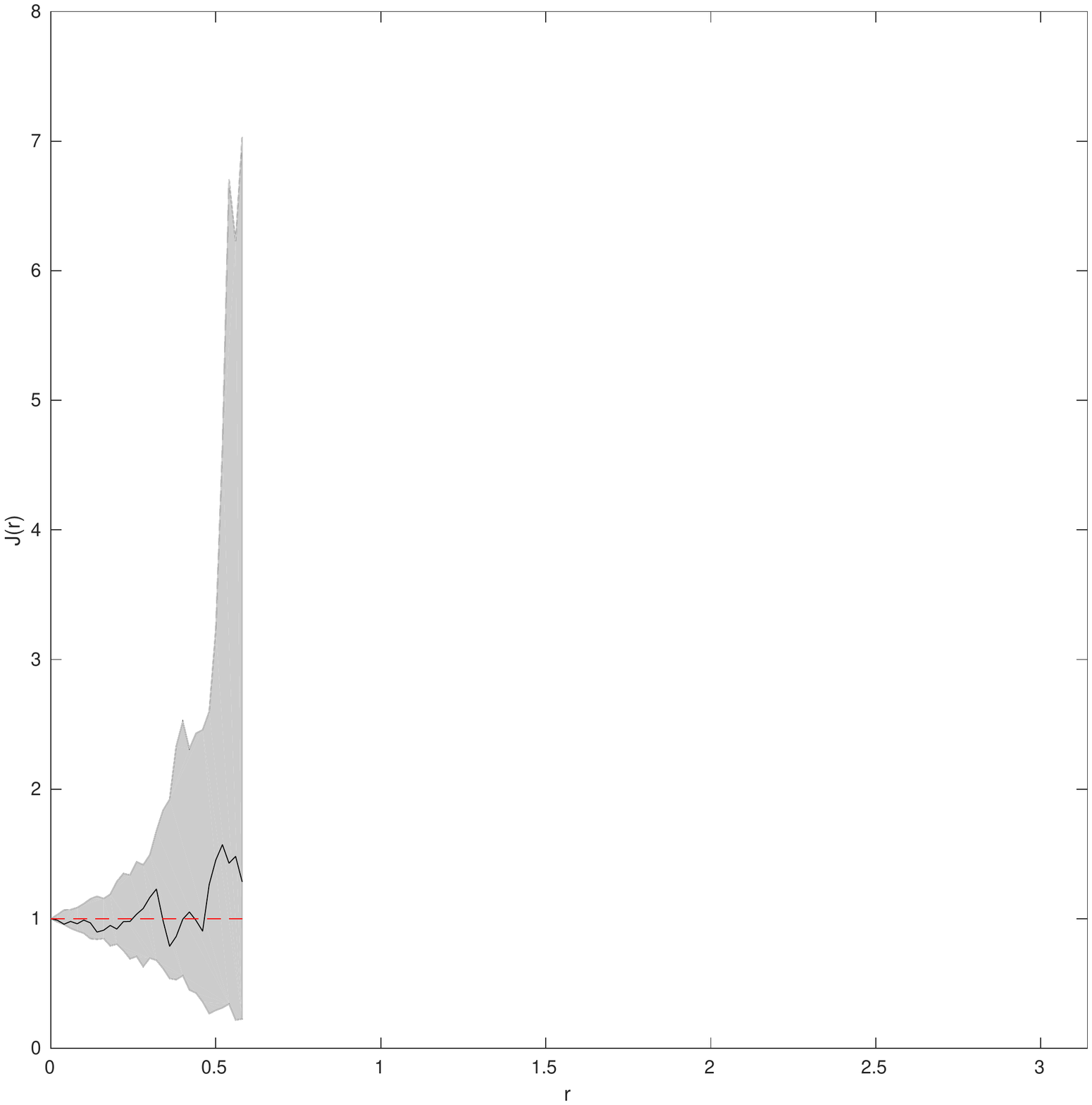}
\end{subfigure}
\caption{Examples of $K_{\text{inhom}}$- (\emph{top left}), $F_{\text{inhom}}$- (\emph{top right}), $H_{\text{inhom}}$- (\emph{bottom left}), and $J_{\text{inhom}}$- (\emph{bottom right}) functions for CSR patterns on a cube with $l=1$ and $\rho=5$. Black line is the estimated functional summary statistic for our observed data, dashed red line is the theoretical functional summary statistic for a Poisson process, and the grey shaded area represents the simulation envelope from 99 Monte Carlo simulations of Poisson processes fitted to the observed data.}
\label{fig:cube:functions}
\hrulefill
\end{figure}

In order to be able to construct the inhomogeneous $F$-, and $H$-function we need to determine $\inf_{\mathbf{x}\in\mathbb{S}^2}\rho^*(\mathbf{x})$. By the nature of the function $f(\mathbf{x})=\mathbf{x}/||\mathbf{x}||$ and assuming that $l\geq 1$, then mapping events from the cube to the sphere causes events to be more concentrated on the sphere compared to the cube, thus increasing the corresponding intensity on the sphere. Therefore, the lowest achievable intensity occurs at the centre of each face of the cube, i.e. for the bottom face it occurs when $x_1=x_2=0$, giving $\inf_{\mathbf{x}\in\mathbb{S}^2}\rho^*(\mathbf{x})=\rho$. Figure \ref{fig:cube:functions} gives examples of the inhomogeneous $K$-, $F$-, $H$-, and $J$-functions where $l=1$ and $\rho=5$ and are typical when the observed process is CSR.

\subsection{Ellipsoid}

\begin{figure}[t]
\begin{subfigure}{0.32\textwidth}
\centering
\includegraphics[trim={8em 20em 6em 18em},clip,width=\textwidth,height=!]{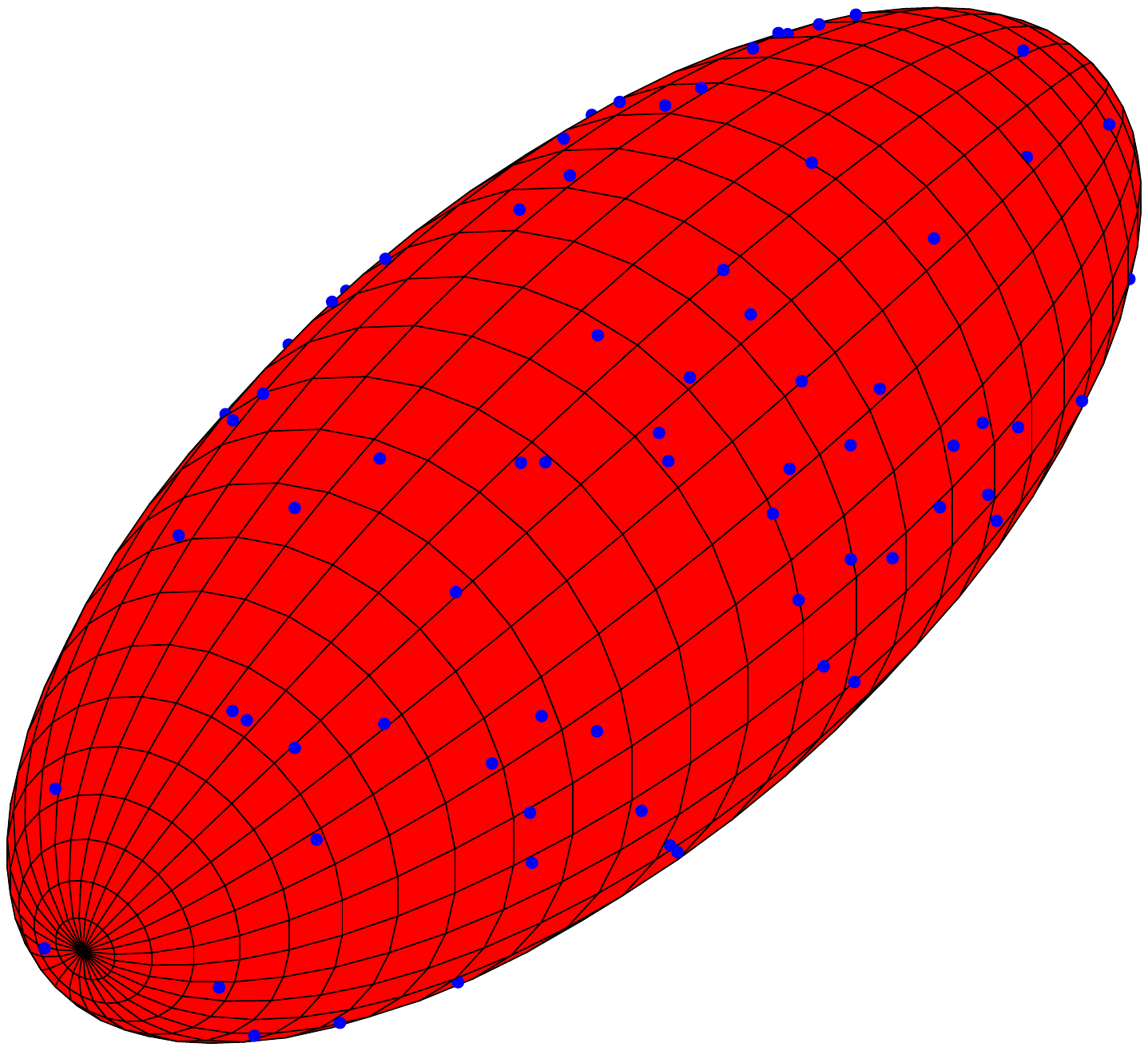}
\end{subfigure}
\begin{subfigure}{0.32\textwidth}
\centering
\includegraphics[trim={8em 20em 6em 18em},clip,width=\textwidth,height=!]{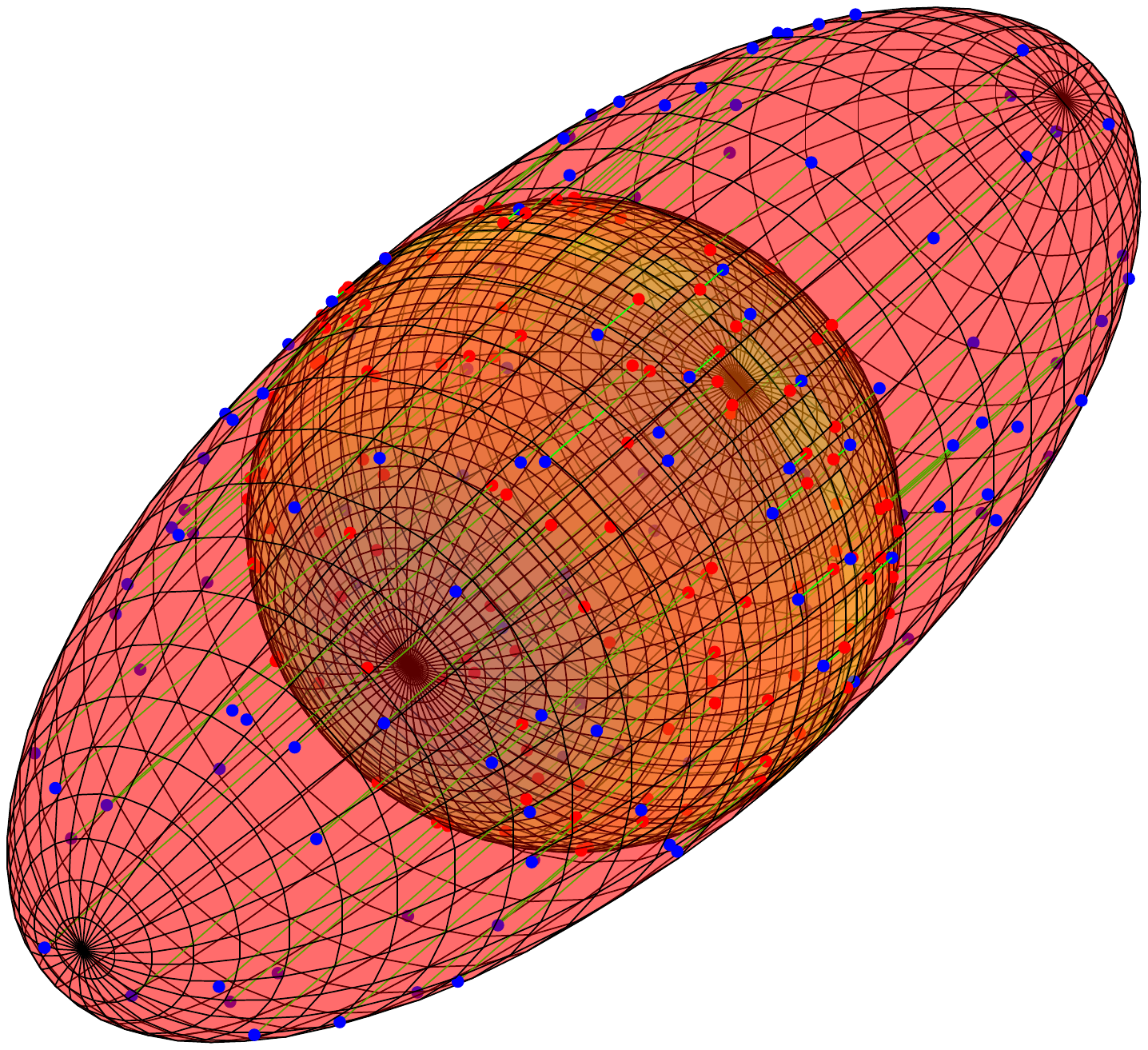}
\end{subfigure}
\begin{subfigure}{0.32\textwidth}
\centering
\includegraphics[trim={8em 20em 6em 18em},clip,width=\textwidth,height=!]{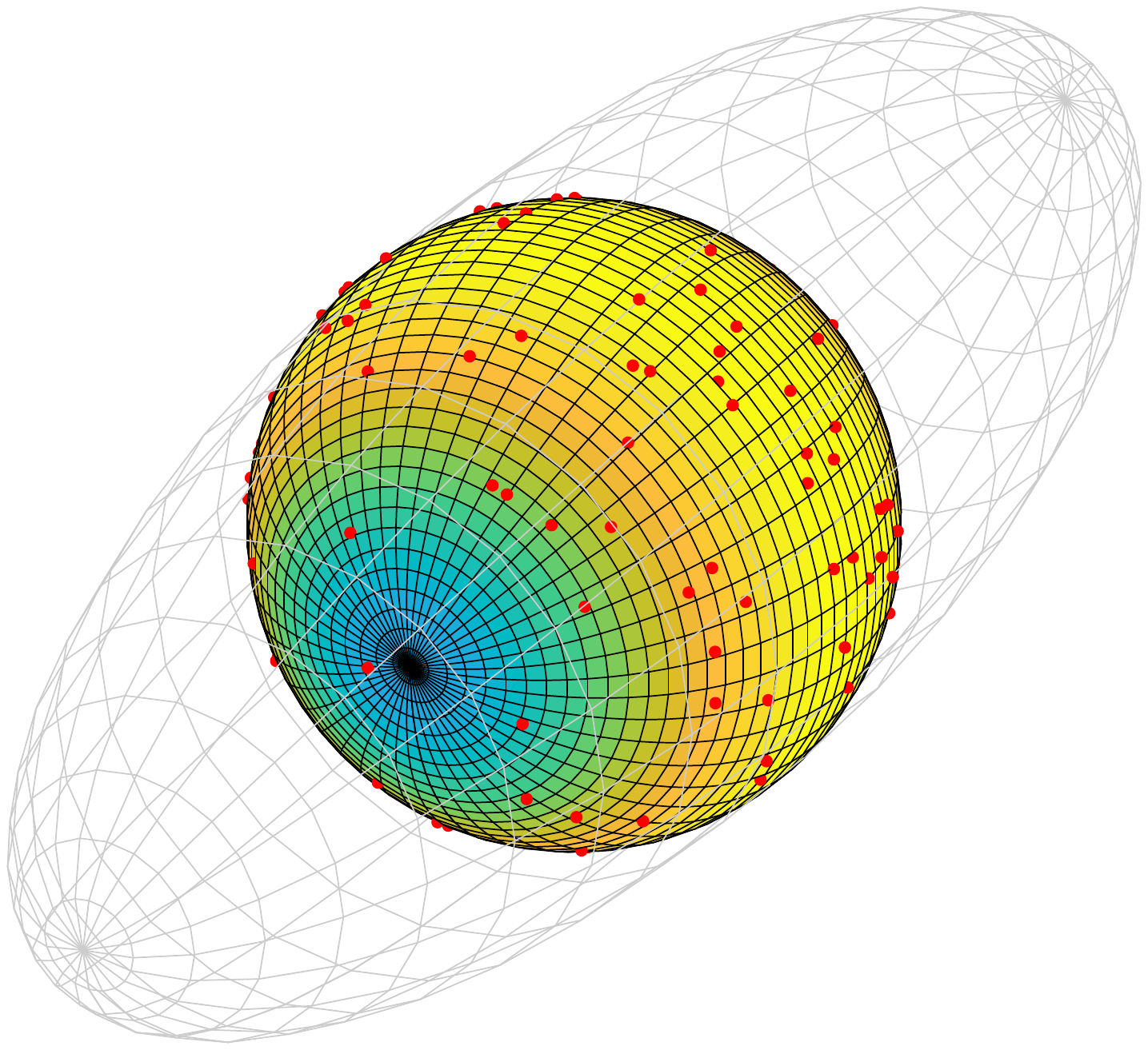}
\end{subfigure}
\caption{Example of simulating and mapping a CSR process on a prolate ellipsoid to the sphere. \emph{Left:} example of a CSR process on a prolate ellipsoid with $a=b=1, c=3$ and $\rho=5$. \emph{Middle:} mapping of points from prolate ellipsoid to the sphere by the function $f(\mathbf{x})=(x_1/a,x_2/b,x_3/c)^T$. \emph{Right:} mapped point pattern on the sphere with the new intensity function indicated by the colour on the sphere. High intensity is indicated in yellow whilst low intensity is indicated in blue.}
\label{fig:ellipsoid:map}
\hrulefill
\end{figure}

An ellipsoid is defined by its semi-major axis lengths $a,b,c\in\mathbb{R}$ along $x$-, $y$-, and $z$-axis respectively. Again we also assume that the ellipsoid is centred at the origin. The level-set function for ellipsoids is given by $g(\mathbf{x})=x_1^2/a^2+x_2^2/b^2+x_3^2/c^2-1$, in this current form $\tilde{g}$ is not well defined in which case we shall use the following equivalent representation
\begin{equation*}
g(\mathbf{x}) =
\begin{cases}
x_1^2/a^2+x_2^2/b^2+x_3^2/c^2-1, \quad \text{for } x_3 \geq 0\\
x_1^2/a^2+x_2^2/b^2+x_3^2/c^2-1, \quad \text{for } x_3 < 0.
\end{cases}
\end{equation*}
This representation then allows for $\tilde{g}$ to be well defined for each partition of the ellipsoid.

\begin{figure}
\begin{subfigure}{0.5\textwidth}
\centering
\includegraphics[trim={0 10em 0 10em},clip,width=0.95\textwidth,height=!]{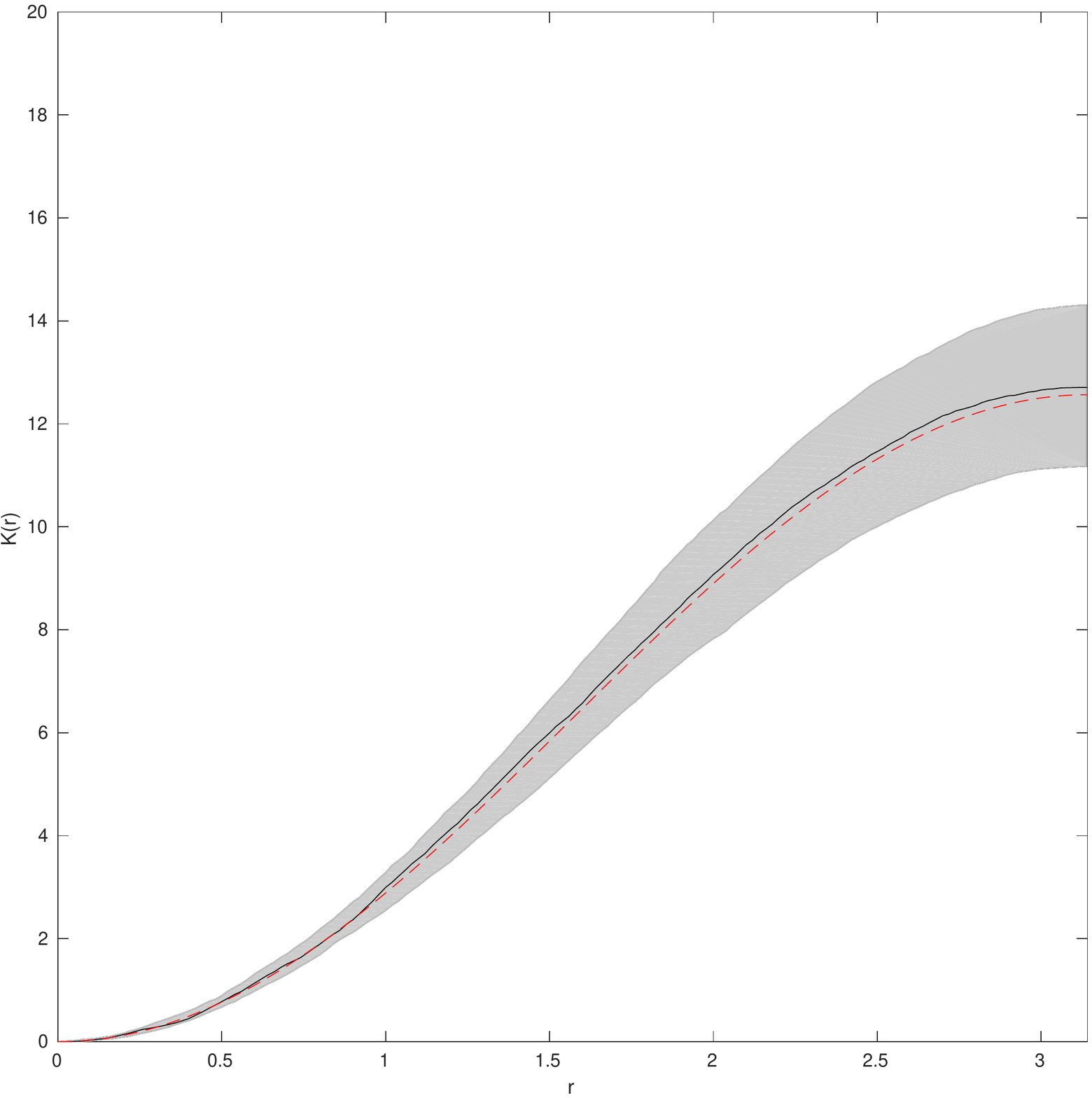}
\end{subfigure}
\begin{subfigure}{0.5\textwidth}
\centering
\includegraphics[trim={0 10em 0 10em},clip,width=0.95\textwidth,height=!]{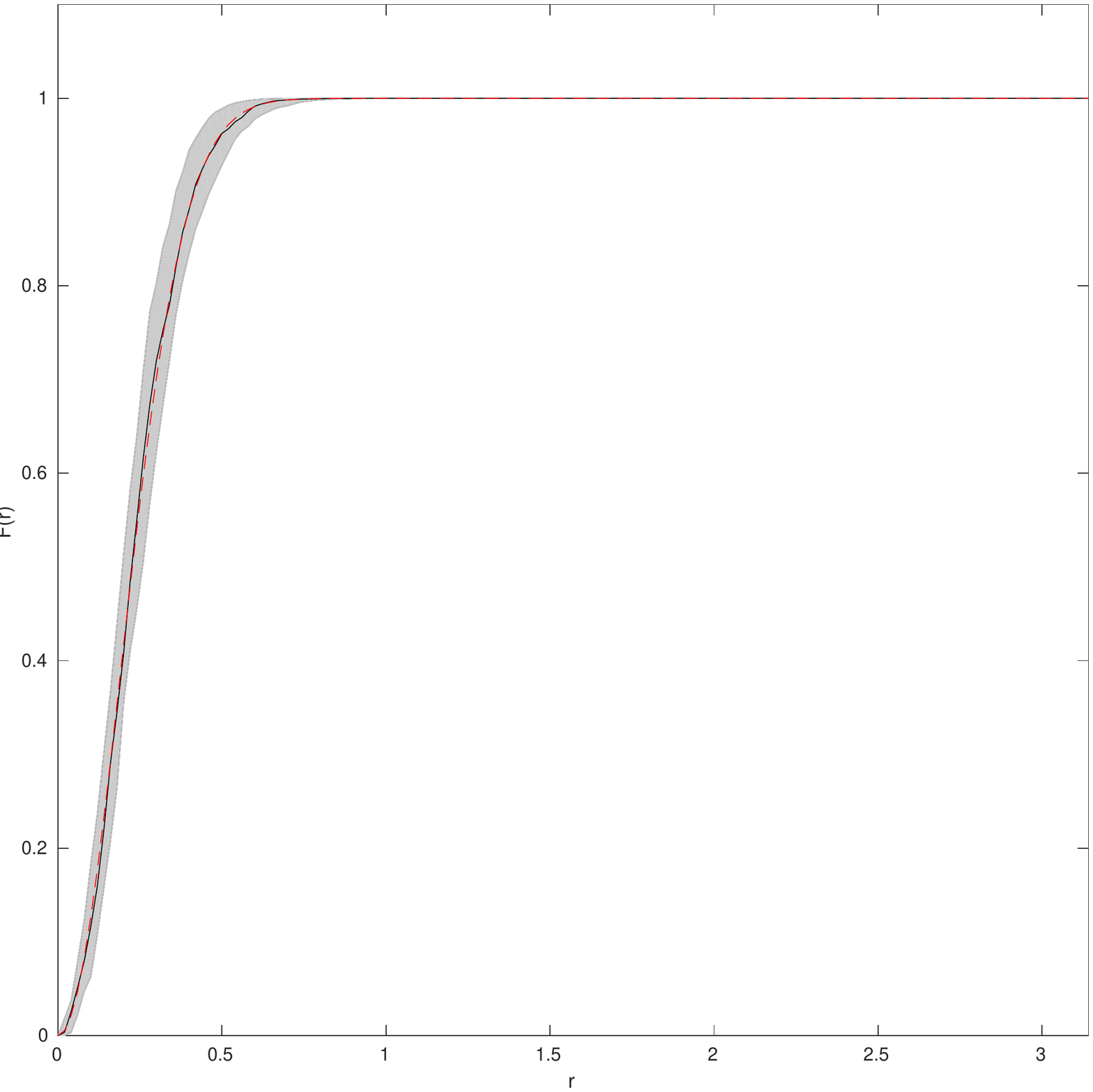}
\end{subfigure}
\begin{subfigure}{0.5\textwidth}
\centering
\includegraphics[trim={0 10em 0 10em},clip,width=0.95\textwidth,height=!]{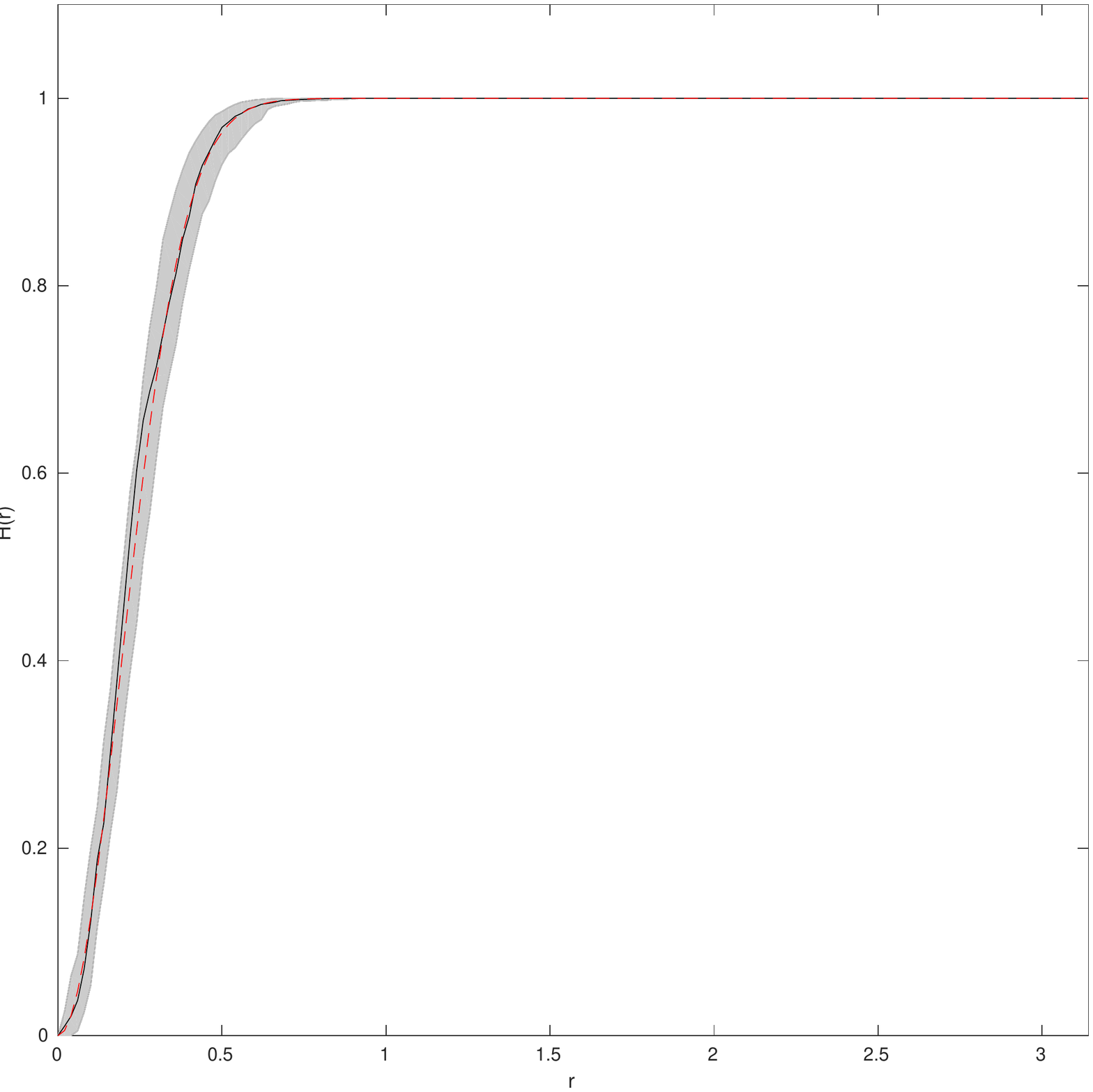}
\end{subfigure}
\begin{subfigure}{0.5\textwidth}
\centering
\includegraphics[trim={0 10em 0 10em},clip,width=0.95\textwidth,height=!]{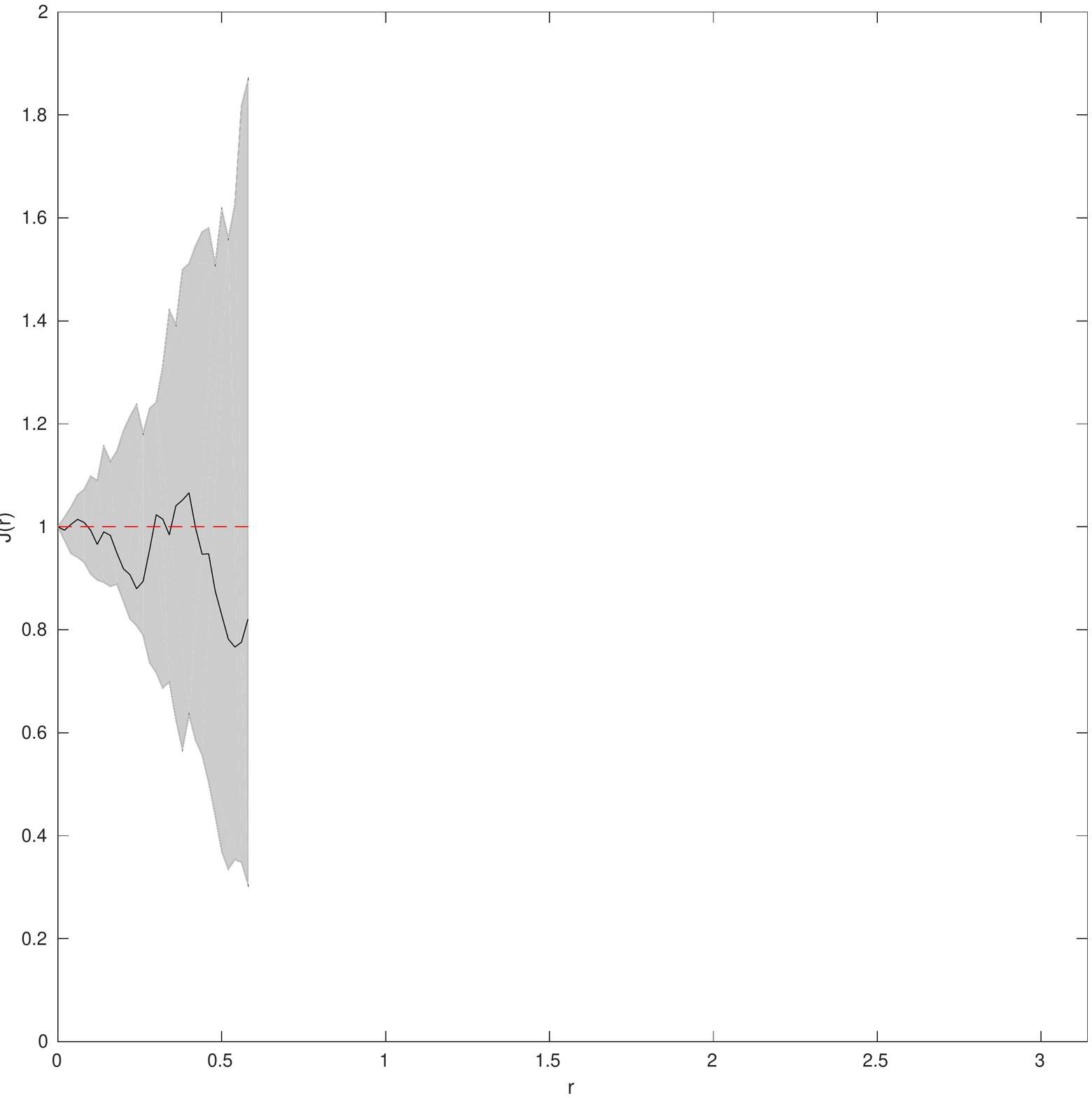}
\end{subfigure}
\caption{Examples of $K_{\text{inhom}}$- (\emph{top left}), $F_{\text{inhom}}$- (\emph{top right}), $H_{\text{inhom}}$- (\emph{bottom left}), and $J_{\text{inhom}}$- (\emph{bottom right}) functions for CSR patterns on a prolate ellipsoid with $a=b=1, c=3,$ and $\rho=5$. Black line is the estimated functional summary statistic for our observed data, dashed red line is the theoretical functional summary statistic for a Poisson process, and the grey shaded area represents the simulation envelope from 99 Monte Carlo simulations of Poisson processes fitted to the observed data.}
\label{fig:ellipsoid:functions}
\hrulefill
\end{figure}

We now demonstrate our methodology on an ellipsoid with semi-major axis lengths $a=1, b=1,$ and $c=3$ along the $x$-, $y$-, and $z$-axis respectively. Instead of using the function $f(\mathbf{x})=\mathbf{x}/||\mathbf{x}||$ to map from the ellipsoid to the sphere, we can use a simpler mapping function which makes calculation of the determinant $J_{(i,f^*)}(\mathbf{x})$ in Theorem \ref{thm:mapping:general:text} significantly easier. We can simply scale along the axis directions, i.e. use the mapping $f(\mathbf{x})=(x_1/a,x_2/b,x_3/c)^T.$ Using this mapping function, as opposed to dividing each vector by the norm of itself, and focusing on the bottom hemiellipsoid (indicated by the minus superscript), then
\begin{equation*}
l_{-}(\mathbf{x})=\sqrt{\frac{1-\left(1-\frac{c^2}{a^2}\right)x_1^2-\left(1-\frac{c^2}{b^2}\right)x_2^2}{1-x_1^2-x_2^2}},\quad J_{(-,f^*)}(\mathbf{x})=ab, 
\end{equation*}
and so on the lower hemisphere the intensity function takes the form
\begin{equation*}
\rho_{-}^*(\mathbf{x})=\rho ab\sqrt{1-\left(1-\frac{c^2}{a^2}\right)x_1^2-\left(1-\frac{c^2}{b^2}\right)x_2^2}.
\end{equation*}
By symmetry the mapped intensity function over the whole sphere is then
\begin{equation*}
\rho^*(\mathbf{x})=\rho ab\sqrt{1-\left(1-\frac{c^2}{a^2}\right)x_1^2-\left(1-\frac{c^2}{b^2}\right)x_2^2}.
\end{equation*}
Again we need to calculate $\inf_{\mathbf{x}\in \mathbb{S}^2}\rho^*(\mathbf{x})$. Noting that $c\geq a=b$, thus $-\left(1-c^2/a^2\right)\geq 0$ and $-\left(1-c^2/b^2\right)\geq 0$, then the square root term is minimised when $x_1$ and $x_2$ are $0$, hence $\inf_{\mathbf{x}\in\mathbb{S}^2}\rho^*(\mathbf{x})=\rho ab$. Using this we can construct the estimators of the inhomogeneous functional summary statistics given by Equations (\ref{F:estimator})-(\ref{K:estimator}). Examples are given in Figure \ref{fig:ellipsoid:functions}. These figures are typical for CSR with the estimated functional summary statistics lying well within the simulation envelopes.

\section{Regular \& cluster processes on $\mathbb{D}$}\label{Section:8}

 We examine some regular and cluster processes on $\mathbb{D}$. In particular, we examine how functional summary statistics constructed under the Poisson hypothesis deviate when the underlying process is in fact not Poisson. We shall be using the Mat\'{e}rn \rom{1} and \rom{2} inhibition processes \cite{chiu2013stochastic} as examples of regular processes, and Thomas processes as a cluster example. Definitions for the Mat\'{e}rn \rom{1}, \rom{2} and Thomas processes on convex shapes will also be presented, whilst properties of such processes are given in Appendix F.

\subsection{Examples of regular and cluster processes on convex shapes}
A common way of defining a regular process is using a minimum distance $R$, known as the hardcore distance, for which no point in the process has a nearest neighbour closer than $R$. In typical applications $R$ is usually the Euclidean distance (in $\mathbb{R}^d$) or the great circle distance (in $\mathbb{S}^2$), but on an arbitrary three dimensional convex shape, $\mathbb{D}$, this distance is taken as the geodesic distance defined by the surface. The following definitions extend the Mat\'{e}rn \rom{1} and \rom{2} processes to a convex shape with geodesic distance $d(\mathbf{x},\mathbf{y})$, $\mathbf{x},\mathbf{y}\in\mathbb{D}$.

\begin{dfn}
Let $X$ be a homogeneous Poisson process on $\mathbb{D}$ with constant intensity function $\rho\in\mathbb{R}_+$. Fix $R\in [0,\pi]$, and thin $X$ according to the following rule: delete events $\mathbf{x}\in X$ if there exists $\mathbf{y}\in X\setminus\{\mathbf{x}\}$ such that $d(\mathbf{x},\mathbf{y})<R$, otherwise retain $\mathbf{x}$. The resulting thinned process is then defined as a Mat\'ern \rom{1} inhibition process on $\mathbb{D}$.
\end{dfn}

\begin{dfn}
Let $X$ be a homogeneous Poisson process on $\mathbb{D}$ with constant intensity function $\rho\in\mathbb{R}_+$. Fix $R\in [0,\pi]$, and let each $\mathbf{x}\in X$ have an associated mark, $M_\mathbf{x}$ drawn from some mark density $P_M$ independently of all other marks and points in $X$.  Thin $X$ according to the following rule: delete the event $\mathbf{x}\in X$ if there exists $\mathbf{y}\in X\setminus\{\mathbf{x}\}$ such that $d(\mathbf{x},\mathbf{y})<R$ and $M_\mathbf{y}<M_\mathbf{x}$, otherwise retain $\mathbf{x}$. The resulting thinned process is then defined as a Mat\'ern \rom{2} inhibition process on $\mathbb{D}$.
\end{dfn}

We also extend the Neyman-Scott process, a class of cluster processes, to arbitrary convex shapes.

\begin{dfn}
Let $X_P$ be a homogeneous Poisson process on $\mathbb{D}$ with constant intensity function $\rho\in\mathbb{R}_+$. Then for each $\mathbf{c}\in X_P$ define $X_{\mathbf{c}}$ to the point process with intensity function $\rho_\mathbf{c}(\mathbf{x})=\alpha k(\mathbf{x},\mathbf{c}),$ where $\alpha>0$ and $k:\mathbb{D}\times\mathbb{D}\mapsto\mathbb{R}$ is a density function and $N_{X_{\mathbf{c}}}(\mathbb{D})$ can be any random counting measure associated to $X_{\mathbf{c}}$. The point process $X=\cup_{\mathbf{c}\in X_p}X_{\mathbf{c}}$ is a Neyman-Scott process.
\end{dfn}

A Thomas process is a specific Neyman-Scott process where the density function $k(\cdot,\cdot)$ has a specific form. In $\mathbb{R}^2$, $k$ is taken to be an isotropic bivariate Gaussian distribution \cite{Moller2004}, whilst on $\mathbb{S}^2$ it is taken as the Von-Mises Fisher distribution \cite{Lawrence2016}. We define a Thomas process on $\mathbb{D}$ to be a Neyman-Scott process with density function $k$ of the form,
\begin{equation*}
k(\mathbf{x},\mathbf{y})=\frac{1}{\chi(\sigma^2)}\exp\left(-\frac{d^2(\mathbf{x},\mathbf{y})}{2\sigma^2}\right),
\end{equation*}
where $\sigma$ is a bandwidth parameter and $\chi(\sigma^2)=\int_{\mathbb{D}}\exp\left(-d(\mathbf{x},\mathbf{y})/2\sigma^2\right)\lambda_{\mathbb{D}}(\mathbf{y})$. This is known as the Riemannian Gaussian distribution \cite{Said2017}, where on the plane this would reduce to an isotropic bivariate Gaussian and on a sphere to the Von-Mises Fisher distribution.

\subsection{Functional summary statistics assuming a homogeneous Poisson process}

\begin{sidewaysfigure}
  \begin{subfigure}{0.24\hsize}\centering
      \includegraphics[trim={2cm 2cm 2cm 2cm},clip=TRUE,width=1\hsize]{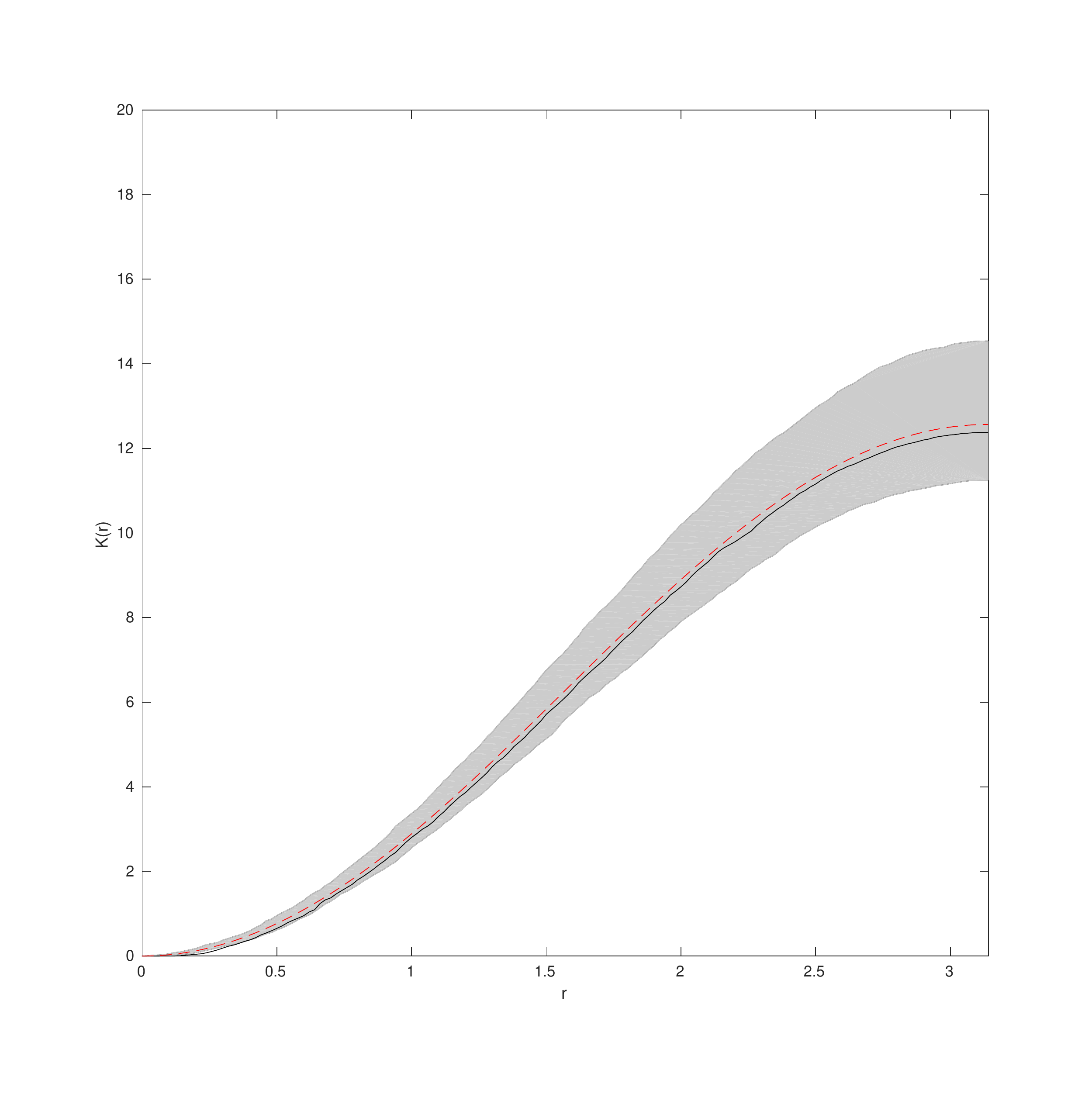}
  \end{subfigure}
  \begin{subfigure}{0.24\hsize}\centering
      \includegraphics[trim={2cm 2cm 2cm 2cm},clip=TRUE,width=1\hsize]{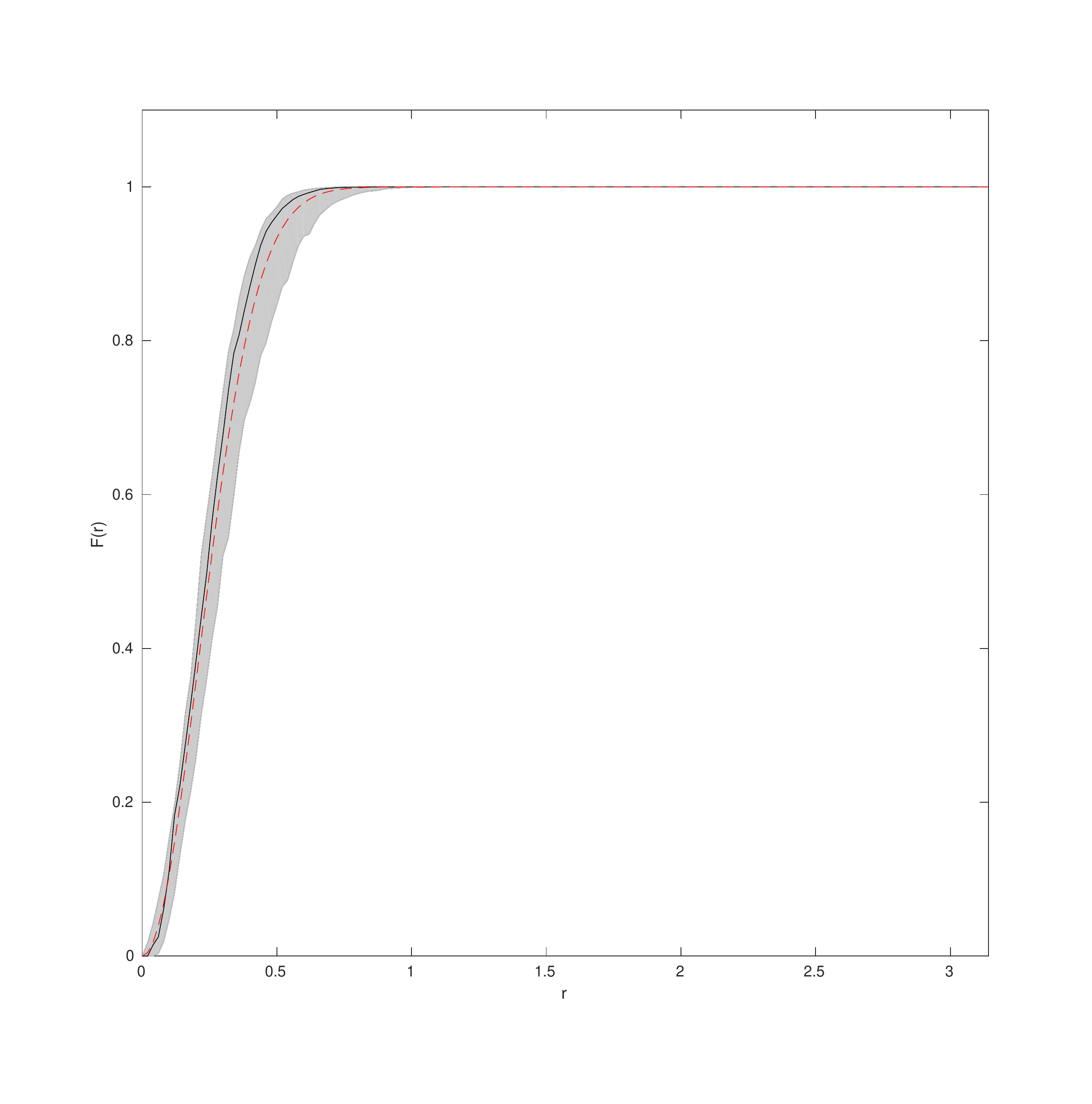}
  \end{subfigure}
  \begin{subfigure}{0.24\hsize}\centering
      \includegraphics[trim={2cm 2cm 2cm 2cm},clip=TRUE,width=1\hsize]{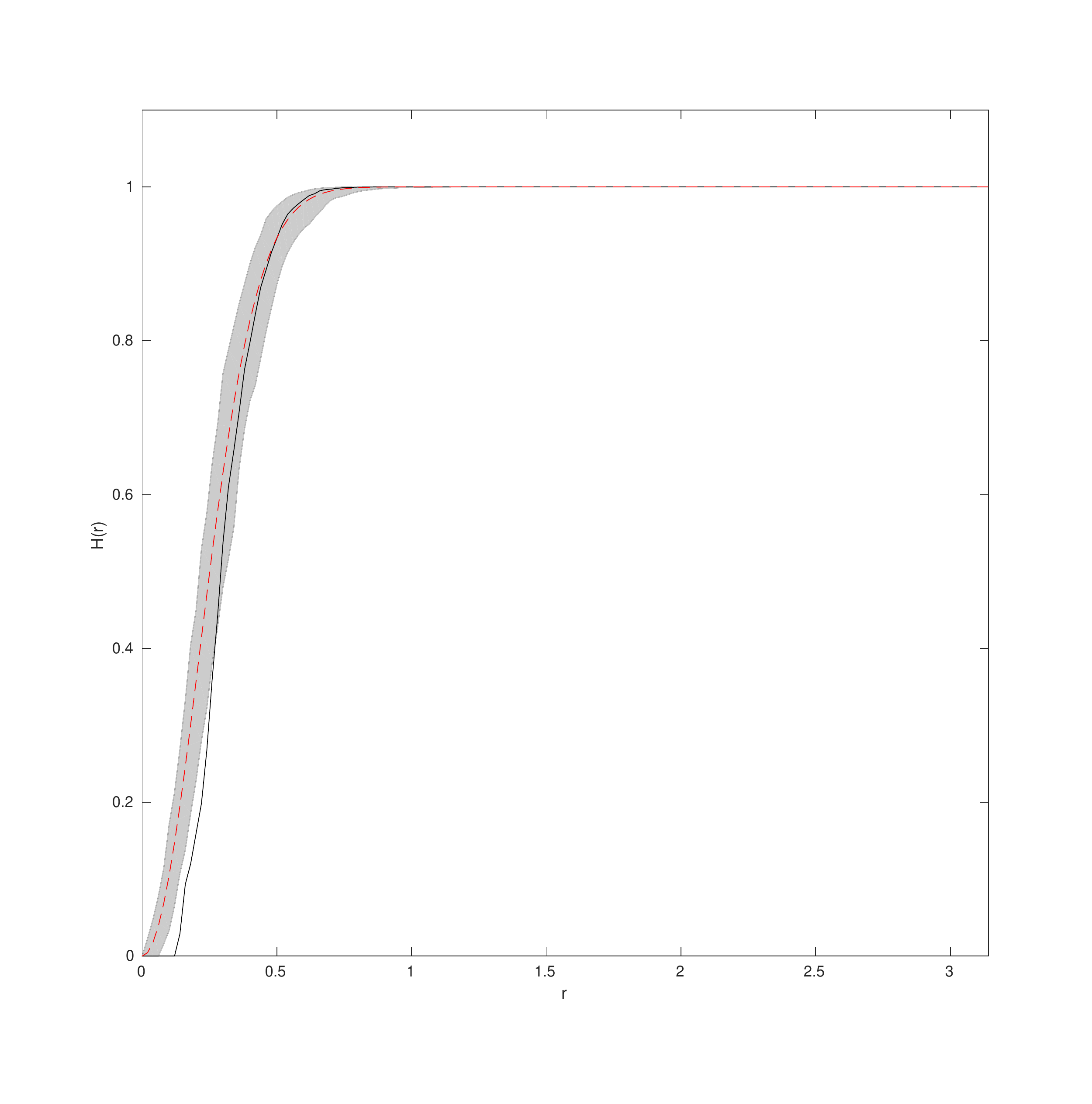}
  \end{subfigure}
  \begin{subfigure}{0.24\hsize}\centering
      \includegraphics[trim={2cm 2cm 2cm 2cm},clip=TRUE,width=1\hsize]{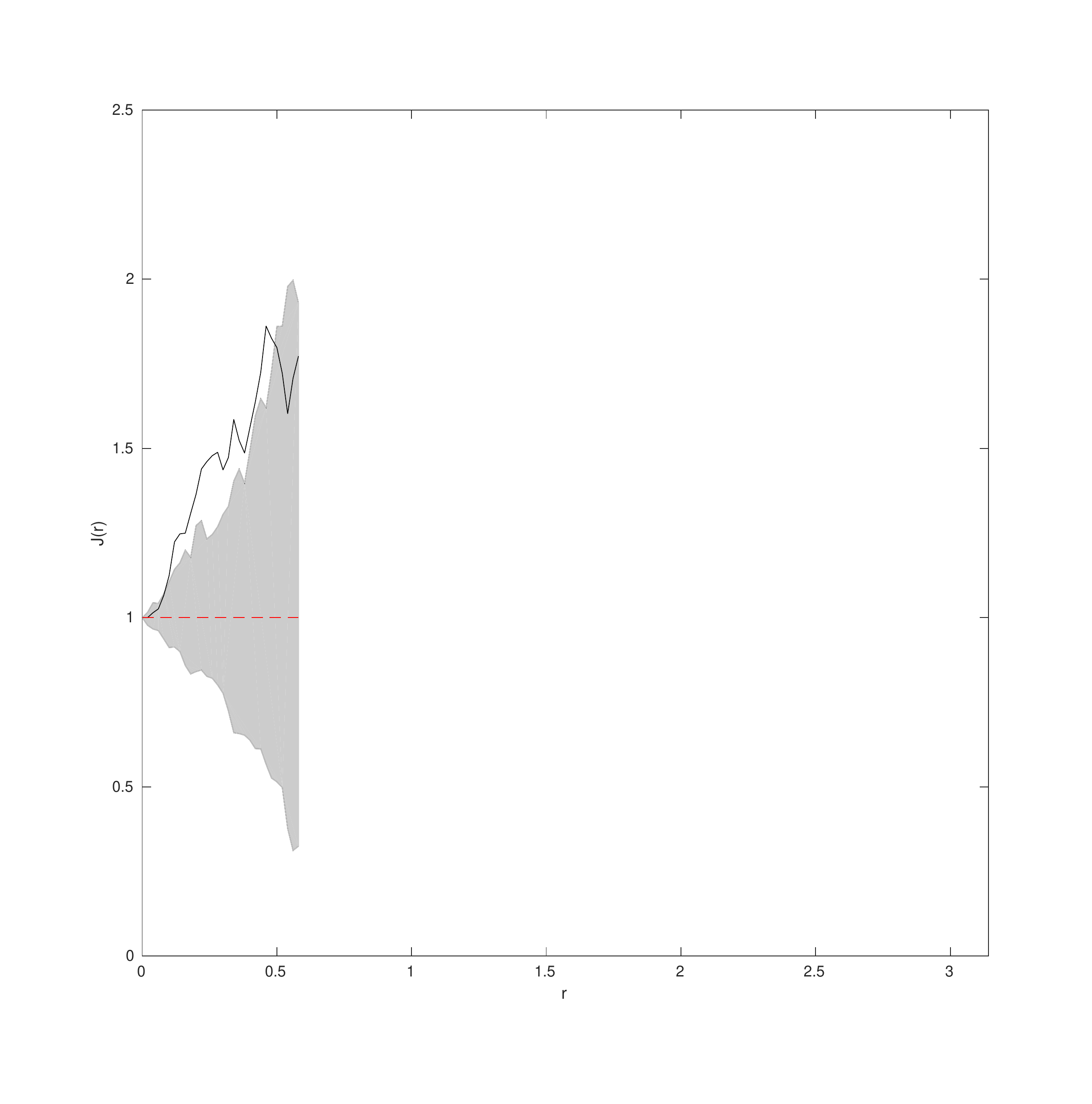}
  \end{subfigure}

  \begin{subfigure}{0.24\hsize}\centering
      \includegraphics[trim={2cm 2cm 2cm 2cm},clip=TRUE,width=1\hsize]{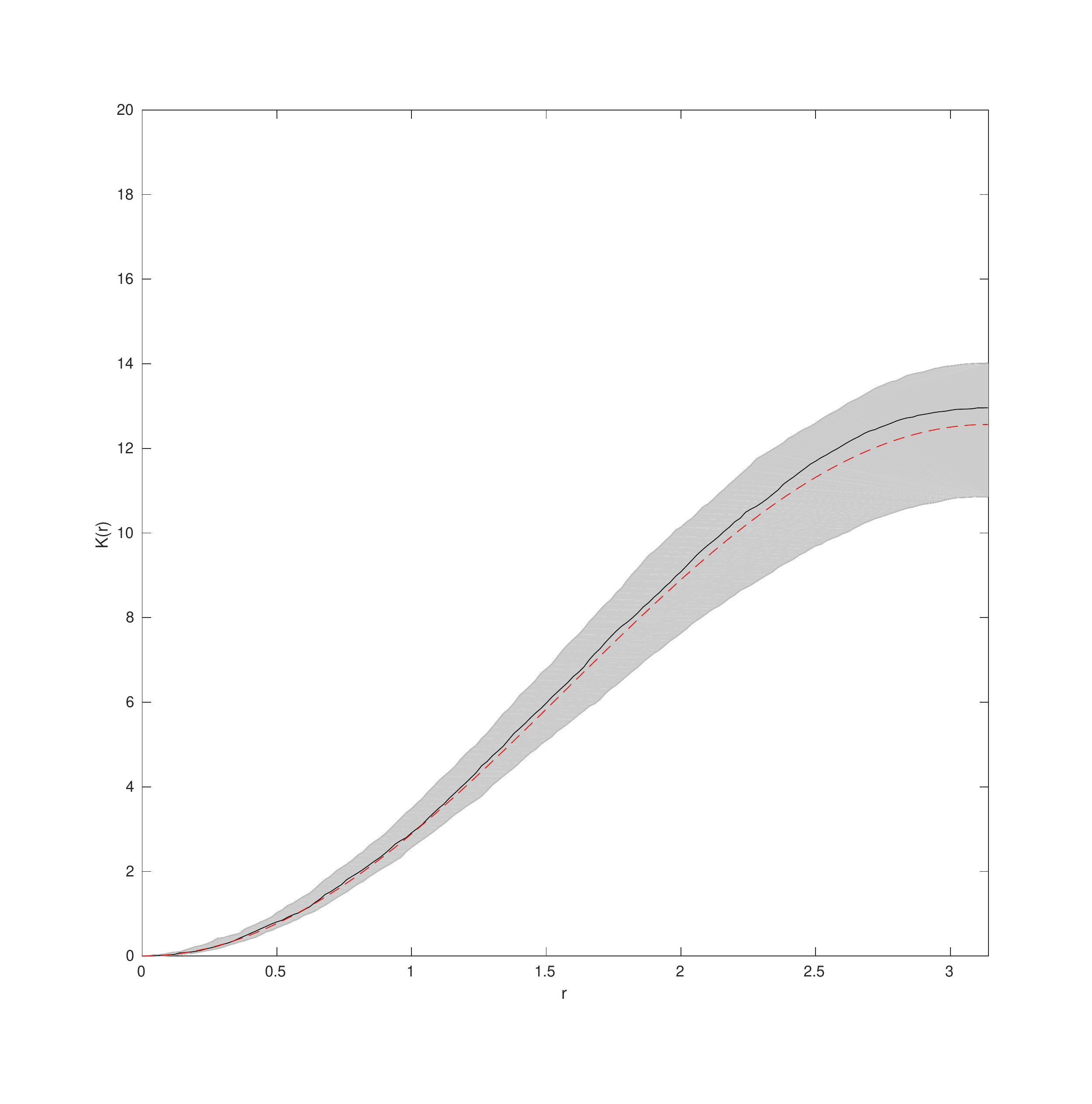}
  \end{subfigure}
  \begin{subfigure}{0.24\hsize}\centering
      \includegraphics[trim={2cm 2cm 2cm 2cm},clip=TRUE,width=1\hsize]{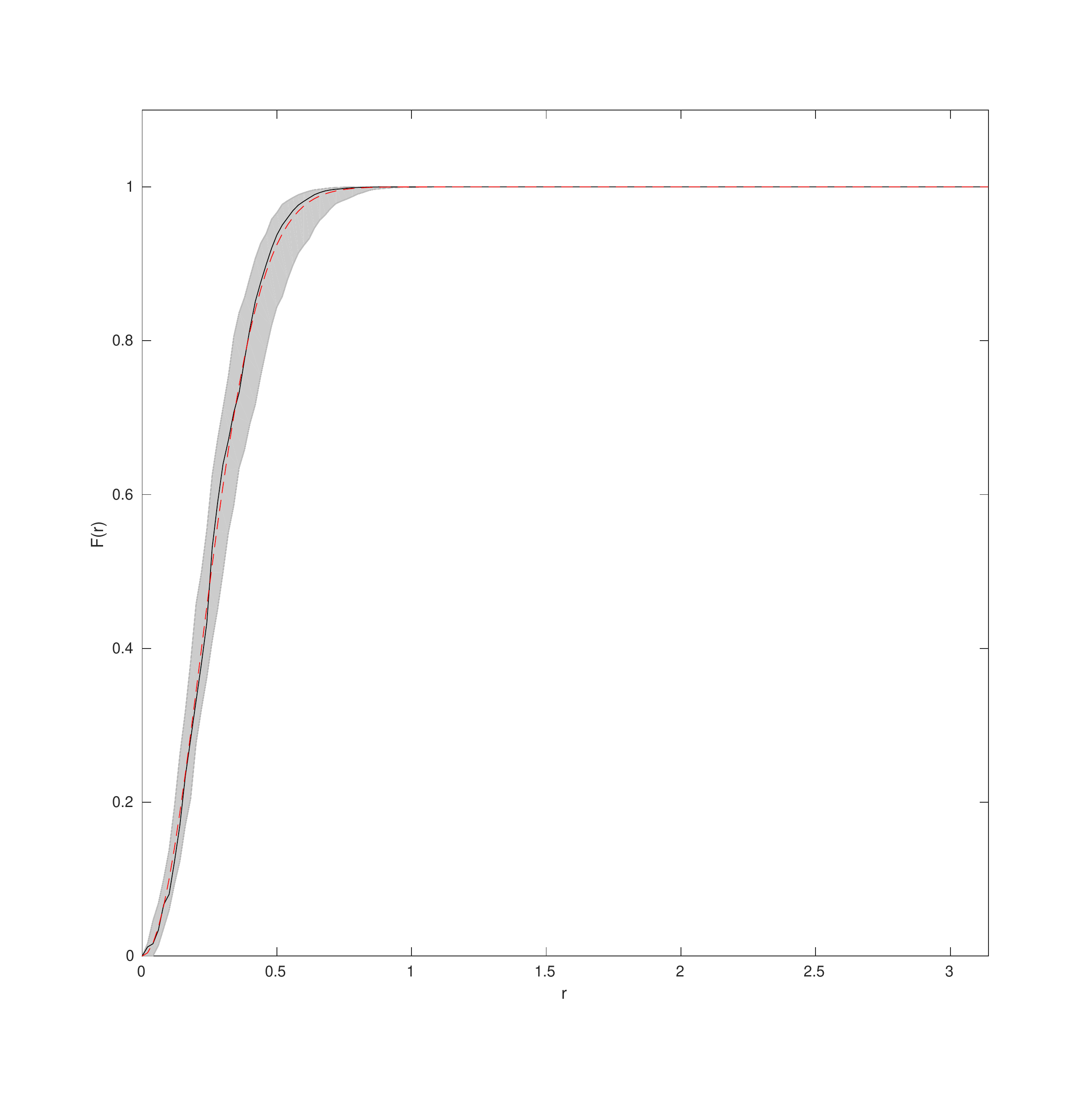}
  \end{subfigure}
  \begin{subfigure}{0.24\hsize}\centering
      \includegraphics[trim={2cm 2cm 2cm 2cm},clip=TRUE,width=1\hsize]{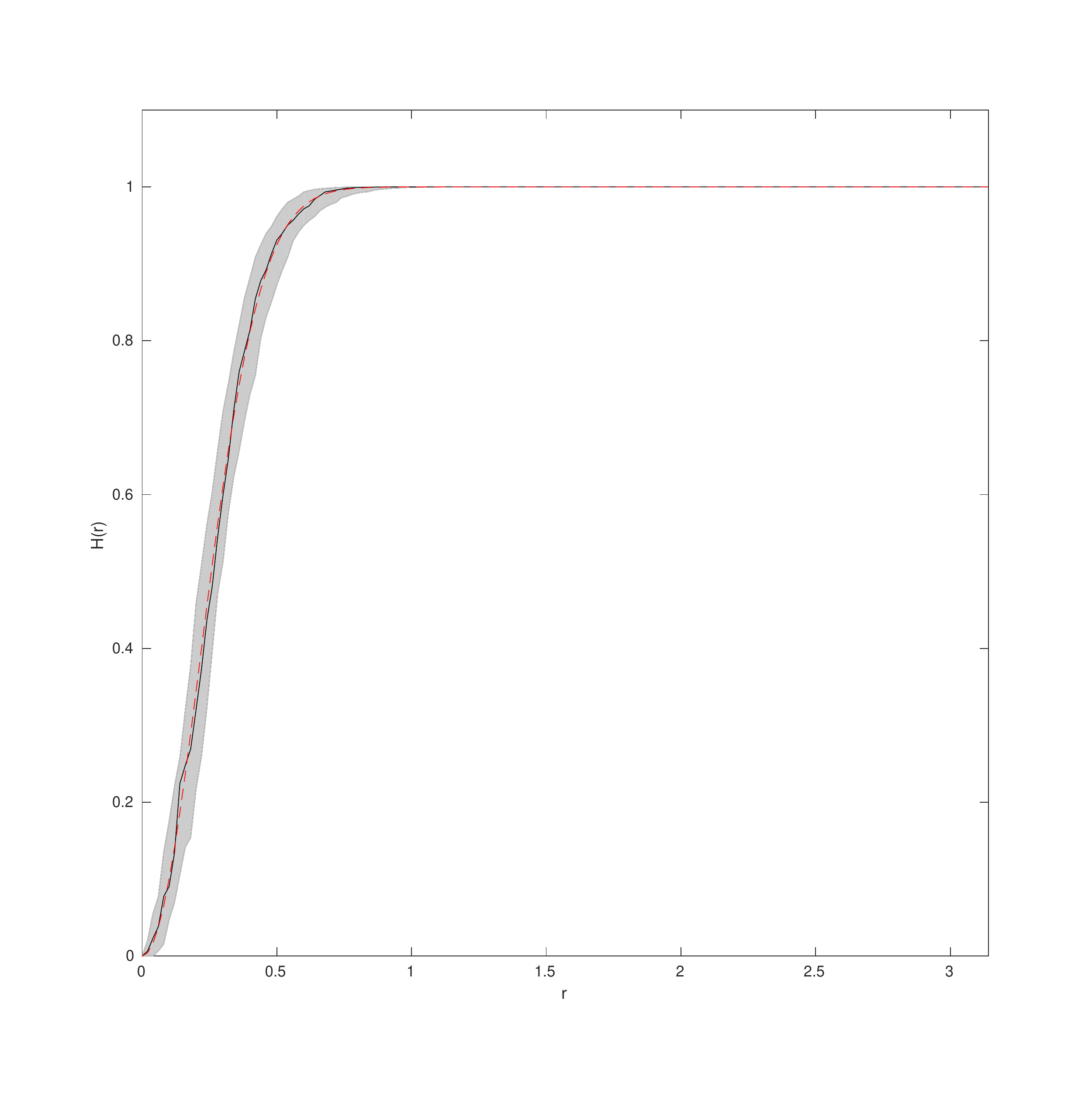}
  \end{subfigure}
  \begin{subfigure}{0.24\hsize}\centering
      \includegraphics[trim={2cm 2cm 2cm 2cm},clip=TRUE,width=1\hsize]{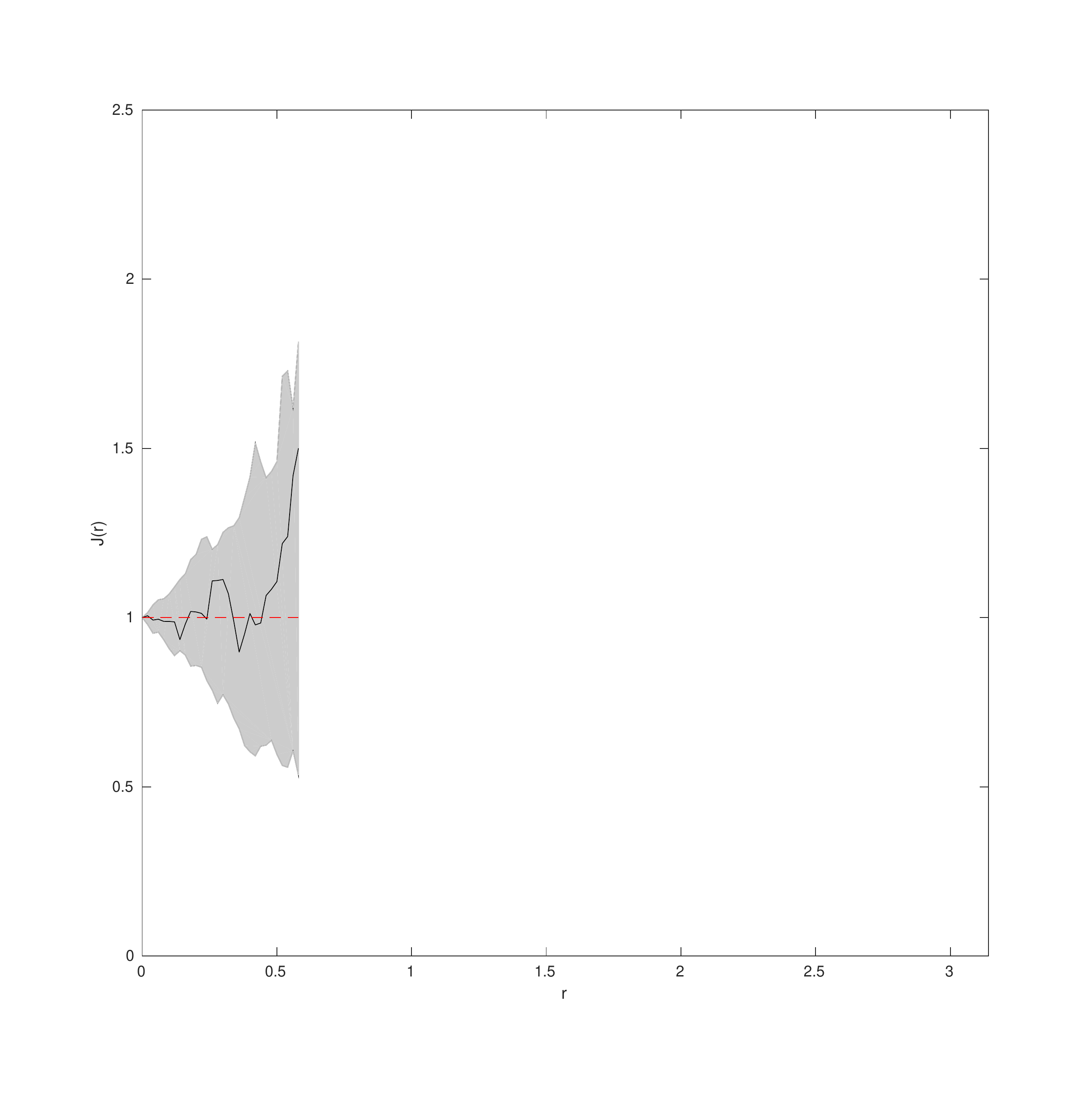}
  \end{subfigure}

  \begin{subfigure}{0.24\hsize}\centering
      \includegraphics[trim={2cm 2cm 2cm 2cm},clip=TRUE,width=1\hsize]{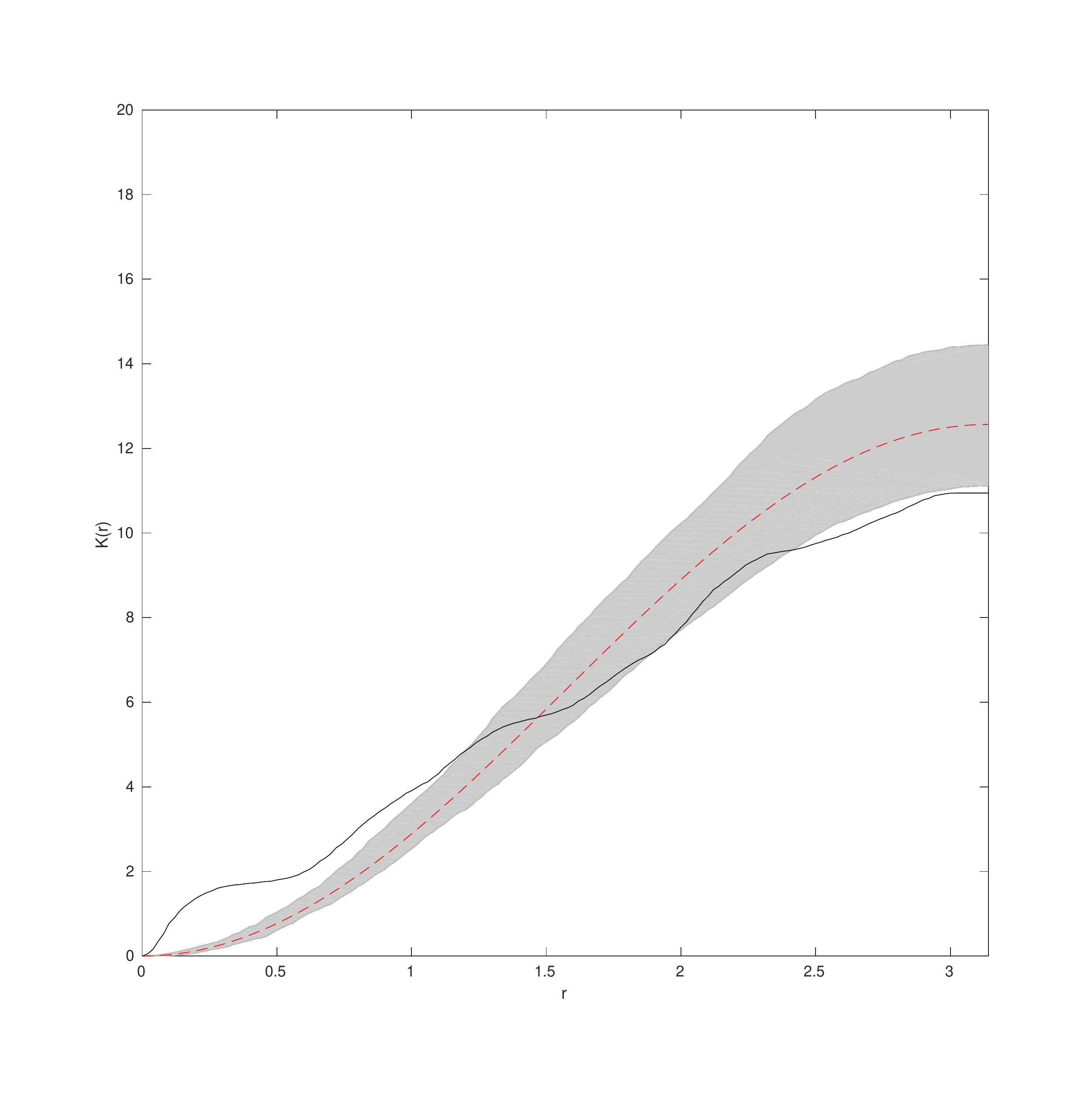}
  \end{subfigure}
  \begin{subfigure}{0.24\hsize}\centering
      \includegraphics[trim={2cm 2cm 2cm 2cm},clip=TRUE,width=1\hsize]{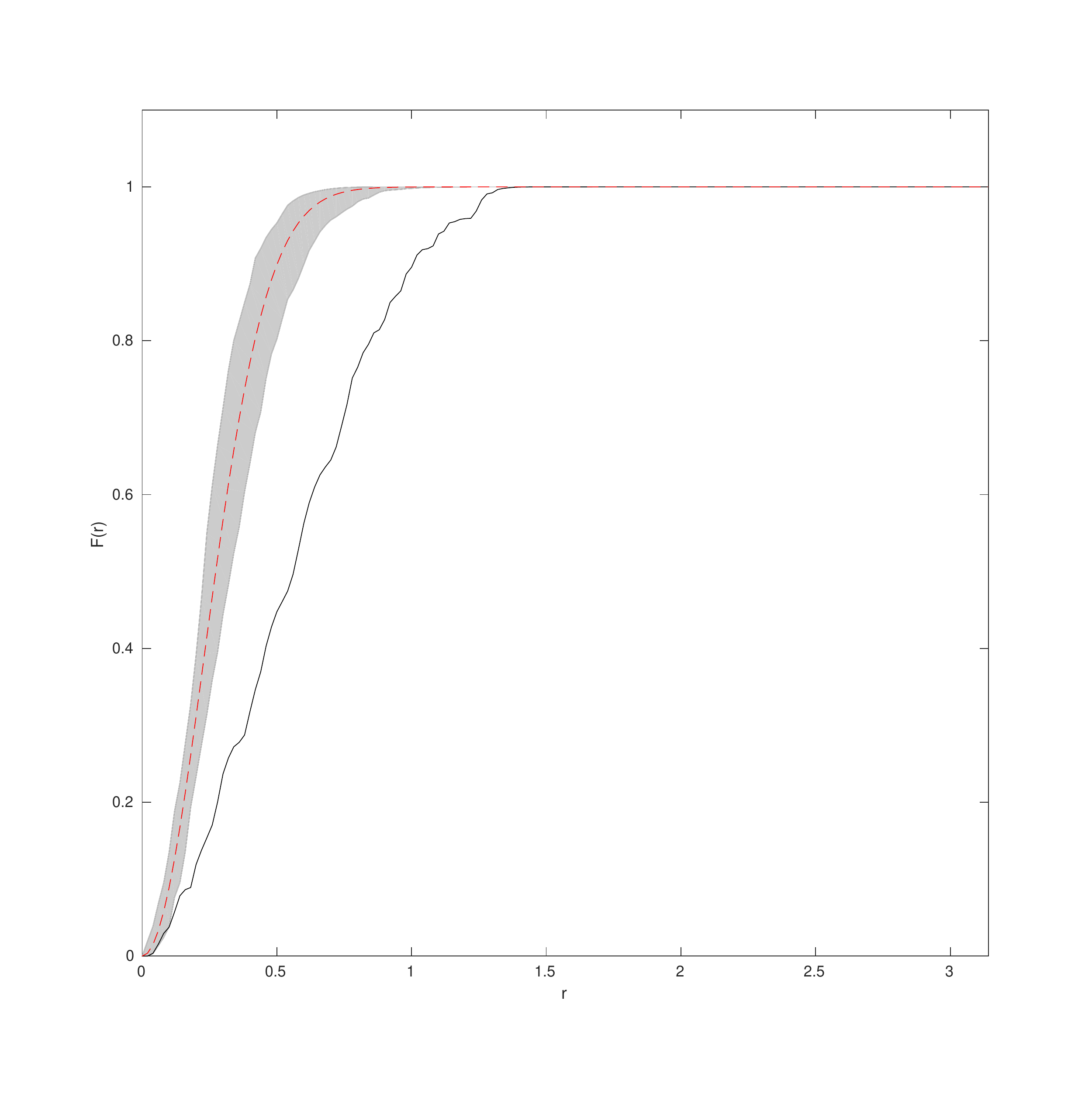}
  \end{subfigure}
  \begin{subfigure}{0.24\hsize}\centering
      \includegraphics[trim={2cm 2cm 2cm 2cm},clip=TRUE,width=1\hsize]{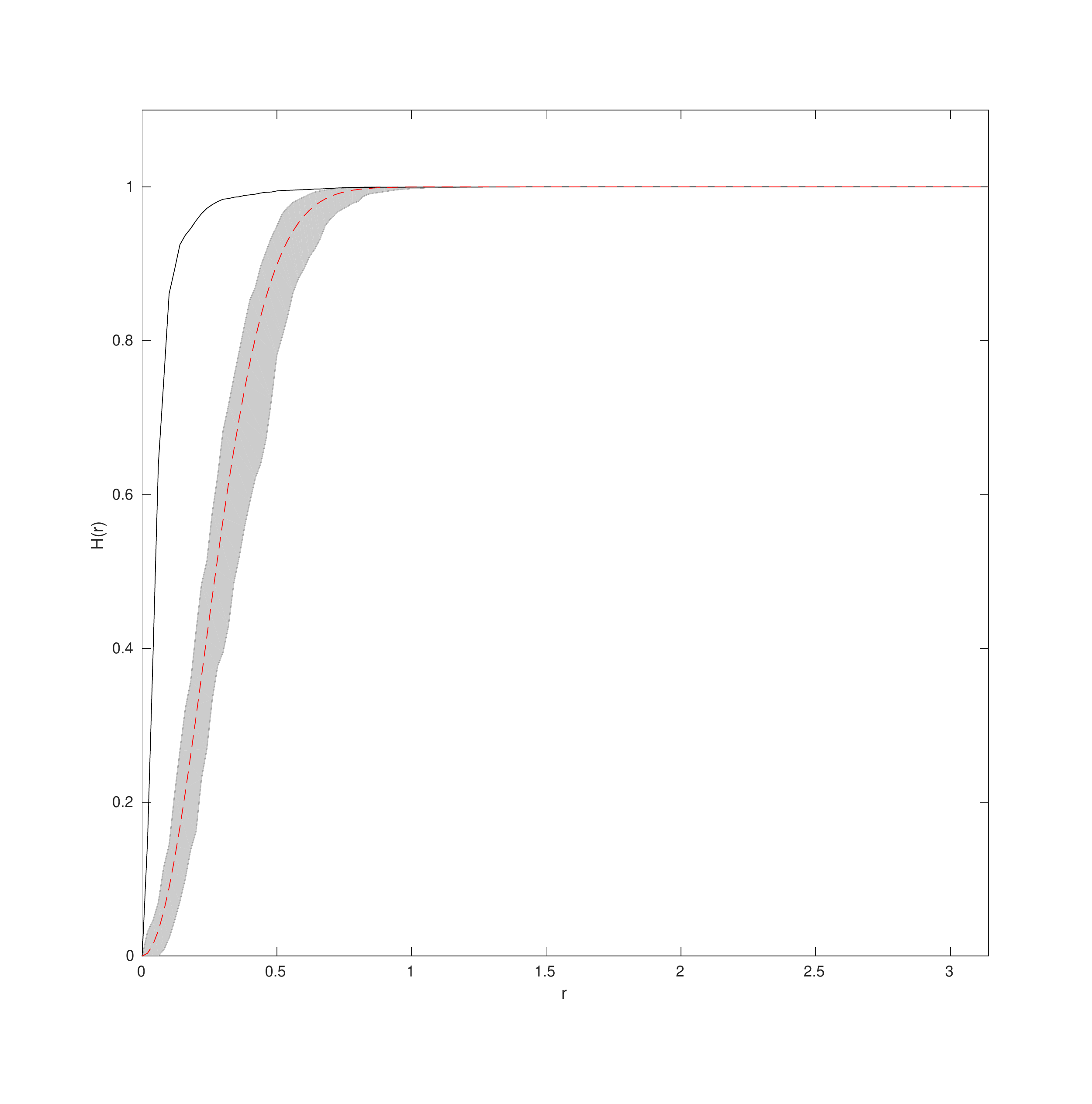}
  \end{subfigure}
  \begin{subfigure}{0.24\hsize}\centering
      \includegraphics[trim={2cm 2cm 2cm 2cm},clip=TRUE,width=1\hsize]{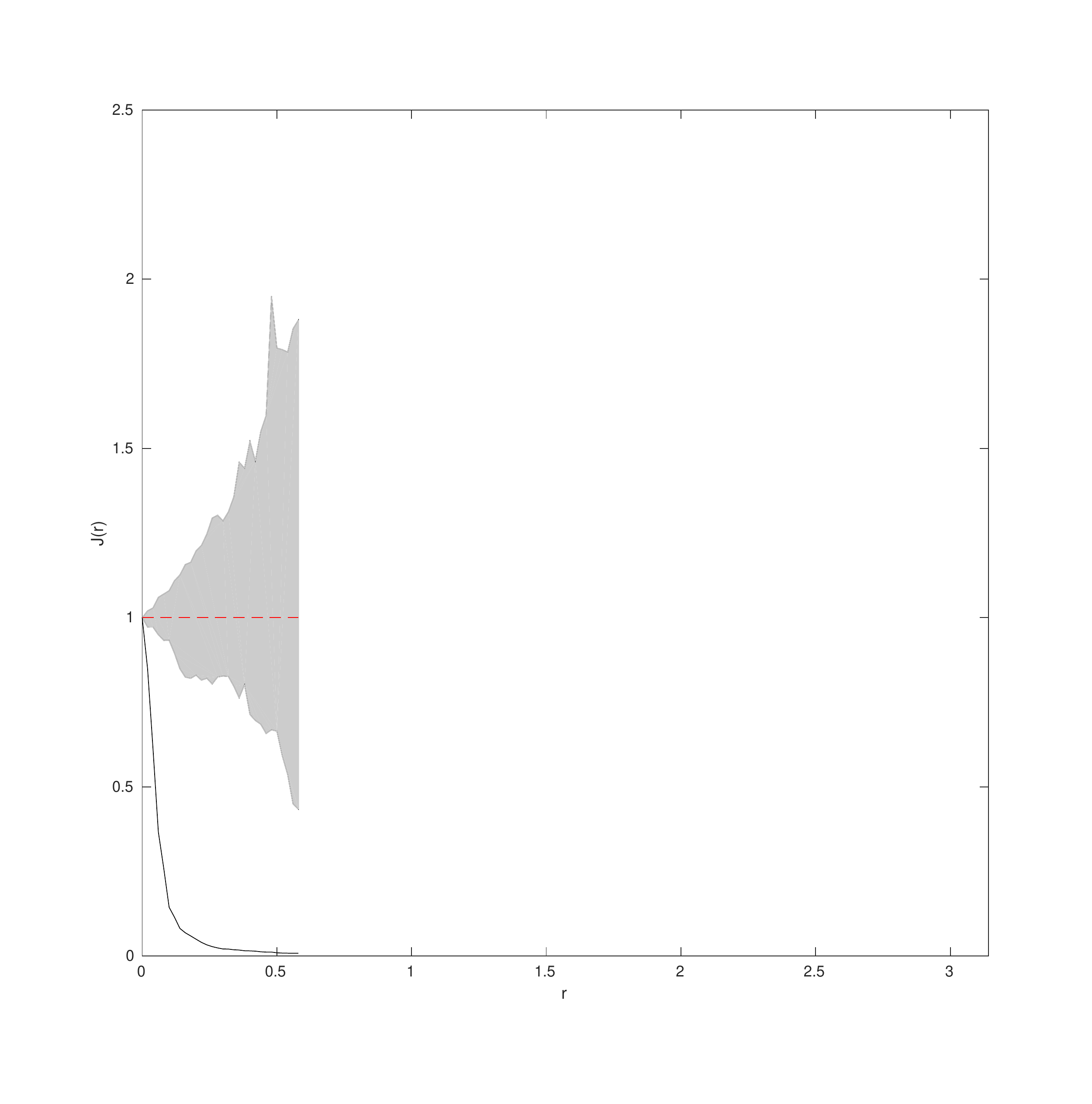}
  \end{subfigure}
  \vspace{0cm}

\caption{Example of (from left to right) $K_{\text{inhom}}$-, $J_{\text{inhom}}$-, $F_{\text{inhom}}$-, and $H_{\text{inhom}}$-functions for a Mat\'{e}rn \rom{2} with parameters $R=0.3$ (\emph{top row}), Poisson process, and Thomas process with parameters $\kappa =0.1$, exponential mark distribution with rate $\lambda =1$ and offspring expectation 15 (\emph{bottom row}) on a prolate spheroid with dimensions $(a,b,c)=(1,1,3)$ all with expectation 100. Black line is the estimated functional summary statistics for our observed data, dashed red line is the theoretical functional summary statistic for a Poisson process, and the grey shaded area is the simulation envelopes from 99 Monte Carlo simulations of Poisson processes fitted to the observed data.}\label{fig:sub2}
\hrulefill
\end{sidewaysfigure}

We simulate Mat\'{e}rn \rom{2}, and Thomas processes and construct estimates of their functional summary statistics under the assumption that they are CSR. The inhomogenous functional summary statistics are displayed in Figure \ref{fig:sub2}. Comparing Figure \ref{fig:sub2} to typical functional summary statistics for regular and cluster processes in $\mathbb{R}^2$, we see the same types of deviations away from CSR. In particular, we see for regular processes with small $r$ that there are negative deviations, whilst the cluster process has large positive deviations for the $\tilde{K}_{\text{inhom}}$-function. Furthermore, the $\hat{J}_{\text{inhom}}$-function shows significant positive deviations for regular processes whilst negative ones are observed for cluster processes.

\section{Testing for CSR on convex shapes in $\mathbb{R}^3$}\label{Section:6}

Exploratory data analysis for spatial point patterns in $\mathbb{R}^2$ typically begins with testing whether the observed point pattern exhibits CSR where test statistics are frequently based on the $L$-function, $L(r)=\sqrt{K(r)/\pi}$. On $\mathbb{R}^2$ and under CSR the $L$-function is linear in $r$ and variance stabilised \cite{Besag1977} whilst \cite{Lawrence2018} discusses the analogue $L$-function in $\mathbb{S}^2$ where again it is variance stabilised when the underlying process is CSR. As we are working with inhomogeneous Poisson processes on $\mathbb{S}^2$ an equivalent transformation for the $L$ function has not been discussed previously and instead we propose using test statistics derived from standardisations of the functional summary statistics \cite{Lagache2013}. In order to construct such test statistics we must derive first and second order properties of the estimated functional summary statistics. Section \ref{Section:8} discusses derivations for any spherical Poisson process when $\rho$ is known. In this section we consider the scenario when we have a homogeneous Poisson process on $\mathbb{D}$ with \emph{unknown}, constant intensity function $\rho\in\mathbb{R}_+$. Furthermore we shall only focus on the inhomogeneous $K$-function as standardisation of the remaining functional summary statistics follow identically.

\subsection{Test statistic for CSR}

Given a homogeneous Poisson process on $\mathbb{D}$ with constant intensity function $\rho\in\mathbb{R}_+$, we map this to $\mathbb{S}^2$ giving a new Poisson process on the sphere with inhomogeneous intensity function given by Theorem \ref{thm:mapping:general:text} as
\begin{equation}\label{rho:CSR:mapping}
\rho^*(\mathbf{x})=
\left\{\begin{alignedat}{2}
\rho l_1(f^{-1}(\mathbf{x})) &J_{(1,f^*)}(\mathbf{x})\sqrt{1-x_1^2-x_2^2}, \quad \mathbf{x}\in f(\mathbb{D}_1)\\
&\vdots\\
\rho l_n(f^{-1}(\mathbf{x})) &J_{(n,f^*)}(\mathbf{x})\sqrt{1-x_1^2-x_2^2}, \quad \mathbf{x}\in f(\mathbb{D}_n).
  \end{alignedat}\right. 
\end{equation}

Using Theorems \ref{mean:functional:summary:statistics} and \ref{thm:variance:K:inhom:text} we can calculate the mean and variances of the inhomogeneous $K$-function when $\rho$ is known. When $\rho$ is unknown we use estimators of $\rho$ when constructing functional summary statistics. In particular we use
\begin{equation*}
\hat{\rho}=\frac{N_X({\mathbb{D}})}{\lambda_{\mathbb{D}}(\mathbb{D})}, \quad \hat{\rho^2}=\frac{N_X{(\mathbb{D})}(N_X{(\mathbb{D})}-1)}{\lambda_{\mathbb{D}}^2(\mathbb{D})},
\end{equation*}
which are both unbiased for $\rho$ and $\rho^2$ respectively by application of the Campbell-Mecke Theorem \cite{Moller2004}. Thus our estimator for $K_\text{inhom}(r)$ when $\rho$ is unknown takes the following form,
\begin{equation}\label{K:unknown:rho}
\tilde{K}_{\text{inhom}}(r)=
\begin{cases}
\frac{\lambda_{\mathbb{D}}^2(\mathbb{D})}{4\pi N_Y(\mathbb{S}^2)(N_Y(\mathbb{S}^2)-1)}\sum_{\mathbf{x}\in Y}\sum_{\mathbf{y}\in Y\setminus \{\mathbf{x}\}}\frac{\mathbbm{1}[d(\mathbf{x},\mathbf{y})\leq r]}{\tilde{\rho}(\mathbf{x})\tilde{\rho}(\mathbf{y})}, &\text{if } N_Y(\mathbb{S}^2)>1\\
0, &\text{otherwise},
\end{cases}
\end{equation}
where $Y=f(X)$, $f$ is our mapping from the ellipsoid to the sphere, and $\tilde{\rho}(\mathbf{x})$ is given by,
\begin{equation}\label{rho:tilde}
\tilde{\rho}(\mathbf{x})=
\left\{\begin{alignedat}{2}
l_1(f^{-1}(\mathbf{x})) &J_{(1,f^*)}(\mathbf{x})\sqrt{1-x_1^2-x_2^2}, \quad \mathbf{x}\in f(\mathbb{D}_1)\\
&\vdots\\
l_n(f^{-1}(\mathbf{x})) &J_{(n,f^*)}(\mathbf{x})\sqrt{1-x_1^2-x_2^2}, \quad \mathbf{x}\in f(\mathbb{D}_n).
\end{alignedat}\right.
\end{equation} Note that $N_{Y}(\mathbb{S}^2)=N_X(\mathbb{D}^2)$. 

\cite{Diggle2003} proposes using the maximum absolute value between the theoretical and the estimated functional summary statistics to test for CSR. Based on this we follow the work of \cite{Lagache2013} and propose the test statistic
\begin{equation}\label{K:test:stat}
T = \sup_{r\in [0,\pi]}\left|\frac{\tilde{K}_{\text{inhom}}(r)-2\pi(1-\cos(r))}{\widehat{\text{Var}}(\tilde{K}_{\text{inhom}}(r))}\right|.
\end{equation}

In order to be able to construct the test statistic $T$, an estimate of the variance of the empirical functional summary statistics are required. Further, we need show that the bias of $\tilde{K}_{\text{inhom}}(r)$ is negligible and hence $\mathbb{E}[\tilde{K}_{\text{inhom}}(r)]\approx 2\pi(1-\cos(r))$ for Poisson processes, validating its use in (\ref{K:test:stat}). By using estimators for $\rho$ and $\rho^2$ we alter the first and second order properties given by Theorems (\ref{mean:functional:summary:statistics}) and (\ref{thm:variance:K:inhom:text}). In the following we consider the first and second order moments of $\tilde{K}_{\text{inhom}}$.

\subsection{Estimating moments of $\tilde{K}_{\text{inhom}}(r)$ on $\mathbb{S}^2$ for CSR process on $\mathbb{D}$}

\begin{thm}\label{thm:variance:K:unknown:rho:text}
The bias and variance of $\tilde{K}_{\text{inhom}}(r)$ are,
\begin{align*}
\text{Bias}(\tilde{K}_{\text{inhom}}(r))=-P(N_Y(\mathbb{S}^2)\leq 1) 2\pi(1-\cos r),
\end{align*}
and,
\begin{equation}\label{variance:K:unknown:rho:2}
\begin{split}
&\text{Var}(\tilde{K}_{\text{inhom}}(r))=4\pi^2(1-\cos r)^2(1-P(N_Y(\mathbb{S}^2)\leq 1))P(N_Y(\mathbb{S}^2)\leq 1)\\
&\phantom{AAAA}+\rho^3\lambda_{\mathbb{D}}^4(\mathbb{D})(1-\cos r)^2\left(\int_{\mathbb{S}^2}\frac{1}{\tilde{\rho}(\mathbf{x})}\lambda_{\mathbb{S}^2}(d\mathbf{x})-\frac{16\pi^2}{\lambda_{\mathbb{D}}(\mathbb{D})}\right)\mathbb{E}\left[\frac{1}{(N_Y(\mathbb{S}^2)+3)^2(N_Y(\mathbb{S}^2)+2)^2}\right]\\
&\phantom{AAAA}+\frac{\rho^2\lambda_{\mathbb{D}}^4(\mathbb{D})}{8\pi^2}\left(\int_{\mathbb{S}^2}\int_{\mathbb{S}^2} \frac{\mathbbm{1}[d(\mathbf{x}_1,\mathbf{x}_2)\leq r]}{\tilde{\rho}(\mathbf{x}_1)\tilde{\rho}(\mathbf{x}_2)} \lambda_{\mathbb{S}^2}(d\mathbf{x}_1) \lambda_{\mathbb{S}^2}(d\mathbf{x}_2) - \frac{64\pi^4(1-\cos r)^2}{\lambda_{\mathbb{D}}^2(\mathbb{D})} \right)\mathbb{E}\left[\frac{1}{(N_Y(\mathbb{S}^2)+2)^2(N_Y(\mathbb{S}^2)+1)^2}\right],
\end{split}
\end{equation}
where $\tilde{\rho}(\mathbf{x})$ is given by Equation (\ref{rho:tilde}).
\end{thm}
\begin{proof}
See Theorem 6 in Appendix G.
\end{proof}

The form of the variance derived in Theorem \ref{thm:variance:K:unknown:rho:text} is near identical to that derived by \cite{Lang2010} except that our derivations considers inhomogeneous Poisson processes, does not require corrections for edge effects, and the space is $\mathbb{S}^2$ instead of $\mathbb{R}^2$. Further we can bound the absolute value of the bias as follows
\begin{align}
|\text{Bias}(\tilde{K}_{\text{inhom}}(r))|&=P(N_Y(\mathbb{S}^2)\leq 1)2\pi(1-\cos r)\nonumber\\
&=\exp(-\mu_Y(\mathbb{S}^2))(1+\mu_Y(\mathbb{S}^2))2\pi(1-\cos r)\nonumber\\
&\leq 4\pi(1+\mu_Y(\mathbb{S}^2))\exp(-\mu_Y(\mathbb{S}^2))\label{eq:ineq:bias:1}\\
&\leq 4\pi(1+\mu_Y(\mathbb{S}^2))\mu_Y(\mathbb{S}^2)^{-\mathrm{e}}=\text{\large{\emph{O}}}\left(\mu_Y^{1-\mathrm{e}}(\mathbb{S}^2)\right),\label{eq:ineq:bias:2}
\end{align}
where $\mu_Y$ is the intensity measure of $Y=f(X)$, the inequality in (\ref{eq:ineq:bias:1}) is attained by setting $r=\pi$, and (\ref{eq:ineq:bias:2}) follows from $\mathrm{e}^{x}\geq x^{\mathrm{e}}$. Thus, for shapes considered in this work, the bias will be negligible.

From Theorem \ref{thm:variance:K:unknown:rho:text} it is possible to construct a ratio-unbiased estimator for the variance. In particular by the Campbell-Mecke Theorem, and defining the estimator $\hat{\rho}_k=N_Y(\mathbb{S}^2)(N_Y(\mathbb{S}^2)-1)\cdots(N_Y(\mathbb{S}^2)-k-1)/\lambda_{\mathbb{D}}^k(\mathbb{D})$, then $\mathbb{E}[\hat{\rho}_k]=\rho^k$ and so $\hat{\rho}_k$ is unbiased for $\rho^k$. We can substitute the expectations in (\ref{variance:K:unknown:rho:2}) with their corresponding observed values, for example we substitute $(N_Y(\mathbb{S}^2)+3)^{-2}(N_Y(\mathbb{S}^2)+2)^{-2}$ for $\mathbb{E}[(N_Y(\mathbb{S}^2)+3)^{-2}(N_Y(\mathbb{S}^2)+2)^{-2}]$. Additionally, the following lemma helps derive a ratio unbiased estimator for $P(N_Y(\mathbb{S}^2)<1)$.
\begin{lemma}\label{lemma:ratio:unbiased:prob:text}
Let $N\sim \text{Poisson}(\lambda)$, $k\in\mathbb{N}$ and $p\in \mathbb{R}_+$. Define the following random variable,
\begin{equation}
R = \frac{N!e^{N-k}}{(N-k)!(e+p)^N}.
\end{equation}
Then $R$ is ratio-unbiased for $\lambda^ke^{-p\lambda}$.
\end{lemma}{}
\begin{proof}
See Lemma 2 in Appendix G.
\end{proof}
Using Lemma \ref{lemma:ratio:unbiased:prob:text} we can construct a ratio-unbiased estimator for $(1-P(N_Y(\mathbb{S}^2)<1))P(N_Y(\mathbb{S}^2)<1)$. Defining $\lambda=\rho\lambda_{\mathbb{D}}(\mathbb{L})$,
\begin{align*}
(1-P(N_Y(\mathbb{S}^2)<1))P(N_Y(\mathbb{S}^2)<1)&=(1-e^{-\lambda}-\lambda e^{-\lambda})(e^{-\lambda}+\lambda e^{-\lambda})\\
&=e^{-\lambda}+\lambda e^{-\lambda}-e^{-2\lambda}-2\lambda e^{-2\lambda}-\lambda^2e^{-2\lambda},
\end{align*}
and so a ratio-unbiased estimator for $(1-P(N_Y(\mathbb{S}^2)<1))P(N_Y(\mathbb{S}^2)<1)$ is
\begin{equation*}
\frac{e^{N_{Y}(\mathbb{S}^2)}}{(e+1)^{N_{Y}(\mathbb{S}^2)}}+\frac{{N_{Y}(\mathbb{S}^2)}e^{{N_{Y}(\mathbb{S}^2)}-1}}{(e+1)^{N_{Y}(\mathbb{S}^2)}}+\frac{e^{{N_{Y}(\mathbb{S}^2)}}}{(e+2)^{N_{Y}(\mathbb{S}^2)}}-\frac{2{N_{Y}(\mathbb{S}^2)}e^{{N_{Y}(\mathbb{S}^2)}-1}}{(e+2)^{N_{Y}(\mathbb{S}^2)}}-\frac{{N_{Y}(\mathbb{S}^2)}({N_{Y}(\mathbb{S}^2)}-1)e^{{N_{Y}(\mathbb{S}^2)}-2}}{(e+2)^{N_{Y}(\mathbb{S}^2)}}.
\end{equation*}
Plugging the given estimators for $(1-P(N_Y(\mathbb{S}^2)<1))P(N_Y(\mathbb{S}^2)<1)$, $\rho^2$, $\rho^3$, $\mathbb{E}\left[(N_Y(\mathbb{S}^2)+3)^{-2}(N_Y(\mathbb{S}^2)+2)^{-2}\right]$ and $\mathbb{E}\left[(N_Y(\mathbb{S}^2)+2)^{-2}(N_Y(\mathbb{S}^2)+1)^{-2}\right]$ into (\ref{variance:K:unknown:rho:2}) gives a ratio unbiased estimator for $\text{Var}(\tilde{K}_{\text{inhom}}(r))$, which in turn allows for the construction of the test statistic $T$ in (\ref{K:test:stat}).

\subsection{Standardised inhomogeneous $K$-function plots}

\begin{figure}

  \begin{subfigure}{0.49\hsize}\centering
      \includegraphics[trim={2cm 2cm 2cm 2cm},clip=TRUE,width=\hsize,height=6cm]{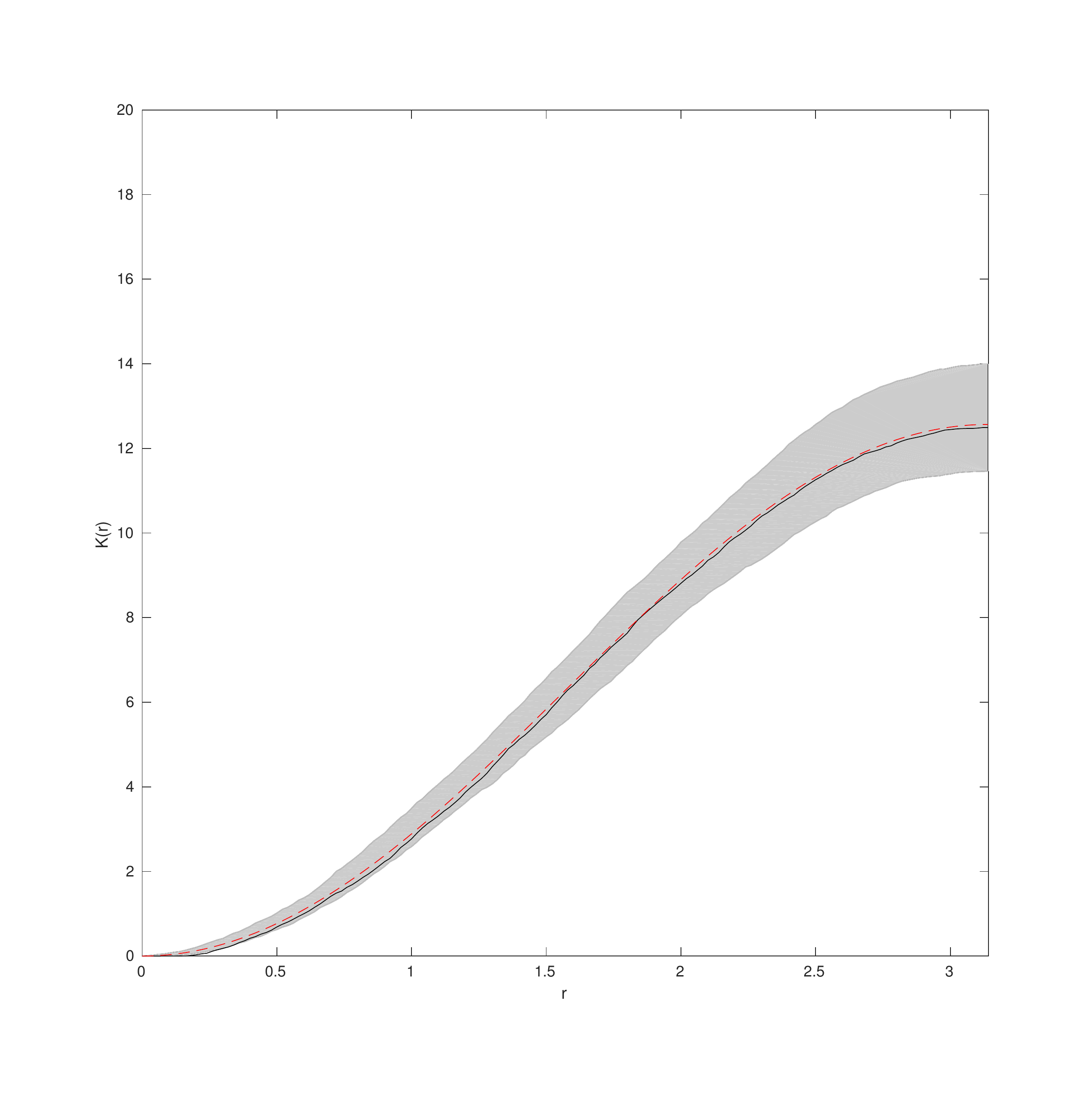}
  \end{subfigure}
  \begin{subfigure}{0.49\hsize}\centering
      \includegraphics[trim={2cm 2cm 2cm 2cm},clip=TRUE,width=\hsize,height=6cm]{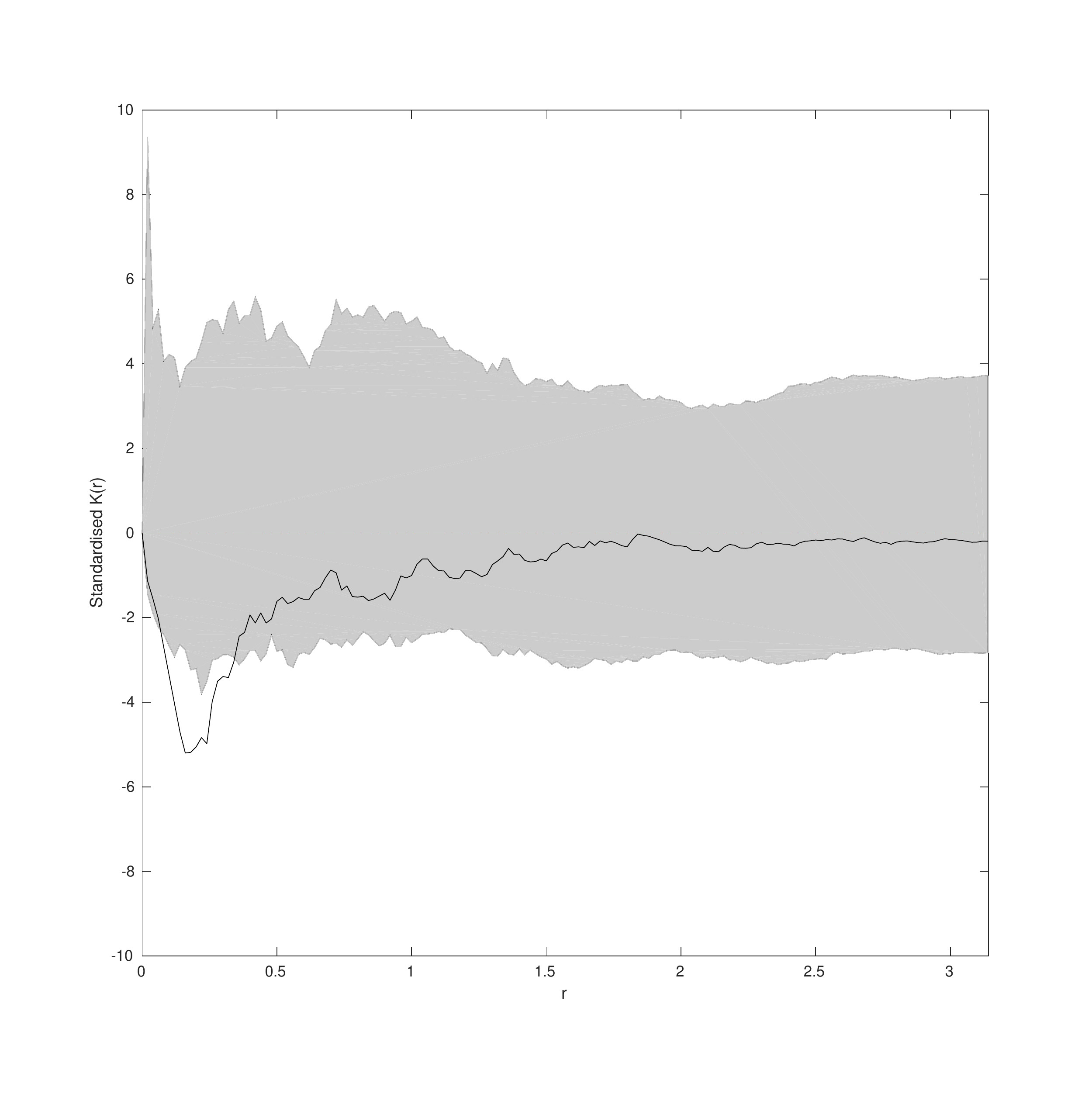}
  \end{subfigure}

  \begin{subfigure}{0.49\hsize}\centering
      \includegraphics[trim={2cm 2cm 2cm 2cm},clip=TRUE,width=\hsize,height=6cm]{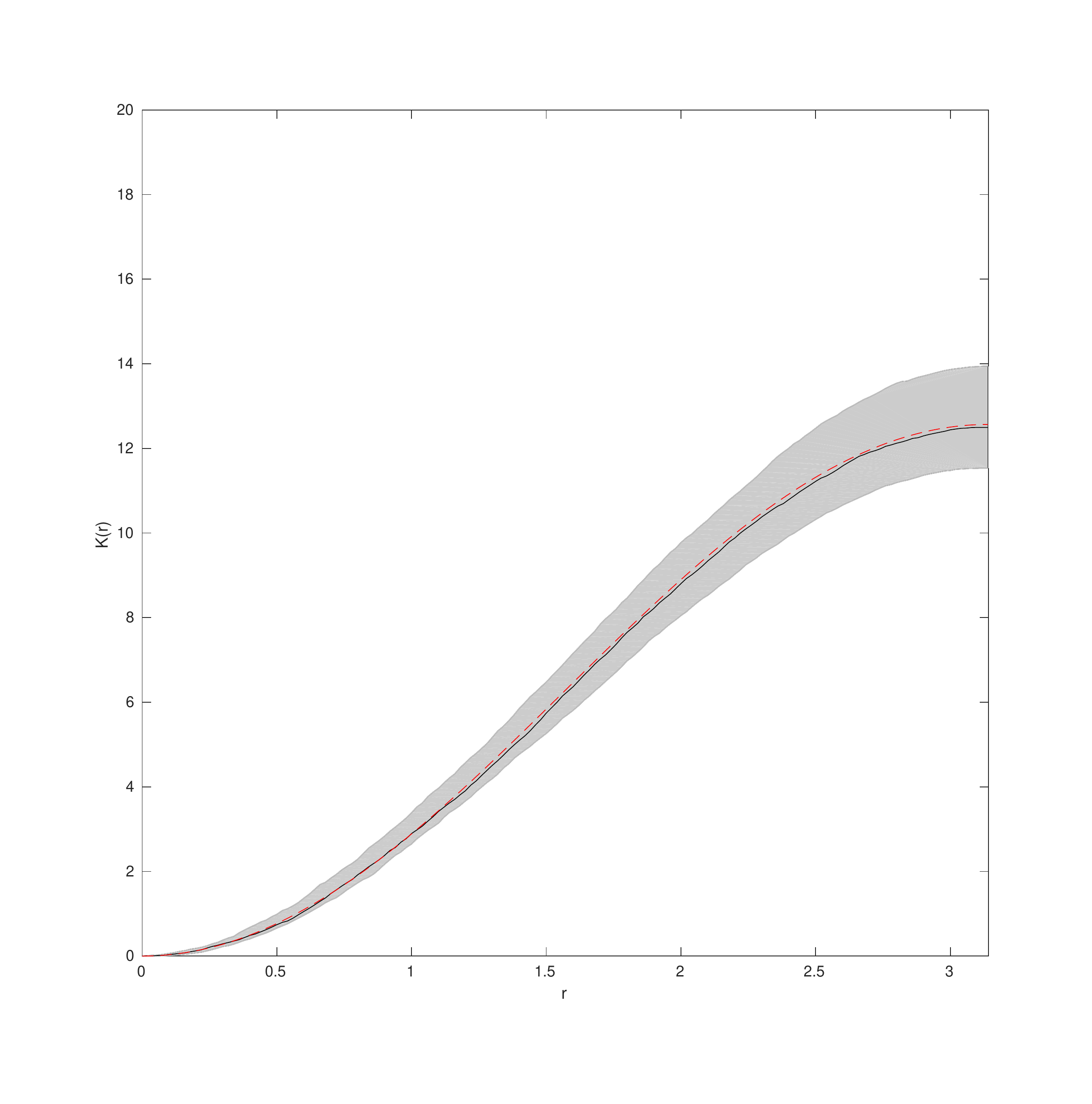}
  \end{subfigure}%
  \begin{subfigure}{0.49\hsize}\centering
      \includegraphics[trim={2cm 2cm 2cm 2cm},clip=TRUE,width=\hsize,height=6cm]{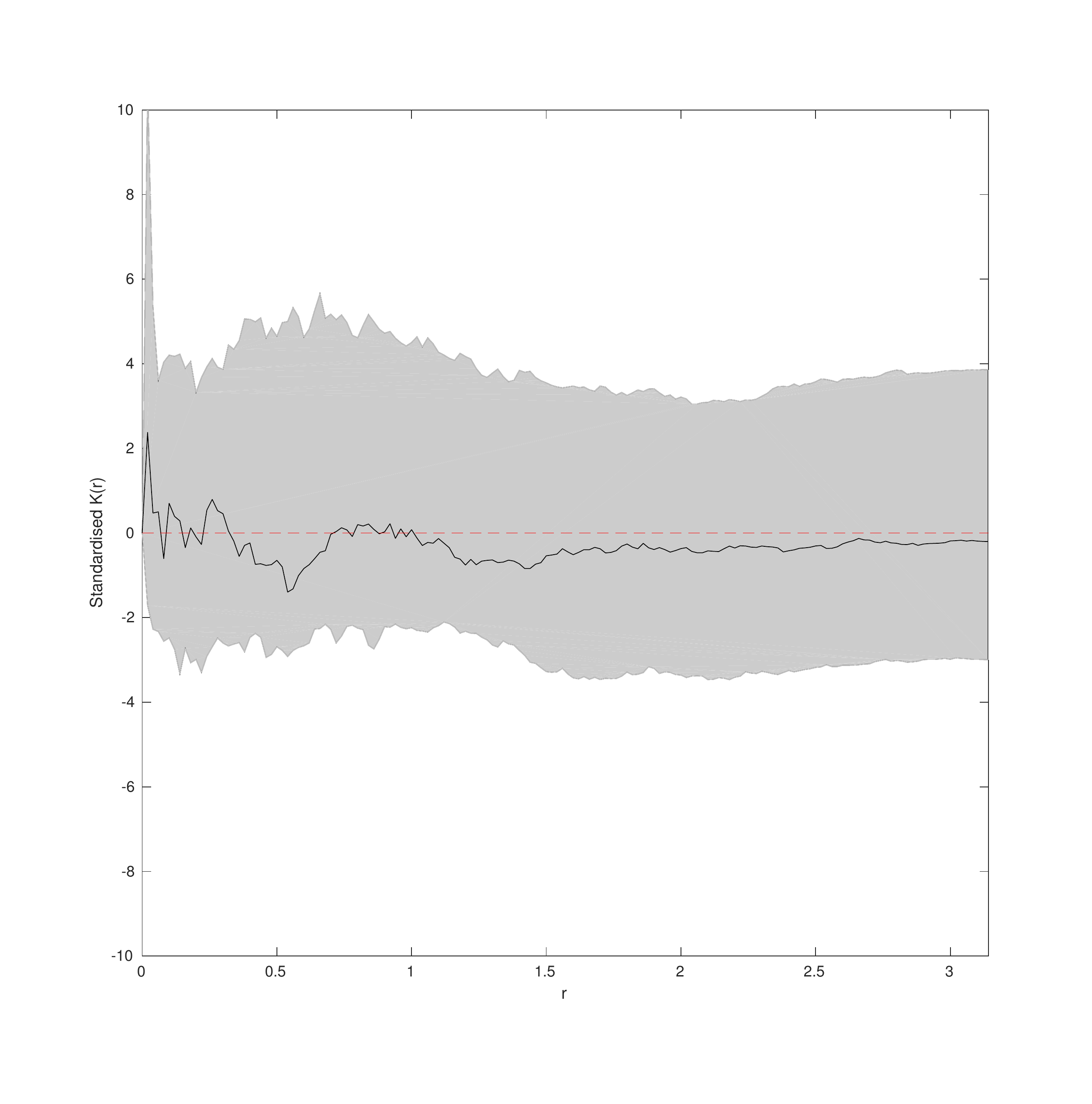}
  \end{subfigure} 

  \begin{subfigure}{0.49\hsize}\centering
      \includegraphics[trim={2cm 2cm 2cm 2cm},clip=TRUE,width=\hsize,height=6cm]{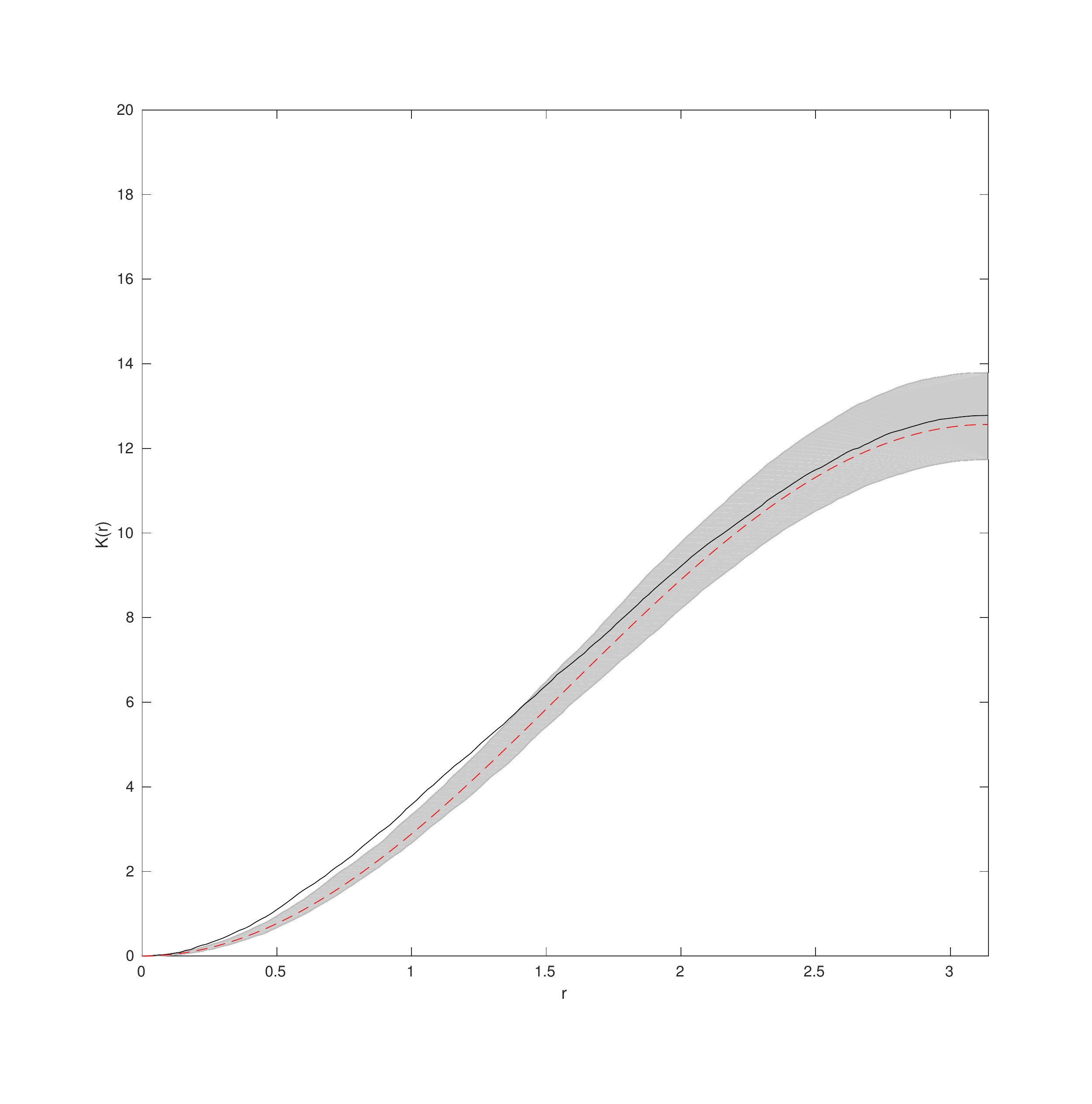}
  \end{subfigure}
  \begin{subfigure}{0.49\hsize}\centering
      \includegraphics[trim={2cm 2cm 2cm 2cm},clip=TRUE,width=\hsize,height=6cm]{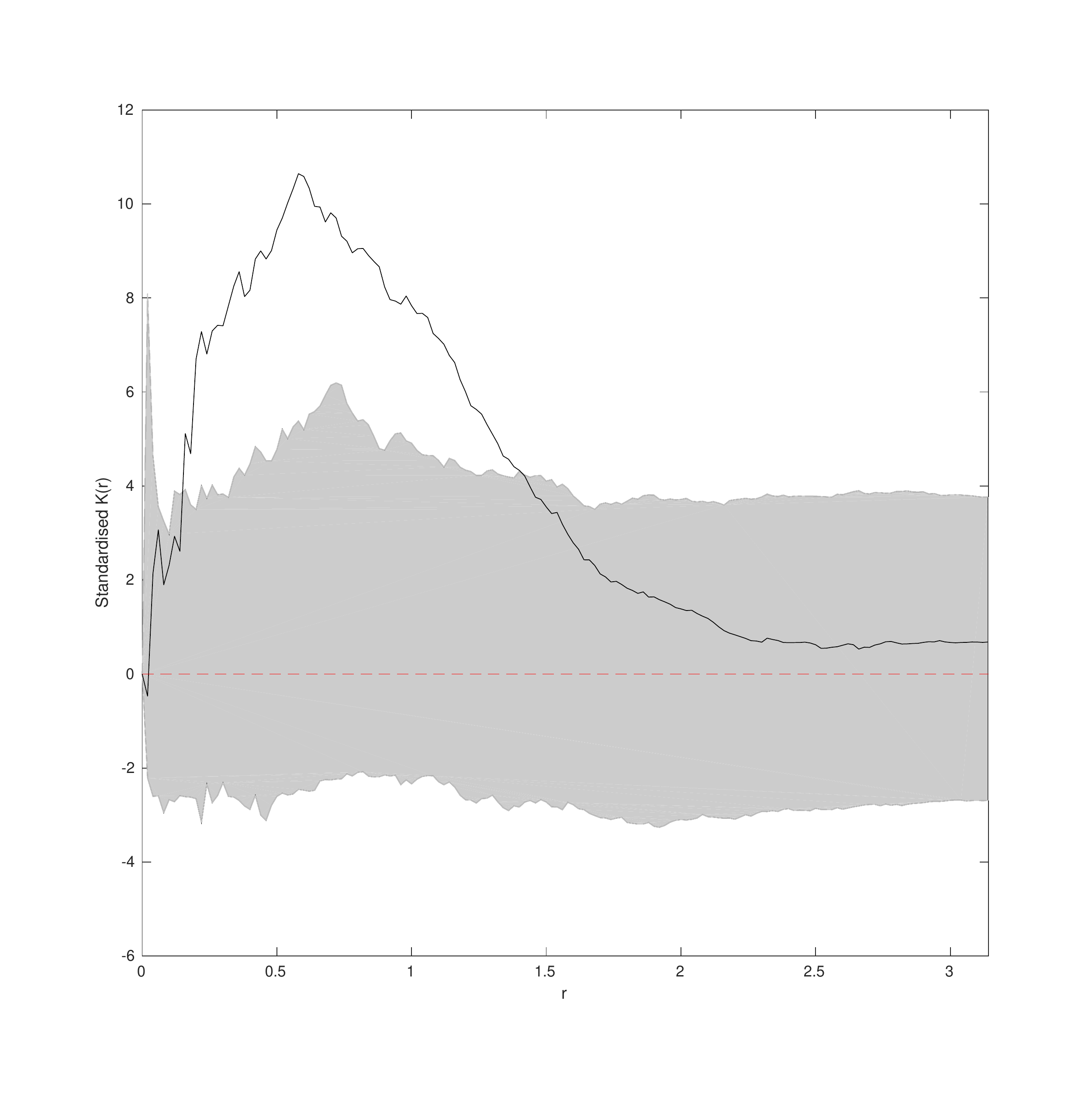}
  \end{subfigure}

\caption{Example of $K_{\text{inhom}}(r)$- (\emph{left}) and standardised $K_{\text{inhom}}(r)$- (\emph{right}) functions for a Mat\'{e}rn \rom{2} with parameters $R=0.2$ and expectation 100 (\emph{top}), Poisson with expectation $40\pi$ (\emph{middle}) and Thomas processes with parameters $\kappa =0.5$, expectation of 150, an exponential mark distribution with rate $\lambda=1$ and offspring mean of 20 on a prolate spheroid with $a=b=0.8000,c=1.43983$ (dimensions chosen so that the area of the ellipsoid is $4\pi$). Black line is the estimated functional summary statistics for our observed data, dazshed red line is the theoretical functional summary statistic for a Poisson process, and the grey shaded area is the simulation envelopes from 999 Monte Carlo simulations of Poisson processes fitted to the observed data.}
\label{fig:matern:functions}
\hrulefill
\end{figure}

Figure \ref{fig:matern:functions} highlights how the empirical $K$-function estimates deviate when the underlying process is not CSR. For the regular processes we notice considerable negative deviations for small $r$ whilst for cluster processes positive deviations are observed, highlighted in the right column of Figure \ref{fig:matern:functions}. 

Intuitively, this is to be expected, with a near identical reasoning to what is observed for the $K$-function in $\mathbb{R}^{2,3}$. Since the regular process has a hard-core distance between events, we observe estimates for $K_{\text{inhom}}(r)$ that are close to zero for small $r$, thus resulting in the large negative deviation observed in Figure \ref{fig:matern:functions}. On the other hand, for the Thomas cluster process, we observe events in closer proximity than would be expected for a CSR process, thus the estimated $K_{\text{inhom}}(r)$ function has large positive deviations away from CSR.

\section{Simulation study}\label{Section:9}

We conduct empirical Type \rom{1} and \rom{2} error studies to evaluate the effectiveness of the proposed test statistic in determining whether or not a point process on a convex shape exhibits CSR. We consider different prolate ellipsoids such that the area of the ellipsoid is the same across differing semi-major axis lengths. This will allow us to determine how the power of our test changes as the space under consideration deforms further away from the unit sphere.

\subsection{Design of simulations}

In order to best understand the properties of our testing procedure we will consider CSR, Mat\'{e}rn \rom{2} and Thomas processes on different prolate spheroids. We design the experiments such that the expected number of events is similar across all experiments. For both the CSR and Thomas process simulations this is easily controlled. For a Poisson process, the expected number on $\mathbb{D}$ is $\rho\lambda_{\mathbb{D}}(\mathbb{D})$ whilst for a Thomas process it is given by Proposition S3 in Appendix F. 

On the other hand, the Mat\'{e}rn \rom{2} process requires a little more attention since Corollary S4 (see Appendix F) limits the maximum expected number of possible events for a given space $\mathbb{D}$. Thus, for a given expected number $\mu$ that is less than or equal to the one prescribed by  Corollary S4 (see Appendix F) we fix the hard-core distance $R$ and solve the following equation for $\rho$
\begin{equation}\label{mu:matern:2}
\int_B \frac{1-\mathrm{e}^{-\rho\lambda_{\mathbb{D}}(B_{\mathbb{D}}(\mathbf{x},R))}}{\lambda_{\mathbb{D}}(B_{\mathbb{D}}(\mathbf{x},R))} d\mathbf{x}=\mu.
\end{equation}

A full outline of all the experiments and the parameters chosen are given in Tables \ref{CSR:table:results}, \ref{Regular:table:results}, and \ref{Cluster:table:results} for CSR, regular, and cluster process simulations respectively. Note that when $R=0$ for the Mat\'{e}rn \rom{2} process and when $\kappa=\infty$ for the Thomas process both processes are CSR.

\subsection{Test Statistics}
Due to the computational intensity of optimising (\ref{K:test:stat}) we instead calculate $|(\tilde{K}_{\text{inhom}}(r)-$ $2\pi(1-\cos r))/\widehat{\text{Var}}(\tilde{K}_{\text{inhom}}(r))|$ for $r\in\mathcal{R}=\{r_1,\dots,r_m\},\; m\in\mathbb{N}$, where $\mathcal{R}$ is a finite set of distinct, evenly spaced points such that $r_{i}\in [0,\pi], i=1,\dots,m$, for the purposes of our simulation studies. We then take our test statistic as
\begin{equation*}
T = \max_{r\in \mathcal{R}}\left|\frac{\tilde{K}_{\text{inhom}}(r)-2\pi(1-\cos r)}{\widehat{\text{Var}}(\tilde{K}_{\text{inhom}}(r))}\right|,
\end{equation*}
where for these simulation studies we set $\mathcal{R}=\{0,0.02,0.04,\dots,\pi\}$. These simulations are tested at a $5\%$ significance level. Each experiment is repeated 1000 times, and for each experiment we simulate 999 Poisson processes to approximate the critical values of the hypothesis test.

\subsection{Results}

\begin{table}[p]
\centering
\begin{tabular}{x{0.14\textwidth}x{0.14\textwidth}x{0.14\textwidth}x{0.14\textwidth}x{0.14\textwidth}x{0.14\textwidth}}
\hline
Experiment No. & Expectation & $a$ & $\rho$ & Accept $H_0$ & Reject $H_0$ \tabularnewline\hline

1a & $40\pi$ & 1 & 10 & 0.9520 & 0.0480  \tabularnewline
1b & $40\pi$ & 0.8 & 10 &  0.9610 & 0.0390  \tabularnewline
1c  & $40\pi$ & 0.6 &10 &  0.9570   & 0.0430   \tabularnewline
1d & $40\pi$ & 0.4 & 10 &  0.9440 &  0.0560 \tabularnewline\hline
\end{tabular}
\vspace{1em}
\caption{Results when the \emph{observed} data is CSR. The semi-major axis length along the $x$-axis, $a$, and $y$-axis, $b$, are equivalent and the semi-major axis length along the $z$-axis is determined such that the area of the ellipsoid is $4\pi$. }
\label{CSR:table:results}
\end{table}

\begin{table}[p]
\centering
\begin{tabular}{x{0.14\textwidth}x{0.14\textwidth}x{0.14\textwidth}x{0.14\textwidth}x{0.14\textwidth}x{0.14\textwidth}}
\hline
Experiment No. & Expectation & $a$ & $R$ & Accept $H_0$ & Reject $H_0$ \tabularnewline\hline
2ai & 100 & 1 & 0  & 0.9250 & 0.0750  \tabularnewline
2aii & 100 & 1 & 0.05 & 0.9730 & 0.0270   \tabularnewline
2aiii & 100 & 1 & 0.1 & 0.5450  & 0.4550   \tabularnewline
2aiv & 100 & 1 & 0.2 & 0.0000 & 1.0000  \tabularnewline\hline
2bi & 100 & 0.8 & 0 & 0.9450 & 0.0550  \tabularnewline
2bii & 100 & 0.8 & 0.05 & 0.9970  & 0.0030  \tabularnewline
2biii & 100 & 0.8 & 0.1 & 0.9630  & 0.0370   \tabularnewline
2biv & 100 & 0.8 & 0.2 & 0.0000 & 1.0000  \tabularnewline\hline
2ci & 100 & 0.6 & 0 & 0.9490  & 0.0510  \tabularnewline
2cii & 100 & 0.6 & 0.05 & 0.9940  & 0.0060   \tabularnewline
2ciii & 100 & 0.6 & 0.1 & 0.9990 & 0.0001
   \tabularnewline
2civ & 100 & 0.6 & 0.2 & 0.2210 & 0.7790  \tabularnewline\hline
2di & 100 & 0.4 & 0 & 0.9590 & 0.0410  \tabularnewline
2dii & 100 & 0.4 & 0.05 & 0.9900 & 0.0100  \tabularnewline
2diii & 100 & 0.4 & 0.1 & 1.0000  &  0.0000  \tabularnewline
2div & 100 & 0.4 & 0.2 & 0.9980 & 0.0020   \tabularnewline\hline
\end{tabular}
\vspace{1em}
\caption{Results when the \emph{observed} data is a Mat\'{e}rn \rom{2} process, with independent mark being exponential with rate $1$. The semi-major axis length along the $x$-axis, $a$, and $y$-axis, $b$, are equivalent and the semi-major axis length along the $z$-axis is determined such that the area of the ellipsoid is $4\pi$. Fixing the expectation, $\mu$, and hard-core distance, $R$, we use Equation \ref{mu:matern:2} to calculate $\rho$ for the underlying constant Poisson process intensity function. When $R=0$ a Mat\'{e}rn \rom{2} process collapses to a CSR process.}
\label{Regular:table:results}
\end{table}

\begin{table}[p]
\centering
\begin{tabular}{x{0.14\textwidth}x{0.14\textwidth}x{0.14\textwidth}x{0.14\textwidth}x{0.14\textwidth}x{0.14\textwidth}}
\hline
Experiment No. & Expectation & $a$ & $\kappa$ & Accept $H_0$ & Reject $H_0$ \tabularnewline\hline
3ai & 150 & 1 & $\infty$  & 0.9560  & 0.0440 \tabularnewline
3aii & 150 & 1 & 5 & 0.9530  & 0.0470
   \tabularnewline
3aiii & 150 & 1 & 1 & 0.4370  & 0.5630   \tabularnewline
3aiv & 150 & 1 & 0.5 & 0.0170& 0.9830  \tabularnewline\hline
3bi & 150 & 0.8 & $\infty$ & 0.9460 & 0.0540  \tabularnewline
3bii & 150 & 0.8 & 5 &  0.9430 & 0.0570   \tabularnewline
3biii & 150 & 0.8 & 1   & 0.7880 &  0.2120  \tabularnewline
3biv & 150 & 0.8 & 0.5 &  0.0660 & 0.9340  \tabularnewline\hline
3ci & 150 & 0.6 & $\infty$ & 0.9540 & 0.0460   \tabularnewline
3cii & 150 & 0.6 & 5 & 0.9390 & 0.0610   \tabularnewline
3ciii & 150 & 0.6 & 1 &  0.8600 & 0.1400   \tabularnewline
3civ & 150 & 0.6 & 0.5 & 0.2200 & 0.7800   \tabularnewline\hline
3di & 150 & 0.4 &  $\infty$ & 0.9400 &  0.0600 \tabularnewline
3dii & 150 & 0.4 & 5 & 0.9640 & 0.0360  \tabularnewline
3diii & 150 & 0.4 & 1 &  0.7980 &  0.2020  \tabularnewline
3div & 150 & 0.4 & 0.5 & 0.3650 &  0.6350 \tabularnewline\hline
\end{tabular}
\vspace{1em}
\caption{Results when the \emph{observed} data is an ellipsoidal Thomas process. The expected number of offspring per parent is $\lambda=20$ and the underlying Poisson parent process has constant intensity function $\rho=\mu/(4\pi\lambda)$, where $\mu$ is the expectation. The semi-major axis length along the $x$-axis, $a$, and $y$-axis, $b$, are equivalent and the semi-major axis length along the $z$-axis is determined such that the area of the ellipsoid is $4\pi$. When $\kappa=\infty$ an ellipsoidal Thomas process collapses to a CSR process.}
\label{Cluster:table:results}
\end{table}

Tables \ref{CSR:table:results}, \ref{Regular:table:results} and \ref{Cluster:table:results} outline the parameter selection and results of our simulations. By the nature of Monte Carlo simulations the CSR results given in Table \ref{CSR:table:results} are as to be expected with an empirical rejection rate close to 0.05. Expectedly, we see that for the same ellipsoid i.e. $a$ kept constant, that when the Mat\'{e}rn \rom{2} parameter $R$ increases and the Thomas process parameter $\kappa$ decreases (each representing an increased departure from CSR), the power of our test improves. In the Appendix F.3 we discuss a potential reason for the power of our test decreasing as $a$ decreases (hence $c$ increases), for both regular and cluster processes, for the same $R$ and $\kappa$ respectively. Additionally, Figure \ref{fig:sub2} suggests we may gain power by considering a two sided test.

\section{Discussion \& Conclusion}\label{Section:10}

In this work we have discussed point patterns observed on arbitrary, bounded convex shapes in $\mathbb{R}^3$, motivated by the need for such exploratory analyses in the area of microbiology. We have highlighted the challenge of handling such spaces due to the lack of isometries for such objects. Using the invariance of Poisson processes \cite{Kingman1993}, we can circumvent this lack of isometries in the original space by mapping to the sphere which has rotational symmetries. By doing so we propose a set of functional summary statistics for the class of Poisson processes. Further to this we have also proposed functional summary statistics for CSR processes on the convex space and explored their properties. Using this we have, in turn, been able to construct test statistics which can be used to reject the hypothesis of CSR for observed point patterns. We have also conducted simulation studies to investigate the effectiveness off the proposed test statistics in rejecting the null hypothesis when the observed data is either regular or clustered.

Interesting extensions to this work would include relaxing the need for convexity of the shape of interest. This presents a significant challenge as how one constructs the required mapping is not obvious. Another consideration is how to construct an estimator of the intensity function on $\mathbb{D}$. One approach might be to construct it on $\mathbb{S}^2$ and inverse map to $\mathbb{D}$. There is, of course, the open question of how one forms summary statistics for multivariate point processes on convex shapes. Answering this would have immediate impact in bioimaging applications where experimentalists are regularly interested in spatial dependencies that exist between two or more different types of molecules.

\section*{Acknowledgements}
Scott Ward is funded by a Wellcome Trust grant (grant number: 210298/Z/18/Z).

\clearpage
\bibliographystyle{unsrt}  
\bibliography{ref3}

\end{document}


\maketitle

\begin{adjustwidth}{1cm}{1cm}
\begin{center}
\textsc{Overview}
\end{center}
\begin{description}
  \setlength\itemsep{0.1cm}
  \item[\ref{Supplementary:H:F:infinite:series}] An infinite series representation for the $F$-, and $H$-functions similar to that of \cite{White1979}.
  \item[\ref{Supplemetary:general:map:theorem}] Poissonial invariance for Poisson processes on a convex shape in $\mathbb{R}^3$ mapped to $\mathbb{S}^2$.
  \item[\ref{Supplementary:FS:means}] Derivation of the expectations of the functional summary statistics: $F_{\text{inhom}}$-, $H_{\text{inhom}}$-, and $K_{\text{inhom}}$-functions.
  \item[\ref{Supplementary:FS:variances}] Derivation of the variances of the functional summary statistics: $F_{\text{inhom}}$-, $H_{\text{inhom}}$-, and $K_{\text{inhom}}$-functions.
  \item[\ref{Supplementary:J:moments}] Existence and approximations to the first and second order moments of the $\hat{J}_{inhom}$-function are given.
  \item[\ref{Supplementary:reg:cluster}] Discussion of regular and cluster processes on convex three dimensional shapes.
  \item[\ref{Supplementary:K:unknown:rho}] Derivation of the moments for the estimator of the $K_{\text{inhom}}$-function when $\rho$ is unknown.
\end{description}
\end{adjustwidth}

\appendix

\section{Infinite series representation of the $F$-, and $H$-functions}
\label{Supplementary:H:F:infinite:series}
In this section we derive an infinite series representation for the $F$-, and $H$-function. \cite{White1979} was the first to derive these in $\mathbb{R}^{d}$. We provide a derivation on $\mathbb{S}^2$ assuming the existence of the $n^{th}$-order factorial moment measures for all $n\in\mathbb{N}$. We also provide an identity when the $n^{th}$-order product densities do exists for all $n\in\mathbb{N}$. Before giving the results we first introduce the $n^{th}$-order factorial moment measures, $\alpha^{(n)}:\mathbb{S}^2\times\cdots\times\mathbb{S}^2\mapsto \mathbb{R}$, as,
\begin{equation*}
\alpha^{(n)}(B_1\times\cdots\times B_n)=\mathbb{E}\sum_{\mathbf{x}_1,\dots,\mathbf{x}_n\in X}^{\neq}\mathbbm{1}[\mathbf{x}_1\in B_1,\dots,\mathbf{x}_n\in B_n],
\end{equation*}
for $B_i\subseteq\mathbb{S}^2,\;i=1,\dots,n$. Then if $\alpha^{(n)}$ is absolutely continuous with respects to the Lebesgue measure, there then exists $\rho^{(n)}$ such that,
\begin{equation}\label{alpha:rho:relation:app}
\alpha^{(n)}(B_1\times\cdots\times B_n)=\int_{B_1}\cdots\int_{B_n}\rho^{(n)}(\mathbf{x}_1,\dots,\mathbf{x}_n)\lambda_{\mathbb{S}^2}(d\mathbf{x}_1)\cdots \lambda_{\mathbb{S}^2}(d\mathbf{x}_n),
\end{equation}
then $\rho^{(n)}$ is the $n^{th}$-order product intensity function \cite{Coeurjolly2015}.

\begin{theorem}\label{thm:infinite:series:f:h:app}
Let $X$ be an isotropic spheroidal point process with constant intensity function $\rho$. Further we assume the existence of all $n^{th}$-order factorial moment measures for both $X$ and its reduced Palm process, $X^!_{\mathbf{x}}$. Then the $F$- and $H$-functions have the following series representation,
\begin{align*}
F(r)&=-\sum_{n=1}^\infty\frac{(-1)^n}{n!}\alpha^{(n)}(B_{\mathbb{S}^2}(\mathbf{o},r),\dots,B_{\mathbb{S}^2}(\mathbf{o},r))\\
H(r)&=-\sum_{n=1}^\infty\frac{(-1)^n}{n!}\alpha^{!(n)}_{\mathbf{o}}(B_{\mathbb{S}^2}(\mathbf{o},r)\dots,B_{\mathbb{S}^2}(\mathbf{o},r))
\end{align*}
where $\alpha^{(n)}$ and $\alpha^{!(n)}_{\mathbf{x}}$ are the factorial moment measure for $X$ and $X^!_{\mathbf{x}}$ and $B_{\mathbb{S}^2}(\mathbf{o},r)$ is the spherical cap of radius $r$ at the origin $\mathbf{o}\in \mathbb{S}^2$. These representations hold provided the series is absolutely convergent, that is if $\lim_{n\rightarrow\infty}|a_{n+1}/a_n| <1$ or $\limsup_{n\rightarrow\infty} (|a_n|)^{1/n} < 1$, where $a_n=((-1)^n/n!)\alpha^{(n)}(B_{\mathbb{S}^2}(\mathbf{o},r),\dots,B_{\mathbb{S}^2}(\mathbf{o},r))$ for the $F$-function or $a_n=((-1)^n/n!)$ $\alpha^{!(n)}_{\mathbf{o}}(B_{\mathbb{S}^2}(\mathbf{o},r)\dots,B_{\mathbb{S}^2}(\mathbf{o},r))$ for the $H$-function.
\end{theorem}
\begin{proof}
Let $\{B_1,\dots,B_n\}$ be a partition of $B_{\mathbb{S}^2}(\mathbf{o},r)$ into n sets, such that as $n$ increases the area of each $B_i$ decreases. Then,
\begin{align}
F(r) &= 1 - P_{X\cap B_{\mathbb{S}^2}(\mathbf{o},r)}(\emptyset)\nonumber\\
&=1-P((X_{B_1}=\emptyset)\cap\cdots\cap(X_{B_n}=\emptyset))\nonumber\\
&=1-\left(1-\sum_{i=1}^n P(X_{B_i}\neq\emptyset) + \sum_{\substack{i,j=1\\i<j}}^n P((X_{B_i}\neq\emptyset)\cap(X_{B_j}\neq\emptyset)) - \cdots \right)\label{inclusion:exclusion}\\
&=\sum_{i=1}^n P(X_{B_i}\neq\emptyset) - \sum_{\substack{i,j=1\\i<j}}^n P((X_{B_i}\neq\emptyset)\cap(X_{B_j}\neq\emptyset)) + \cdots, \label{in:ex:eq}
\end{align}
where the (\ref{inclusion:exclusion}) follows from the inclusion-exclusion principle. Next we shall consider the first term. Define $a_n=\sum_{i=1}^n\mathbbm{1}[X_{B_i}\neq\emptyset]$ and $a=\sum_{\mathbf{x}\in X}\mathbbm{1}[\mathbf{x}\in B_{\mathbb{S}^2}(\mathbf{o},r)]$. Then it can easily be seen that $a_n$ is a monotonically increasing sequence since as the number of partitions increases the number of partitions containing more than one point of $X$ decreases. Further the maxmimum of $a_n$ is $a$. To see this consider small neighbourhoods of each $\mathbf{x}\in X$ such that these neighbourhoods are all disjoint (this is possible since $X$ is a simple process), label these $B_{\mathbf{x}}$, then $\mathbbm{1}[X_{B_{\mathbf{x}}}\neq\emptyset]=1,\; \forall
\mathbf{x}\in X$ and so for this partition $a_n=a$ and cannot increase as this would require at least one point of $X$ to be in two different disjoint $B_{\mathbf{x}}$: a contradiction. Hence as $n\rightarrow \infty$, $a_n\rightarrow a$. Therefore,
\begin{align*}
\sum_{i=1}^n P(X_{B_i}\neq \emptyset) &= \mathbb{E}\sum_{i=1}^n\mathbbm{1}[X_{B_i}\neq \emptyset]\\
&=\int_{N_{lf}}\sum_{i=1}^n\mathbbm{1}[x\neq \emptyset] dP(x).
\end{align*}
Since $a_n\leq a$ and by assumption $\alpha^{(1)}\equiv\alpha$ exists $\mathbb{E}[a]=\mathbb{E}\sum_{\mathbf{x}\in X}\mathbbm{1}[\mathbf{x}\in B_{\mathbb{S}^2}(\mathbf{o},r)]=\alpha(B_{\mathbb{S}^2}(\mathbf{o},r))<\infty$, we can therefore apply the dominated convergence theorem when taking the limit as $n$ increases and the volumes of $B_i$ decrease, i.e.
\begin{align*}
\lim_{ \substack{n\rightarrow\infty\\ |B_i|\rightarrow 0}}\sum_{i=1}^n P(X_{B_i}\neq \emptyset) &= \lim_{ \substack{n\rightarrow\infty\\ |B_i|\rightarrow 0}}\int_{N_{lf}}\sum_{i=1}^n\mathbbm{1}[x\neq \emptyset] dP(x)\\
&= \int_{N_{lf}}\lim_{ \substack{n\rightarrow\infty\\ |B_i|\rightarrow 0}} \sum_{i=1}^n\mathbbm{1}[x\neq \emptyset] dP(x)\\
&= \int_{N_{lf}} \sum_{\mathbf{y}\in x}\mathbbm{1}[y\in B_{\mathbb{S}^2}(\mathbf{o},r)] dP(x)\\
&= \mathbb{E} \sum_{\mathbf{x}\in X}\mathbbm{1}[x\in B_{\mathbb{S}^2}(\mathbf{o},r)]\\
&=\alpha(B_{\mathbb{S}^2}(\mathbf{o},r))
\end{align*}
An identical approach can be used for the remaining terms of \ref{in:ex:eq}. By considering the $k^{th}$ term of \ref{in:ex:eq} we have,
\begin{equation*}
\begin{split}
(-1)^{k+1}\sum_{\substack{i_1,\dots,i_k=1\\i_1<\cdots<i_k}}^n &P((X_{B_{i_1}}\neq\emptyset)\cap\cdots\cap(X_{B_{i_k}}\neq\emptyset))\\
&=\frac{(-1)^{k+1}}{k!}\sum_{i_1=1}^n\cdots\sum_{\substack{i_k=1\\i_k\notin\{i_1,\dots,i_{k-1}\}}}^nP((X_{B_{i_1}}\neq\emptyset)\cap\cdots\cap(X_{B_{i_k}}\neq\emptyset)).
\end{split}
\end{equation*}
By defining $a_n=\sum_{i_1,\dots,i_k\in\{1,\dots,n\}}^{\neq} \mathbbm{1}[X_{B_{i_1}}=\emptyset,\dots,X_{B_{i_k}}=\emptyset]$ and $a=\sum_{\mathbf{x}_1,\dots,\mathbf{x}_k\in X}^{\neq}\mathbbm{1}[\mathbf{x}_1\in B_{\mathbb{S}^2}(\mathbf{o},r),\dots,\mathbf{x}_k\in B_{\mathbb{S}^2}(\mathbf{o},r)]$ identical arguments can be made as in the case for $k=1$ giving,
\begin{equation*}
\lim_{ \substack{n\rightarrow\infty\\ |B_i|\rightarrow 0}} (-1)^{k+1}\sum_{\substack{i_1,\dots,i_k=1\\i_1<\cdots<i_k}}^n P((X_{B_{i_1}}\neq\emptyset)\cap\cdots\cap(X_{B_{i_k}}\neq\emptyset))=\frac{-(-1)^{k}}{k!}\alpha^{(k)}(B_{\mathbb{S}^2}(\mathbf{o},r),\dots,B_{\mathbb{S}^2}(\mathbf{o},r)),
\end{equation*}
and so gives the series for the $F$-function.

The series representation for the $H$-function follows an identical argument to that of the $F$-function, instead using the factorial moment measure for the reduced Palm point process, $X^!_{\mathbf{x}}$.
\end{proof}

The following corollary shows that the infinite series for the $F$-, and $H$-functions can also be represented using the $n^{th}$-order product densities, used by \cite{vanLieshout2011}.

\begin{corollary}
\cite{White1979} Under the same assumptions as Theorem \ref{thm:infinite:series:f:h:app}, let $X$ be an isotropic spheroidal point process with constant intensity function $\rho$. Further we assume the existence of all $n^{th}$-order product intensities for both $X$ and its reduced Palm process, $X^!_{\mathbf{x}}$. Then the $F$- and $H$-functions have the following series representation,
\begin{align*}
F(r)&=-\sum_{n=1}^\infty\frac{(-1)^n}{n!}\int_{B_{\mathbb{S}^2}(\mathbf{o},r)}\cdots\int_{B_{\mathbb{S}^2}(\mathbf{o},r)}\rho^{(n)}(\mathbf{x}_1,\dots,\mathbf{x}_n)\lambda_{\mathbb{S}^2}(d\mathbf{x}_1)\cdots \lambda_{\mathbb{S}^2}(d\mathbf{x}_n)\\
H(r)&=-\sum_{n=1}^\infty\frac{(-1)^n}{n!}\int_{B_{\mathbb{S}^2}(\mathbf{o},r)}\cdots\int_{B_{\mathbb{S}^2}(\mathbf{o},r)}\frac{\rho^{(n+1)}(\mathbf{o},\mathbf{x}_1,\dots,\mathbf{x}_n)}{\rho}\lambda_{\mathbb{S}^2}(d\mathbf{x}_1)\cdots \lambda_{\mathbb{S}^2}(d\mathbf{x}_n)
\end{align*}
provided the series is absolutely convergent, where $B_{\mathbb{S}^2}(\mathbf{o},r)$ is the spherical cap of radius $r$ at the origin $\mathbf{o}\in \mathbb{S}^2$.
\end{corollary}
\begin{proof}
The infinite series for $F$ follows immediately from Theorem \ref{thm:infinite:series:f:h:app} and  Equation \ref{alpha:rho:relation:app}. Further, by assumption we also know that,
\begin{equation*}
\alpha^{!(n)}_{\mathbf{o}}(B_1\times\cdots\times B_n)=\int_{B_1}\cdots\int_{B_n}\rho^{!(n)}_{\mathbf{o}}(\mathbf{x}_1,\dots,\mathbf{x}_n)\lambda_{\mathbb{S}^2}(d\mathbf{x}_1)\cdots \lambda_{\mathbb{S}^2}(d\mathbf{x}_n).
\end{equation*}
Then we have the following relationship between the $n^{th}$-order product intensities of $X$ and $X^!_{\mathbf{o}}$ (see for example \cite{Coeurjolly2015}),
\begin{equation*}
\rho^{!(n)}_{\mathbf{x}}(\mathbf{x}_1,\dots,\mathbf{x}_m)=\frac{\rho^{(n+1)}(\mathbf{x},\mathbf{x}_1,\dots,\mathbf{x}_n)}{\rho(\mathbf{x})},
\end{equation*}
and so for our point process $X$,
\begin{equation*}
\alpha^{!(n)}_{\mathbf{o}}(B_1\times\cdots\times B_n)=\int_{B_1}\cdots\int_{B_n}\frac{\rho^{(n+1)}_{\mathbf{o}}(\mathbf{o},\mathbf{x}_1,\dots,\mathbf{x}_n)}{\rho}\lambda_{\mathbb{S}^2}(d\mathbf{x}_1)\cdots \lambda_{\mathbb{S}^2}(d\mathbf{x}_n),
\end{equation*}
and then by Theorem \ref{thm:infinite:series:f:h:app} we have the infinite series for the $H$-function in terms of the $n^{th}$-order product intensities.
\end{proof}

\section{Invariance of Poisson process between metric spaces}\label{Supplemetary:general:map:theorem}

In this section we show that a Poisson process lying on an arbitrary bounded convex space $\mathbb{D}\subset \mathbb{R}^3$ can be mapped to another Poisson process on a sphere, known as the Mapping Theorem \cite{Kingman1993}. We also show that no two different Poisson processes on $\mathbb{D}$ map to the same Poisson process on the sphere under the same mapping. 

Before starting the main theorem of this section we first introduce a lemma (see for example \cite[Proposition 3.1, pp 15-16]{Moller2004}) which is an expansion of the probability measure for a Poisson process on an arbitrary metric space. We shall state it in the context of an arbitrary convex shape in $\mathbb{R}^3$, represented by $\mathbb{D}$.

\begin{customlemma}{S1}\label{Poisson:expansion:lemma}
(M\o ller et al. \cite{Moller2004}) $X$ is a Poisson process with intensity function $\rho:\mathbb{D}\mapsto\mathbb{R}$ on $\mathbb{D}$ if and only if for all $B\subseteq\mathbb{D}$ with $\mu(B)=\int_B\rho(\mathbf{x})d\mathbf{x}<\infty$ and all $F\subseteq N_{\text{lf}}$,
\begin{equation}\label{Poisson:expansion}
P(X_B\in F) = \sum_{n=0}^\infty\frac{\exp(-\mu(B))}{n!}\int_B\cdots\int_B\mathbbm{1}[\{\mathbf{x}_1,\dots,\mathbf{x}_n\}\in F]\prod_{i=1}^n\rho(\mathbf{x}_i)\lambda_{\mathbb{D}}(d\mathbf{x}_1)\cdots\lambda_{\mathbb{D}}(d\mathbf{x}_n),
\end{equation} 
where the integral for $n=0$ is read as $\mathbbm{1}[\emptyset\in F]$.
\end{customlemma}
\begin{proof}
See Proposition 3.1 of \cite{Moller2004}.
\end{proof}

The following lemma shows that the function $f(\mathbf{x})=\mathbf{x}/||\mathbf{x}||$ is bijective and hence has a well-defined inverse.

\begin{lemmaapp}\label{lemma:f:bijective:app}
Let $\mathbb{D}$ be a convex subspace of $\mathbb{R}^3$ such that the origin in $\mathbb{R}^3$ is in the interior of $\mathbb{D}$, i.e. $\mathbf{o}\in\mathbb{D}_{int}$. Then the function $f(\mathbf{x})=\mathbf{x}/||\mathbf{x}||, f:\mathbb{D}\mapsto\mathbb{S}^2$ is bijective.
\end{lemmaapp}   
\begin{proof}
To show that $f$ is bijective we need to show that it is both injective and surjective. For surjectivity we need to show that for any $\mathbf{x}'\in\mathbb{S}^2$, there exists $\mathbf{x}\in\mathbb{D}$ such that $f(\mathbf{x})=\mathbf{x}'$. To do this first fix any $\mathbf{x}'\in\mathbb{S}^2$ and define the half line, $L_{\mathbf{x}'}=\{\mathbf{y}\in\mathbb{R}^3:\mathbf{y}=\lambda\mathbf{x}',\;\lambda\in\mathbb{R}^+\}$, where $\mathbb{R}^+$ is the positive real line including $0$. Then since $\mathbb{D}$ is compact (i.e. closed and bounded) then the half line must intersect the $\mathbb{D}$ and so there exists $\mathbf{x}\in\mathbb{D}$ such that $\mathbf{x}\in L_{\mathbf{x}'}$. Therefore there must exists $\lambda\in\mathbb{R}^+, \mathbf{x}=\lambda\mathbf{x}'$. Taking norms of both sides and noting that since $\mathbf{x}'\in\mathbb{S}^2$ meaning $||\mathbf{x}'||=1$ then $\lambda=||\mathbf{x}||$ and so, $\mathbf{x}/||\mathbf{x}||=\mathbf{x}'$ and so $f$ is surjective.

For injectivity we need to show that for any $\mathbf{x},\mathbf{y}\in\mathbb{D}$ if $f(\mathbf{x})=f(\mathbf{y})$ then $\mathbf{x}=\mathbf{y}$. Fix $\mathbf{x},\mathbf{y}\in\mathbb{D}$ such that $f(\mathbf{x})=f(\mathbf{y})$ and define $\mathbf{x}'=f(\mathbf{x})=f(\mathbf{y})$. Again define the line $L_{\mathbf{x}'}$ as previous and by convexity of $\mathbb{D}$, the fact that $\mathbf{o}\in\mathbb{D}$ and since $L_{\mathbf{x}'}$ is a half line then there is precisely only one intersection of $L_{\mathbf{x}'}$ with $\mathbb{D}$. Therefore both $\mathbf{x}$ and $\mathbf{y}$ must be this point of intersection meaning $\mathbf{x}=\mathbf{y}$. Hence $f$ is bijective.
\end{proof}

Now we give the main theorem of our work which shows that a Poisson process on $\mathbb{D}$ is mapped to a Poisson process on $\mathbb{S}^2$. This is known as the Mapping Theorem \cite{Kingman1993}. Here we use Lemma \ref{Poisson:expansion:lemma}.

Before beginning this theorem we lay down a little notation in order to avoid confusion. $\mathbf{x}=(x,y,z)$ will be an element of $\mathbb{R}^3$ where we may subscript with an $n\in\mathbb{N}$ when we are referring to a single vector within a set. $x$ will be an element of $N_{lf}$ and may also be subscripted with $n\in\mathbb{N}$ when we are referring to a single element in a set of finite point configurations. Notice that we are using $x$ to be both the first element of $\mathbf{x}$ and an element in $N_{lf}$, based on context it will be clear to which we are referring too.

\begin{theorem}\label{thm:mapping:general:app}
Let $X$ be a Poisson process on an arbitrary bounded convex shape $\mathbb{D}\subset\mathbb{R}^3$ with intensity function $\rho:\mathbb{D}\mapsto \mathbb{R}$. We assume that $\mathbb{D}=\{\mathbf{x}\in\mathbb{R}^3: g(\mathbf{x})=0\}$ where $g(\mathbf{x})=0$ is the zero set function and is defined as,
\begin{equation*}
g(\mathbf{x})=
\left\{\begin{alignedat}{2}
    g_1(\mathbf{x})=0,&\quad \mathbf{x}\in\mathbb{D}_1\\
    &\vdots \\
    g_n(\mathbf{x})=0,&\quad \mathbf{x}\in\mathbb{D}_n
  \end{alignedat}\right.
\end{equation*}
such that $\cup_{i=1}^n\mathbb{D}_i=\mathbb{D}$ and $\mathbb{D}_i\cap\mathbb{D}_j=\emptyset,\; \forall i\neq j$. Then define $Y=f(X)$, where $f(\mathbf{x})=\mathbf{x}/||\mathbf{x}||$ and we have taken the convention that $f(X)=\{\mathbf{y}\in\mathbb{S}^2: \mathbf{y}=\mathbf{x}/||\mathbf{x}||, \mathbf{x}\in X\}$. Then $Y$ is a Poisson process on $\mathbb{S}^2$, with intensity function,
\begin{equation*}
\rho^*(\mathbf{x})=
\left\{\begin{alignedat}{2}
\rho(f^{-1}(\mathbf{x}))l_1(f^{-1}(\mathbf{x})) &J_{(1,f^*)}(\mathbf{x})\sqrt{1-x^2-y^2}, \quad \mathbf{x}\in f(\mathbb{D}_1)\\
&\vdots\\
\rho(f^{-1}(\mathbf{x}))l_n(f^{-1}(\mathbf{x})) &J_{(n,f^*)}(\mathbf{x})\sqrt{1-x^2-y^2}, \quad \mathbf{x}\in f(\mathbb{D}_n)
  \end{alignedat}\right.
\end{equation*}
where,
\begin{align*}
z&=\tilde{g}_i(x,y)\\
l_i(\mathbf{x}) &=\left[1+\left(\frac{\partial \tilde{g}_i}{\partial x}\right)^2+\left(\frac{\partial \tilde{g}_i}{\partial  y}\right)^2\right]^{\frac{1}{2}}\\
J_{(i,f^{*-1})}(\mathbf{x})&=\frac{1}{( x^2+ y^2+\tilde{g}_i^2( x, y))^3}\\
&\det\left[
\begin{pmatrix}
 y^2+\tilde{g}_i^2( x, y)- x\tilde{g}_i( x, y)\frac{\partial \tilde{g}_i}{\partial  x} & - x\left( y+\tilde{g}_i( x, y)\frac{\partial \tilde{g}_i}{\partial  y}\right) \\
- y\left( x+\tilde{g}_i( x, y)\frac{\partial \tilde{g}_i}{\partial  x}\right) &  x^2+\tilde{g}_i^2( x, y)- y\tilde{g}_i( x, y)\frac{\partial \tilde{g}_i}{\partial  y}
\end{pmatrix}\right]\\
J_{(i,f^*)}(\boldsymbol{\mathbf{x}}) &= \frac{1}{J_{(i,f^{*-1})}(f^{-1}(\boldsymbol{\mathbf{x}}))},
\end{align*}
where $\mathbf{x}=(x,y,z)^T$, $f^{-1}$ is the inverse of $f$, $\det(\cdot)$ is the determinant operator, and $f^*:\mathbb{R}^2\mapsto\mathbb{R}^2$ is the function which maps $ x\mapsto x/||\mathbf{x}||$ and $ y\mapsto  y/||\mathbf{x}||$.
\end{theorem}

\begin{proof}
In order to show that $Y\equiv f(X)$ is a Poisson process we show that its distribution function can be expanded as given by Equation \ref{Poisson:expansion}. Then $\forall B \subseteq \mathbb{D}$ and $\forall F \subseteq N_{lf}$ and noting that $f$ is a measurable map (since the map is bijective and hence an inverse exists) we have that,
\begin{align*}
P(Y_B\in F) &= P(f^{-1}(Y_B)\in f^{-1}(F))\\
&= P(X_{f^{-1}(B)}\in f^{-1}(F))
\end{align*}
Now we define $f^{-1}_{i}(F),$ for $i=1,\dots,n$ where $F\subseteq N_{lf}$,
\begin{align*}
F&=&&\{x\in N_{lf}: x=\{\mathbf{x}_1,\dots,\mathbf{x}_m\}, m\in\mathbb{N}, x\in F\}\\
f^{-1}(F) &= &&\{f^{-1}(x): f^{-1}(x)=\{f^{-1}(\mathbf{x}_1),\dots,f^{-1}(\mathbf{x}_m)\},m\in\mathbb{N}, x\in F\}\\
\intertext{we want to partition $f^{-1}(F)$ over each $f^{-1}(\mathbb{D}_i), i = 1,\dots,n$}
f^{-1}(F) &= &&\left\{\cup_{i=1}^n f^{-1}(x_i): f^{-1}(x_i)=\{f^{-1}(\mathbf{x}_{(i,1)}),\dots,f^{-1}(\mathbf{x}_{(i,m_i)})\},\right.\\
& && \left.f^{-1}(\mathbf{x}_{(i,j)})\in f^{-1}(\mathbb{D}_i), j=1,\dots,m_i,\;m_i\in\mathbb{N},\;i=1,\dots,n ,\cup_{i=1}^n x_i\in F\right\}.
\end{align*}
To understand the notation $\mathbf{x}_{(i,j)}$ consider first a single element $x\in F$. Then since F is a subset of $N_{lf}$ this means that $|x|\in\mathbb{N}$. Then define $m_i=|x\cap \mathbb{D}_i|,$ hence $\sum_{i=1}^n m_i=|x|$. Then $\mathbf{x}_{(i,j)}$ is $j^{th}$ element of $x\cap\mathbb{D}_i$ such that $j=1,\dots,m_i$. We define $f^{-1}_{i}(F)\equiv \{f^{-1}(x): f^{-1}(x)\equiv\{f^{-1}(\mathbf{x}_1),\dots,f^{-1}(\mathbf{x}_n)\},\;f^{-1}(\mathbf{x}_i)\in f^{-1}(\mathbb{D}_i),n\in\mathbb{N}, \exists y\in F \text{ such that } x\subseteq y\}$. Then,
\begin{align*}
P(X_{f^{-1}(B)}\in f^{-1}(F))&= P(\{X_{f^{-1}(B)\cap \mathbb{D}_1},\dots,X_{f^{-1}(B)\cap \mathbb{D}_n}\}\in f^{-1}(F))\\
&= P(X_{f^{-1}(B)\cap \mathbb{D}_1}\in f^{-1}_1(F),\dots,X_{f^{-1}(B)\cap \mathbb{D}_n}\in f^{-1}_n(F))\\
&= P(X_{f_1^{-1}(B)}\in f_1^{-1}(F),\dots,X_{f^{-1}_n(B)}\in f_n^{-1}(F))\\
&= \prod_{i=1}^n P(X_{f_i^{-1}(B)}\in f_i^{-1}(F))
\end{align*}
where $f^{-1}_i(A)=\{(x,y,z)^T\in A: (x,y,z)^T\in\mathbb{D}_i\}$ if $A\subseteq\mathbb{D}$. We emphasize the dual meaning of $f^{-1}_i(A)$ where the defintion depends on the nature of $A$, i.e. if $A\subseteq\mathbb{D}$ or if $A\subseteq N_{lf}$. Without loss of generality lets suppose that all projections of each $\mathbb{D}_i$ onto $\mathbb{R}^2$ are invertible, if not then we can divide $\mathbb{D}_i$ into further subsets $\cup_{j=1}^m \mathbb{D}_{i,j}$ such that the projection of each $\mathbb{D}_{i,j}$ is then invertible. For example an ellipsoid is defined by the zero-set equation $x^2/a^2+y^2/b^2+z^2/c^2=1$, but if the entire space were projected down to $\mathbb{R}^2$ its inverse does not exists, instead we divide the ellipsoid into the upper and lower hemiellipsoids and then the projections restricted to these segments of the ellipsoids are then invertible. Let us also define the projection of $\mathbb{D}_i$ to $\mathbb{R}^2$ as $P_{\mathbb{D}_i}$. Further, for $\mathbf{x}=(x,y,z)^T$, then $\mathbf{x}$ lies on $\mathbb{D}$ if $g(\mathbf{x})=0$, we define $\tilde{g}$ to be the rearrangement of $g$ such that $z$ is a function of $x$ and $y$, i.e. $z=\tilde{g}(x,y)$. Then by Lemma \ref{Poisson:expansion:lemma},
\begin{align*}
&P(X_{f_i^{-1}(B)}\in f_i^{-1}(F))\\
&= \sum_{n=0}^\infty\frac{\exp(-\mu(f^{-1}_i(B)))}{n!}\\
&\phantom{AAAA}\int_{f^{-1}_i(B)}\cdots\int_{f^{-1}_i(B)}\mathbbm{1}[\{\mathbf{x}_1,\dots,\mathbf{x}_n\}\in f^{-1}_i(F)]\prod_{i=1}^n \rho(\mathbf{x}_i) \lambda_{\mathbb{D}_i}(d\mathbf{x}_1)\cdots \lambda_{\mathbb{D}_i}(d\mathbf{x}_n)\\
&= \sum_{n=0}^\infty\frac{\exp(-\mu(f^{-1}_i(B))}{n!}\\
&\phantom{AAAA}\iint\limits_{P_{\mathbb{D}_i}[f^{-1}_i(B)]}\cdots\iint\limits_{P_{\mathbb{D}_i}[f^{-1}_i(B)]}\mathbbm{1}[\{(x_1,y_1,\tilde{g}_1(x_1,y_1))^T,\dots,(x_n,y_n,\tilde{g}_1(x_n,y_n))^T\}\in f^{-1}_i(F)]\\
&\phantom{AAAAAAAA}\prod_{i=1}^n\rho\left(x_i,y_i,\tilde{g}_i(x_i,y_i)\right)\sqrt{1+\left(\frac{\partial\tilde{g}_i}{\partial x_i}\right)^2+\left(\frac{\partial \tilde{g}_i}{\partial y_i}\right)^2}dx_idy_i\\
&= \sum_{n=0}^\infty\frac{\exp(-\mu(f^{-1}_i(B))}{n!}\\
&\phantom{AAAA}\iint\limits_{P_{\mathbb{D}_i}[f^{-1}_i(B)]}\cdots\iint\limits_{P_{\mathbb{D}_i}[f^{-1}_i(B)]}\mathbbm{1}[\{(x_1,y_1,\tilde{g}_1(x_1,y_1))^T,\dots,(x_n,y_n,\tilde{g}_1(x_n,y_n))^T\}\in f^{-1}_i(F)]\\
&\phantom{AAAAAAAA}\prod_{i=1}^n\rho\left(x_i,y_i,\tilde{g}_i(x_i,y_i)\right)l_i(x_i,y_i)dx_idy_i,
\end{align*}
where $l_i(x_i,y_i)=\sqrt{1+\left(\partial\tilde{g}_i/\partial x_i \right)^2+\left(\partial \tilde{g}_i/\partial y_i\right)^2}$. Now consider the indicator term, we need to show that when $x_i\mapsto x_i/||\mathbf{x}||$ and $y_i\mapsto y_i/||\mathbf{x}||$ then $\mathbbm{1}[\{\mathbf{x}_1,\dots,\mathbf{x}_n\}\in f^{-1}_i(F)]\mapsto\mathbbm{1}[\{\mathbf{y}_1,\dots,\mathbf{y}_n\}\in F_i]$, where $F_i=\{x: x=\{\mathbf{x}_1,\dots,\mathbf{x}_m\}, \mathbf{x}_j\in f(\mathbb{D}_i), j=1,\dots,m, m\in\mathbb{N}, \exists y \in F \text{ such that } x\subseteq y\}$. Let us consider first an individual point $\mathbf{x}\in\mathbb{D}$. Further let us define $r=|\mathbf{x}|=\sqrt{x^2+y^2+z^2}$. Then since $z=\tilde{g}_i(x,y),$ $r$ is thus a function of $x$ and $y$, let us write $r(x,y)=||\mathbf{x}||$. Thus we can rewrite $z=\tilde{g}_i(x,y)=\sqrt{r^2(x,y)-x^2-y^2}$. Then suppose we apply the transformations $x'=x/r(x,y)$ and $y'=y/r(x,y)$, we have  $z=\sqrt{r^2(x,y)-r^2(x,y)x'^2-r^2(x,y)y'^2}\Rightarrow z=r(x,y)\sqrt{1-x'^2-y'^2}$. Therefore,
\begin{align*}
&\mathbbm{1}[\{\mathbf{x}_1,\dots,\mathbf{x}_n\}\in f^{-1}_i(F)]=\mathbbm{1}[\{(x_1,y_1,\tilde{g}_1(x_1,y_1))^T,\dots,(x_n,y_n,\tilde{g}_1(x_n,y_n))^T\}\in f^{-1}_i(F)]\\
\intertext{apply transformations $x'=x/r(x,y)$ and $y'=y/r(x,y)$,}
&=\mathbbm{1}\left[\left\{\left(r(x_1,y_1)x'_1,r(x_1,y_1)y'_1,r(x_1,y_1)\sqrt{1-x_1'^2-y_1'^2}\right)^T,\dots,\right.\right.\\
&\phantom{AAAA}\left.\left.\left(r(x_n,y_n)x'_n,r(x_n,y_n)y'_n,r(x_n,y_n)\sqrt{1-x_n'^2-y_n'^2}\right)^T\right\}\in f^{-1}_i(F)\right]\\
&=\mathbbm{1}\left[\left\{\left(x'_1,y'_1,\sqrt{1-x_1'^2-y_1'^2}\right)^T,\dots,\left(x'_n,y'_n,\sqrt{1-x_n'^2-y_n'^2}\right)^T\right\}\in F_i\right]\\
&=\mathbbm{1}[\{\mathbf{y}_1,\dots,\mathbf{y}_n\}\in F_i]
\end{align*}
Let us return to $P(X_{f_i^{-1}(B)}\in f_i^{-1}(F))$. We apply the transformation $f^*(x,y)=\left(x/r(x,y),\right.$ $\left.y/r(x,y)\right)^T$ and define the inverse as $f^{*-1}(x,y)$. We have,
\begin{align*}
&P(X_{f_i^{-1}(B)}\in f_i^{-1}(F)) \\
&=\sum_{n=0}^\infty\frac{\exp(-\mu(f^{-1}_i(B))}{n!}\iint\limits_{f^*(P_{\mathbb{D}_i}[f^{-1}_i(B)])}\cdots\iint\limits_{f^*(P_{\mathbb{D}_i}[f^{-1}_i(B)])}\mathbbm{1}[\{\mathbf{y},\dots,\mathbf{y}_n\}\in F_i]\\
&\phantom{AAAA}\prod_{i=1}^n\left[\rho\left(f^{*-1}_1(x'_i,y'_i),f^{*-1}_2(x'_i,y'_i),\tilde{g}_i(f^{*-1}_1(x'_i,y'_i),f^{*-1}_2(x'_i,y'_i))\right)\right.\\
&\phantom{AAAAAAAA}\left.l_i(f^{*-1}_1(x'_i,y'_i),f^{*-1}_2(x'_i,y'_i)) J_{(i,f^*)}(x'_i,y'_i)dx'_idy'_i\right],
\end{align*}
where $J_{(i,f^*)}(x'_i,y'_i)$ is the Jacobian of the transformation $f^*$. The Jacobian of the transformation is defined as,
\begin{align*}
J_{(i,f^*)}(\boldsymbol{\mathbf{x}}) &= \frac{1}{J_{(i,f)}(f^{*-1}(\boldsymbol{\mathbf{x}}))}\\
\intertext{where we can use the inverse property to obtain $J_{(i,f^{*-1})}$,}
J_{(i,f^{*-1})}(\mathbf{x})&=\det\left[
\begin{pmatrix}
\frac{\partial x'}{\partial x} & \frac{\partial x'}{\partial y}\\
\frac{\partial y'}{\partial x} & \frac{\partial y'}{\partial y}
\end{pmatrix}\right]
\end{align*}
The entries of $J_{(i,f^{*-1})}(\mathbf{x})$ are given as follows,
\begin{align*}
\frac{\partial x'}{\partial x}&=\frac{y^2+\tilde{g}_i^2(x,y)-x\tilde{g}(x,y)\frac{\partial \tilde{g}}{\partial x}}{r^3(x,y)}\\
\frac{\partial x'}{\partial y}&=\frac{-x\left(y+\tilde{g}(x,y)\frac{\partial \tilde{g}}{\partial y}\right)}{r^3(x,y)}\\
\frac{\partial y'}{\partial x}&=\frac{-y\left(x+\tilde{g}(x,y)\frac{\partial \tilde{g}}{\partial x}\right)}{r^3(x,y)}\\
\frac{\partial y'}{\partial y}&=\frac{x^2+\tilde{g}_i^2(x,y)-y\tilde{g}(x,y)\frac{\partial \tilde{g}}{\partial y}}{r^3(x,y)},
\end{align*}
where $r(x,y)=||\mathbf{x}||,$ and so 
\begin{align*}
J_{(i,f^{*-1})}(\mathbf{x}) &= J_{(i,f^{*-1})}(x,y)\\
&=\frac{1}{r^6(x,y)}\det\left[
\begin{pmatrix}
y^2+\tilde{g}_i^2(x,y)-x\tilde{g}(x,y)\frac{\partial \tilde{g}}{\partial x} & -x\left(y+\tilde{g}(x,y)\frac{\partial \tilde{g}}{\partial y}\right)\\
-y\left(x+\tilde{g}(x,y)\frac{\partial \tilde{g}}{\partial x}\right) & x^2+\tilde{g}_i^2(x,y)-y\tilde{g}(x,y)\frac{\partial \tilde{g}}{\partial y}
\end{pmatrix}
\right].
\end{align*}
Therefore,
\begin{equation*}
J_{(i,f^*)}(\boldsymbol{\mathbf{x}'}) = J_{(i,f^*)}(x',y')   = \frac{1}{J_{(i,f)}(f_1^{*-1}(x',y'),f_2^{*-1}(x',y'))}
\end{equation*}
Projecting onto the sphere,
\begin{align}
&P(X_{f_i^{-1}(B)}\in f_i^{-1}(F))= \sum_{n=0}^\infty\frac{\exp(-\mu(f^{-1}_i(B))}{n!}\nonumber\\
&\phantom{AAAA}\int\limits_{P^{-1}_{\mathbb{S}^2}[f^*(P_{\mathbb{D}_i}[f^{-1}_i(B)])]}\cdots\int\limits_{P^{-1}_{\mathbb{S}^2}[f^*(P_{\mathbb{D}_i}[f^{-1}_i(B)])]}\mathbbm{1}[\{\mathbf{y},\dots,\mathbf{y}_n\}\in F_i]\prod_{i=1}^n\rho^*_i(\mathbf{y}_i)\lambda_{\mathbb{S}^2}(d\mathbf{y}_i),\label{multiplicand}
\end{align}
where $\rho^*_i(\mathbf{y})=\rho(f^{*-1}_1(x)$ $,f^{*-1}_2(y),$ $\tilde{g}_i(f^{*-1}_1(x),$ $f^{*-1}_2(y)))$ $l_i(f^{*-1}_1(x),f^{*-1}_2(y))$ $J_{i,f^{-1}}(x,y)$ $\sqrt{1-x^2-y^2},$ $\mathbf{y}=(x,y,z)^T\in\mathbb{S}^2$. 

We now need to show that $P^{-1}_{\mathbb{S}^2}[f^*(P_{\mathbb{D}_i}[f^{-1}_i(B)])]=B\cap f(\mathbb{D}_i)$. Equivalently, we can show that $f^*(P_{\mathbb{D}_i}[f^{-1}_i(B)])=P_{\mathbb{S}^2}[B\cap f(\mathbb{D}_i)]$. It is easy to see that,
\begin{align*}
f^*(P_{\mathbb{D}_i}[f^{-1}_i(B)])&=\{(x/||\mathbf{x}||,y/||\mathbf{x}||)^2\in\mathbb{R}^2: f(\mathbf{x})\in B\cap f(\mathbb{D}_i), \mathbf{x}=(x,y,z)^T\in\mathbb{D}_i\}\\
P_{\mathbb{S}^2}[B\cap f(\mathbb{D}_i)]&=\{(x,y)\in\mathbb{R}^2:(x,y,z)^T\in B\cap f(\mathbb{D}_i)\}.
\end{align*}{}
Then since $f$ is bijective (see Lemma \ref{lemma:f:bijective:app}) this means that for all $\mathbf{x}\in \mathbb{D},$ there exists $\mathbf{y}\in\mathbb{S}^2$ such that $\mathbf{y}=f(\mathbf{x})$. Further since $\mathbb{D}_i,\;i=1,\dots,n$ partition $\mathbb{D}$ this means that $f(\mathbb{D}_i),\; i=1,\dots,n$ partitions $\mathbb{S}^2$ and so $\mathbf{x}\in f^{-1}(B)\cap\mathbb{D}_i\Rightarrow \mathbf{y}=f(\mathbf{x})\in B\cap f(\mathbb{D}_i)$. Hence taking the set $P_{\mathbb{S}^2}[B\cap f(\mathbb{D}_i)]$,
\begin{align*}
P_{\mathbb{S}^2}[B\cap f(\mathbb{D}_i)]&=\{(x,y)\in\mathbb{R}^2:(x,y,z)^T\in B\cap f(\mathbb{D}_i)\}\\
&=\{(x,y)\in\mathbb{R}^2: (x,y,z)^T=f(\mathbf{x}')\in B\cap f(\mathbb{D}_i), \mathbf{x}'\in\mathbb{D}_i\}\\
&=\{(x,y)\in\mathbb{R}^2: (x,y,z)^T=(x'/||\mathbf{x}'||,y'/||\mathbf{x}'||,z'/||\mathbf{x}'||)^T\in B\cap f(\mathbb{D}_i), \\
&\phantom{AAAA}\mathbf{x}'=(x',y',z')^T\in\mathbb{D}_i\}\\
&=\{(x'/||\mathbf{x}'||,y'/||\mathbf{x}'||)^T\in\mathbb{R}^2: (x'/||\mathbf{x}'||,y'/||\mathbf{x}'||,z'/||\mathbf{x}'||)^T\in B\cap f(\mathbb{D}_i), \\
&\phantom{AAAA}\mathbf{x}'=(x',y',z')^T\in\mathbb{D}_i\}\\
&=f^*(P_{\mathbb{D}_i}[f^{-1}_i(B)]).
\end{align*}
Therefore,
\begin{equation*}
\begin{split}
P(X_{f_i^{-1}(B)}\in f_i^{-1}(F)) &= \sum_{n=0}^\infty\frac{\exp(-\mu(f^{-1}_i(B))}{n!}\\
&\int\limits_{B\cap f(\mathbb{D}_i)}\cdots\int\limits_{B\cap f(\mathbb{D}_i)}\mathbbm{1}[\{\mathbf{y},\dots,\mathbf{y}_n\}\in F_i]\prod_{i=1}^n\rho^*_i(\mathbf{y}_i)\lambda_{\mathbb{S}^2}(d\mathbf{y}_i).
\end{split}
\end{equation*}
With this identity we thus have,
\begin{align*}
P(Y_B\in F) &= \prod_{i=1}^n P(X_{f_i^{-1}(B)}\in f_i^{-1}(F))\\
&= \prod_{i=1}^n \left(\sum_{n=0}^\infty\frac{\exp(-\mu(f^{-1}_i(B))}{n!}\right.\\
&\phantom{AAAA}\left.\int\limits_{B\cap f(\mathbb{D}_i)}\cdots\int\limits_{B\cap f(\mathbb{D}_i)}\mathbbm{1}[\{\mathbf{y},\dots,\mathbf{y}_n\}\in F_i]\prod_{i=1}^n\rho^*_i(\mathbf{y}_i)\lambda_{\mathbb{S}^2}(d\mathbf{y}_i)\right).
\end{align*}
Consider the multiplication of two of the multiplicands $i\neq j$. Define $\rho_{i,j}^*(\mathbf{x})=\rho^*_i(\mathbf{x})$ if $\mathbf{x}\in\mathbb{S}^2\cap f(\mathbb{D}_i)$ and $\rho^*_{i,j}(\mathbf{x})=\rho^*_j(\mathbf{x})$ if $\mathbf{x}\in\mathbb{S}^2\cap f(\mathbb{D}_j)$ and also let $k_{p,n}(\{\mathbf{y}_i\}_{i=1}^n)=\mathbbm{1}[\{\mathbf{y}_1,\dots,\mathbf{y}_n\}\in F_p]\prod_{i=1}^n\rho^*_i(\mathbf{y}_i)$, for $p=i,j$. Then we have,
\begin{align*}
&P(X_{f_i^{-1}(B)}\in f_i^{-1}(F))P(X_{f_j^{-1}(B)}\in f_j^{-1}(F))\\
&=\left(\sum_{n=0}^\infty\frac{\exp(-\mu(f^{-1}_i(B)))}{n!}\int_{B\cap f(\mathbb{D}_i)}\cdots\int_{B\cap f(\mathbb{D}_i)}k_{i,n}(\{\mathbf{y}_i\}_{i=1}^n)\lambda_{\mathbb{S}^2}(d\mathbf{y}_1)\cdots\lambda_{\mathbb{S}^2}(d\mathbf{y}_n) \right)\\
&\phantom{AAAA}\times \left(\sum_{m=0}^\infty\frac{\exp(-\mu(f^{-1}_j(B)))}{m!}\int_{B\cap f(\mathbb{D}_j)}\cdots\int_{B\cap f(\mathbb{D}_j)}k_{j,m}(\{\mathbf{y}_i\}_{i=1}^m)\lambda_{\mathbb{S}^2}(d\mathbf{y}_1)\cdots\lambda_{\mathbb{S}^2}(d\mathbf{y}_m) \right)\\
&=\sum_{n=0}^\infty \sum_{m=0}^{\infty} \exp(-\mu(f^{-1}_i(B)))\exp(-\mu(f^{-1}_j(B)))\\
&\phantom{AAAA}\times\left(\frac{1}{n!}\int_{B\cap f(\mathbb{D}_i)}\cdots\int_{B\cap f(\mathbb{D}_i)}k_{i,n}(\{\mathbf{y}_i\}_{i=1}^n)\lambda_{\mathbb{S}^2}(d\mathbf{y}_1)\cdots\lambda_{\mathbb{S}^2}(d\mathbf{y}_n)\right)\\
&\phantom{AAAAAAAA}\times\left(\frac{1}{m!}\int_{B\cap f(\mathbb{D}_j)}\cdots\int_{B\cap f(\mathbb{D}_j)}k_{j,m}(\{\mathbf{y}_i\}_{i=1}^m)\lambda_{\mathbb{S}^2}(d\mathbf{y}_1)\cdots\lambda_{\mathbb{S}^2}(d\mathbf{y}_m)\right)\\
\intertext{we then use the substitution $m'=m+n$,}
&=\sum_{n=0}^\infty \sum_{m=n}^{\infty} \exp(-\mu(f^{-1}_i(B)))\exp(-\mu(f^{-1}_j(B)))\\
&\phantom{AAAA}\times\left(\frac{1}{n!}\int_{B\cap f(\mathbb{D}_i)}\cdots\int_{B\cap f(\mathbb{D}_i)}k_{i,n}(\{\mathbf{y}_i\}_{i=1}^n)\lambda_{\mathbb{S}^2}(d\mathbf{y}_1)\cdots\lambda_{\mathbb{S}^2}(d\mathbf{y}_n)\right)\\
&\phantom{AAAAAAAA}\times\left(\frac{1}{(m-n)!}\int_{B\cap f(\mathbb{D}_j)}\cdots\int_{B\cap f(\mathbb{D}_j)}k_{j,m-n}(\{\mathbf{y}_i\}_{i=1}^{m-n})\lambda_{\mathbb{S}^2}(d\mathbf{y}_1)\cdots\lambda_{\mathbb{S}^2}(d\mathbf{y}_{m-n})\right)\\
&=\sum_{m=0}^\infty \sum_{n=0}^{m} \exp(-\mu(f^{-1}_i(B)))\exp(-\mu(f^{-1}_j(B)))\frac{1}{m!}\frac{m!}{n!(m-n)!}\\
&\phantom{AAAA}\times\left(\int_{B\cap f(\mathbb{D}_i)}\cdots\int_{B\cap f(\mathbb{D}_i)}k_{i,n}(\{\mathbf{y}_i\}_{i=1}^n)\lambda_{\mathbb{S}^2}(d\mathbf{y}_1)\cdots\lambda_{\mathbb{S}^2}(d\mathbf{y}_n)\right)\\
&\phantom{AAAAAAAA}\times\left(\int_{B\cap f(\mathbb{D}_j)}\cdots\int_{B\cap f(\mathbb{D}_j)}k_{j,m-n}(\{\mathbf{y}_i\}_{i=1}^{m-n})\lambda_{\mathbb{S}^2}(d\mathbf{y}_1)\cdots\lambda_{\mathbb{S}^2}(d\mathbf{y}_{m-n})\right)\\
&=\sum_{m=0}^\infty \frac{\exp(-\mu(f^{-1}_i(B)))\exp(-\mu(f^{-1}_j(B)))}{m!} \\
&\phantom{AAAA}\sum_{n=0}^{m} {m\choose n}\left(\int_{B\cap f(\mathbb{D}_i)}\cdots\int_{B\cap f(\mathbb{D}_i)}k_{i,n}(\{\mathbf{y}_i\}_{i=1}^n)\lambda_{\mathbb{S}^2}(d\mathbf{y}_1)\cdots\lambda_{\mathbb{S}^2}(d\mathbf{y}_n)\right)\\
&\phantom{AAAAAAAA}\times\left(\int_{B\cap f(\mathbb{D}_j)}\cdots\int_{B\cap f(\mathbb{D}_j)}k_{j,m-n}(\{\mathbf{y}_i\}_{i=1}^{m-n})\lambda_{\mathbb{S}^2}(d\mathbf{y}_1)\cdots\lambda_{\mathbb{S}^2}(d\mathbf{y}_{m-n})\right)\\
&=\sum_{m=0}^\infty \frac{\exp(-\mu(f^{-1}_i(B)))\exp(-\mu(f^{-1}_j(B)))}{m!}\\
&\phantom{AAAA}\sum_{n=0}^{m}{m\choose n}\left(\int_{B\cap f(\mathbb{D}_i)}\cdots\int_{B\cap f(\mathbb{D}_i)}k_{i,n}(\{\mathbf{y}_i\}_{i=1}^n)\lambda_{\mathbb{S}^2}(d\mathbf{y}_1)\cdots\lambda_{\mathbb{S}^2}(d\mathbf{y}_n)\right)\\
&\phantom{AAAAAAAA}\times\left(\int_{B\cap f(\mathbb{D}_j)}\cdots\int_{B\cap f(\mathbb{D}_j)}k_{j,m-n}(\{\mathbf{y}_i\}_{i=n+1}^{m})\lambda_{\mathbb{S}^2}(d\mathbf{y}_{n+1})\cdots\lambda_{\mathbb{S}^2}(d\mathbf{y}_{m})\right)\\
&=\sum_{m=0}^\infty \frac{\exp(-\mu(f^{-1}_i(B)))\exp(-\mu(f^{-1}_j(B)))}{m!}\\
&\phantom{AAAA}\times\left(\int_{B\cap f(\mathbb{D}_i)} + \int_{B\cap f(\mathbb{D}_j)}\right)\cdots\\
&\phantom{AAAAAAAA}\left(\int_{B\cap f(\mathbb{D}_i)} + \int_{B\cap f(\mathbb{D}_j)}\right) k_{i,n}(\{\mathbf{y}\}_{i=1}^n)k_{j,m-n}(\{\mathbf{y}\}_{i=n+1}^m) \lambda_{\mathbb{S}^2}(d\mathbf{y}_1)\cdots \lambda_{\mathbb{S}^2}(d\mathbf{y}_m)\\
\intertext{Define $k_{i,j,m}(\{\mathbf{y}_i\}_{i=1}^m)=\mathbbm{1}[\{\mathbf{y}_1,\dots,\mathbf{y}_n\}\in F_{i,j}]\prod_{i=1}^n\rho^*_{i,j}(\mathbf{y}_i)$, where $F_{i,j}=\{x:x=\{\mathbf{x}_1,$ $\dots,\mathbf{x}_n\},n\in\mathbb{N},\mathbf{x}_i\in \mathbb{D}_i\cup\mathbb{D}_j,\exists y\in F \text{ such that } x\subseteq y\}$. Then it can be seen that $k_{i,n}(\{\mathbf{y}\}_{i=1}^n)k_{j,m-n}(\{\mathbf{y}\}_{i=n+1}^m)$ $=k_{i,j,m}(\{\mathbf{y}_i\}_{i=1}^m)$ and so,}
&=\sum_{m=0}^\infty \frac{\exp(-\mu(f^{-1}_i(B)))\exp(-\mu(f^{-1}_j(B)))}{m!}\\
&\phantom{AAAA}\times\left(\int_{B\cap f(\mathbb{D}_i)} + \int_{B\cap f(\mathbb{D}_j)}\right)\cdots\\
&\phantom{AAAAAAAA}\left(\int_{B\cap f(\mathbb{D}_i)} + \int_{B\cap f(\mathbb{D}_j)}\right) k_{i,j,m}(\{\mathbf{y}_i\}_{i=1}^{m}) \lambda_{\mathbb{S}^2}(d\mathbf{y}_1)\cdots \lambda_{\mathbb{S}^2}(d\mathbf{y}_m)\\
&=\sum_{m=0}^\infty \frac{\exp(-\mu(f^{-1}_i(B)))\exp(-\mu(f^{-1}_j(B)))}{m!}\\
&\phantom{AAAA}\times\int_{B\cap \left(f(\mathbb{D}_i)\cup f(\mathbb{D}_j)\right)} \cdots \int_{B\cap \left(f(\mathbb{D}_i)\cup f(\mathbb{D}_j)\right)} k_{i,j,m}(\{\mathbf{y}_i\}_{i=1}^{m}) \lambda_{\mathbb{S}^2}(d\mathbf{y}_1)\cdots \lambda_{\mathbb{S}^2}(d\mathbf{y}_m)\\
\intertext{Define $B_{i,j}=B\cap \left(f(\mathbb{D}_i)\cup f(\mathbb{D}_j)\right)$ and $f^{-1}_{i,j}(B)=f^{-1}_i(B)\cup f^{-1}_{j}(B)$ and noting that $f^{-1}_i(B)\cap f^{-1}_{j}(B)=\emptyset$,}
&=\sum_{m=0}^\infty \frac{\exp(-\mu(f_{i,j}^{-1}(B)))}{m!}\int_{B_{i,j}} \cdots \int_{B_{i,j}} k_{i,j,m}(\{\mathbf{y}_i\}_{i=1}^{m}) \lambda_{\mathbb{S}^2}(d\mathbf{y}_1)\cdots \lambda_{\mathbb{S}^2}(d\mathbf{y}_m)\\
&=\sum_{m=0}^\infty \frac{\exp(-\mu(f_{i,j}^{-1}(B)))}{m!}\int_{B_{i,j}} \cdots \int_{B_{i,j}} \mathbbm{1}[\{\mathbf{y}_1,\dots,\mathbf{y}_m\}\in F_{i,j}]\prod_{i=1}^m\rho_{i,j}^*(\mathbf{y}_i) \lambda_{\mathbb{S}^2}(d\mathbf{y}_1)\cdots \lambda_{\mathbb{S}^2}(d\mathbf{y}_m)
\end{align*}
We can then repeat this argument for each multiplicand in Equation \ref{multiplicand} and reduce the probability measure $P(Y_B\in F)$ to,
\begin{equation*}
P(Y_B\in F) = \sum_{n=0}\frac{\exp(-\mu(f^{-1}(B))}{n!}\int\limits_{B}\cdots\int\limits_{B}\mathbbm{1}[\{\mathbf{y},\dots,\mathbf{y}_n\}\in F]\prod_{i=1}^n\rho^*(\mathbf{y}_i)d\mathbf{y}_i,
\end{equation*}
where,
\begin{equation}
\rho^*(\mathbf{x})=
\begin{cases}
\rho^*_1(\mathbf{x}),&\quad\mathbf{x}\in f(\mathbb{D}_1)\\
& \vdots\\
\rho^*_n(\mathbf{x}),&\quad\mathbf{x}\in f(\mathbb{D}_n).
\end{cases}
\end{equation}
Then defining $\mu^*(B)=\mu(f^{-1}(B))$ we need to show that $\mu^*(B)=\int_B\rho^*(\mathbf{x})d\mathbf{x}$. This follows since,
\begin{align*}
\mu^*(B)&=\mu(f^{-1}(B))\\
&=\int_{f^{-1}(B)}\rho(\mathbf{x})d\mathbf{x}\\
&=\sum_{i=1}^n\int_{f^{-1}_i(B)}\rho(\mathbf{x})d\mathbf{x}\\
&=\sum_{i=1}^n \int_{B\cap f(\mathbb{D}_j)}\rho^*_i(\mathbf{x})d\mathbf{x}\\
&=\int_{B}\rho^*(\mathbf{x})d\mathbf{x},
\end{align*}  
where the penultimate line follows by a the same argument used before, first projecting from $\mathbb{D}_i$ to $\mathbb{R}^2$, then applying the transformation $x'\mapsto x/||\mathbf{x}||$ and $y'\mapsto y/||\mathbf{x}||$, then doing the inverse projection back to the sphere. This finishes the proof by again applying Lemma \ref{Poisson:expansion:lemma}.
\end{proof}

The following corollary shows that two Poisson processes lying on $\mathbb{D}$ with different intensity functions maps to distinct Poisson processes on $\mathbb
{S}^2$ under the same transformation, $f(\mathbf{x})=\mathbf{x}/||\mathbf{x}||$.
\begin{customcorollary}{S1}
Suppose that $X_1$ and $X_2$ are two Poisson processes on $\mathbb{D}$ with intensity functions $\rho_1$ and $\rho_2$ respectively, such that $\rho_1(\mathbf{x})\neq \rho_2(\mathbf{x})$ for at least one $\mathbf{x}\in \mathbb{D}$. Define the transformed processes $Y_1=f(X_1)$ and $Y_2=f(X_2)$ from $\mathbb{D}$ to $\mathbb{S}^2$ where $f(\mathbf{x})=\mathbf{x}/||\mathbf{x}||$. Then the $Y_1$ and $Y_2$ are Poisson processes with intensity functions $\rho^*_1$ and $\rho^*_2$ respectively such that $\rho^*_1(\mathbf{x})\neq\rho^*_2(\mathbf{x})$ for at least one $\mathbf{x}\in\mathbb{S}^2$. 
\end{customcorollary}
\begin{proof}
By assumption we suppose there is at least one point in $\mathbb{D}$ such that $\rho_1(\mathbf{x})\neq \rho_2(\mathbf{x})$. Denote this point $\mathbf{x}^*$ and define $\mathbf{y}^*\equiv f(\mathbf{x}^*)\in f(\mathbb{D}_j)$ for some $j\in\{1,\dots,n\}$. Then at the point $\mathbf{y}^*$ the intensities of $Y_1$ and $Y_2$ are,
\begin{align*}
\rho^*_1(\mathbf{y}^*)&=\rho_1(f^{-1}(\mathbf{y}^*)) l_1(f^{-1}(\mathbf{y}^*)) J_{(1,f^{-1})}(\mathbf{y}^*)\sqrt{1-y^{*2}_1-y^{*2}_2}\\
\rho^*_2(\mathbf{y}^*)&=\rho_1(f^{-1}(\mathbf{y}^*)) l_1(f^{-1}(\mathbf{y}^*)) J_{(1,f^{-1})}(\mathbf{y}^*)\sqrt{1-y^{*2}_1-y^{*2}_2},
\end{align*}
respectively. Then $\rho^*_1(\mathbf{y}^*)\neq \rho^*_2(\mathbf{y}^*)$, since $\rho_1(f^{-1}(\mathbf{y}^*))=\rho_1(\mathbf{x}^*)\neq\rho_2(\mathbf{x}^*)=\rho_2(f^{-1}(\mathbf{y}^*))$.
\end{proof}

The following theorem is the inverse of Theorem \ref{thm:mapping:general:app}, where we take $\mathbf{y}=(x,y,z)\in\mathbb{R}^3$.

\begin{customthm}{S1}\label{thm:inverse:mapping:general:text}
Let $X$ be a Poisson process on $\mathbb{S}^2$ with intensity function $\rho:\mathbb{S}^2\mapsto \mathbb{R}$. Further let $\mathbb{D}\subset\mathbb{R}^3$ be an arbitrary bounded convex shape that $\mathbb{D}=\{\mathbf{y}\in\mathbb{R}^3: g(\mathbf{y})=0\}$ where $g(\mathbf{y})=0$ is the zero set function and is defined as,
\begin{equation*}
g(\mathbf{y})=
\left\{\begin{alignedat}{2}
    g_1(\mathbf{y})=0,&\quad \mathbf{y}\in\mathbb{D}_1\\
    &\vdots \\
    g_n(\mathbf{y})=0,&\quad \mathbf{y}\in\mathbb{D}_n
  \end{alignedat}\right.
\end{equation*}
such that $\cup_{i=1}^n\mathbb{D}_i=\mathbb{D}$ and $\mathbb{D}_i\cap\mathbb{D}_j=\emptyset,\; \forall i\neq j$. Then define $Y=f^{-1}(X)$, where $f(\mathbf{y})=\mathbf{y}/||\mathbf{y}||$ and we have taken the convention that $f^{-1}(X)=\{\mathbf{y}\in\mathbb{D}: f(\mathbf{y})\in X\}$. Then $Y$ is a Poisson process on $\mathbb{D}$, with intensity function,
\begin{equation*}
\rho^*(\mathbf{y})=
\begin{cases}
\frac{\rho(\mathbf{y})J_{1,f^{*-1}}(\mathbf{y})}{l_1(\mathbf{y})\sqrt{1-||\mathbf{y}||^2 x^2-||\mathbf{y}||^2 y^2}},&\quad\mathbf{y}\in\mathbb{D}_1\\
&\vdots\\
\frac{\rho(\mathbf{y})J_{n,f^{*-1}}(\mathbf{y})}{l_n(\mathbf{y})\sqrt{1-||\mathbf{y}||^2 x^2-||\mathbf{y}||^2 y^2}},&\quad\mathbf{y}\in\mathbb{D}_n,
\end{cases}
\end{equation*}
where,
\begin{align*}
z&=\tilde{g}_i(x,y)\\
l_i(\mathbf{y}) &=\left[1+\left(\frac{\partial \tilde{g}_i}{\partial x}\right)^2+\left(\frac{\partial \tilde{g}_i}{\partial y}\right)^2\right]^{\frac{1}{2}}\\
J_{(i,f^{*-1})}(\mathbf{y})&=\frac{1}{(x^2+y^2+\tilde{g}_i^2(x,y))^6}\\
&\det\left[
\begin{pmatrix}
y^2+\tilde{g}_i^2(x,y)-x\tilde{g}_i(x,y)\frac{\partial \tilde{g}_i}{\partial x} & -x\left(y+\tilde{g}_i(x,y)\frac{\partial \tilde{g}_i}{\partial y}\right) \\
-y\left(x+\tilde{g}_i(x,y)\frac{\partial \tilde{g}_i}{\partial x}\right) & x^2+\tilde{g}_i^2(x,y)-y\tilde{g}_i(x,y)\frac{\partial \tilde{g}_i}{\partial y}
\end{pmatrix}\right],
\end{align*}
with $f^{-1}$ is the inverse of $f$, $\det(\cdot)$ is the determinant operator, and $f^*:\mathbb{R}^2\mapsto\mathbb{R}^2$ is the function which maps $x\mapsto x/||\mathbf{y}||$ and $y\mapsto y/||\mathbf{y}||$.
\end{customthm}
\begin{proof}
Proof follows identically to the proof of Theorem \ref{thm:mapping:general:app} but following the argument in reverse.
\end{proof}

\section{Expectation of functional summary statistics}\label{Supplementary:FS:means}
In this section we derive the means for the estimators of $F_{\text{inhom}}$-, $H_{\text{inhom}}$-, and $K_{\text{inhom}}$-functions. To do this we will make use of the Campbell-Mecke Theorem \cite{Moller2004} on $\mathbb{S}^2$,
\begin{customthm}{S2}
Let $X$ be a point process on $\mathbb{S}^2$ and $h$ be any non-negative, measurable function such that $h:\mathbb{S}^2\times N_{lf}\mapsto\mathbb{R}$. Then,
\begin{equation*}
\mathbb{E}\sum_{\mathbf{x}\in X}h(\mathbf{x},X\setminus\{\mathbf{x}\})=\int_{\mathbb{S}^2}\mathbb{E}[h(\mathbf{x},X^!_{\mathbf{x}})] \mu(d\mathbf{x}),
\end{equation*}
and if the intensity function exists then,
\begin{equation*}
\mathbb{E}\sum_{\mathbf{x}\in X}h(\mathbf{x},X\setminus\{\mathbf{x}\})=\int_{\mathbb{S}^2}\mathbb{E}[h(\mathbf{x},X^!_{\mathbf{x}})] \rho(\mathbf{x}) \lambda_{\mathbb{S}^2}(d\mathbf{x}).
\end{equation*}
\end{customthm}
\begin{proof}
See \cite[Appendix C, pp. 248-249]{Moller2004}.
\end{proof}
Noting that the reduced Palm process, $X^!_\mathbf{x}$, for a Poisson process, $X$, is again the same Poisson process \cite{Coeurjolly2015} we get the Slivnyak-Mecke Theorem \cite[Theorem 3.2, p. 21]{Moller2004},
\begin{equation*}
\mathbb{E}\sum_{\mathbf{x}\in X}h(\mathbf{x},X\setminus\{\mathbf{x}\})=\int_{\mathbb{S}^2}\mathbb{E}[h(\mathbf{x},X)] \rho(\mathbf{x}) \lambda_{\mathbb{S}^2}(d\mathbf{x}).
\end{equation*}
We now prove give the expectation of estimates of the inhomogeneous functional summary statistics,

\begin{theorem}\label{expect:FS}
Let $X$ be a spherical Poisson process on $\mathbb{S}^2$ with known intensity function $\rho:\mathbb{S}^2\mapsto \mathbb{R}_+$, such that $\bar{\rho}=\inf_{\mathbf{x}\in\mathbb{S}^2}\rho(\mathbf{x})>0$. Then the estimators for $\hat{F}_{\text{inhom}}(r)$, and $\hat K_{\text{inhom}}(r)$ are unbiased whilst $\hat{H}_{\text{inhom}}(r)$ is ratio-unbiased. More precisely, 
\begin{align*}
\mathbb{E}[\hat{F}_{\text{inhom}}(r)]&=1-\exp(-\bar{\rho}2\pi(1-\cos r))\\
\mathbb{E}[\hat{H}_{\text{inhom}}(r)]&=1-\frac{\mathrm{exp}(-\bar{\rho}2\pi(1-\cos r))-\mathrm{exp}(-\mu(\mathbb{S}^2))}{1-\frac{\bar{\rho}2\pi(1-\cos r)}{\mu(\mathbb{S}^2)}}\\
\mathbb{E}[\hat{K}_{\text{inhom}}(r)]&=2\pi (1-\cos r),
\end{align*}
where $r\in [0,\pi]$, and $\bar{\rho}=\inf_{\mathbf{x}\in\mathbb{S}^2}\rho(\mathbf{x})>0$. Further by unbiasedness and ratio-unbiasedness of $\hat{F}_{\text{inhom}}(r)$ and $\hat{H}_{\text{inhom}}(r)$, respectively, we immediately have ratio-unbiasedness of $\hat{J}_{\text{inhom}}(r)$. 
\end{theorem}
\begin{proof}
Proofs for the expectation of $\hat{F}_{\text{inhom}}$-, and $\hat{K}_{\text{inhom}}$-functions are found in \cite{vanLieshout2011} and \cite{Lawrence2016,Moller2016} respectively, whilst ratio-unbiasedness of the $\hat{H}_{\text{inhom}}$-function is also found in \cite{vanLieshout2011}. The proofs found in \cite{vanLieshout2011} are in $\mathbb{R}^n$ but can be extended easily to $\mathbb{S}^2$. Then for the expectation of $\hat{H}_{\text{inhom}}(r)$,
\begin{align}
&\mathbb{E}\left[\frac{1}{N_X(\mathbb{S}^2)}\sum_{\mathbf{x}\in X}\prod_{\mathbf{y}\in X\setminus\{\mathbf{x}\}}\left(1-\frac{\bar{\rho}\mathbbm{1}[\mathbf{y}\in B_{\mathbb{S}^2}(\mathbf{x},r)]}{\rho(\mathbf{y})}\right)\right]\nonumber\\
&\phantom{AAAA}=\mathbb{E}\left[\sum_{\mathbf{x}\in X}\frac{1}{N_X(\mathbb{S}^2)}\prod_{\mathbf{y}\in X\setminus\{\mathbf{x}\}}\left(1-\frac{\bar{\rho}\mathbbm{1}[\mathbf{y}\in B_{\mathbb{S}^2}(\mathbf{x},r)]}{\rho(\mathbf{y})}\right)\right]\nonumber\\
&\phantom{AAAA}=\mathbb{E}\left[\sum_{\mathbf{x}\in X}\frac{1}{|X|}\prod_{\mathbf{y}\in X\setminus\{\mathbf{x}\}}\left(1-\frac{\bar{\rho}\mathbbm{1}[\mathbf{y}\in B_{\mathbb{S}^2}(\mathbf{x},r)]}{\rho(\mathbf{y})}\right)\right]\nonumber\\
&\phantom{AAAA}=\mathbb{E}\left[\sum_{\mathbf{x}\in X}\frac{1}{|X\setminus\{\mathbf{x}\}\cup\{\mathbf{x}\}|}\prod_{\mathbf{y}\in X\setminus\{\mathbf{x}\}}\left(1-\frac{\bar{\rho}\mathbbm{1}[\mathbf{y}\in B_{\mathbb{S}^2}(\mathbf{x},r)]}{\rho(\mathbf{y})}\right)\right],\nonumber\\
\intertext{applying the Slivnyak-Mecke Theorem,}
&\phantom{AAAA}=\int_{\mathbb{S}^2}\mathbb{E}\left[\frac{1}{|X \cup\{\mathbf{x}\}|}\prod_{\mathbf{y}\in X}\left(1-\frac{\bar{\rho}\mathbbm{1}[\mathbf{y}\in B_{\mathbb{S}^2}(\mathbf{x},r)]}{\rho(\mathbf{y})}\right)\right]\rho(\mathbf{x})\lambda_{\mathbb{S}^2}(d\mathbf{x})\nonumber\\
&\phantom{AAAA}=\int_{\mathbb{S}^2}\mathbb{E}\left[\frac{1}{N_{X}(\mathbb{S}^2)+1}\prod_{\mathbf{y}\in X}\left(1-\frac{\bar{\rho}\mathbbm{1}[\mathbf{y}\in B_{\mathbb{S}^2}(\mathbf{x},r)]}{\rho(\mathbf{y})}\right)\right]\rho(\mathbf{x})\lambda_{\mathbb{S}^2}(d\mathbf{x}).\label{int:x:1}
\end{align} 
We then take the expectation in the integrand and handle it separately, using the definition of a Poisson process being the independent distribution of a Poisson number of points.
\begin{align*}
&\mathbb{E}\left[\frac{1}{N_{X}(\mathbb{S}^2)+1}\prod_{\mathbf{y}\in X}\left(1-\frac{\bar{\rho}\mathbbm{1}[\mathbf{y}\in B_{\mathbb{S}^2}(\mathbf{x},r)]}{\rho(\mathbf{y})}\right)\right]\\
&\phantom{AAAA}=\mathbb{E}\left[\frac{1}{N_{X}(\mathbb{S}^2)+1}\mathbb{E}\left[\prod_{\mathbf{y}\in X}\left(1-\frac{\bar{\rho}\mathbbm{1}[\mathbf{y}\in B_{\mathbb{S}^2}(\mathbf{x},r)]}{\rho(\mathbf{y})}\right)\Bigg|N_{X}(\mathbb{S}^2)=n\right]\right]\\
&\phantom{AAAA}=\mathbb{E}\left[\frac{1}{N_{X}(\mathbb{S}^2)+1}\mathbb{E}\left[\prod_{i=1}^n\left(1-\frac{\bar{\rho}\mathbbm{1}[\mathbf{X}_i\in B_{\mathbb{S}^2}(\mathbf{x},r)]}{\rho(\mathbf{X}_i)}\right)\right]\right],
\end{align*}
where $\mathbf{X}_i$ are independently distributed across $\mathbb{S}^2$ with density $\frac{\rho(\mathbf{x})}{\mu(\mathbb{S}^2)}$. Then taking the first expectation,
\begin{align*}
\mathbb{E}\left[\prod_{i=1}^n\left(1-\frac{\bar{\rho}\mathbbm{1}[\mathbf{X}_i\in B_{\mathbb{S}^2}(\mathbf{x},r)]}{\rho(\mathbf{X}_i)}\right)\right]&=\overbrace{\int_{\mathbb{S}^2}\cdots\int_{\mathbb{S}^2}}^n \prod_{i=1}^n\left(1-\frac{\bar{\rho}\mathbbm{1}[\mathbf{x}_i\in B_{\mathbb{S}^2}(\mathbf{x},r)]}{\rho(\mathbf{x}_i)}\right)\frac{\rho(\mathbf{x}_i)}{\mu(\mathbb{S}^2)}\lambda_{\mathbb{S}^2}(d\mathbf{x})_i\\
&=\left(\int_{\mathbb{S}^2} \frac{\rho(\mathbf{y})}{\mu(\mathbb{S}^2)} - \bar{\rho}\mathbbm{1}[\mathbf{y}\in B_{\mathbb{S}^2}(\mathbf{x},r)] d\mathbf{y}\right)^n\\
&=\left(1-\frac{\bar{\rho}}{\mu(\mathbb{S}^2)}2\pi (1-\cos r)\right)^n
\end{align*}
Returning to the expectation over $N_X(\mathbb{S}^2)$ we have,
\begin{align*}
&\mathbb{E}\left[\frac{1}{N_{X}(\mathbb{S}^2)+1}\prod_{\mathbf{y}\in X}\left(1-\frac{\bar{\rho}\mathbbm{1}[\mathbf{y}\in B_{\mathbb{S}^2}(\mathbf{x},r)]}{\rho(\mathbf{y})}\right)\right]\\
&\phantom{AAAA}=\mathbb{E}\left[\frac{1}{N_{X}(\mathbb{S}^2)+1}\left(1-\frac{\bar{\rho}}{\mu(\mathbb{S}^2)}2\pi (1-\cos r)\right)^{N_X(\mathbb{S}^2)}\right],\\
\intertext{define $A\equiv 1-\frac{\bar{\rho}}{\mu(\mathbb{S}^2)}2\pi (1-\cos r)$, then}
&\phantom{AAAA}=\mathbb{E}\left[\frac{1}{N_{X}(\mathbb{S}^2)+1}A^{N_X(\mathbb{S}^2)}\right],\\
&\phantom{AAAA}=\sum_{n=0}^\infty \frac{A^n}{n+1}\frac{\lambda^n\mathrm{e}^{-\lambda}}{n!},\\
\intertext{where $\lambda\equiv \mu(\mathbb{S}^2)$,}
&\phantom{AAAA}=\frac{\mathrm{e}^{-\lambda}}{A\lambda}\sum_{n=0}^\infty \frac{(A\lambda)^{n+1}}{(n+1)!}\\
&\phantom{AAAA}=\frac{\mathrm{e}^{-\lambda}}{A\lambda}\sum_{n=1}^\infty \frac{(A\lambda)^{n}}{n!}\\
&\phantom{AAAA}=\frac{\mathrm{e}^{-\lambda}}{A\lambda}\left(\sum_{n=0}^\infty \frac{(A\lambda)^{n}}{n!}-1\right)\\
&\phantom{AAAA}=\frac{\mathrm{e}^{-\lambda}}{A\lambda}\left(\mathrm{e}^{A\lambda}-1\right)\\
&\phantom{AAAA}=\frac{\mathrm{e}^{-\bar{\rho}2\pi(1-\cos r)}-\mathrm{e}^{-\mu(\mathbb{S}^2)}}{\mu(\mathbb{S}^2)-\bar{\rho}2\pi(1-\cos r)},
\end{align*}
plugging this into Equation \ref{int:x:1},
\begin{align*}
\int_{\mathbb{S}^2}\frac{\mathrm{e}^{-\bar{\rho}2\pi(1-\cos r)}-\mathrm{e}^{-\mu(\mathbb{S}^2)}}{\mu(\mathbb{S}^2)-\bar{\rho}2\pi(1-\cos r)}\rho(\mathbf{x})\lambda_{\mathbb{S}^2}(d\mathbf{x}) = \frac{\mathrm{e}^{-\bar{\rho}2\pi(1-\cos r)}-\mathrm{e}^{-\mu(\mathbb{S}^2)}}{1-\frac{\bar{\rho}2\pi(1-\cos r)}{\mu(\mathbb{S}^2)}},
\end{align*}
and so,
\begin{equation*}
\mathbb{E}[\hat{H}_{\text{inhom}}(r)]=1-\frac{\mathrm{e}^{-\bar{\rho}2\pi(1-\cos r)}-\mathrm{e}^{-\mu(\mathbb{S}^2)}}{1-\frac{\bar{\rho}2\pi(1-\cos r)}{\mu(\mathbb{S}^2)}}.
\end{equation*}
\end{proof}

The following corollary puts a bound on the bias of $\hat{H}_{\text{inhom}}(r)$ when we consider the scenario of CSR on $\mathbb{D}$ mapped on $\mathbb{S}^2$.

\begin{corollary}
With the same assumptions as Theorem \ref{expect:FS}, let $X$ be a spherical Poisson process on $\mathbb{S}^2$ with intensity function $\rho:\mathbb{S}^2\mapsto \mathbb{R}_+$. Defining $\bar{\rho}=\inf_{\mathbf{x}\in\mathbb{S}^2}\rho(\mathbf{x})$, the bias of the estimator $\hat{H}_{\text{inhom}}(r)$ is bounded by
\begin{equation*}
|\text{Bias}(\hat{H}_{\text{inhom}}(r))|\leq \exp(-\mu(\mathbb{S}^2))\leq \exp(-4\pi\bar{\rho}),
\end{equation*}
for all $r\in[0,\pi]$.
\end{corollary}
\begin{proof}
\begin{align*}
\text{Bias}(\hat{H}_{\text{inhom}}(r)) &= \left(1-\frac{\mathrm{e}^{-\bar{\rho}2\pi(1-\cos r)}-\mathrm{e}^{-\mu(\mathbb{S}^2)}}{1-\frac{\bar{\rho}2\pi(1-\cos r)}{\mu(\mathbb{S}^2)}}\right) - \left(1-\mathrm{e}^{-\bar{\rho}2\pi(1-\cos r)}\right)\\
&=\frac{\exp(\mu(\mathbb{S}^2))-\frac{\bar{\rho}2\pi(1-\cos r)}{\mu(\mathbb{S}^2)}\exp(-2\pi(1-\cos r)\bar{\rho})}{1-\frac{\bar{\rho}2\pi(1-\cos r)}{\mu(\mathbb{S}^2)}}\\
&=\frac{\mu(\mathbb{S}^2)\exp(\mu(\mathbb{S}^2))-\bar{\rho}2\pi(1-\cos r)\exp(-2\pi(1-\cos r)\bar{\rho})}{\mu(\mathbb{S}^2)-\bar{\rho}2\pi(1-\cos r)}
\end{align*}
Taking the absolute value of the bias and the numerator is bounded above by,
\begin{align*}
&|\mu(\mathbb{S}^2)\exp(\mu(\mathbb{S}^2))-\bar{\rho}2\pi(1-\cos r)\exp(-2\pi(1-\cos r)\bar{\rho})|\\
&\phantom{AAAAAAAA}\leq \mu(\mathbb{S}^2)\exp(\mu(\mathbb{S}^2))+\bar{\rho}2\pi(1-\cos r)\exp(-2\pi(1-\cos r)\bar{\rho})\\
&\phantom{AAAAAAAA}\leq \mu(\mathbb{S}^2)\exp(\mu(\mathbb{S}^2)),
\end{align*}
where the fist line follows from the triangle inequality. The denominator is bounded below,
\begin{align*}
|\mu(\mathbb{S}^2)-\bar{\rho}2\pi(1-\cos r)|&\geq \mu(\mathbb{S}^2)+\bar{\rho}2\pi(1-\cos r)\\
&\geq\mu(\mathbb{S}^2),
\end{align*}
where again the first line follows from the triangle inequality, since $|a-b+b|\leq |a-b|+|b|\Rightarrow |a|-|b|\leq |a-b|,$ for $a,b\in\mathbb{R}$. Hence the absolute of the bias is bounded by,
\begin{equation*}
|\text{Bias}(\hat{H}_{\text{inhom}}(r))|\leq \exp(\mu(\mathbb{S}^2)).
\end{equation*}
The final inequality follows by noting that $\mu(\mathbb{S}^2)=\int_{\mathbb{S}^2}\rho(\mathbf{x})\lambda_{\mathbb{S}^2}(d\mathbf{x})\geq \int_{\mathbb{S}^2}\bar{\rho}\lambda_{\mathbb{S}^2}(d\mathbf{x})=4\pi\bar{\rho}$.
\end{proof}

\section{Variance of functional summary statistics}\label{Supplementary:FS:variances}

In this section we derive the variance of the functional summary statistics. Throughout this section we will frequently refer to the area of a spherical cap, at any point $\mathbf{o}\in\mathbb{S}^2$ with geodesic distance $r$ by $B_{\mathbb{S}^2}(\mathbf{o},r)$. We will also make use of the extended Campbell-Mecke Theorem \cite{Moller2004} throughout,
\begin{customthm}{S3}
Let $X$ be a point process on $\mathbb{S}^2$ and $h$ be any non-negative, measurable function such that $h:\left(\times_{i=1}^n\mathbb{S}^2\right)\times N_{lf}\mapsto\mathbb{R}$. Then,
\begin{align*}
&\mathbb{E}\sum_{\mathbf{x}_1,\dots,\mathbf{x}_n\in X}^{\neq}h(\mathbf{x}_1,\dots,\mathbf{x}_n,X\setminus\{\mathbf{x}_1,\dots,\mathbf{x}_n\})\\
&\phantom{AAAA}=\int_{\mathbb{S}^2}\cdots\int_{\mathbb{S}^2}\mathbb{E}[h(\mathbf{x}_1,\dots,\mathbf{x}_n,X^!_{\mathbf{x}_1,\dots,\mathbf{x}_n})] \alpha(\lambda_{\mathbb{S}^2}(d\mathbf{x})_1,\dots,\lambda_{\mathbb{S}^2}(d\mathbf{x})_n)\\
\intertext{and if the intensity function of order $n$ exists then,}
&\mathbb{E}\sum_{\mathbf{x}_1,\dots,\mathbf{x}_n\in X}^{\neq}h(\mathbf{x}_1,\dots,\mathbf{x}_n,X\setminus\{\mathbf{x}_1,\dots,\mathbf{x}_n\})\\
&\phantom{AAAA}=\int_{\mathbb{S}^2}\cdots\int_{\mathbb{S}^2}\mathbb{E}[h(\mathbf{x}_1,\dots,\mathbf{x}_n,X^!_{\mathbf{x}_1,\dots,\mathbf{x}_n})] \rho^{(n)}(\mathbf{x}_1,\dots,\mathbf{x}_n)\lambda_{\mathbb{S}^2}(\lambda_{\mathbb{S}^2}(d\mathbf{x})_1)\cdots\lambda_{\mathbb{S}^2}(\lambda_{\mathbb{S}^2}(d\mathbf{x})_n),
\end{align*}
\end{customthm}
where $X^!_{\mathbf{x}_1,\dots,\mathbf{x}_n}$ is the $n^{th}$-order reduced Palm process of $X$.
\begin{proof}
See \cite[Appendix C, pp. 248-249]{Moller2004}.
\end{proof}
Again, as in the case for when $n=1$, for any order $n$ the reduced Palm process of order $n$, $X^!_{\mathbf{x}_1,\dots,\mathbf{x}_n}$, for a Poisson process, $X$, is again the same Poisson process. Hence we get the extended Slivnyak-Mecke Theorem \cite[Theorem 3.3, p. 22]{Moller2004},
\begin{equation}\label{extended:SM}
\begin{split}
&\mathbb{E}\sum_{\mathbf{x}_1,\dots,\mathbf{x}_n\in X}^{\neq}h(\mathbf{x}_1,\dots,\mathbf{x}_n,X\setminus\{\mathbf{x}_1,\dots,\mathbf{x}_n\})\\
&\phantom{AAAA}=\int_{\mathbb{S}^2}\cdots\int_{\mathbb{S}^2}\mathbb{E}[h(\mathbf{x}_1,\dots,\mathbf{x}_n,X)] \prod_{i=1}^n \rho(\mathbf{x}_i) \lambda_{\mathbb{S}^2}(\lambda_{\mathbb{S}^2}(d\mathbf{x})_1)\cdots\lambda_{\mathbb{S}^2}(\lambda_{\mathbb{S}^2}(d\mathbf{x})_n),
\end{split}
\end{equation}
To derive the variance of $\hat{K}_{\text{inhom}}(r)$, before which we require the following lemma,
\begin{customlemma}{S2}\label{set:partition:lemma}
Let $X$ be a finite set and define $X^n$ as the Cartesian product of $X$ $n$ times, i.e. $X^n=X\times\cdots\times X$, and  the following sets,
\begin{align*}
Y & = \{(\mathbf{x}_1,\mathbf{x}_2,\mathbf{x}_3,\mathbf{x}_4)^T\in X^4 : \mathbf{x}_1\in X, \mathbf{x}_2\in X\setminus\{\mathbf{x}_1\}, \mathbf{x}_3\in X, \mathbf{x}_4\in X\setminus\{\mathbf{x}_3\}\}\\
Y_1 &=\{(\mathbf{x}_1,\mathbf{x}_2,\mathbf{x}_3,\mathbf{x}_4)^T\in X^4 : \mathbf{x}_1\in X, \mathbf{x}_2\in X\setminus\{\mathbf{x}_1\}, \mathbf{x}_3\in X\setminus\{\mathbf{x}_1,\mathbf{x}_2\}, \mathbf{x}_4\in X\setminus\{\mathbf{x}_1,\mathbf{x}_2,\mathbf{x}_3\}\}\\
Y_2 &=\{(\mathbf{x}_1,\mathbf{x}_2,\mathbf{x}_3,\mathbf{x}_4)^T\in X^4 : \mathbf{x}_1\in X, \mathbf{x}_2\in X\setminus\{\mathbf{x}_1\}, \mathbf{x}_3\in \{\mathbf{x}_1,\mathbf{x}_2\}, \mathbf{x}_4\in \{\mathbf{x}_1,\mathbf{x}_2\}\setminus\{\mathbf{x}_3\}\}\\
Y_3 &=\{(\mathbf{x}_1,\mathbf{x}_2,\mathbf{x}_3,\mathbf{x}_4)^T\in X^4 : \mathbf{x}_1\in X, \mathbf{x}_2\in X\setminus\{\mathbf{x}_1\}, \mathbf{x}_3\in X\setminus\{\mathbf{x}_1,\mathbf{x}_2\}, \mathbf{x}_4\in \{\mathbf{x}_1,\mathbf{x}_2\}\setminus\{\mathbf{x}_3\}\}\\
Y_4 &=\{(\mathbf{x}_1,\mathbf{x}_2,\mathbf{x}_3,\mathbf{x}_4)^T\in X^4 : \mathbf{x}_1\in X, \mathbf{x}_2\in X\setminus\{\mathbf{x}_1\}, \mathbf{x}_3\in \{\mathbf{x}_1,\mathbf{x}_2\}, \mathbf{x}_4\in X\setminus\{\mathbf{x}_1,\mathbf{x}_2,\mathbf{x}_3\}\},
\end{align*}
then $Y=\cup_{i=1,\dots,4}Y_i$ and $Y_i, i=1,\dots,4$ are pairwise disjoint.
\end{customlemma}
\begin{proof}
Pairwise disjointness of the sets $Y_i, i=1\dots4$ follows from the definitions of the sets. To prove equality we will show the following,
\begin{align}
\mathbf{r}\in Y&\Rightarrow \mathbf{r}\in\cup_{i=1,\dots,4}Y_i, \quad \text{ and}\label{set:identity:1}\\
\mathbf{r}\in \cup_{i=1,\dots,4}Y_i&\Rightarrow \mathbf{r}\in Y\label{set:identity:2}
\end{align}
Statement \ref{set:identity:2}, holds by considering each set $Y_i$ in turn. In particular it is clear by the definitions of the sets that for each $i=1,\dots,4, Y_i\subseteq Y$ and so Statement \ref{set:identity:2} holds. 

To show Statement \ref{set:identity:1} let $\mathbf{x}\in Y$ and fix $\mathbf{x}_1\in X$ and $\mathbf{x}_2\in X\setminus\{\mathbf{x}_1\}$. Then there are two possibilities for $\mathbf{x}_3$, either $\mathbf{x}_3\in X\setminus\{\mathbf{x}_1,\mathbf{x}_2\}$ or $\mathbf{x}_3\in \{\mathbf{x}_1,\mathbf{x}_2\}$. If the former holds then $\mathbf{x}_4$ can either be in $X\setminus\{\mathbf{x}_1,\mathbf{x}_2,\mathbf{x}_3\}$ or $\{\mathbf{x}_1,\mathbf{x}_2\}\setminus\{\mathbf{x}_3\}$.  If the first holds then $\mathbf{x}\in Y_1$ and if the second holds then $\mathbf{x}\in Y_3$. Considering all possible combinations proves Statement \ref{set:identity:1}. Hence it follows that $Y=\cup_{i=1,\dots,4}Y_i$.
\end{proof}

We now proceed with the proof for the variance of estimates of the inhomogeneous functional summary statistics.

\begin{theorem}\label{thm:variance:K:inhom:app}
Let $X$ be a spherical Poisson process on $\mathbb{S}^2$ with known intensity function $\rho:\mathbb{S}^2\mapsto \mathbb{R}_+$, such that $\bar{\rho}=\inf_{\mathbf{x}\in\mathbb{S}^2}\rho(\mathbf{x})>0$. Then the estimators $\hat{K}_{\rm{inhom}}(r)$, $\hat{F}_{\rm{inhom}}(r)$, and $\hat{H}_{\rm{inhom}}(r)$ have variance,
\begin{align*}
&\emph{Var}(\hat{K}_{\rm{inhom}}(r)) = \frac{1}{8\pi^2}\int_{\mathbb{S}^2}\int_{\mathbb{S}^2}\frac{\mathbbm{1}[d(\mathbf{x},\mathbf{y})\leq r]}{\rho(\mathbf{x})\rho(\mathbf{y})}\lambda_{\mathbb{S}^2}(d\mathbf{x})\lambda_{\mathbb{S}^2}(d\mathbf{y}) +(1-\cos r)^2\int_{\mathbb{S}^2}\frac{1}{\rho(\mathbf{x})}\lambda_{\mathbb{S}^2}(d\mathbf{x}),\\
&\emph{Var}(\hat{F}_{\rm{inhom}}(r)) = \frac{\exp\left(-2\bar{\rho}\lambda_{\mathbb{S}^2}(B_{\mathbb{S}^2}(\mathbf{o},r))\right)}{|P|^2}\\
&\phantom{AAAAAAAA}\sum_{\mathbf{p}\in P}\sum_{\mathbf{p}'\in P} \exp\left(\int_{B_{\mathbb{S}^2}(\mathbf{p},r)\cap B_{\mathbb{S}^2}(\mathbf{p}',r)}\frac{\bar{\rho}^2}{\rho(\mathbf{x})}\lambda_{\mathbb{S}^2}(d\mathbf{x})\right)-\exp\left(-2\bar{\rho}\lambda_{\mathbb{S}^2}(B_{\mathbb{S}^2}(\mathbf{o},r))\right),\\
&\emph{Var}(\hat{H}_{\rm{inhom}}(r))\\
&=\frac{1}{\mu^2(\mathbb{S}^2)}\int_{\mathbb{S}^2}\int_{\mathbb{S}^2}\left(\rho(\mathbf{x})-\bar{\rho}\mathbbm{1}[\mathbf{x}\in B_{\mathbb{S}^2}(\mathbf{y},r)]\right)\left(\rho(\mathbf{y})-\bar{\rho}\mathbbm{1}[\mathbf{y}\in B_{\mathbb{S}^2}(\mathbf{x},r)]\right)\\
&\phantom{=+}\frac{\mathrm{e}^{-\mu(\mathbb{S}^2)}}{A_1^{2}(\mathbf{x},\mathbf{y})}\left(\mathrm{e}^{\mu(\mathbb{S}^2) A_1(\mathbf{x},\mathbf{y})}-1 -\rm{Ei}(\mu(\mathbb{S}^2) A_1(\mathbf{x},\mathbf{y}))+\gamma +\log(\mu(\mathbb{S}^2) A_1(\mathbf{x},\mathbf{y}))\right) \lambda_{\mathbb{S}^2}(d\mathbf{x}) \lambda_{\mathbb{S}^2}(d\mathbf{y})\\
&\phantom{=}+\frac{1}{\mu(\mathbb{S}^2)}\int_{\mathbb{S}^2}\frac{\mathrm{e}^{-\mu(\mathbb{S}^2)}}{A_2(\mathbf{x})}\left(\gamma +\log(\mu(\mathbb{S}^2) A_2(\mathbf{x}))-\rm{Ei}(\mu(\mathbb{S}^2) A_2(\mathbf{x}))\right)\rho(\mathbf{y})\lambda_{\mathbb{S}^2}(d\mathbf{y})\\
&\phantom{=}-\frac{\mathrm{e}^{-2\mu(\mathbb{S}^2)}}{\left(1-\frac{\bar{\rho}}{\mu(\mathbb{S}^2)}2\pi(1-\cos r)\right)^2}\left(\mathrm{e}^{\mu(\mathbb{S}^2) \left(1-\frac{\bar{\rho}}{\mu(\mathbb{S}^2)}2\pi(1-\cos r)\right)}-1\right)^2
\end{align*}
where,
\begin{align*}
A_1(\mathbf{x},\mathbf{y})&= 1-\frac{2\bar{\rho}}{\mu(\mathbb{S}^2)}2\pi(1-\cos r)+\frac{\bar{\rho}^2}{\mu(\mathbb{S}^2)}\int_{B_{\mathbb{S}^2}(\mathbf{x},r)\cap B_{\mathbb{S}^2}(\mathbf{y},r)}\frac{1}{\rho(\mathbf{z})}d\mathbf{z}\\
A_2(\mathbf{x})&= 1-\frac{2\bar{\rho}}{\mu(\mathbb{S}^2)}2\pi(1-\cos r)+\frac{\bar{\rho}^2}{\mu(\mathbb{S}^2)}\int_{B_{\mathbb{S}^2}(\mathbf{x},r)}\frac{1}{\rho(\mathbf{y})}\lambda_{\mathbb{S}^2}(d\mathbf{y})\\
\rm{Ei}(x) &= -\int_{-x}^{\infty}\frac{\mathrm{e}^{-t}}{t} dt
\end{align*}
and $\rm{Ei}(x)$ is the exponential integral and $r\in [0,\pi]$.
\end{theorem}
The proof is spread over three parts, one for each functional summary statistic.
\begin{proof}
\subsection*{Variance of $\hat{K}_{\text{inhom}}(r)$}
We expand the variance as $\text{Var}(X)=\mathbb{E}[X^2]-\mathbb{E}^2[X]$,
\begin{align}
\text{Var}(\tilde{K}_{\text{inhom}}(r)) &= \text{Var}\left(\frac{1}{4\pi}\sum_{\mathbf{x}\in X}\sum_{\mathbf{y}\in X\setminus\{\mathbf{x}\}} \frac{\mathbbm{1}[d(\mathbf{x},\mathbf{y})\leq r]}{\rho(\mathbf{x})\rho(\mathbf{y})}\right)\nonumber\\
&= \frac{1}{16\pi^2}\Bigg[\underbrace{\mathbb{E}\Bigg(\sum_{\mathbf{x}\in X}\sum_{\mathbf{y}\in X\setminus\{\mathbf{x}\}} \frac{\mathbbm{1}[d(\mathbf{x},\mathbf{y})\leq r]}{\rho(\mathbf{x})\rho(\mathbf{y})}\Bigg)^2}_{(1)}-\underbrace{\mathbb{E}^2\sum_{\mathbf{x}\in X}\sum_{\mathbf{y}\in X\setminus\{\mathbf{x}\}} \frac{\mathbbm{1}(d(\mathbf{x},\mathbf{y})\leq r)}{\rho(\mathbf{x})\rho(\mathbf{y})}}_{(2)}\Bigg],\label{variance:E:X2}
\end{align}
We deal with each of the terms individually. First consider term $(2)$ of the previous equation, this is simply the inhomogeneous $K$-function for a Poisson process,
\begin{align*}
\frac{1}{16\pi^2}\mathbb{E}^2\sum_{\mathbf{x}\in X}\sum_{\mathbf{y}\in X\setminus\{\mathbf{x}\}} \frac{\mathbbm{1}(d(\mathbf{x},\mathbf{y})\leq r)}{\rho(\mathbf{x})\rho(\mathbf{y})} &= \mathbb{E}^2\frac{1}{4\pi}\sum_{\mathbf{x}\in X}\sum_{\mathbf{y}\in X\setminus\{\mathbf{x}\}} \frac{\mathbbm{1}(d(\mathbf{x},\mathbf{y})\leq r)}{\rho(\mathbf{x})\rho(\mathbf{y})}\\
&= K_{\text{inhom}}^2(r)\\ 
&= 4\pi^2(1-\cos r)^2,
\end{align*}
where the penultimate equality follows from the definition of $K_{\text{inhom}}$ (taking the arbitrary area $B$ to be $\mathbb{S}^2$) and the final equality follows from Proposition \ref{expect:FS}. To handle term (1) we first expand the square,
\begin{align*}
&\frac{1}{16\pi^2}\mathbb{E}\Bigg(\sum_{\mathbf{x}\in X}\sum_{\mathbf{y}\in X\setminus\{\mathbf{x}\}} \frac{\mathbbm{1}[d(\mathbf{x},\mathbf{y})\leq r]}{\rho(\mathbf{x})\rho(\mathbf{y})}\Bigg)^2\\
&\phantom{AAAA}=\frac{1}{16\pi^2}\mathbb{E}\sum_{\mathbf{x}\in X}\sum_{\mathbf{y}\in X\setminus\{\mathbf{x}\}} \frac{\mathbbm{1}[d(\mathbf{x},\mathbf{y})\leq r]}{\rho(\mathbf{x})\rho(\mathbf{y})}\sum_{\mathbf{x}'\in X}\sum_{\mathbf{y}'\in X\setminus\{\mathbf{x}'\}} \frac{\mathbbm{1}[d(\mathbf{x}',\mathbf{y}')\leq r]}{\rho(\mathbf{x}')\rho(\mathbf{y}')}\\
&\phantom{AAAA}=\frac{1}{16\pi^2}\mathbb{E}\sum_{\mathbf{x}\in X}\sum_{\mathbf{y}\in X\setminus\{\mathbf{x}\}}\sum_{\mathbf{x}'\in X}\sum_{\mathbf{y}'\in X\setminus\{\mathbf{x}'\}} \frac{\mathbbm{1}[d(\mathbf{x},\mathbf{y})\leq r]\mathbbm{1}[d(\mathbf{x}',\mathbf{y}')\leq r]}{\rho(\mathbf{x})\rho(\mathbf{y})\rho(\mathbf{x}')\rho(\mathbf{y}')}
\end{align*}
Let us define the summand as,
\begin{equation*}
f_r(\mathbf{x},\mathbf{y},\mathbf{x}',\mathbf{y}') = \frac{\mathbbm{1}[d(\mathbf{x},\mathbf{y})\leq r]\mathbbm{1}[d(\mathbf{x}',\mathbf{y}')\leq r]}{\rho(\mathbf{x})\rho(\mathbf{y})\rho(\mathbf{x}')\rho(\mathbf{y}')},
\end{equation*}
and then by Lemma \ref{set:partition:lemma} we can divide the sum in the expectation into 4 terms,
\begin{align*}
&\underbrace{\frac{1}{16\pi^2}\mathbb{E}\sum_{\mathbf{x}\in X}\sum_{\mathbf{y}\in X\setminus\{\mathbf{x}\}} \sum_{\mathbf{x}'\in X\setminus\{\mathbf{x},\mathbf{y}\}}\sum_{\mathbf{y}'\in X\setminus\{\mathbf{x}',\mathbf{x},\mathbf{y}\}} f_r(\mathbf{x},\mathbf{y},\mathbf{x}',\mathbf{y}')}_{(a)} \\
&\phantom{AAAA}+\underbrace{\frac{1}{16\pi^2}\mathbb{E}\sum_{\mathbf{x}\in X}\sum_{\mathbf{y}\in X\setminus\{\mathbf{x}\}}  \sum_{\mathbf{x}'\in \{\mathbf{x},\mathbf{y}\}}\sum_{\mathbf{y}'\in \{\mathbf{x},\mathbf{y}\}\setminus\{\mathbf{x}'\}} f_r(\mathbf{x},\mathbf{y},\mathbf{x}',\mathbf{y}')}_{(b)}\\
&\phantom{AAAA}+\underbrace{\frac{1}{16\pi^2}\mathbb{E}\sum_{\mathbf{x}\in X}\sum_{\mathbf{y}\in X\setminus\{\mathbf{x}\}}  \sum_{\mathbf{x}'\in X\setminus\{\mathbf{x},\mathbf{y}\}}\sum_{\mathbf{y}'\in \{\mathbf{x},\mathbf{y}\}\setminus\{\mathbf{x}'\}} f_r(\mathbf{x},\mathbf{y},\mathbf{x}',\mathbf{y}')}_{(c)}\\
&\phantom{AAAA}+\underbrace{\frac{1}{16\pi^2}\mathbb{E}\sum_{\mathbf{x}\in X}\sum_{\mathbf{y}\in X\setminus\{\mathbf{x}\}}  \sum_{\mathbf{x}'\in \{\mathbf{x},\mathbf{y}\}}\sum_{\mathbf{y}'\in X\setminus\{\mathbf{x},\mathbf{y},\mathbf{x}'\}} f_r(\mathbf{x},\mathbf{y},\mathbf{x}',\mathbf{y}')}_{(d)}.
\end{align*}
We handle these terms independently. For term $(a)$ we can directly apply the extended Slivnyak-Mecke Theorem given by Equation \ref{extended:SM},
\begin{align*}
\frac{1}{16\pi^2}\mathbb{E}\sum_{\mathbf{x}\in X}&\sum_{\mathbf{y}\in X\setminus\{\mathbf{x}\}} \sum_{\mathbf{x}'\in X\setminus\{\mathbf{x},\mathbf{y}\}}\sum_{\mathbf{y}'\in X\setminus\{\mathbf{x}',\mathbf{x},\mathbf{y}\}} f_r(\mathbf{x},\mathbf{y},\mathbf{x}',\mathbf{y}')\\
&= \frac{1}{16\pi^2}\int_{\mathbb{S}^2}\cdots\int_{\mathbb{S}^2} f_r(\mathbf{x},\mathbf{y},\mathbf{x}',\mathbf{y}')\rho(\mathbf{x})\rho(\mathbf{y})\rho(\mathbf{x}')\rho(\mathbf{y}')\lambda_{\mathbb{S}^2}(d\mathbf{x}) \lambda_{\mathbb{S}^2}(d\mathbf{y}) \lambda_{\mathbb{S}^2}(d\mathbf{x})'\lambda_{\mathbb{S}^2}(d\mathbf{y})'\\
&= \frac{1}{16\pi^2}\int_{\mathbb{S}^2}\cdots\int_{\mathbb{S}^2} \mathbbm{1}[d(\mathbf{x},\mathbf{y})\leq r] \mathbbm{1}[d(\mathbf{x}',\mathbf{y}')\leq r]\lambda_{\mathbb{S}^2}(d\mathbf{x}) \lambda_{\mathbb{S}^2}(d\mathbf{y}) \lambda_{\mathbb{S}^2}(d\mathbf{x})'\lambda_{\mathbb{S}^2}(d\mathbf{y})'\\
&= \frac{1}{16\pi^2}\left(\int_{\mathbb{S}^2}\int_{\mathbb{S}^2}\mathbbm{1}[d(\mathbf{x},\mathbf{y})\leq r]\lambda_{\mathbb{S}^2}(d\mathbf{x}) \lambda_{\mathbb{S}^2}(d\mathbf{y})\right)^2\\
&= \frac{1}{16\pi^2}\left(\int_{\mathbb{S}^2}\lambda_{\mathbb{S}^2}(B_{\mathbb{S}^2}(\mathbf{o},r))\lambda_{\mathbb{S}^2}(d\mathbf{y})\right)^2\\
&= \frac{1}{16\pi^2}\left(4\pi\lambda_{\mathbb{S}^2}(B_{\mathbb{S}^2}(\mathbf{o},r))\right)^2\\
&= \lambda_{\mathbb{S}^2}(B_{\mathbb{S}^2}(\mathbf{o},r))^2,
\end{align*}
where the penultimate equality follows since the area of the spherical cap is constant for a fixed geodesic radius for any centre, $\mathbf{0}$ indicates any arbitrary point in $\mathbb{S}^2$. Term $(b)$ can be handled in a similar manner as term $(a)$,
\begin{align*}
\frac{1}{16\pi^2}\mathbb{E}\sum_{\mathbf{x}\in X}&\sum_{\mathbf{y}\in X\setminus\{\mathbf{x}\}}  \sum_{\mathbf{x}'\in \{\mathbf{x},\mathbf{y}\}}\sum_{\mathbf{y}'\in \{\mathbf{x},\mathbf{y}\}\setminus\{\mathbf{x}'\}} f_r(\mathbf{x},\mathbf{y},\mathbf{x}',\mathbf{y}')\\
&=\frac{1}{16\pi^2}\mathbb{E}\sum_{\mathbf{x}\in X}\sum_{\mathbf{y}\in X\setminus\{\mathbf{x}\}} f_r(\mathbf{x},\mathbf{y},\mathbf{x},\mathbf{y}) + f_r(\mathbf{x},\mathbf{y},\mathbf{y},\mathbf{x})\\
&=\frac{1}{16\pi^2}\mathbb{E}\sum_{\mathbf{x}\in X}\sum_{\mathbf{y}\in X\setminus\{\mathbf{x}\}}\frac{\mathbbm{1}[d(\mathbf{x},\mathbf{y})\leq r]\mathbbm{1}[d(\mathbf{x},\mathbf{y})\leq r]}{\rho(\mathbf{x})\rho(\mathbf{y})\rho(\mathbf{x})\rho(\mathbf{y})} \\
&\phantom{\sum\sum\sum\sum\sum\sum\sum\sum}+ \frac{\mathbbm{1}[d(\mathbf{x},\mathbf{y})\leq r]\mathbbm{1}[d(\mathbf{y},\mathbf{x})\leq r]}{\rho(\mathbf{x})\rho(\mathbf{y})\rho(\mathbf{y})\rho(\mathbf{x})}\\
&=\frac{1}{8\pi^2} \mathbb{E}\sum_{\mathbf{x}\in X}\sum_{\mathbf{y}\in X\setminus\{\mathbf{x}\}}\frac{\mathbbm{1}[d(\mathbf{x},\mathbf{y})\leq r]}{\rho^2(\mathbf{x})\rho^2(\mathbf{y})}
\end{align*}
Hence by the extended Slivnyak-Mecke Theorem again we have,
\begin{align*}
\frac{1}{8\pi^2} \mathbb{E}\sum_{\mathbf{x}\in X}\sum_{\mathbf{y}\in X\setminus\{\mathbf{x}\}}\frac{\mathbbm{1}[d(\mathbf{x},\mathbf{y})\leq r]}{\rho^2(\mathbf{x})\rho^2(\mathbf{y})} &= \frac{1}{8\pi^2}\int_{\mathbb{S}^2}\int_{\mathbb{S}^2}\frac{\mathbbm{1}[d(\mathbf{x},\mathbf{y})\leq r]}{\rho^2(\mathbf{x})\rho^2(\mathbf{y})}\rho(\mathbf{x})\rho(\mathbf{y})\lambda_{\mathbb{S}^2}(d\mathbf{x}) \lambda_{\mathbb{S}^2}(d\mathbf{y})\\
&=\frac{1}{8\pi^2}\int_{\mathbb{S}^2}\int_{\mathbb{S}^2}\frac{\mathbbm{1}[d(\mathbf{x},\mathbf{y})\leq r]}{\rho(\mathbf{x})\rho(\mathbf{y})}\lambda_{\mathbb{S}^2}(d\mathbf{x}) \lambda_{\mathbb{S}^2}(d\mathbf{y})
\end{align*}
We now consider term $(c)$, 
\begin{align*}
\frac{1}{16\pi^2}\mathbb{E}\sum_{\mathbf{x}\in X}&\sum_{\mathbf{y}\in X\setminus\{\mathbf{x}\}}  \sum_{\mathbf{x}'\in X\setminus\{\mathbf{x},\mathbf{y}\}}\sum_{\mathbf{y}'\in \{\mathbf{x},\mathbf{y}\}\setminus\{\mathbf{x}'\}} f(\mathbf{x},\mathbf{y},\mathbf{x}',\mathbf{y}')\\
&=\frac{1}{16\pi^2}\mathbb{E}\sum_{\mathbf{x}\in X}\sum_{\mathbf{y}\in X\setminus\{\mathbf{x}\}}  \sum_{\mathbf{x}'\in X\setminus\{\mathbf{x},\mathbf{y}\}}\underbrace{f_r(\mathbf{x},\mathbf{y},\mathbf{x}',\mathbf{x})}_{(\rom{1})}+\underbrace{f_r(\mathbf{x},\mathbf{y},\mathbf{x}',\mathbf{y})}_{(\rom{2})}.
\end{align*}
We can handle these terms independently.
\begin{align*}
\frac{1}{16\pi^2}\mathbb{E}\sum_{\mathbf{x}\in X}&\sum_{\mathbf{y}\in X\setminus\{\mathbf{x}\}}  \sum_{\mathbf{x}'\in X\setminus\{\mathbf{x},\mathbf{y}\}} f_r(\mathbf{x},\mathbf{y},\mathbf{x}',\mathbf{x})\\
&= \frac{1}{16\pi^2}\mathbb{E}\sum_{\mathbf{x}\in X}\sum_{\mathbf{y}\in X\setminus\{\mathbf{x}\}}  \sum_{\mathbf{x}'\in X\setminus\{\mathbf{x},\mathbf{y}\}} \frac{\mathbbm{1}[d(\mathbf{x},\mathbf{y})\leq r]\mathbbm{1}[d(\mathbf{x}',\mathbf{x})\leq r]}{\rho(\mathbf{x})\rho(\mathbf{y})\rho(\mathbf{x}')\rho(\mathbf{x})}\\
&= \frac{1}{16\pi^2}\int_{\mathbb{S}^2}\int_{\mathbb{S}^2}\int_{\mathbb{S}^2}\frac{\mathbbm{1}[d(\mathbf{x},\mathbf{y})\leq r]\mathbbm{1}[d(\mathbf{x}',\mathbf{x})\leq r]}{\rho^2(\mathbf{x})\rho(\mathbf{y})\rho(\mathbf{x}')}\rho(\mathbf{x})\rho(\mathbf{y})\rho(\mathbf{x}')\lambda_{\mathbb{S}^2}(d\mathbf{x}) \lambda_{\mathbb{S}^2}(d\mathbf{y}) \lambda_{\mathbb{S}^2}(d\mathbf{x})'\\
&= \frac{1}{16\pi^2} \int_{\mathbb{S}^2}\int_{\mathbb{S}^2}\int_{\mathbb{S}^2}\frac{\mathbbm{1}[d(\mathbf{x},\mathbf{y})\leq r]\mathbbm{1}[d(\mathbf{x}',\mathbf{x})\leq r]}{\rho(\mathbf{x})}\lambda_{\mathbb{S}^2}(d\mathbf{x}) \lambda_{\mathbb{S}^2}(d\mathbf{y}) \lambda_{\mathbb{S}^2}(d\mathbf{x})'\\
&= \frac{1}{16\pi^2} \int_{\mathbb{S}^2}\int_{\mathbb{S}^2}\frac{\mathbbm{1}[d(\mathbf{x},\mathbf{y})\leq r]}{\rho(\mathbf{x})}\left(\int_{\mathbb{S}^2}\mathbbm{1}[d(\mathbf{x}',\mathbf{x})\leq r]\lambda_{\mathbb{S}^2}(d\mathbf{x})' \right)\lambda_{\mathbb{S}^2}(d\mathbf{x}) \lambda_{\mathbb{S}^2}(d\mathbf{y})\\
&= \frac{\lambda_{\mathbb{S}^2}(B_{\mathbb{S}^2}(\mathbf{o},r))}{16\pi^2}\int_{\mathbb{S}^2}\frac{1}{\rho(\mathbf{x})}\left(\int_{\mathbb{S}^2}\mathbbm{1}[d(\mathbf{x},\mathbf{y})\leq r] \lambda_{\mathbb{S}^2}(d\mathbf{y})\right) \lambda_{\mathbb{S}^2}(d\mathbf{x})\\
& = \frac{\lambda_{\mathbb{S}^2}(B_{\mathbb{S}^2}(\mathbf{o},r))^2}{16\pi^2}\int_{\mathbb{S}^2}\frac{1}{\rho(\mathbf{x})}\lambda_{\mathbb{S}^2}(d\mathbf{x}),
\end{align*}
where the second equality follows by the Slivnyak-Mecke Theorem. By an identical argument term \rom{2} is given by,
\begin{equation*}
\frac{1}{16\pi^2}\mathbb{E}\sum_{\mathbf{x}\in X}\sum_{\mathbf{y}\in X\setminus\{\mathbf{x}\}}  \sum_{\mathbf{x}'\in X\setminus\{\mathbf{x},\mathbf{y}\}} f_r(\mathbf{x},\mathbf{y},\mathbf{x}',\mathbf{y}) = \frac{\lambda_{\mathbb{S}^2}(B_{\mathbb{S}^2}(\mathbf{o},r))^2}{16\pi^2}\int_{\mathbb{S}^2}\frac{1}{\rho(\mathbf{y})}\lambda_{\mathbb{S}^2}(d\mathbf{y})
\end{equation*}
Term $(d)$ can be handled in an identical manner as term $(c)$. To see this consider the summation over $\mathbf{y}'$ in term $(d)$. Since $\mathbf{x}'\in\{\mathbf{x},\mathbf{y}\}$ this means that the set $X\setminus\{\mathbf{x},\mathbf{y},\mathbf{x}'\}$ is identical to $X\setminus\{\mathbf{x},\mathbf{y}\}$ and so the summations over $\mathbf{x}'$ and $\mathbf{y}'$ can be interchanged. Therefore,
\begin{equation*}
\frac{1}{16\pi^2}\mathbb{E}\sum_{\mathbf{x}\in X}\sum_{\mathbf{y}\in X\setminus\{\mathbf{x}\}}  \sum_{\mathbf{x}'\in \{\mathbf{x},\mathbf{y}\}}\sum_{\mathbf{y}'\in X\setminus\{\mathbf{x},\mathbf{y},\mathbf{x}'\}} f_r(\mathbf{x},\mathbf{y},\mathbf{x}',\mathbf{y}') = \frac{2\lambda_{\mathbb{S}^2}(B_{\mathbb{S}^2}(\mathbf{o},r))^2}{16\pi^2}\int_{\mathbb{S}^2}\frac{1}{\rho(\mathbf{x})}\lambda_{\mathbb{S}^2}(d\mathbf{x}).
\end{equation*}
Further we note that the area of spherical cap is given by $\lambda_{\mathbb{S}^2}(B_{\mathbb{S}^2}(\mathbf{o},r))=2\pi(1-\cos r)$. Thus collecting all the terms gives the form of the variance,
\begin{equation}
\text{Var}(\hat{K}_{\text{inhom}}(r)) = \frac{1}{8\pi^2}\int_{\mathbb{S}^2}\int_{\mathbb{S}^2}\frac{\mathbbm{1}[d(\mathbf{x},\mathbf{y})\leq r]}{\rho(\mathbf{x})\rho(\mathbf{y})}\lambda_{\mathbb{S}^2}(d\mathbf{x})\lambda_{\mathbb{S}^2}(d\mathbf{y}) +(1-\cos r)^2\int_{\mathbb{S}^2}\frac{1}{\rho(\mathbf{x})}\lambda_{\mathbb{S}^2}(d\mathbf{x}).
\end{equation}
\end{proof}

\begin{proof}
\subsection*{Variance of $\hat{F}_{\text{inhom}}(r)$}
Restating the estimator for $F_{\text{inhom}(r)}$,
\begin{equation*}
\hat{F}_{\text{inhom}}(r)=1-\frac{\sum_{\mathbf{p}\in P}\prod_{\mathbf{x}\in X\cap B_{\mathbb{S}^2}(\mathbf{p},r)}\left(1-\frac{\bar{\rho}}{\rho(\mathbf{x})}\right)}{|P|}.
\end{equation*}
Taking the variance using $\text{Var}(X)=\mathbb{E}[X^2]-\mathbb{E}^2[X]$,
\begin{align}
&\text{Var}(\hat{F}_{\text{inhom}}(r))=\frac{1}{|P|^2}\text{Var}\sum_{\mathbf{p}\in P}\prod_{\mathbf{x}\in X}\left(1-\frac{\bar{\rho}\mathbbm{1}[\mathbf{x}\in B_{\mathbb{S}^2}(\mathbf{p},r)]}{\rho(\mathbf{x})}\right)\nonumber\\
&=\frac{1}{|P|^2}\mathbb{E}\left(\sum_{\mathbf{p}\in P}\prod_{\mathbf{x}\in X}\left(1-\frac{\bar{\rho}\mathbbm{1}[\mathbf{x}\in B_{\mathbb{S}^2}(\mathbf{p},r)]}{\rho(\mathbf{x})}\right)\right)^2-\frac{1}{|P|^2}\mathbb{E}^2\sum_{\mathbf{p}\in P}\prod_{\mathbf{x}\in X}\left(1-\frac{\bar{\rho}\mathbbm{1}[\mathbf{x}\in B_{\mathbb{S}^2}(\mathbf{p},r)]}{\rho(\mathbf{x})}\right)\label{variance:F:inhom}
\end{align}
Dealing with each term independently, we have
\begin{align}
&\mathbb{E}\left(\sum_{\mathbf{p}\in P}\prod_{\mathbf{x}\in X}\left(1-\frac{\bar{\rho}\mathbbm{1}[\mathbf{x}\in B_{\mathbb{S}^2}(\mathbf{p},r)]}{\rho(\mathbf{x})}\right)\right)^2\nonumber\\
&\phantom{AAAA}=\mathbb{E}\sum_{\mathbf{p}\in P}\prod_{\mathbf{x}\in X}\left(1-\frac{\bar{\rho}\mathbbm{1}[\mathbf{x}\in B_{\mathbb{S}^2}(\mathbf{p},r)]}{\rho(\mathbf{x})}\right) \sum_{\mathbf{p}'\in P}\prod_{\mathbf{y}\in X}\left(1-\frac{\bar{\rho}\mathbbm{1}[\mathbf{y}\in B_{\mathbb{S}^2}(\mathbf{p}',r)]}{\rho(\mathbf{y})}\right)\nonumber\\
&\phantom{AAAA}=\sum_{\mathbf{p}\in P}\sum_{\mathbf{p}'\in P}\mathbb{E}\prod_{\mathbf{x}\in X}\left(1-\frac{\bar{\rho}\mathbbm{1}[\mathbf{x}\in B_{\mathbb{S}^2}(\mathbf{p},r)]}{\rho(\mathbf{x})}\right)\prod_{\mathbf{y}\in X}\left(1-\frac{\bar{\rho}\mathbbm{1}[\mathbf{y}\in B_{\mathbb{S}^2}(\mathbf{p}',r)]}{\rho(\mathbf{y})}\right)\label{product:expectation}
\end{align}
From the proof of Theorem 1 given by \cite{vanLieshout2011}, we have the following identity,
\begin{equation*}
\prod_{\mathbf{y}\in X}\left(1-\frac{\bar{\rho}\mathbbm{1}[\mathbf{y}\in B_{\mathbb{S}^2}(\mathbf{x},r)]}{\rho(\mathbf{y})}\right)= 1+\sum_{n=1}^{\infty}\frac{(-\bar{\rho})^n}{n!}\sum^{\neq}_{\mathbf{x}_1,\dots,\mathbf{x}_n\in X}\prod_{i=1}^n\frac{\mathbbm{1}[\mathbf{x}_i\in B_{\mathbb{S}^2}(\mathbf{x},r)]}{\rho(\mathbf{x}_i)},
\end{equation*}
and using the convention that a sum over an emptyset is $0$, i.e. $\sum_{k=1}^0=\sum_{x\subseteq\emptyset} = \sum_{\emptyset\in x}= 1$,
\begin{equation}\label{prod:identity}
\prod_{\mathbf{y}\in X}\left(1-\frac{\bar{\rho}\mathbbm{1}[\mathbf{y}\in B_{\mathbb{S}^2}(\mathbf{x},r)]}{\rho(\mathbf{y})}\right)= \sum_{n=0}^{\infty}\frac{(-\bar{\rho})^n}{n!}\sum^{\neq}_{\mathbf{x}_1,\dots,\mathbf{x}_n\in X}\prod_{i=1}^n\frac{\mathbbm{1}[\mathbf{x}_i\in B_{\mathbb{S}^2}(\mathbf{x},r)]}{\rho(\mathbf{x}_i)}.
\end{equation}
Taking just the expectation from Equation \ref{product:expectation} we expand the first product over $\mathbf{x}$ using the previous identity to give,
\begin{align}
&\mathbb{E}\prod_{\mathbf{x}\in X}\left(1-\frac{\bar{\rho}\mathbbm{1}[\mathbf{x}\in B_{\mathbb{S}^2}(\mathbf{p},r)]}{\rho(\mathbf{x})}\right)\prod_{\mathbf{y}\in X}\left(1-\frac{\bar{\rho}\mathbbm{1}[\mathbf{y}\in B_{\mathbb{S}^2}(\mathbf{p}',r)]}{\rho(\mathbf{y})}\right)\nonumber\\
&=\mathbb{E}\left(\sum_{n=0}^{\infty}\frac{(-\bar{\rho})^n}{n!}\sum^{\neq}_{\mathbf{x}_1,\dots,\mathbf{x}_n\in X}\prod_{i=1}^n\frac{\mathbbm{1}[\mathbf{x}_i\in B_{\mathbb{S}^2}(\mathbf{p},r)]}{\rho(\mathbf{x}_i)}\right)\prod_{\mathbf{y}\in X}\left(1-\frac{\bar{\rho}\mathbbm{1}[\mathbf{y}\in B_{\mathbb{S}^2}(\mathbf{p}',r)]}{\rho(\mathbf{y})}\right)\nonumber\\
&=\mathbb{E}\sum_{n=0}^{\infty}\frac{(-\bar{\rho})^n}{n!}\sum^{\neq}_{\mathbf{x}_1,\dots,\mathbf{x}_n\in X}\left(\prod_{i=1}^n\frac{\mathbbm{1}[\mathbf{x}_i\in B_{\mathbb{S}^2}(\mathbf{p},r)]}{\rho(\mathbf{x}_i)}\right.\nonumber
\\&\phantom{AAAA}\left.\prod_{\mathbf{y}\in \{X\setminus\{\mathbf{x}_1,\dots,\mathbf{x}_n\},\{\mathbf{x}_1,\dots,\mathbf{x}_n\}\}}\left(1-\frac{\bar{\rho}\mathbbm{1}[\mathbf{y}\in B_{\mathbb{S}^2}(\mathbf{p}',r)]}{\rho(\mathbf{y})}\right)\right)\nonumber\\
&=\mathbb{E}\sum_{n=0}^{\infty}\frac{(-\bar{\rho})^n}{n!}\sum^{\neq}_{\mathbf{x}_1,\dots,\mathbf{x}_n\in X}\left(\prod_{i=1}^n\frac{\mathbbm{1}[\mathbf{x}_i\in B_{\mathbb{S}^2}(\mathbf{p},r)]}{\rho(\mathbf{x}_i)}\prod_{\mathbf{y}\in \{\mathbf{x}_1,\dots,\mathbf{x}_n\}}\left(1-\frac{\bar{\rho}\mathbbm{1}[\mathbf{y}\in B_{\mathbb{S}^2}(\mathbf{p}',r)]}{\rho(\mathbf{y})}\right)\right.\nonumber\\
&\phantom{=AAAA}\left.\prod_{\mathbf{y}\in X\setminus\{\mathbf{x}_1,\dots,\mathbf{x}_n\}}\left(1-\frac{\bar{\rho}\mathbbm{1}[\mathbf{y}\in B_{\mathbb{S}^2}(\mathbf{p}',r)]}{\rho(\mathbf{y})}\right)\right)\nonumber\\
&=\mathbb{E}\sum_{n=0}^{\infty}\frac{(-\bar{\rho})^n}{n!}\sum^{\neq}_{\mathbf{x}_1,\dots,\mathbf{x}_n\in X}\left(\prod_{i=1}^n\frac{\mathbbm{1}[\mathbf{x}_i\in B_{\mathbb{S}^2}(\mathbf{p},r)]}{\rho(\mathbf{x}_i)}\left(1-\frac{\bar{\rho}\mathbbm{1}[\mathbf{x}_i\in B_{\mathbb{S}^2}(\mathbf{p}',r)]}{\rho(\mathbf{x}_i)}\right)\right.\nonumber\\
&\phantom{=AAAA}\left.\prod_{\mathbf{y}\in X\setminus\{\mathbf{x}_1,\dots,\mathbf{x}_n\}}\left(1-\frac{\bar{\rho}\mathbbm{1}[\mathbf{y}\in B_{\mathbb{S}^2}(\mathbf{p}',r)]}{\rho(\mathbf{y})}\right)\right)\nonumber\\
&=\sum_{n=0}^{\infty}\frac{(-\bar{\rho})^n}{n!}\mathbb{E}\sum^{\neq}_{\mathbf{x}_1,\dots,\mathbf{x}_n\in X}\left(\prod_{i=1}^n\frac{\mathbbm{1}[\mathbf{x}_i\in B_{\mathbb{S}^2}(\mathbf{p},r)]}{\rho(\mathbf{x}_i)}\left(1-\frac{\bar{\rho}\mathbbm{1}[\mathbf{x}_i\in B_{\mathbb{S}^2}(\mathbf{p}',r)]}{\rho(\mathbf{x}_i)}\right)\right.\nonumber\\
&\phantom{=AAAA}\left.\prod_{\mathbf{y}\in X\setminus\{\mathbf{x}_1,\dots,\mathbf{x}_n\}}\left(1-\frac{\bar{\rho}\mathbbm{1}[\mathbf{y}\in B_{\mathbb{S}^2}(\mathbf{p}',r)]}{\rho(\mathbf{y})}\right)\right)\nonumber\\
\intertext{Using the extended Slivnyak-Mecke Theorem,}
&=\sum_{n=0}^{\infty}\frac{(-\bar{\rho})^n}{n!}\overbrace{\int_{\mathbb{S}^2}\cdots\int_{\mathbb{S}^2}}^{n}\nonumber\\
&\phantom{AAAA}\mathbb{E}\left(\prod_{i=1}^n\frac{\mathbbm{1}[\mathbf{x}_i\in B_{\mathbb{S}^2}(\mathbf{p},r)]}{\rho(\mathbf{x}_i)}\left(1-\frac{\bar{\rho}\mathbbm{1}[\mathbf{x}_i\in B_{\mathbb{S}^2}(\mathbf{p}',r)]}{\rho(\mathbf{x}_i)}\right)\right.\nonumber\\
&\phantom{AAAAAAAA}\left.\prod_{\mathbf{y}\in X}\left(1-\frac{\bar{\rho}\mathbbm{1}[\mathbf{y}\in B_{\mathbb{S}^2}(\mathbf{p}',r)]}{\rho(\mathbf{y})}\right)\right)\prod_{i=1}^n\rho(\mathbf{x}_i)\lambda_{\mathbb{S}^2}(d\mathbf{x}_i)\nonumber\\
&=\sum_{n=0}^{\infty}\frac{(-\bar{\rho})^n}{n!}\overbrace{\int_{\mathbb{S}^2}\cdots\int_{\mathbb{S}^2}}^{n}\nonumber\\
&\phantom{AAAA}\prod_{i=1}^n\frac{\mathbbm{1}[\mathbf{x}_i\in B_{\mathbb{S}^2}(\mathbf{p},r)]}{\rho(\mathbf{x}_i)}\left(1-\frac{\bar{\rho}\mathbbm{1}[\mathbf{x}_i\in B_{\mathbb{S}^2}(\mathbf{p}',r)]}{\rho(\mathbf{x}_i)}\right)\nonumber\\
&\phantom{AAAAAAAA}\mathbb{E}\prod_{\mathbf{y}\in X}\left(1-\frac{\bar{\rho}\mathbbm{1}[\mathbf{y}\in B_{\mathbb{S}^2}(\mathbf{p}',r)]}{\rho(\mathbf{y})}\right)\prod_{i=1}^n\rho(\mathbf{x}_i)\lambda_{\mathbb{S}^2}(d\mathbf{x}_i)\nonumber\\
&=\sum_{n=0}^{\infty}\frac{(-\bar{\rho})^n}{n!}\mathbb{E}\prod_{\mathbf{y}\in X}\left(1-\frac{\bar{\rho}\mathbbm{1}[\mathbf{y}\in B_{\mathbb{S}^2}(\mathbf{p}',r)]}{\rho(\mathbf{y})}\right)\nonumber\\
&\phantom{AAAA}\overbrace{\int_{\mathbb{S}^2}\cdots\int_{\mathbb{S}^2}}^{n}\prod_{i=1}^n\frac{\mathbbm{1}[\mathbf{x}_i\in B_{\mathbb{S}^2}(\mathbf{p},r)]}{\rho(\mathbf{x}_i)}\left(1-\frac{\bar{\rho}\mathbbm{1}[\mathbf{x}_i\in B_{\mathbb{S}^2}(\mathbf{p}',r)]}{\rho(\mathbf{x}_i)}\right)\prod_{i=1}^n\rho(\mathbf{x}_i)\lambda_{\mathbb{S}^2}(d\mathbf{x}_i)\nonumber\\
\intertext{Note the expectation is just the definition of the generating functional of X which does not depend on the point $\mathbf{p}'$ (see proof of Theorem 1 by \cite{vanLieshout2011}),}
&=\sum_{n=0}^{\infty}\frac{(-\bar{\rho})^n}{n!}G(1-u^{\mathbf{y}}_r)\overbrace{\int_{\mathbb{S}^2}\cdots\int_{\mathbb{S}^2}}^{n}\prod_{i=1}^n\frac{\mathbbm{1}[\mathbf{x}_i\in B_{\mathbb{S}^2}(\mathbf{p},r)]}{\rho(\mathbf{x}_i)}\nonumber\\
&\phantom{AAAA}\left(1-\frac{\bar{\rho}\mathbbm{1}[\mathbf{x}_i\in B_{\mathbb{S}^2}(\mathbf{p}',r)]}{\rho(\mathbf{x}_i)}\right)\prod_{i=1}^n\rho(\mathbf{x}_i)\lambda_{\mathbb{S}^2}(d\mathbf{x}_i)\nonumber\\
&=\sum_{n=0}^{\infty}\frac{(-\bar{\rho})^n}{n!}G(1-u^{\mathbf{y}}_r)\nonumber\\
&\phantom{AAAA}\overbrace{\int_{\mathbb{S}^2\cap B_{\mathbb{S}^2}(\mathbf{p},r)}\cdots\int_{\mathbb{S}^2\cap B_{\mathbb{S}^2}(\mathbf{p},r)}}^{n}\prod_{i=1}^n\left(1-\frac{\bar{\rho}\mathbbm{1}[\mathbf{x}_i\in B_{\mathbb{S}^2}(\mathbf{p}',r)]}{\rho(\mathbf{x}_i)}\right)\lambda_{\mathbb{S}^2}(d\mathbf{x})_1\cdots \lambda_{\mathbb{S}^2}(d\mathbf{x})_n\nonumber\\
&=\sum_{n=0}^{\infty}\frac{(-\bar{\rho})^n}{n!}G(1-u^{\mathbf{y}}_r)\left(\int_{\mathbb{S}^2\cap B_{\mathbb{S}^2}(\mathbf{p},r)} \left(1-\frac{\bar{\rho}\mathbbm{1}[\mathbf{x}\in B_{\mathbb{S}^2}(\mathbf{p}',r)]}{\rho(\mathbf{x})}\right)\lambda_{\mathbb{S}^2}(d\mathbf{x}) \right)^n\nonumber\\
&=G(1-u^{\mathbf{y}}_r) \sum_{n=0}^{\infty}\frac{(-\bar{\rho})^n}{n!} \left(\lambda_{\mathbb{S}^2}(\mathbb{S}^2_{B_{\mathbb{S}^2}(\mathbf{p},r)})-\int_{\mathbb{S}^2_{B_{\mathbb{S}^2}(\mathbf{p},r)\cap B_{\mathbb{S}^2}(\mathbf{p}',r)}}\frac{\bar{\rho}}{\rho(\mathbf{x})}\lambda_{\mathbb{S}^2}(d\mathbf{x})\right)^n\nonumber\\
\intertext{Using the series definition for the exponential function,}
&=G(1-u^{\mathbf{y}}_r) \exp\left(-\bar{\rho}\left(\lambda_{\mathbb{S}^2}(\mathbb{S}^2_{B_{\mathbb{S}^2}(\mathbf{p},r)})-\int_{\mathbb{S}^2_{B_{\mathbb{S}^2}(\mathbf{o},r)\cap B_{\mathbb{S}^2}(\mathbf{p}',r)}}\frac{\bar{\rho}}{\rho(\mathbf{x})}\lambda_{\mathbb{S}^2}(d\mathbf{x})\right)\right)\nonumber\\
&=G(1-u^{\mathbf{y}}_r) \exp\left(\bar{\rho}\int_{\mathbb{S}^2_{B_{\mathbb{S}^2}(\mathbf{p},r)\cap B_{\mathbb{S}^2}(\mathbf{p}',r)}}\frac{\bar{\rho}}{\rho(\mathbf{x})}\lambda_{\mathbb{S}^2}(d\mathbf{x}) -\bar{\rho}\lambda_{\mathbb{S}^2}(B_{\mathbb{S}^2}(\mathbf{o},r))\right)\nonumber\\
&=G(1-u^{\mathbf{y}}_r)\exp\left(-\bar{\rho}\lambda_{\mathbb{S}^2}(B_{\mathbb{S}^2}(\mathbf{o},r))\right) \exp\left(\int_{\mathbb{S}^2_{B_{\mathbb{S}^2}(\mathbf{p},r)\cap B_{\mathbb{S}^2}(\mathbf{p}',r)}}\frac{\bar{\rho}^2}{\rho(\mathbf{x})}\lambda_{\mathbb{S}^2}(d\mathbf{x})\right)\nonumber
\end{align}
Substituting this into the first term of Equation \ref{variance:F:inhom} gives,
\begin{align*}
&\frac{1}{|P|^2}\mathbb{E}\left(\sum_{\mathbf{p}\in P}\prod_{\mathbf{x}\in X}\left(1-\frac{\bar{\rho}\mathbbm{1}[\mathbf{x}\in B_{\mathbb{S}^2}(\mathbf{p},r)]}{\rho(\mathbf{x})}\right)\right)^2\\
&= G(1-u^{\mathbf{y}}_r)\exp\left(-\bar{\rho}\lambda_{\mathbb{S}^2}(B_{\mathbb{S}^2}(\mathbf{o},r))\right)\frac{1}{|P|^2}\sum_{\mathbf{p}\in P}\sum_{\mathbf{p}'\in P} \exp\left(\int_{\mathbb{S}^2_{B_{\mathbb{S}^2}(\mathbf{p},r)\cap B_{\mathbb{S}^2}(\mathbf{p}',r)}}\frac{\bar{\rho}^2}{\rho(\mathbf{x})}\lambda_{\mathbb{S}^2}(d\mathbf{x})\right)
\end{align*}
The second term of Equation \ref{variance:F:inhom}, by unbiasedness of $\hat{F}_{\text{inhom}}$ shown in Theorem \ref{expect:FS}, gives,
\begin{align*}
\mathbb{E}^2\frac{1}{|P|}\sum_{\mathbf{p}\in P}\prod_{\mathbf{x}\in X}\left(1-\frac{\bar{\rho}\mathbbm{1}[\mathbf{x}\in B_{\mathbb{S}^2}(\mathbf{p},r)]}{\rho(\mathbf{x})}\right)&=\mathbb{E}^2\left[1-\hat{F}_{\text{inhom}}(r)\right]\\
&=\left(1-\mathbb{E}[\hat{F}_{\text{inhom}}(r)]\right)^2\\
&=G^2(1-u^{\mathbf{y}}_r),
\end{align*}
where the generating functional does not depend on $\mathbf{y}$. Thus the variance is,
\begin{equation*}
\begin{split}
&G(1-u^{\mathbf{y}}_r)\exp\left(-\bar{\rho}\lambda_{\mathbb{S}^2}(B_{\mathbb{S}^2}(\mathbf{o},r))\right)\\
&\phantom{AAAA}\times\frac{1}{|P|^2}\sum_{\mathbf{p}\in P}\sum_{\mathbf{p}'\in P} \exp\left(\int_{\mathbb{S}^2_{B_{\mathbb{S}^2}(\mathbf{p},r)\cap B_{\mathbb{S}^2}(\mathbf{p}',r)}}\frac{\bar{\rho}^2}{\rho(\mathbf{x})}\lambda_{\mathbb{S}^2}(d\mathbf{x})\right)-G^2(1-u^{\mathbf{y}}_r),
\end{split}
\end{equation*}
For a Poisson process, $G(1-u^{\mathbf{y}}_r)=\exp\left(-\bar{\rho}\lambda_{\mathbb{S}^2}(\mathcal{S}^2_{B_{\mathbb{S}^2}(\mathbf{o},r)})\right)$, thus the variance is,
\begin{equation*}
\begin{split}
&\exp\left(-2\bar{\rho}\lambda_{\mathbb{S}^2}(B_{\mathbb{S}^2}(\mathbf{o},r))\right)\frac{1}{|P|^2}\sum_{\mathbf{p}\in P}\sum_{\mathbf{p}'\in P} \exp\left(\int_{\mathbb{S}^2_{B_{\mathbb{S}^2}(\mathbf{p},r)\cap B_{\mathbb{S}^2}(\mathbf{p}',r)}}\frac{\bar{\rho}^2}{\rho(\mathbf{x})}\lambda_{\mathbb{S}^2}(d\mathbf{x})\right)\\
&\phantom{AAAA}-\exp\left(-2\bar{\rho}\lambda_{\mathbb{S}^2}(\mathcal{S}^2_{B_{\mathbb{S}^2}(\mathbf{o},r)})\right).
\end{split}
\end{equation*}
\end{proof}
Before the proof for the variance of $\hat{H}_{\text{inhom}}(r)$ we introduce the exponential integral. The exponential integral, denoted $\text{Ei}(x)$, is defined as the following integral,
\begin{equation*}
\text{Ei}(x) = -\int_{-x}^{\infty}\frac{\mathrm{e}^{-t}}{t} dt.
\end{equation*}
It can then be shown that the exponential integral has the following infinite series representation,
\begin{equation}\label{Ei:series}
\text{Ei}(x) = \gamma +\log(x) +\sum_{k=1}^{\infty} \frac{x^k}{k\cdot k!},
\end{equation}
where $\gamma$ is known as the Euler-Mascheroni constant and defined as,
\begin{equation*}
\gamma = \lim_{n\rightarrow \infty} \left(-\log(n)+\sum_{k=1}^n\frac{1}{k}\right).
\end{equation*}
The variance for $\hat{H}_{\text{inhom}}(r)$ will be given in terms of $\text{Ei}(x)$.
\begin{proof}
\subsection*{Variance of $\hat{H}_{\text{inhom}}(r)$}
Restating the estimator for $\hat{H}_{\text{inhom}}(r),$
\begin{equation*}
\hat{H}_{\text{inhom}}(r) = 1 -\frac{\sum_{\mathbf{x}\in X \prod_{\mathbf{y}\in X\setminus \{\mathbf{x}\}}}\left(1-\frac{\bar{\rho}\mathbbm{1}[\mathbf{y}\in B_{\mathbb{S}^2}(\mathbf{x},r)]}{\rho(\mathbf{y})}\right)}{N_X(\mathbb{S}^2)},
\end{equation*}
we note that this estimator is only well defined when $\mathbbm{1}[N_X(\mathbb{S}^2)>0]$. In the event that $N_X(\mathbb{S}^2)=0$ we shall define $\hat{H}_{\text{inhom}}(r) = 0,$ in which case we can write our estimator as,
\begin{equation*}
\hat{H}_{\text{inhom}}(r) = \mathbbm{1}[N_X(\mathbb{S}^2)>0]\left(1 -\frac{\sum_{\mathbf{x}\in X \prod_{\mathbf{y}\in X\setminus \{\mathbf{x}\}}}\left(1-\frac{\bar{\rho}\mathbbm{1}[\mathbf{y}\in B_{\mathbb{S}^2}(\mathbf{x},r)]}{\rho(\mathbf{y})}\right)}{N_X(\mathbb{S}^2)}\right).
\end{equation*}
By using the law of total variance, conditioning on $N_X(\mathbb{S}^2)=n$,
\begin{equation}
\text{Var}(\hat{H}_{\text{inhom}}(r)) = \underbrace{\mathbb{E}[\text{Var}(\hat{H}_{\text{inhom}}(r)|N_X(\mathbb{S}^2=n)]}_{(1)} + \underbrace{\text{Var}(\mathbb{E}[\hat{H}_{\text{inhom}}(r)|N_X(\mathbb{S}^2)=n])}_{(2)}.\label{variance:total:law:1}
\end{equation}
The variance in term $(1)$ is,
\begin{align}
\text{Var}(\hat{H}_{\text{inhom}}(r)|N_X(\mathbb{S}^2=n)
&=\mathbbm{1}[n>0]\text{Var}\left(1+\frac{1}{n}\sum_{i=1}^n \prod_{\substack{j=1\\j\neq i}}^n\left(1-\frac{\bar{\rho}\mathbbm{1}[\mathbf{X}_j\in B(\mathbf{X}_i,r)]}{\rho(\mathbf{X}_j)}\right)\right),\nonumber\\
&=\frac{\mathbbm{1}[n>0]}{n^2}\text{Var}\left(\sum_{i=1}^n \prod_{\substack{j=1\\j\neq i}}^n\left(1-\frac{\bar{\rho}\mathbbm{1}[\mathbf{X}_j\in B(\mathbf{X}_i,r)]}{\rho(\mathbf{X}_j)}\right)\right),\nonumber\\
\intertext{where $\mathbf{X}_k,\;k=1,\dots,n$ are independently distributed points with density proportional to $\rho(\mathbf{x}_k)$, by definition of a Poisson process. We use the identity $\text{Var}(\sum_{i=1}^n X_i)=\sum_{i=1}^n$ $\sum_{j=1}^n\text{Cov}(X_i,X_j)$ and define $L_i=\prod_{\substack{j=1\\j\neq i}}^n\left(1-\frac{\bar{\rho}\mathbbm{1}[\mathbf{X}_j\in B(\mathbf{X}_i,r)]}{\rho(\mathbf{X}_j)}\right)$,}
&=\frac{\mathbbm{1}[n>0]}{n^2}\sum_{p=1}^n\sum_{q=1}^n \text{Cov}\left(L_p,L_q\right)\nonumber\\
&=\frac{\mathbbm{1}[n>0]}{n^2}\sum_{p=1}^n\sum_{q=1}^n \Big(\mathbb{E}[L_pL_q]-\mathbb{E}[L_p]\mathbb{E}[L_q]\Big)\nonumber\\
&=\frac{\mathbbm{1}[n>0]}{n^2}\sum_{p,q\in\{1,\dots,n\}}^{\neq}\left(\underbrace{\mathbb{E}[L_pL_q]}_{(a)}-\underbrace{\mathbb{E}[L_p]\mathbb{E}[L_q]}_{(b)}\right)\nonumber\\
&\phantom{AAAA}+\frac{\mathbbm{1}[n>0]}{n^2}\sum_{p=1}^n\left(\underbrace{\mathbb{E}[L_p^2]}_{(c)}-\underbrace{\mathbb{E}^2[L_p]}_{(d)}\right)\label{covariance:expansion}
\end{align}
Taking term $(a)$,
\begin{align*}
&\mathbb{E}[L_pL_q]=\mathbb{E}\left[\prod_{\substack{j=1\\j\neq p}}^n\left(1-\frac{\bar{\rho}\mathbbm{1}[\mathbf{X}_j\in B_{\mathbb{S}^2}(\mathbf{X}_p,r)]}{\rho(\mathbf{X}_j)}\right)\prod_{\substack{j=1\\j\neq q}}^n\left(1-\frac{\bar{\rho}\mathbbm{1}[\mathbf{X}_j\in B_{\mathbb{S}^2}(\mathbf{X}_q,r)]}{\rho(\mathbf{X}_j)}\right)\right]\\
\intertext{using iterated expectation, conditioning on $\mathbf{X_p}=\mathbf{x}_p,\mathbf{X_q}=\mathbf{x}_q$,}
&=\mathbb{E}\left[\mathbb{E}\left[\prod_{\substack{j=1\\j\neq p}}^n\left(1-\frac{\bar{\rho}\mathbbm{1}[\mathbf{X}_j\in B_{\mathbb{S}^2}(\mathbf{X}_p,r)]}{\rho(\mathbf{X}_j)}\right)\right.\right.\\
&\phantom{AAAAAAAAAAAA}\left.\left.\prod_{\substack{j=1\\j\neq q}}^n\left(1-\frac{\bar{\rho}\mathbbm{1}[\mathbf{X}_j\in B_{\mathbb{S}^2}(\mathbf{X}_q,r)]}{\rho(\mathbf{X}_j)}\right)\Bigg|\mathbf{X_p}=\mathbf{x}_p,\mathbf{X_q}=\mathbf{x}_q\right]\right]\\
&=\mathbb{E}\left[\left(1-\frac{\bar{\rho}\mathbbm{1}[\mathbf{X}_p\in B_{\mathbb{S}^2}(\mathbf{X}_q,r)]}{\rho(\mathbf{X}_p)}\right)\left(1-\frac{\bar{\rho}\mathbbm{1}[\mathbf{X}_q\in B_{\mathbb{S}^2}(\mathbf{X}_p,r)]}{\rho(\mathbf{X}_q)}\right)\right.\\
&\times\left.\mathbb{E}\left[\prod_{\substack{j=1\\j\neq p,q}}^n \left(1-\frac{\bar{\rho}\mathbbm{1{}}[\mathbf{X}_j\in B_{\mathbb{S}^2}(\mathbf{X}_p,r)]}{\rho(\mathbf{X}_j)}\right)\left(1-\frac{\bar{\rho}\mathbbm{1}[\mathbf{X}_j\in B_{\mathbb{S}^2}(\mathbf{X}_q,r)]}{\rho(\mathbf{X}_j)}\right) \Bigg|\mathbf{X_p}=\mathbf{x}_p,\mathbf{X_q}=\mathbf{x}_q\right]\right]
\end{align*}
The expectation conditioned on $(\mathbf{X}_q,\mathbf{X}_p)$ is,
\begin{align*}
&\mathbb{E}\left[\prod_{\substack{j=1\\j\neq p,q}}^n \left(1-\frac{\bar{\rho}\mathbbm{1}[\mathbf{X}_j\in B_{\mathbb{S}^2}(\mathbf{X}_p,r)]}{\rho(\mathbf{X}_j)}\right)\left(1-\frac{\bar{\rho}\mathbbm{1}[\mathbf{X}_j\in B_{\mathbb{S}^2}(\mathbf{X}_q,r)]}{\rho(\mathbf{X}_j)}\right) \Bigg|\mathbf{X_p}=\mathbf{x}_p,\mathbf{X_q}=\mathbf{x}_q\right]\\
&\phantom{AAAA}=\overbrace{\int_{\mathbb{S}^2}\cdots\int_{\mathbb{S}^2}}^{n-2}\prod_{\substack{j=1\\j\neq p,q}}^n \left(1-\frac{\bar{\rho}\mathbbm{1}[\mathbf{x}_j\in B_{\mathbb{S}^2}(\mathbf{x}_p,r)]}{\rho(\mathbf{x}_j)}\right)\left(1-\frac{\bar{\rho}\mathbbm{1}[\mathbf{x}_j\in B_{\mathbb{S}^2}(\mathbf{x}_q,r)]}{\rho(\mathbf{x}_j)}\right)\frac{\rho(\mathbf{x}_j)}{\mu(\mathbb{S}^2)}\lambda_{\mathbb{S}^2}(d\mathbf{x}_j)\\
&\phantom{AAAA}=\left(\int_{\mathbb{S}^2} \left(1-\frac{\bar{\rho}\mathbbm{1}[\mathbf{x}\in B_{\mathbb{S}^2}(\mathbf{x}_p,r)]}{\rho(\mathbf{x})}\right)\left(1-\frac{\bar{\rho}\mathbbm{1}[\mathbf{x}\in B_{\mathbb{S}^2}(\mathbf{x}_q,r)]}{\rho(\mathbf{x})}\right)\frac{\rho(\mathbf{x})}{\mu(\mathbb{S}^2)}\lambda_{\mathbb{S}^2}(d\mathbf{x})\right)^{n-2}\\
&\phantom{AAAA}=\left(\int_{\mathbb{S}^2} \frac{\rho(\mathbf{x})}{\mu(\mathbb{S}^2)} -\frac{\bar{\rho}}{\mu(\mathbb{S}^2)}\mathbbm{1}[\mathbf{x}\in B_{\mathbb{S}^2}(\mathbf{x}_p,r)] -\frac{\bar{\rho}}{\mu(\mathbb{S}^2)}\mathbbm{1}[\mathbf{x}\in B_{\mathbb{S}^2}(\mathbf{x}_q,r)]\right.\\
&\phantom{AAAAAAAA}\left.+\frac{\bar{\rho}^2}{\mu(\mathbb{S}^2)}\frac{\mathbbm{1}[\mathbf{x}\in B_{\mathbb{S}^2}(\mathbf{x}_p,r),\mathbf{x}\in B_{\mathbb{S}^2}(\mathbf{x}_q,r)]}{\rho(\mathbf{x})}  \lambda_{\mathbb{S}^2}(d\mathbf{x}) \right)^{n-2}\\
&\phantom{AAAA}=\left(1-\frac{2\bar{\rho}}{\mu(\mathbb{S}^2)}2\pi(1-\cos r)+\frac{\bar{\rho}^2}{\mu(\mathbb{S}^2)}\int_{B_{\mathbb{S}^2}(\mathbf{x}_p,r)\cap B_{\mathbb{S}^2}(\mathbf{x}_q,r)}\frac{1}{\rho(\mathbf{x})}\lambda_{\mathbb{S}^2}(d\mathbf{x})\right)^{n-2}
\end{align*}
We then define $A_1(\mathbf{x}_p,\mathbf{x}_q)\equiv 1-\frac{2\bar{\rho}}{\mu(\mathbb{S}^2)}2\pi(1-\cos r)+\frac{\bar{\rho}^2}{\mu(\mathbb{S}^2)}\int_{  B_{\mathbb{S}^2}(\mathbf{x}_p,r)\cap B_{\mathbb{S}^2}(\mathbf{x}_q,r)}\frac{1}{\rho(\mathbf{x})}\lambda_{\mathbb{S}^2}(d\mathbf{x})$ and returning to $\mathbb{E}[L_pL_q],$
\begin{align*}
\mathbb{E}[L_pL_q] &= \mathbb{E}\left[\left(1-\frac{\bar{\rho}\mathbbm{1}[\mathbf{X}_p\in B_{\mathbb{S}^2}(\mathbf{X}_q,r)]}{\rho(\mathbf{X}_p)}\right)\left(1-\frac{\bar{\rho}\mathbbm{1}[\mathbf{X}_q\in B_{\mathbb{S}^2}(\mathbf{X}_p,r)]}{\rho(\mathbf{X}_q)}\right)A_1^{n-2}(\mathbf{X}_p,\mathbf{X}_q)\right]\\
&=\int_{\mathbb{S}^2}\int_{\mathbb{S}^2}\left(1-\frac{\bar{\rho}\mathbbm{1}[\mathbf{x}_p\in B_{\mathbb{S}^2}(\mathbf{x}_q,r)]}{\rho(\mathbf{x}_p)}\right)\left(1-\frac{\bar{\rho}\mathbbm{1}[\mathbf{x}_q\in B_{\mathbb{S}^2}(\mathbf{x}_p,r)]}{\rho(\mathbf{x}_q)}\right)\\
&\phantom{AAAAAAAA}A_1^{n-2}(\mathbf{x}_p,\mathbf{x}_q)\frac{\rho(\mathbf{x}_p)\rho(\mathbf{x}_p)}{\mu^2(\mathbb{S}^2)}\lambda_{\mathbb{S}^2}(d\mathbf{x}_p) \lambda_{\mathbb{S}^2}(d\mathbf{x}_q)\\
&=\frac{1}{\mu^2(\mathbb{S}^2)}\int_{\mathbb{S}^2}\int_{\mathbb{S}^2}\left(\rho(\mathbf{x}_p)-\bar{\rho}\mathbbm{1}[\mathbf{x}_p\in B_{\mathbb{S}^2}(\mathbf{x}_q,r)]\right)\left(\rho(\mathbf{x}_q)-\bar{\rho}\mathbbm{1}[\mathbf{x}_q\in B_{\mathbb{S}^2}(\mathbf{x}_p,r)]\right)\\
&\phantom{AAAAAAAA}A_1^{n-2}(\mathbf{x}_p,\mathbf{x}_q)\lambda_{\mathbb{S}^2}(d\mathbf{x}_p) \lambda_{\mathbb{S}^2}(d\mathbf{x}_q)\\
&=\frac{1}{\mu^2(\mathbb{S}^2)}\int_{\mathbb{S}^2}\int_{\mathbb{S}^2}\left(\rho(\mathbf{x})-\bar{\rho}\mathbbm{1}[\mathbf{x}\in B_{\mathbb{S}^2}(\mathbf{y},r)]\right)\left(\rho(\mathbf{y})-\bar{\rho}\mathbbm{1}[\mathbf{y}\in B_{\mathbb{S}^2}(\mathbf{x},r)]\right)\\
&\phantom{AAAAAAAA}A_1^{n-2}(\mathbf{x},\mathbf{y})\lambda_{\mathbb{S}^2}(d\mathbf{x}) \lambda_{\mathbb{S}^2}(d\mathbf{y})
\end{align*}
Then taking term $(b)$ of Equation \ref{covariance:expansion} and examining one of the expectations,
\begin{align*}
\mathbb{E}[L_p]&=\mathbb{E}\left[\prod_{\substack{j=1\\j\neq p}}^n\left(1-\frac{\bar{\rho}\mathbbm{1}[\mathbf{X}_j\in B_{\mathbb{S}^2}(\mathbf{X}_p,r)]}{\rho(\mathbf{X}_j)}\right)\right]\\
\intertext{using iterated expectation conditioning on $\mathbf{X}_p=\mathbf{x}_p$,}
&=\mathbb{E}\left[\mathbb{E}\left[\prod_{\substack{j=1\\j\neq p}}^n\left(1-\frac{\bar{\rho}\mathbbm{1}[\mathbf{X}_j\in B_{\mathbb{S}^2}(\mathbf{X}_p,r)]}{\rho(\mathbf{X}_j)}\right)\Bigg| \mathbf{X}_p=\mathbf{x}_p\right]\right]\\
&=\mathbb{E}\left[\mathbb{E}\left[\prod_{\substack{j=1\\j\neq p}}^n\left(1-\frac{\bar{\rho}\mathbbm{1}[\mathbf{X}_j\in B_{\mathbb{S}^2}(\mathbf{x}_p,r)]}{\rho(\mathbf{X}_j)}\right)\right]\right].
\end{align*}
Taking the condition expectation,
\begin{align*}
\mathbb{E}\left[\prod_{\substack{j=1\\j\neq p}}^n\left(1-\frac{\bar{\rho}\mathbbm{1}[\mathbf{X}_j\in B_{\mathbb{S}^2}(\mathbf{x}_p,r)]}{\rho(\mathbf{X}_j)}\right)\right]&=\overbrace{\int_{\mathbb{S}^2}\cdots\int_{\mathbb{S}^2}}^{n-1}\prod_{\substack{j=1\\j\neq p}}^n\left(1-\frac{\bar{\rho}\mathbbm{1}[\mathbf{x}_j\in B_{\mathbb{S}^2}(\mathbf{x}_p,r)]}{\rho(\mathbf{x}_j)}\right)\frac{\rho(\mathbf{x}_j)}{\mu(\mathbb{S}^2)}\lambda_{\mathbb{S}^2}(d\mathbf{x}_j)\\
&=\left(\int_{\mathbb{S}^2} \left(1-\frac{\bar{\rho}\mathbbm{1}[\mathbf{x}\in B_{\mathbb{S}^2}(\mathbf{x}_p,r)]}{\rho(\mathbf{x})}\right)\frac{\rho(\mathbf{x})}{\mu(\mathbb{S}^2)}\lambda_{\mathbb{S}^2}(d\mathbf{x})\right)^{n-1}\\
&=\left(1-\frac{\bar{\rho}}{\mu(\mathbb{S}^2)}2\pi(1-\cos r)\right)^{n-1}.
\end{align*}
Hence,
\begin{equation*}
\mathbb{E}[L_p]=\left(1-\frac{\bar{\rho}}{\mu(\mathbb{S}^2)}2\pi(1-\cos r)\right)^{n-1}
\end{equation*}
Next we will deal with term $(c)$.
\begin{align*}
\mathbb{E}[L_p^2]&= \mathbb{E}\left[\left(\prod_{\substack{j=1\\j\neq p}}^n\left(1-\frac{\bar{\rho}\mathbbm{1}[\mathbf{X}_j\in B_{\mathbb{S}^2}(\mathbf{X}_p,r)]}{\rho(\mathbf{X}_j)}\right)\right)^2\right]\\
&=\mathbb{E}\left[\prod_{\substack{j=1\\j\neq p}}^n\left(1-\frac{\bar{\rho}\mathbbm{1}[\mathbf{X}_j\in B_{\mathbb{S}^2}(\mathbf{X}_p,r)]}{\rho(\mathbf{X}_j)}\right)^2\right]\\
&=\mathbb{E}\left[\mathbb{E}\left[\prod_{\substack{j=1\\j\neq p}}^n\left(1-\frac{\bar{\rho}\mathbbm{1}[\mathbf{X}_j\in B_{\mathbb{S}^2}(\mathbf{x}_p,r)]}{\rho(\mathbf{X}_j)}\right)^2\Bigg|\mathbf{X}_p=\mathbf{x}_p\right]\right]\\
&=\mathbb{E}\left[\prod_{\substack{j=1\\j\neq p}}^n\mathbb{E}\left[\left(1-\frac{\bar{\rho}\mathbbm{1}[\mathbf{X}_j\in B_{\mathbb{S}^2}(\mathbf{x}_p,r)]}{\rho(\mathbf{X}_j)}\right)^2\Bigg|\mathbf{X}_p=\mathbf{x}_p\right]\right]\\
&=\mathbb{E}\left[\prod_{\substack{j=1\\j\neq p}}^n\int_{\mathbb{S}^2}\left(1-\frac{\bar{\rho}\mathbbm{1}[\mathbf{x}\in B_{\mathbb{S}^2}(\mathbf{x}_p,r)]}{\rho(\mathbf{x})}\right)^2\frac{\rho(\mathbf{x})}{\mu(\mathbb{S}^2)}\lambda_{\mathbb{S}^2}(d\mathbf{x})\right]\\
&=\mathbb{E}\left[\prod_{\substack{j=1\\j\neq p}}^n\int_{\mathbb{S}^2}\frac{\rho(\mathbf{x})}{\mu(\mathbb{S}^2)}-\frac{2\bar{\rho}}{\mu(\mathbb{S}^2)}\mathbbm{1}[\mathbf{x}\in B_{\mathbb{S}^2}(\mathbf{x}_q,r)] + \frac{\bar{\rho}^2}{\mu(\mathbb{S}^2)}\frac{\mathbbm{1}[\mathbf{x}\in B_{\mathbb{S}^2}(\mathbf{x}_p,r)]}{\rho(\mathbf{x})} \lambda_{\mathbb{S}^2}(d\mathbf{x})\right]\\
&=\mathbb{E}\left[\prod_{\substack{j=1\\j\neq p}}^n\left(1-\frac{2\bar{\rho}}{\mu(\mathbb{S}^2)}2\pi(1-\cos r) + \frac{\bar{\rho}^2}{\mu(\mathbb{S}^2)}\int_{\mathbb{S}^2\cap B_{\mathbb{S}^2}(\mathbf{x}_p,r)}\frac{1}{\rho(\mathbf{x})}\lambda_{\mathbb{S}^2}(d\mathbf{x})\right)\right]\\
&=\mathbb{E}\left[\left(1-\frac{2\bar{\rho}}{\mu(\mathbb{S}^2)}2\pi(1-\cos r) + \frac{\bar{\rho}^2}{\mu(\mathbb{S}^2)}\int_{\mathbb{S}^2\cap B_{\mathbb{S}^2}(\mathbf{X}_p,r)}\frac{1}{\rho(\mathbf{x})}\lambda_{\mathbb{S}^2}(d\mathbf{x})\right)^{n-1}\right]\\
&= \int_{\mathbb{S}^2}\left(1-\frac{2\bar{\rho}}{\mu(\mathbb{S}^2)}2\pi(1-\cos r) + \frac{\bar{\rho}^2}{\mu(\mathbb{S}^2)}\int_{\mathbb{S}^2\cap B_{\mathbb{S}^2}(\mathbf{y},r)}\frac{1}{\rho(\mathbf{x})}\lambda_{\mathbb{S}^2}(d\mathbf{x})\right)^{n-1}\frac{\rho(\mathbf{y})}{\mu(\mathbb{S}^2)}\lambda_{\mathbb{S}^2}(d\mathbf{y})\\
\intertext{Let us define $A_2(\mathbf{y})=1-\frac{2\bar{\rho}}{\mu(\mathbb{S}^2)}2\pi(1-\cos r) + \frac{\bar{\rho}^2}{\mu(\mathbb{S}^2)}\int_{\mathbb{S}^2\cap B_{\mathbb{S}^2}(\mathbf{y},r)}\frac{1}{\rho(\mathbf{x})}\lambda_{\mathbb{S}^2}(d\mathbf{x})$ then,}
&=\frac{1}{\mu(\mathbb{S}^2)}\int_{\mathbb{S}^2}A_2^{n-1}(\mathbf{y})\rho(\mathbf{y})\lambda_{\mathbb{S}^2}(d\mathbf{y})
\end{align*}
The final term $(d)$ is identical to that of $(b)$. So plugging into Equation \ref{covariance:expansion} gives,
\begin{align*}
&\text{Var}(\hat{H}_{\text{inhom}}(r)|N_X(\mathbb{S}^2=n)=\frac{\mathbbm{1}[n>0]}{n^2}\sum_{p,q\in\{1,\dots,n\}}^{\neq}\\
&\Bigg(\frac{1}{\mu^2(\mathbb{S}^2)}\int_{\mathbb{S}^2}\int_{\mathbb{S}^2}\left(\rho(\mathbf{x})-\bar{\rho}\mathbbm{1}[\mathbf{x}\in B_{\mathbb{S}^2}(\mathbf{y},r)]\right)\left(\rho(\mathbf{y})-\bar{\rho}\mathbbm{1}[\mathbf{y}\in B_{\mathbb{S}^2}(\mathbf{x},r)]\right)A_1^{n-2}(\mathbf{x},\mathbf{y})\lambda_{\mathbb{S}^2}(d\mathbf{x}) \lambda_{\mathbb{S}^2}(d\mathbf{y})\\
&\phantom{AAAA}\left.-\left(1-\frac{\bar{\rho}}{\mu(\mathbb{S}^2)}2\pi(1-\cos r)\right)^{n-1}\left(1-\frac{\bar{\rho}}{\mu(\mathbb{S}^2)}2\pi(1-\cos r)\right)^{n-1}\right)\\
&\phantom{AAAA}+\frac{\mathbbm{1}[n>0]}{n^2}\sum_{p=1}^n \left(\frac{1}{\mu(\mathbb{S}^2)}\int_{\mathbb{S}^2}A_2^{n-1}(\mathbf{y})\rho(\mathbf{y})\lambda_{\mathbb{S}^2}(d\mathbf{y})-\left(1-\frac{\bar{\rho}}{\mu(\mathbb{S}^2)}2\pi(1-\cos r)\right)^{2n-2}\right)\\
&=\frac{1}{\mu^2(\mathbb{S}^2)}\\
&\phantom{=}\int_{\mathbb{S}^2}\int_{\mathbb{S}^2}\left(\rho(\mathbf{x})-\bar{\rho}\mathbbm{1}[\mathbf{x}\in B_{\mathbb{S}^2}(\mathbf{y},r)]\right)\left(\rho(\mathbf{y})-\bar{\rho}\mathbbm{1}[\mathbf{y}\in B_{\mathbb{S}^2}(\mathbf{x},r)]\right)\\
&\phantom{AAAA}\left(\frac{\mathbbm{1}[n>0](n-1)}{n}A_1^{n-2}(\mathbf{x},\mathbf{y})\right)\lambda_{\mathbb{S}^2}(d\mathbf{x}) \lambda_{\mathbb{S}^2}(d\mathbf{y})\\
&\phantom{AAAA}+\frac{1}{\mu(\mathbb{S}^2)}\int_{\mathbb{S}^2}\left(\frac{\mathbbm{1}[n>0]}{n}A_2^{n-1}(\mathbf{y})\right)\rho(\mathbf{y})\lambda_{\mathbb{S}^2}(d\mathbf{y})-\mathbbm{1}[n>0]\left(1-\frac{\bar{\rho}}{\mu(\mathbb{S}^2)}2\pi(1-\cos r)\right)^{2n-2}
\end{align*}
We need to then take the expectation of this variance over $N_X(\mathbb{S}^2)$. By application of Tonelli's Theorem we can interchange the expectation over $N_X(\mathbb{S}^2)$ with the integrals over $\mathbf{x}$ and $\mathbf{y}$. This comes by showing that both $A_1(\mathbf{x},\mathbf{y})$ and $A_2(\mathbf{x})$ are non-negative for all $\mathbf{x},\mathbf{y}\in\mathbb{S}^2$. Obviously $\rho(\mathbf{x})-\bar{\rho}\mathbbm{1}[\mathbf{x}\in B_{\mathbb{S}^2}(\mathbf{y},r)]$ and $\rho(\mathbf{y})-\bar{\rho}\mathbbm{1}[\mathbf{y}\in B_{\mathbb{S}^2}(\mathbf{x},r)]$ are greater than or equal to 0 since $\bar{\rho}=\inf_{\mathbf{x}\in\mathbb{S}^2}\rho(\mathbf{x})$. It then follows since the integrand of $A_1(\mathbf{x}_p,\mathbf{x}_q)$, $\left(1-\bar{\rho}\mathbbm{1}[\mathbf{x}\in B_{\mathbb{S}^2}(\mathbf{x}_p,r)]/\rho(\mathbf{x})\right)\left(1-\bar{\rho}\mathbbm{1}[\mathbf{x}\in B_{\mathbb{S}^2}(\mathbf{x}_q,r)]/\rho(\mathbf{x})\right)\rho(\mathbf{x})/\mu(\mathbb{S}^2)$, is then non-negative for all $\mathbf{x}$ and so the integral over $\mathbf{f}$ is non-negative and thus Tonelli's Theorem can be applied. A near identitical argument can be applied to $A_2(\mathbf{x})$ to show that it is always non-negative and hence Tonelli's Theorem can then be applied. We then calculate the following expectations,
\begin{enumerate}[a)]
\item $\mathbb{E}\left[\frac{\mathbbm{1}[N_X(\mathbb{S}^2)>0](N_X(\mathbb{S}^2)-1)}{N_X(\mathbb{S}^2)}A_1^{N_X(\mathbb{S}^2)-2}(\mathbf{x},\mathbf{y})\right]$

\item $\mathbb{E}\left[\frac{\mathbbm{1}[N_X(\mathbb{S}^2)>0]}{N_X(\mathbb{S}^2)}A_2^{N_X(\mathbb{S}^2)-1}(\mathbf{x})\right]$
\item $\mathbb{E}\left[\mathbbm{1}[N_X(\mathbb{S}^2)>0]\left(1-\frac{\bar{\rho}}{\mu(\mathbb{S}^2)}2\pi(1-\cos r)\right)^{2N_X(\mathbb{S}^2)-2}\right]$
\end{enumerate}
Setting $\Lambda =\mu(\mathbb{S}^2)$, expectation (a) is,
\begin{align*}
&\mathbb{E}\left[\frac{\mathbbm{1}[N_X(\mathbb{S}^2)>0](N_X(\mathbb{S}^2)-1)}{N_X(\mathbb{S}^2)}A_1^{N_X(\mathbb{S}^2)-2}(\mathbf{x},\mathbf{y})\right]\\
&\phantom{AAAA}=\mathbb{E}\left[\mathbbm{1}[N_X(\mathbb{S}^2)>0]A_1^{N_X(\mathbb{S}^2)-2}(\mathbf{x},\mathbf{y})\right]-\mathbb{E}\left[\frac{\mathbbm{1}[N_X(\mathbb{S}^2)>0]}{N_X(\mathbb{S}^2)}A_1^{N_X(\mathbb{S}^2)-2}(\mathbf{x},\mathbf{y})\right]\\
&\phantom{AAAA}=\sum_{n=0}^{\infty}\mathbbm{1}[n>0]A_1^{n-2}(\mathbf{x},\mathbf{y})\frac{\Lambda^n\mathrm{e}^{-\Lambda}}{n!}-\sum_{n=0}^{\infty}\frac{\mathbbm{1}[n>0]}{n}A_1^{n-2}(\mathbf{x},\mathbf{y})\frac{\Lambda^n\mathrm{e}^{-\Lambda}}{n!}\\
&\phantom{AAAA}=\frac{\mathrm{e}^{-\Lambda}}{A_1^{2}(\mathbf{x},\mathbf{y})}\sum_{n=1}^{\infty}\frac{(\Lambda A_1(\mathbf{x},\mathbf{y}))^{n}}{n!}-\frac{\mathrm{e}^{-\Lambda}}{A_1^{2}(\mathbf{x},\mathbf{y})}\sum_{n=1}^{\infty}\frac{1}{n}\frac{(\Lambda A_1(\mathbf{x},\mathbf{y}))^n}{n!}\\
&\phantom{AAAA}=\frac{\mathrm{e}^{-\Lambda}}{A_1^{2}(\mathbf{x},\mathbf{y})}\left(\sum_{n=0}^{\infty}\frac{(\Lambda A_1(\mathbf{x},\mathbf{y}))^{n}}{n!}-1\right)-\frac{\mathrm{e}^{-\Lambda}}{A_1^{2}(\mathbf{x},\mathbf{y})}(\text{Ei}(\Lambda A_1(\mathbf{x},\mathbf{y}))-\gamma-\log(\Lambda A_1(\mathbf{x},\mathbf{y})))\\
&\phantom{AAAA}=\frac{\mathrm{e}^{-\Lambda}}{A_1^{2}(\mathbf{x},\mathbf{y})}\left(\mathrm{e}^{\Lambda A_1(\mathbf{x},\mathbf{y})}-1 -\text{Ei}(\Lambda A_1(\mathbf{x},\mathbf{y}))+\gamma +\log(\Lambda A_1(\mathbf{x},\mathbf{y}))\right)
\end{align*}
where the penultimate line follows from Equation \ref{Ei:series}. Similarly for $(b)$ and $(c)$,
\begin{align*}
\mathbb{E}\left[\frac{\mathbbm{1}[N_X(\mathbb{S}^2)>0]}{N_X(\mathbb{S}^2)}A_2^{N_X(\mathbb{S}^2)-1}(\mathbf{x})\right]&= \frac{\mathrm{e}^{-\Lambda}}{A_2(\mathbf{x})}\left(\gamma +\log(\Lambda A_2(\mathbf{x}))-\text{Ei}(\Lambda A_2(\mathbf{x}))\right)\\
\mathbb{E}\left[\mathbbm{1}[N_X(\mathbb{S}^2)>0]C^{2N_X(\mathbb{S}^2)-2}\right]&= \frac{\mathrm{e}^{-\Lambda}}{C^2}\left(\mathrm{e}^{\Lambda C^2}-1\right),
\end{align*}
where $C=\left(1-\frac{\bar{\rho}}{\mu(\mathbb{S}^2)}2\pi(1-\cos r)\right)$. Then taking expectations of $\text{Var}(\hat{H}_{\text{inhom}}(r)|N_X(\mathbb{S}^2)=n)$ over $N_X(\mathbb{S}^2)$ gives,
\begin{align*}
&\mathbb{E}[\text{Var}(\hat{H}_{\text{inhom}}(r)|N_X(\mathbb{S}^2)=n)]\\
&=\frac{1}{\mu^2(\mathbb{S}^2)}\\
&\phantom{=}\int_{\mathbb{S}^2}\int_{\mathbb{S}^2}\left(\rho(\mathbf{x})-\bar{\rho}\mathbbm{1}[\mathbf{x}\in B_{\mathbb{S}^2}(\mathbf{y},r)]\right)\left(\rho(\mathbf{y})-\bar{\rho}\mathbbm{1}[\mathbf{y}\in B_{\mathbb{S}^2}(\mathbf{x},r)]\right)\\
&\phantom{AAAAAAAA}\frac{\mathrm{e}^{-\Lambda}}{A_1^{2}(\mathbf{x},\mathbf{y})}\left(\mathrm{e}^{\Lambda A_1(\mathbf{x},\mathbf{y})}-1 -\text{Ei}(\Lambda A_1(\mathbf{x},\mathbf{y}))+\gamma +\log(\Lambda A_1(\mathbf{x},\mathbf{y}))\right) \lambda_{\mathbb{S}^2}(d\mathbf{x}) \lambda_{\mathbb{S}^2}(d\mathbf{y})\\
&\phantom{=}+\frac{1}{\mu(\mathbb{S}^2)}\int_{\mathbb{S}^2}\frac{\mathrm{e}^{-\Lambda}}{A_2(\mathbf{x})}\left(\gamma +\log(\Lambda A_2(\mathbf{x}))-\text{Ei}(\Lambda A_2(\mathbf{x}))\right)\rho(\mathbf{y})\lambda_{\mathbb{S}^2}(d\mathbf{y})\\
&\phantom{=}-\frac{\mathrm{e}^{-\Lambda}}{\left(1-\frac{\bar{\rho}}{\mu(\mathbb{S}^2)}2\pi(1-\cos r)\right)^2}\left(\mathrm{e}^{\Lambda \left(1-\frac{\bar{\rho}}{\mu(\mathbb{S}^2)}2\pi(1-\cos r)\right)^2}-1\right)
\end{align*}
Term $(2)$ of Equation \ref{variance:total:law:1}, the expectation conditioned on $N_X(\mathbb{S}^2)=n$,
\begin{align*}
&\mathbb{E}[\hat{H}_{\text{inhom}}(r)|N_X(\mathbb{S}^2)=n]=\\
&\phantom{AAAA}\mathbb{E}\left[1-\frac{\mathbbm{1}[N_X(\mathbb{S}^2)>0]}{N_X(\mathbb{S}^2)}\sum_{\mathbf{x}\in X}\prod_{\mathbf{y}\in X\setminus\{\mathbf{x}\}}\left(1-\frac{\bar{\rho}\mathbbm{1}[\mathbf{y}\in B_{\mathbb{S}^2}(\mathbf{x},r)]}{\rho(\mathbf{y})}\right)\Bigg|N_X(\mathbb{S}^2)=n\right]\\
&\phantom{AAAA}=1-\frac{\mathbbm{1}[n>0]}{n}\sum_{i=1}^n\mathbb{E}\left[\prod_{\substack{j=1\\j\neq i}}^n\left(1-\frac{\bar{\rho}\mathbbm{1}[\mathbf{x}_j\in B(\mathbf{x}_i,r)]}{\rho(\mathbf{x}_j)}\right)\right]\\
&\phantom{AAAA}=1-\frac{\mathbbm{1}[n>0]}{n}\sum_{i=1}^n\mathbb{E}\left[L_i\right],
\intertext{where $L_i$ is as defined in Equation \ref{covariance:expansion},}
&\phantom{AAAA}=1-\frac{\mathbbm{1}[n>0]}{n}\sum_{i=1}^n \left(1-\frac{\bar{\rho}}{\mu(\mathbb{S}^2)}2\pi(1-\cos r)\right)^{n-1}\\
&\phantom{AAAA}=1-\mathbbm{1}[n>0]\left(1-\frac{\bar{\rho}}{\mu(\mathbb{S}^2)}2\pi(1-\cos r)\right)^{n-1}\\
&\phantom{AAAA}=1-\mathbbm{1}[n>0]C^{n-1},
\end{align*}
where $C=1-\frac{\bar{\rho}}{\mu(\mathbb{S}^2)}2\pi(1-\cos r)$. Then taking the variance over $N_X(\mathbb{S}^2)$ gives,
\begin{align*}
&\text{Var}(1-\mathbbm{1}[N_X(\mathbb{S}^2)>0]C^{N_X(\mathbb{S}^2)-1})\\
&\phantom{AAAA}=\text{Var}(\mathbbm{1}[N_X(\mathbb{S}^2)>0]C^{N_X(\mathbb{S}^2)-1})\\
&\phantom{AAAA}=\mathbb{E}\left[\left(\mathbbm{1}[N_X(\mathbb{S}^2)>0]C^{N_X(\mathbb{S}^2)-1}\right)^2\right]-\mathbb{E}^2\left[\mathbbm{1}[N_X(\mathbb{S}^2)>0]C^{N_X(\mathbb{S}^2)-1}\right]\\
&\phantom{AAAA}=\mathbb{E}\left[\mathbbm{1}[N_X(\mathbb{S}^2)>0]C^{2N_X(\mathbb{S}^2)-2}\right]-\mathbb{E}^2\left[\mathbbm{1}[N_X(\mathbb{S}^2)>0]C^{N_X(\mathbb{S}^2)-1}\right]\\
&\phantom{AAAA}=\frac{\mathrm{e}^{-\Lambda}}{C^2}\left(\mathrm{e}^{\Lambda C^2}-1\right)-\frac{\mathrm{e}^{-2\Lambda}}{C^2}\left(\mathrm{e}^{\Lambda C}-1\right)^2,
\end{align*}
where $\Lambda=\mu(\mathbb{S}^2)$. Therefore the variance of $\hat{H}_{\text{inhom}}(r)$ is,
\begin{equation}\label{var:sub}
\begin{split}
&\text{Var}(\hat{H}_{\text{inhom}}(r))\\
&=\frac{1}{\mu^2(\mathbb{S}^2)}\int_{\mathbb{S}^2}\int_{\mathbb{S}^2}\left(\rho(\mathbf{x})-\bar{\rho}\mathbbm{1}[\mathbf{x}\in B_{\mathbb{S}^2}(\mathbf{y},r)]\right)\left(\rho(\mathbf{y})-\bar{\rho}\mathbbm{1}[\mathbf{y}\in B_{\mathbb{S}^2}(\mathbf{x},r)]\right)\\
&\frac{\mathrm{e}^{-\mu(\mathbb{S}^2)}}{A_1^{2}(\mathbf{x},\mathbf{y})}\left(\mathrm{e}^{\mu(\mathbb{S}^2) A_1(\mathbf{x},\mathbf{y})}-1 -\text{Ei}(\mu(\mathbb{S}^2) A_1(\mathbf{x},\mathbf{y}))+\gamma +\log(\mu(\mathbb{S}^2) A_1(\mathbf{x},\mathbf{y}))\right) \lambda_{\mathbb{S}^2}(d\mathbf{x}) \lambda_{\mathbb{S}^2}(d\mathbf{y})\\
&\phantom{=}+\frac{1}{\mu(\mathbb{S}^2)}\int_{\mathbb{S}^2}\frac{\mathrm{e}^{-\mu(\mathbb{S}^2)}}{A_2(\mathbf{x})}\left(\gamma +\log(\mu(\mathbb{S}^2) A_2(\mathbf{x}))-\text{Ei}(\mu(\mathbb{S}^2) A_2(\mathbf{x}))\right)\rho(\mathbf{y})\lambda_{\mathbb{S}^2}(d\mathbf{y})\\
&\phantom{=}-\frac{\mathrm{e}^{-2\mu(\mathbb{S}^2)}}{\left(1-\frac{\bar{\rho}}{\mu(\mathbb{S}^2)}2\pi(1-\cos r)\right)^2}\left(\mathrm{e}^{\mu(\mathbb{S}^2) \left(1-\frac{\bar{\rho}}{\mu(\mathbb{S}^2)}2\pi(1-\cos r)\right)}-1\right)^2
\end{split}
\end{equation}
The final part of the proof is to ensure that the integrands are truly Lebesgue integrable, that is the integrals given in the previous equation are finite. We shall work with the second, $\int_{\mathbb{S}^2}(\mathrm{e}^{-\mu(\mathbb{S}^2)}/A_2(\mathbf{x}))\left(\gamma +\log(\mu(\mathbb{S}^2) A_2(\mathbf{x}))-\text{Ei}(\mu(\mathbb{S}^2) A_2(\mathbf{x}))\right)\rho(\mathbf{y})\lambda_{\mathbb{S}^2}(d\mathbf{y})$, the first then follows by a similar argument. By using the series expansion we know that $\text{Ei}(\mu(\mathbb{S}^2) A_2(\mathbf{x}))- \gamma - \log(\mu(\mathbb{S}^2) A_2(\mathbf{x}))=\sum_{n=1}^\infty (\left(\mu(\mathbb{S}^2) A_2(\mathbf{x})\right)^n)/(n\cdot n!)$ then it can be bound from above as,
\begin{align*}
\text{Ei}(\mu(\mathbb{S}^2) A_2(\mathbf{x})) - \gamma - \log(\mu(\mathbb{S}^2) A_2(\mathbf{x}))&\leq \sum_{n=0}^\infty \frac{\left(\mu(\mathbb{S}^2) A_2(\mathbf{x})\right)^n}{n\cdot n!}\\
& \sum_{n=0}^\infty \frac{\left(\mu(\mathbb{S}^2) A_2(\mathbf{x})\right)^n}{n!}\\
& = \mathrm{e}^{\mu(\mathbb{S}^2) A_2(\mathbf{x})},
\end{align*}
and bound from below as,
\begin{align*}
\text{Ei}(\mu(\mathbb{S}^2) A_2(\mathbf{x})) - \gamma - \log(\mu(\mathbb{S}^2) A_2(\mathbf{x}))&\geq \sum_{n=1}^\infty \frac{\left(\mu(\mathbb{S}^2) A_2(\mathbf{x})\right)^n}{(n+1)!}\\
&=\frac{1}{\mu(\mathbb{S}^2) A_2(\mathbf{x})}\left(\sum_{n=0}^{\infty}\frac{\left(\mu(\mathbb{S}^2) A_2(\mathbf{x})\right)^{n+1}}{(n+1)!}-1\right)\\
&=\frac{1}{\mu(\mathbb{S}^2) A_2(\mathbf{x})}\left(\mathrm{e}^{\mu(\mathbb{S}^2) A_2(\mathbf{x})}-1\right)
\end{align*}
Hence the integrand is bounded,
\begin{equation}\label{bound:integrand}
\begin{split}
\frac{\mathrm{e}^{-\mu(\mathbb{S}^2)}}{\mu(\mathbb{S}^2) A_2^2(\mathbf{x})}&\left(\mathrm{e}^{\mu(\mathbb{S}^2) A_2(\mathbf{x})}-1\right)\leq\\
&\frac{\mathrm{e}^{-\mu(\mathbb{S}^2)}}{A_2(\mathbf{x})}\left(\gamma +\log(\mu(\mathbb{S}^2) A_2(\mathbf{x}))-\text{Ei}(\mu(\mathbb{S}^2) A_2(\mathbf{x}))\right)\leq\frac{\mathrm{e}^{-\mu(\mathbb{S}^2)}}{A_2(\mathbf{x})}\mathrm{e}^{\mu(\mathbb{S}^2) A_2(\mathbf{x})}.
\end{split}
\end{equation}
We show that the lower and upper bounds can be bounded further such that they do not depend on $r$ or $\mathbf{x}$. First we bound $\mathrm{e}^{-\mu(\mathbb{S}^2)}\mathrm{e}^{\mu(\mathbb{S}^2) A_2(\mathbf{x})}$ below and above,
\begin{align*}
\mathrm{e}^{-\mu(\mathbb{S}^2)}\mathrm{e}^{\mu(\mathbb{S}^2) A_2(\mathbf{x})}&=\exp\left(-4\pi\bar{\rho}(1-\cos r)+\frac{\bar{\rho}^2}{\mu(\mathbb{S}^2)}\int_{\mathbb{S}^2\cap B_{\mathbb{S}^2}(\mathbf{x},r)}\frac{1}{\rho(\mathbf{y})}\lambda_{\mathbb{S}^2}(d\mathbf{y})\right)\\
&\leq\exp\left(\frac{\bar{\rho}^2}{\mu(\mathbb{S}^2)}\int_{\mathbb{S}^2\cap B_{\mathbb{S}^2}(\mathbf{x},r)}\frac{1}{\rho(\mathbf{y})}\lambda_{\mathbb{S}^2}(d\mathbf{y})\right)\\
&\leq\exp\left(\frac{\bar{\rho}^2}{\mu(\mathbb{S}^2)}\int_{\mathbb{S}^2\cap B_{\mathbb{S}^2}(\mathbf{x},r)}\frac{1}{\bar{\rho}}\lambda_{\mathbb{S}^2}(d\mathbf{y})\right)\\
&=\exp\left(\frac{4\pi\bar{\rho}}{\mu(\mathbb{S}^2)}\right)\\
&\leq\mathrm{e},
\end{align*}
where the final inequality follows from $\mu(\mathbb{S}^2)\geq 4\pi\bar{\rho}$. The lower bound is,
\begin{align*}
\mathrm{e}^{-\mu(\mathbb{S}^2)}\mathrm{e}^{\mu(\mathbb{S}^2) A_2(\mathbf{x})}&=\exp\left(-4\pi\bar{\rho}(1-\cos r)+\frac{\bar{\rho}^2}{\mu(\mathbb{S}^2)}\int_{\mathbb{S}^2\cap B_{\mathbb{S}^2}(\mathbf{x},r)}\frac{1}{\rho(\mathbf{y})}\lambda_{\mathbb{S}^2}(d\mathbf{y})\right)\\
&\geq \exp\left(-4\pi\bar{\rho}(1-\cos r)\right)\\
&\geq \exp\left(-8\pi\bar{\rho}\right).
\end{align*}
Next we need to show that $A_2(\mathbf{x})$ is strictly greater than 0,
\begin{align*}
A_2(\mathbf{x})&=1-\frac{4\pi\bar{\rho}}{\mu(\mathbb{S}^2)}(1-\cos r) +\frac{\bar{\rho}^2}{\mu(\mathbb{S}^2)}\int_{\mathbb{S}^2\cap B_{\mathbb{S}^2}(\mathbf{x},r)}\frac{1}{\rho(\mathbf{y})}\lambda_{\mathbb{S}^2}(d\mathbf{y})\\
&=\int_{\mathbb{S}^2}\left(1-\frac{\bar{\rho}\mathbbm{1}[\mathbf{y}\in B_{\mathbb{S}^2}(\mathbf{x},r)]}{\rho(\mathbf{y})}\right)^2\frac{\rho(\mathbf{y})}{\mu(\mathbb{S}^2)}\lambda_{\mathbb{S}^2}(d\mathbf{y})\\
&=\int_{\mathbb{S}^2\cap B_{\mathbb{S}^2}(\mathbf{x},r)}\left(1-\frac{\bar{\rho}}{\rho(\mathbf{y})}\right)^2\frac{\rho(\mathbf{y})}{\mu(\mathbb{S}^2)}\lambda_{\mathbb{S}^2}(d\mathbf{y})+\int_{\mathbb{S}^2\setminus B_{\mathbb{S}^2}(\mathbf{x},r)}\frac{\rho(\mathbf{y})}{\mu(\mathbb{S}^2)}\lambda_{\mathbb{S}^2}(d\mathbf{y}).
\end{align*}
Then for $r\in [0,\pi)$ $\mu_L(\mathbb{S}^2\setminus B_{\mathbb{S}^2}(\mathbf{x},r))>0$ and since $\rho(\mathbf{y})\geq \bar{\rho}>0$ this means that the second term is strictly greater than $0$. Further the first term is always non-negative since $\rho(\mathbf{y})\geq \bar{\rho}>0$. Then if $r=\pi$ we have that,
\begin{align*}
A_2(\mathbf{x})&=\int_{\mathbb{S}^2}\left(1-\frac{\bar{\rho}}{\rho(\mathbf{y})}\right)^2\frac{\rho(\mathbf{y})}{\mu(\mathbb{S}^2)}\lambda_{\mathbb{S}^2}(d\mathbf{y})\\
&=\int_{\mathbb{S}^2\cap E}\left(1-\frac{\bar{\rho}}{\rho(\mathbf{y})}\right)^2\frac{\rho(\mathbf{y})}{\mu(\mathbb{S}^2)}\lambda_{\mathbb{S}^2}(d\mathbf{y}) + \int_{\mathbb{S}^2\setminus E}\left(1-\frac{\bar{\rho}}{\rho(\mathbf{y})}\right)^2\frac{\rho(\mathbf{y})}{\mu(\mathbb{S}^2)}\lambda_{\mathbb{S}^2}(d\mathbf{y}), 
\end{align*}
then by assumption $\rho(\mathbf{y})>\bar{\rho}$ for all $\mathbf{y}\in E\subset \mathbb{S}^2$, hence the first term is strictly greater than 0 whilst the second term is non-negative and so we have shown that for any $r\in [0,\pi]$ $A_2(\mathbf{x})>0$, therefore when taking the reciprocal of $A_2(\mathbf{x})$ it is well defined for all $\mathbf{x}\in\mathbb{S}^2$. Then we can bound the absolute value of the reciprocal of $A_2(\mathbb{S}^2)$. Further, in the case when $r\in [0,\pi),$ $\int_{\mathbb{S}^2\setminus B_{\mathbb{S}^2}(\mathbf{x},r)}\frac{\rho(\mathbf{y})}{\mu(\mathbb{S}^2)}\lambda_{\mathbb{S}^2}(d\mathbf{y}) \geq \int_{\mathbb{S}^2\setminus B_{\mathbb{S}^2}(\mathbf{x},r)}\frac{\bar{\rho}}{\mu(\mathbb{S}^2)}\lambda_{\mathbb{S}^2}(d\mathbf{y})=(1-\cos r) \frac{2\pi\bar{\rho}}{\mu(\mathbb{S}^2)}$ then define,
\begin{equation*}
\bar{A}_2=\begin{cases}
(1-\cos r) \frac{2\pi\bar{\rho}}{\mu(\mathbb{S}^2)},\quad &r\in [0,\pi)\\
\int_{\mathbb{S}^2\setminus E}\left(1-\frac{\bar{\rho}}{\rho(\mathbf{y})}\right)^2\frac{\rho(\mathbf{y})}{\mu(\mathbb{S}^2)}\lambda_{\mathbb{S}^2}(d\mathbf{y}),\quad &r=\pi,
\end{cases}
\end{equation*}
which does not depend on $\mathbf{x}$. Then $A_2(\mathbf{x})\geq \bar{A}_2> 0$ and so $0 \leq \frac{1}{A_2(\mathbf{x})}\leq \frac{1}{\bar{A}_2}$, which means we can bound the the absolute value of the reciprocal of $A_2(\mathbb{S}^2)$. Further we can obtain a non-zero lower bound for $\frac{1}{A_2(\mathbf{x})}$ since $A_2(\mathbf{x})=1-\frac{4\pi\bar{\rho}}{\mu(\mathbb{S}^2)}(1-\cos r) +\frac{\bar{\rho}^2}{\mu(\mathbb{S}^2)}\int_{\mathbb{S}^2\cap B_{\mathbb{S}^2}(\mathbf{x},r)}\frac{1}{\rho(\mathbf{y})}\lambda_{\mathbb{S}^2}(d\mathbf{y})\leq 1+\frac{\bar{\rho}^2}{\mu(\mathbb{S}^2)}\int_{\mathbb{S}^2}\frac{1}{\rho(\mathbf{y})}\lambda_{\mathbb{S}^2}(d\mathbf{y})$. Let us define $\tilde{A}_2=1+\frac{\bar{\rho}^2}{\mu(\mathbb{S}^2)}\int_{\mathbb{S}^2}\frac{1}{\rho(\mathbf{y})}\lambda_{\mathbb{S}^2}(d\mathbf{y})$ then returning to Equation \ref{bound:integrand} we have, 
\begin{equation*}
\begin{split}
\frac{1}{\mu(\mathbb{S}^2) \tilde{A}_2^2}\left(\mathrm{e}^{-8\pi\bar{\rho}}-\mathrm{e}^{-\mu(\mathbb{S}^2)}\right)\leq
\frac{\mathrm{e}^{-\mu(\mathbb{S}^2)}}{A_2(\mathbf{x})}\left(\gamma +\log(\mu(\mathbb{S}^2) A_2(\mathbf{x}))-\text{Ei}(\mu(\mathbb{S}^2) A_2(\mathbf{x}))\right)\leq\frac{\mathrm{e}}{\bar{A}_2},
\end{split}
\end{equation*}
and therefore,
\begin{equation*}
\left|\frac{\mathrm{e}^{-\mu(\mathbb{S}^2)}}{A_2(\mathbf{x})}\left(\gamma +\log(\mu(\mathbb{S}^2) A_2(\mathbf{x}))-\text{Ei}(\mu(\mathbb{S}^2) A_2(\mathbf{x}))\right)\right|\leq \max\left\{\left|\frac{1}{\mu(\mathbb{S}^2) \tilde{A}_2^2}\left(\mathrm{e}^{-8\pi\bar{\rho}}-\mathrm{e}^{-\mu(\mathbb{S}^2)}\right)\right|,\frac{\mathrm{e}}{\bar{A}_2}\right\}
\end{equation*}
Since the right hand side of this inequality is a constant this means that it is Lebesgue integrable over $\mathbb{S}^2$ and so by the dominated convergence theorem so is the left hand side, thus showing that the integrands are truly Lebesgue integrable. An identical approach can be used to show that the first term of Equation \ref{var:sub} is also Lebesgue integrable.
\end{proof}

\section{Existence and approximation of first and second order moments of the $\hat{J}_{\text{inhom}}$-function}\label{Supplementary:J:moments}

In this section we discuss the existence of the first two moments of the empirical inhomogeneous $J$-function, outlining sufficient conditions for the moments to exist. We then derive approximations to the first and second order moments of $\hat{J}_{\text{inhom}}(r)$. The following theorem gives conditions under which the first two moments of $\hat{J}_{\text{inhom}}(r)$ exist.

\subsection{Existence of first and second moments of $J_{\text{inhom}(r)}$}

\begin{theorem}
Let $X$ be a spheroidal Poisson process with intensity function $\rho:\mathbb{S}^2\mapsto\mathbb{R}_+$ such that $\bar{\rho}\equiv\inf_{\mathbf{x}\in\mathbb{S}^2}\rho(\mathbf{x})>0$. Let $P$ be any finite grid on $\mathbb{S}^2$ and define $r_{\max}=\sup\{r\in [0,\pi]:\text{ there exists } \mathbf{p}\in P \text{ such that } \rho(\mathbf{x})\neq \bar{\rho} \text{ for all } \mathbf{x} \in B_{\mathbb{S}^2}(\mathbf{p},r)\}$. Then for any given $r\in [0,r_{\max}]$ both $\mathbb{E}[\hat{J}_{\rm{inhom}}(r)]$ and ${\rm{Var}}(\hat{J}_{\rm{inhom}}(r))$ exist.
\end{theorem}
\begin{proof}
Starting with the expectation,
\begin{align*}
\mathbb{E}[\hat{J}_{\text{inhom}}(r)]&=\int_{N_{lf}}\hat{J}_{\text{inhom}}(r)P_X(dx),\\
\intertext{where $P_X(X\in F), F\subseteq N_{lf}$ is the probability measure of $X$. Define $N_{lf,0}=\{x\in N_{lf}| n_x(\mathbb{S}^2)=0\}$ and $N_{lf,1}=\{x\in N_{lf}| n_x(\mathbb{S}^2)>0\}$, then $N_{lf}=N_{lf,0}\cup N_{lf,1}$ and $N_{lf,0}\cap N_{lf,1}=\emptyset$ and so,}
\mathbb{E}[\hat{J}_{\text{inhom}}(r)]&=\int_{N_{lf,0}}\hat{J}_{\text{inhom}}(r)P_X(dx) + \int_{N_{lf,1}}\hat{J}_{\text{inhom}}(r)P_X(dx),\\
\intertext{taking the convention that $\frac{0}{0}=1$, the first term is finite,}
&=\int_{N_{lf,0}}P_X(dx) + \int_{N_{lf,1}}\hat{J}_{\text{inhom}}(r)P_X(dx),\\
&=P_X(X\in N_{lf,0}) + \int_{N_{lf,1}}\hat{J}_{\text{inhom}}(r)P_X(dx),\\
\intertext{it can be shown that $P_X(X\in N_{lf,0})=P_{N_X(\mathbb{S}^2)}(N_X(\mathbb{S}^2)=0)$ then,}
&=P(N_X(\mathbb{S}^2)=0) + \int_{N_{lf,1}}\hat{J}_{\text{inhom}}(r)P_X(dx).
\end{align*}
We now show that the second term can be bounded and hence the expectation is well defined. Taking the integrand and noting that random variables are now deterministic,
\begin{equation*}
\hat{J}_{\text{inhom}}(r)=\frac{|P|\sum_{\mathbf{x}\in x}\prod_{\mathbf{y}\in x\setminus \mathbf{x}}\left(1-\frac{\bar{\rho}\mathbbm{1}[\mathbf{y}\in B_{\mathbb{S}^2}(\mathbf{x},r)]}{\rho(\mathbf{y})}\right)}{n_x(\mathbb{S}^2)\sum_{\mathbf{p}\in P}\prod_{\mathbf{z}\in x}\left(1-\frac{\bar{\rho}\mathbbm{1}[\mathbf{z}\in B_{\mathbb{S}^2}(\mathbf{p},r)]}{\rho(\mathbf{z})}\right)}.
\end{equation*}
First we note that $\hat{J}_{\text{inhom}}(r)\geq 0$ and is so bounded below. Working with the numerator we have that,
\begin{align*}
\sum_{\mathbf{x}\in x}\prod_{\mathbf{y}\in x\setminus \mathbf{x}}\left(1-\frac{\bar{\rho}\mathbbm{1}[\mathbf{y}\in B_{\mathbb{S}^2}(\mathbf{x},r)]}{\rho(\mathbf{y})}\right)&\leq \sum_{\mathbf{x}\in x} \prod_{y\in x\setminus\mathbf{x}} 1\\
&=n_x(\mathbb{S}^2).
\end{align*}
Now working with the denominator showing that it is bounded below and hence its reciprocal bounded above. By the assumption we have that there exists $\tilde{\mathbf{p}}\in P$ such that for any $\mathbf{x}\in B_{\mathbb{S}^2}(\tilde{\mathbf{p}},r), \rho(\mathbf{x})\neq \bar{\rho}$, then
\begin{align*}
&\sum_{\mathbf{p}\in P}\prod_{\mathbf{z}\in x}\left(1-\frac{\bar{\rho}\mathbbm{1}[\mathbf{z}\in B_{\mathbb{S}^2}(\mathbf{p},r)]}{\rho(\mathbf{z})}\right) = \\
&\phantom{AAAA}\prod_{\mathbf{z}\in x}\left(1-\frac{\bar{\rho}\mathbbm{1}[\mathbf{z}\in B_{\mathbb{S}^2}(\tilde{\mathbf{p}},r)]}{\rho(\mathbf{z})}\right) + \sum_{\mathbf{p}\in P\setminus \tilde{p}}\prod_{\mathbf{z}\in x}\left(1-\frac{\bar{\rho}\mathbbm{1}[\mathbf{z}\in B_{\mathbb{S}^2}(\mathbf{p},r)]}{\rho(\mathbf{z})}\right).
\end{align*}
Then the first term is strictly greater than 0 by our assumption. Thus,
\begin{align*}
\sum_{\mathbf{p}\in P}\prod_{\mathbf{z}\in x}\left(1-\frac{\bar{\rho}\mathbbm{1}[\mathbf{z}\in B_{\mathbb{S}^2}(\mathbf{p},r)]}{\rho(\mathbf{z})}\right) \geq \prod_{\mathbf{z}\in x}\left(1-\frac{\bar{\rho}\mathbbm{1}[\mathbf{z}\in B_{\mathbb{S}^2}(\tilde{\mathbf{p}},r)]}{\rho(\mathbf{z})}\right)>0.
\end{align*}
Further by the assumption we can define $\bar{\rho}_{\tilde{\mathbf{p}}}=\inf_{\mathbf{x}\in B_{\mathbb{S}^2}(\tilde{\mathbf{p}},r)}\rho(\mathbf{x})$ and then,
\begin{align*}
\sum_{\mathbf{p}\in P}\prod_{\mathbf{z}\in x}\left(1-\frac{\bar{\rho}\mathbbm{1}[\mathbf{z}\in B_{\mathbb{S}^2}(\mathbf{p},r)]}{\rho(\mathbf{z})}\right) \geq \prod_{\mathbf{z}\in x}\left(1-\frac{\bar{\rho}}{\bar{\rho}_{\tilde{\mathbf{p}}}}\right)=\left(1-\frac{\bar{\rho}}{\bar{\rho}_{\tilde{\mathbf{p}}}}\right)^{n_x(\mathbf{S}^2)},
\end{align*}
and so,
\begin{align*}
\left|\hat{J}_{\text{inhom}}(r)\right| \leq \frac{|P|n_x(\mathbb{S}^2)}{n_x(\mathbb{S}^2)\left(1-\frac{\bar{\rho}}{\bar{\rho}_{\tilde{\mathbf{p}}}}\right)^{n_x(\mathbf{S}^2)}} = |P|\left(1-\frac{\bar{\rho}}{\bar{\rho}_{\tilde{\mathbf{p}}}}\right)^{-n_x(\mathbf{S}^2)}.
\end{align*}
We now show that the right hand side of the above inequality is integrable. Define the sets $N_{lf,1,(i)}=\{x\in N_{lf,1}: n(x)=i\},$ then $N_{lf,1}=\bigcup_{i=1}^\infty N_{lf,1,(i)}$ and $N_{lf,1,(i)}\cap N_{lf,1,(j)}=\emptyset$ for $i\neq j$ and hence we have partitioned the space $N_{lf,1}$. Thus,
\begin{align*}
\int_{N_{lf,1}}|P|\left(1-\frac{\bar{\rho}}{\bar{\rho}_{\tilde{\mathbf{p}}}}\right)^{-n_x(\mathbf{S}^2)}P_X(dx)&=\sum_{i=1}^\infty \int_{N_{lf,1,(i)}}|P|\left(1-\frac{\bar{\rho}}{\bar{\rho}_{\tilde{\mathbf{p}}}}\right)^{-n_x(\mathbb{S}^2)}P_X(dx)\\
\intertext{$n_x(\mathbb{S}^2)=i$ for all $x\in N_{lf,1,(i)}$ and so is constant over each partition of the space,}
&=\sum_{i=1}^\infty |P|\left(1-\frac{\bar{\rho}}{\bar{\rho}_{\tilde{\mathbf{p}}}}\right)^{-i}\int_{N_{lf,1,(i)}}P_X(dx)\\
&=\sum_{i=1}^\infty |P|\left(1-\frac{\bar{\rho}}{\bar{\rho}_{\tilde{\mathbf{p}}}}\right)^{-i} P_X(X\in N_{lf,1,(i)})\\
\intertext{but it is easy to see that $P_X(X\in N_{lf,1,(i)})=P_{N_X(\mathbb{S}^2)}(N_X(\mathbb{S}^2)=i)$, and since $X$ is Poisson and defining $a=(1-\frac{\bar{\rho}}{\bar{\rho}_{\tilde{\mathbf{p}}}})$ and $\lambda =\mu(\mathbb{S}^2)$,}
&=|P|\sum_{i=1}^\infty a^{-i}\frac{\lambda^n\mathrm{e}^{-\lambda}}{i!}\\
&=|P|\mathrm{e}^{-\lambda}\sum_{i=1}^\infty\frac{(\lambda/a)^i}{i!}\\
&=|P|\mathrm{e}^{-\lambda}\left(\sum_{i=0}^\infty\frac{(\lambda/a)^i}{i!}-1\right)\\
&=|P|\mathrm{e}^{-\lambda}\left(\mathrm{e}^{-\lambda/a}-1\right).
\end{align*}
Hence we have shown that,
\begin{equation*}
\int_{N_{lf,1}}\left|\hat{J}_{\text{inhom}}(r)\right|P_X(dx)\leq |P|\mathrm{e}^{-\lambda}\left(\mathrm{e}^{-\lambda/a}-1\right)
\end{equation*}
and so by the dominated convergence theorem the expectation of $\hat{J}_{\text{inhom}}(r)$ exists.

Existence of the variance of $\hat{J}_{\text{inhom}}(r)$ follows simply based on the existence of the expectation,
\begin{align*}
\text{Var}(\hat{J}_{\text{inhom}}(r))=\int_{N_{lf}}\left(\hat{J}_{\text{inhom}}(r)-\mathbb{E}[\hat{J}_{\text{inhom}}(r)]\right)^2 P_X(dx).
\end{align*}
Partitioning the space $N_{lf}$ again into $N_{lf,0}\equiv\{x\in N_{lf}| N(x) = 0\}$ and $N_{lf,1}\equiv\{x\in N_{lf}| N(x) > 0\}$ we have that,
\begin{align}
&\text{Var}(\hat{J}_{\text{inhom}}(r))\\
&\phantom{AAA}=\int_{N_{lf,0}}\left(\hat{J}_{\text{inhom}}(r)-\mathbb{E}[\hat{J}_{\text{inhom}}(r)]\right)^2 P_X(dx)+\int_{N_{lf,1}}\left(\hat{J}_{\text{inhom}}(r)-\mathbb{E}[\hat{J}_{\text{inhom}}(r)]\right)^2 P_X(dx).\label{var:int:bound}
\end{align}
Taking the convention that $\frac{0}{0}=1$ then the first term is simply,
\begin{align*}
\int_{N_{lf,0}}\left(\hat{J}_{\text{inhom}}(r)-\mathbb{E}[\hat{J}_{\text{inhom}}(r)]\right)^2 P_X(dx)&=\mathbb{E}^2[\hat{J}_{\text{inhom}}(r)]\int_{N_{lf,0}}P_X(dx)\\
&=\mathbb{E}^2[\hat{J}_{\text{inhom}}(r)]P_X(X\in N_{lf,0}),\\
\intertext{then we note that $P_X(X\in N_{lf,0})=P(N_X(\mathbb{S}^2)=0)$ and so,}
&=\mathbb{E}^2[\hat{J}_{\text{inhom}}(r)]P(N_X(\mathbb{S}^2)=0),
\end{align*}
which is finite since the expectation exists. The second term of Equation \ref{var:int:bound} over the space $N_{lf,1}$ can also be shown to be finite. Since $\mathbb{E}[\hat{J}_{\text{inhom}}(r)]$ is finite we just need to show that $\hat{J}_{\text{inhom}}(r)$ is bounded in order show that the integrand is bounded and hence integrable. But from proving the expectation exists we have that,
\begin{equation*}
0\leq \hat{J}_{\text{inhom}}(r)\leq |P|\left(1-\frac{\bar{\rho}}{\bar{\rho}_{\tilde{\mathbf{p}}}}\right)^{-n_x(\mathbf{S}^2)},
\end{equation*}
and so the integrand of second term in Equation \ref{var:int:bound} is bounded and thus the variance exists.
\end{proof}
\begin{customremark}{S1}
From this theorem the requirement of the process being Poisson was only needed for a closed form of the distribution for its corresponding counting process $N_X(\mathbb{S}^2)$. We can drop the requirement of $X$ being Poisson but instead require that the probability generating function of $N_X(\mathbb{S}^2)$, $G_{N_X(\mathbb{S}^2)}(s),$ has radius of convergence $|s|\leq \left(1-\frac{\bar{\rho}}{\bar{\rho}_{\tilde{\mathbf{p}}}}\right)^{-1}$. This condition would be sufficient for the theorem to hold true.
\end{customremark}

\subsection{Covariance between $\hat{H}_{\text{inhom}}(r)$ and $\hat{F}_{\text{inhom}}(r)$}

\begin{proposition}\label{prop:J:moment:approx:app}
  Let $X$ be a spheroidal Poisson process with known intensity function $\rho:\mathbb{S}^2\mapsto\mathbb{R}$. Then the covariance between $1-\hat{H}_{\text{inhom}}(r)$ and $1-\hat{F}_{\text{inhom}}(r)$ for $r\in [0,\pi]$ is,
  \begin{equation*}
  \begin{split}
  &\text{Cov}(1-\hat{H}_{\text{inhom}}(r),1-\hat{F}_{\text{inhom}}(r))\\
  &=\frac{1}{|P|}\sum_{\mathbf{p}\in P}\int_{\mathbb{S}^2}\left(1-\frac{\bar{\rho}\mathbbm{1}[\mathbf{x}\in B_{\mathbb{S}^2}(\mathbf{p},r)]}{\rho(\mathbf{x})}\right)\\
  &\phantom{=-}\frac{\mathrm{exp}\left\{-2\bar{\rho}2\pi(1-\cos r)-\int_{B_{\mathbb{S}^2}(\mathbf{x},r)\cap B_{\mathbb{S}^2}(\mathbf{p},r)}\frac{\bar{\rho}^2}{\rho(\mathbf{y})}\lambda_{\mathbb{S}^2}(d\mathbf{y})\right\}}{A(\mathbf{x},\mathbf{p})}\frac{\rho(\mathbf{x})}{\mu(\mathbb{S}^2)}\lambda_{\mathbb{S}^2}(d\mathbf{x})\\
  &\phantom{=}-\exp(-2\pi(1-\cos r)\bar{\rho})\left(\exp(-2\pi(1-\cos r)\bar{\rho}\right)-\exp(-\mu(\mathbb{S}^2))\frac{\mu(\mathbb{S}^2)}{\mu(\mathbb{S}^2)-2\pi(1-\cos r)\bar{\rho}},
  \end{split}
  \end{equation*}
  where $P$ is a finite grid of points on $\mathbb{S}^2$ and,
  \begin{equation*}
  A(\mathbf{x},\mathbf{p})=1-\frac{2\bar{\rho}}{\mu(\mathbb{S}^2)}2\pi(1-\cos r) +\frac{1}{\mu(\mathbb{S}^2)}\int_{ B_{\mathbb{S}^2}(\mathbf{x},r)\cap B_{\mathbb{S}^2}(\mathbf{p},r)}\frac{\bar{\rho}^2}{\rho(\mathbf{y})}\lambda_{\mathbb{S}^2}(d\mathbf{y}).
  \end{equation*}
\end{proposition}

\begin{proof}
Define,
\begin{align*}
X &\equiv 1-\hat{H}_{\text{inhom}}(r) = \frac{1}{N_X(\mathbb{S}^2)}\sum_{\mathbf{x}\in X}\prod_{\mathbf{y}\in X\setminus \{\mathbf{x}\}}\left(1-\frac{\bar{\rho}\mathbbm{1}[\mathbf{y}\in B_{\mathbb{S}^2}(\mathbf{x},r)]}{\rho(\mathbf{y})}\right)\\
Y &\equiv 1-\hat{F}_{\text{inhom}}(r) = \frac{1}{|P|}\sum_{\mathbf{p}\in P}\prod_{\mathbf{y}\in X }\left(1-\frac{\bar{\rho}\mathbbm{1}[\mathbf{y}\in B_{\mathbb{S}^2}(\mathbf{p},r)]}{\rho(\mathbf{y})}\right)
\end{align*}
Then,
\begin{equation*}
\text{Cov}(X,Y)=\underbrace{\mathbb{E}[XY]}_{(a)}-\underbrace{\mathbb{E}[X]\mathbb{E}[Y]}_{(b)}
\end{equation*}
Term $(a)$ is given by,
\begin{align*}
&\mathbb{E}[XY]\\
&=\mathbb{E}\left[\left(\frac{1}{N_X(\mathbb{S}^2)}\sum_{\mathbf{x}\in X}\prod_{\mathbf{y}\in X\setminus \{\mathbf{x}\}}\left(1-\frac{\bar{\rho}\mathbbm{1}[\mathbf{y}\in B_{\mathbb{S}^2}(\mathbf{x},r)]}{\rho(\mathbf{y})}\right)\right)\left(\frac{1}{|P|}\sum_{\mathbf{p}\in P}\prod_{\mathbf{y}\in X}\left(1-\frac{\bar{\rho}\mathbbm{1}[\mathbf{y}\in B_{\mathbb{S}^2}(\mathbf{p},r)]}{\rho(\mathbf{y})}\right)\right)\right]\\
&=\frac{1}{|P|}\sum_{\mathbf{p}\in P}\mathbb{E}\left[\frac{1}{N_X(\mathbb{S}^2)}\sum_{\mathbf{x}\in X}\left(\prod_{\mathbf{y}\in X\setminus \{\mathbf{x}\}}\left(1-\frac{\bar{\rho}\mathbbm{1}[\mathbf{y}\in B_{\mathbb{S}^2}(\mathbf{x},r)]}{\rho(\mathbf{y})}\right)\right)\left(\prod_{\mathbf{y}\in X}\left(1-\frac{\bar{\rho}\mathbbm{1}[\mathbf{y}\in B_{\mathbb{S}^2}(\mathbf{p},r)]}{\rho(\mathbf{y})}\right)\right)\right]\\
&=\frac{1}{|P|}\sum_{\mathbf{p}\in P}\mathbb{E}\left(\frac{1}{N_X(\mathbb{S}^2)}\sum_{i = 1}^{N_X(\mathbb{S}^2)}\right.\\
&\phantom{AAAA}\left.\underbrace{\mathbb{E}\left[\left(\prod_{\substack{j=1\\j\neq i}}^n\left(1-\frac{\bar{\rho}\mathbbm{1}[\mathbf{X}_j\in B_{\mathbb{S}^2}(\mathbf{X}_i,r)]}{\rho(\mathbf{X}_j)}\right)\right)\left(\prod_{j=1}^n\left(1-\frac{\bar{\rho}\mathbbm{1}[\mathbf{X}_j\in B_{\mathbb{S}^2}(\mathbf{p},r)]}{\rho(\mathbf{X}_j)}\right)\right)\Bigg|N_X(\mathbb{S}^2)=n\right]}_{(\star)}\right),
\end{align*}
where $X_k,\; k=1,\dots,n$ are independently distributed on $\mathbb{S}^2$ with density proportional to $\rho(\mathbf{x}_k)$, when conditioned on $N_X(\mathbb{S}^2)$. Taking $(\star)$,
\begin{align}
&\mathbb{E}\left[\left(\prod_{\substack{j=1\\j\neq i}}^n\left(1-\frac{\bar{\rho}\mathbbm{1}[\mathbf{X}_j\in B_{\mathbb{S}^2}(\mathbf{X}_i,r)]}{\rho(\mathbf{X}_j)}\right)\right)\left(\prod_{j=1}^n\left(1-\frac{\bar{\rho}\mathbbm{1}[\mathbf{X}_j\in B_{\mathbb{S}^2}(\mathbf{p},r)]}{\rho(\mathbf{X}_j)}\right)\right)\Bigg|N_X(\mathbb{S}^2)=n\right]\nonumber\\
&=\mathbb{E}\left[\left(1-\frac{\bar{\rho}\mathbbm{1}[\mathbf{X}_i\in B_{\mathbb{S}^2}(\mathbf{p},r)]}{\rho(\mathbf{X}_i)}\right)
\left(\prod_{\substack{j=1\\j\neq i}}^n\left(1-\frac{\bar{\rho}\mathbbm{1}[\mathbf{X}_j\in B_{\mathbb{S}^2}(\mathbf{X}_i,r)]}{\rho(\mathbf{X}_j)}\right)\left(1-\frac{\bar{\rho}\mathbbm{1}[\mathbf{X}_j\in B_{\mathbb{S}^2}(\mathbf{p},r)]}{\rho(\mathbf{X}_j)}\right)\right)\right],\nonumber\\
\intertext{Using iterated expectations and conditioning on $\mathbf{X}_i$,}
&=\mathbb{E}\left[\left(1-\frac{\bar{\rho}\mathbbm{1}[\mathbf{X}_i\in B_{\mathbb{S}^2}(\mathbf{p},r)]}{\rho(\mathbf{X}_i)}\right)
\prod_{\substack{j=1\\j\neq i}}^n\mathbb{E}\left[\left(1-\frac{\bar{\rho}\mathbbm{1}[\mathbf{X}_j\in B_{\mathbb{S}^2}(\mathbf{x},r)]}{\rho(\mathbf{X}_j)}\right)\left(1-\frac{\bar{\rho}\mathbbm{1}[\mathbf{X}_j\in B_{\mathbb{S}^2}(\mathbf{p},r)]}{\rho(\mathbf{X}_j)}\right)\right]\right]\nonumber\\
&=\mathbb{E}\left[\left(1-\frac{\bar{\rho}\mathbbm{1}[\mathbf{X}_i\in B_{\mathbb{S}^2}(\mathbf{p},r)]}{\rho(\mathbf{X}_i)}\right)
\mathbb{E}^{n-1}\left[\left(1-\frac{\bar{\rho}\mathbbm{1}[\mathbf{Y}\in B_{\mathbb{S}^2}(\mathbf{x},r)]}{\rho(\mathbf{Y})}\right)\left(1-\frac{\bar{\rho}\mathbbm{1}[\mathbf{Y}\in B_{\mathbb{S}^2}(\mathbf{p},r)]}{\rho(\mathbf{Y})}\right)\right]\right]\label{iterated:1}
\end{align}
where $\mathbf{Y}$ is distributed with density proportional to $\rho(\mathbf{y})$ on $\mathbb{S}^2$. It is then easy to show that,
\begin{align*}
&\mathbb{E}\left[\left(1-\frac{\bar{\rho}\mathbbm{1}[\mathbf{Y}\in B_{\mathbb{S}^2}(\mathbf{x},r)]}{\rho(\mathbf{Y})}\right)\left(1-\frac{\bar{\rho}\mathbbm{1}[\mathbf{Y}\in B_{\mathbb{S}^2}(\mathbf{p},r)]}{\rho(\mathbf{Y})}\right)\right]\\
&\phantom{AAAA}=1-\frac{2\bar{\rho}}{\mu(\mathbb{S}^2)}2\pi(1-\cos r) +\frac{1}{\mu(\mathbb{S}^2)}\int_{B_{\mathbb{S}^2}(\mathbf{x},r)\cap B_{\mathbb{S}^2}(\mathbf{p},r)}\frac{\bar{\rho}^2}{\rho(\mathbf{y})}\lambda_{\mathbb{S}^2}(d\mathbf{y}).
\end{align*}
Let us define $A(\mathbf{x},\mathbf{p})=\mathbb{E}\left[\left(1-\frac{\bar{\rho}\mathbbm{1}[\mathbf{Y}\in B_{\mathbb{S}^2}(\mathbf{x},r)]}{\rho(\mathbf{Y})}\right)\left(1-\frac{\bar{\rho}\mathbbm{1}[\mathbf{Y}\in B_{\mathbb{S}^2}(\mathbf{p},r)]}{\rho(\mathbf{Y})}\right)\right]$, and so returning to \ref{iterated:1},
\begin{align*}
&\mathbb{E}\left[\left(1-\frac{\bar{\rho}\mathbbm{1}[\mathbf{X}_i\in B_{\mathbb{S}^2}(\mathbf{p},r)]}{\rho(\mathbf{X}_i)}\right)
A^{n-1}(\mathbf{X}_i,\mathbf{p})\right] = \mathbb{E}\left[\left(1-\frac{\bar{\rho}\mathbbm{1}[\mathbf{X}\in B_{\mathbb{S}^2}(\mathbf{p},r)]}{\rho(\mathbf{X})}\right)
A^{n-1}(\mathbf{X},\mathbf{p})\right],\\
\intertext{where $\mathbf{X}$ has density proportional to $\rho(\mathbf{x})$ on $\mathbb{S}^2$,}
&\phantom{AAAAAAAAAAAA}=\int_{\mathbb{S}^2}\left(1-\frac{\bar{\rho}\mathbbm{1}[\mathbf{x}\in B_{\mathbb{S}^2}(\mathbf{p},r)]}{\rho(\mathbf{x})}\right)
A^{n-1}(\mathbf{x},\mathbf{p})\frac{\rho(\mathbf{x})}{\mu(\mathbb{S}^2)}\lambda_{\mathbb{S}^2}(d\mathbf{x}).
\end{align*}
And so,
\begin{align*}
\mathbb{E}[XY]&=\frac{1}{|P|}\sum_{\mathbf{p}\in P}\mathbb{E}\left[\frac{1}{N_X(\mathbb{S}^2)}\sum_{i=1}^{N_X(\mathbb{S}^2)}\int_{\mathbb{S}^2}\left(1-\frac{\bar{\rho}\mathbbm{1}[\mathbf{x}\in B_{\mathbb{S}^2}(\mathbf{p},r)]}{\rho(\mathbf{x})}\right)
A^{N_X(\mathbb{S}^2)-1}(\mathbf{x},\mathbf{p})\frac{\rho(\mathbf{x})}{\mu(\mathbb{S}^2)}\lambda_{\mathbb{S}^2}(d\mathbf{x})\right]\\
&=\frac{1}{|P|}\sum_{\mathbf{p}\in P}\mathbb{E}\left[\int_{\mathbb{S}^2}\left(1-\frac{\bar{\rho}\mathbbm{1}[\mathbf{x}\in B_{\mathbb{S}^2}(\mathbf{p},r)]}{\rho(\mathbf{x})}\right)
A^{N_X(\mathbb{S}^2)-1}(\mathbf{x},\mathbf{p})\frac{\rho(\mathbf{x})}{\mu(\mathbb{S}^2)}\lambda_{\mathbb{S}^2}(d\mathbf{x})\right]\\
&=\frac{1}{|P|}\sum_{\mathbf{p}\in P}\int_{\mathbb{S}^2}\left(1-\frac{\bar{\rho}\mathbbm{1}[\mathbf{x}\in B_{\mathbb{S}^2}(\mathbf{p},r)]}{\rho(\mathbf{x})}\right)
\mathbb{E}\left[A^{N_X(\mathbb{S}^2)-1}(\mathbf{x},\mathbf{p})\right]\frac{\rho(\mathbf{x})}{\mu(\mathbb{S}^2)}\lambda_{\mathbb{S}^2}(d\mathbf{x})
\end{align*}
Then it can easily be shown that,
\begin{align*}
\mathbb{E}\left[A^{N_X(\mathbb{S}^2)-1}(\mathbf{x},\mathbf{p})\right]=\frac{\mathrm{exp}\left\{-2\bar{\rho}2\pi(1-\cos r)-\int_{ B_{\mathbb{S}^2}(\mathbf{x},r)\cap B_{\mathbb{S}^2}(\mathbf{p},r)}\frac{\bar{\rho}^2}{\rho(\mathbf{y})}\lambda_{\mathbb{S}^2}(d\mathbf{y})\right\}}{A(\mathbf{x},\mathbf{p})}.
\end{align*}
Then from Theorem \ref{expect:FS} we have that,
\begin{align*}
\mathbb{E}[X]&=\exp(2\pi(1-\cos r)\bar{\rho})\\
\mathbb{E}[Y]&=\mu^2(\mathbb{S})\frac{\exp(2\pi(1-\cos r)\bar{\rho})-\exp(\mu(\mathbb{S}^2))}{\mu(\mathbb{S}^2)-2\pi(1-\cos r)\bar{\rho}}.
\end{align*}
And so we have the covariance between $\hat{H}_{\text{inhom}}(r)$ and $\hat{F}_{\text{inhom}}(r)$.
\end{proof}

\subsection{Taylor Series Expansion}
We discuss Taylor series expansions in general and their use to approximate moments of random variables \cite{wolter2007}, after which we apply this to the  $\hat{J}_{\text{inhom}}$-function. For any function $f:\mathbb{R}^2\mapsto\mathbb{R}$, we have its Taylor expansion up to second order around $\boldsymbol{\theta}=(\theta_x,\theta_y)$ as,
\begin{equation*}
\begin{split}
f(x,y)& = f(\theta_x,\theta_y)+\frac{\partial f}{\partial  x}\bigg|_{\theta_X,\theta_Y}(x-\theta_x)+\frac{\partial f}{\partial y}\bigg|_{\theta_X,\theta_Y}(y-\theta_y)+\\
&\frac{1}{2}\left[\frac{\partial^2 f}{\partial x^2}\bigg|_{\theta_X,\theta_Y}(x-\theta_x)^2+2\frac{\partial^2 f}{\partial x\partial y}\bigg|_{\theta_X,\theta_Y}(x-\theta_x)(y-\theta_x)+\frac{\partial^2 f}{\partial y^2}\bigg|_{\theta_X,\theta_Y}(y-\theta_y)^2\right] + R(x,y),
\end{split}
\end{equation*}
where $R(x,y)$ is a remainder term. Thus using random variables $X$ and $Y$, with $\boldsymbol{\theta}=(\mathbb{E}[X],\mathbb{E}[Y])\equiv (\mu_X,\mu_Y)$, whilst also assuming that $\mathbb{E}[R(X,Y)]$ is close to 0 then,
\begin{equation*}
\mathbb{E}[f(X,Y)]\approx f(\mu_X,\mu_Y) +\frac{1}{2}\left[\frac{\partial^2 f}{\partial x^2}\text{Var}(X)+2\frac{\partial^2 f}{\partial x\partial y}\text{Cov}(X,Y)+\frac{\partial^2 f}{\partial y^2}\text{Var}(Y)\right].
\end{equation*}
Then to approximate the variance we first note that using the first order Taylor series expansion $\mathbb{E}[f(X,Y)]\approx f(\mu_X,\mu_Y)$. Then,
\begin{align*}
\text{Var}(f(X,Y))&=\mathbb{E}[(f(X,Y)-\mathbb{E}[f(X,Y)])^2]\\
&\approx\mathbb{E}[(f(X,Y)-f(\mu_X,\mu_Y))^2]\\
&\approx\mathbb{E}\left[\left(f(\mu_X,\mu_Y)+\frac{\partial f}{\partial  x}\bigg|_{\mu_X,\mu_Y}(X-\mu_X)+\frac{\partial f}{\partial y}\bigg|_{\mu_X,\mu_Y}(Y-\mu_Y)-f(\mu_X,\mu_Y)\right)^2\right]\\
&=\mathbb{E}\left[\left(\frac{\partial f}{\partial  x}\bigg|_{\mu_X,\mu_Y}(X-\mu_X)+\frac{\partial f}{\partial y}\bigg|_{\mu_X,\mu_Y}(Y-\mu_Y)\right)^2\right]\\
&=\frac{\partial f}{\partial  x}\bigg|_{\mu_X,\mu_Y}^2\text{Var}(X) + 2\frac{\partial f}{\partial  x}\bigg|_{\mu_X,\mu_Y}\frac{\partial f}{\partial y}\bigg|_{\mu_X,\mu_Y}\text{Cov}(X,Y)+\frac{\partial f}{\partial  y}\bigg|_{\mu_X,\mu_Y}^2\text{Var}(Y)
\end{align*}
Hence defining $f(X,Y)=\frac{X}{Y}$, we have the following approximations to the first and second order moments of $f(X,Y),$
\begin{align}
\mathbb{E}\left[\frac{X}{Y}\right] &\approx\frac{\mu_X}{\mu_Y}-\frac{\text{Cov}(X,Y)}{\mu_Y^2}+\frac{\text{Var}(Y)\mu_X}{\mu_Y^3}\label{E:approx}\\
\text{Var}\left(\frac{X}{Y}\right) &\approx \frac{\mu_X}{\mu_Y}\left[\frac{\text{Var}(X)}{\mu_X^2}-2\frac{\text{Cov}(X,Y)}{\mu_X\mu_Y}+\frac{\text{Var}(Y)}{\mu_Y^2}\right].\label{V:approx}
\end{align}

Then since $\hat{J}_{\text{inhom}}(r)$ is defined as the ratio of two random variables, in particular $1-\hat{H}_{\text{inhom}}(r)$ and $1-\hat{F}_{\text{inhom}}(r)$. Thus combining Proposition \ref{prop:J:moment:approx:app} with the Taylor series expansions given in Equations \ref{E:approx} and \ref{V:approx} provides an estimate for the expectation and variance of $\hat{J}_{\text{inhom}}(r)$.

\section{Regular and Cluster Processes}\label{Supplementary:reg:cluster}

In this section we discuss properties of the regular and cluster processes discussed in the paper. These properties help to design the simulations discussed in the paper. We also discuss why we observe decreasing power of our test statistic as $\mathbb{D}$ \emph{deforms} further away from the unit sphere.

\subsection{Properties of Regular and Cluster Processes}

The following proposition gives the expected number of points for Mat\'{e}rn \rom{1} and \rom{2} processes in any subset of $\mathbb{D}$.

\begin{customprop}{S1}
Let $X_1$ and $X_2$ be a Mat\'{e}rn \rom{1} and \rom{2} inhibition processes respectively. In the Mat\'{e}rn \rom{2} case we also define a mark distribution $P_{M_{\mathbf{x}}}$ such that the mark is independent not only of all other points $\mathbf{y}\in X_2\setminus\{\mathbf{x}\}$ and marks $M_\mathbf{y},$ for $\mathbf{y}\in X_2\setminus\{\mathbf{x}\}$ but also of $\mathbf{x}$, the point associated to the mark $M_{\mathbf{x}}$. Define $N_1$ and $N_2$ to be the random counting measures of $X_1$ and $X_2$. Then the expectations of $N_{X_1}$ and $N_{X_2}$ are given as,
\begin{align*}
\mathbb{E}[N_{X_1}(B)]&=\rho\int_B \mathrm{e}^{-\rho\lambda_{\mathbb{D}}(B_{\mathbb{D}}(\mathbf{x},R))}  \lambda_{\mathbb{D}}(d\mathbf{x})\\
\mathbb{E}[N_{X_2}(B)]&=\int_B \frac{1-\mathrm{e}^{-\rho\lambda_{\mathbb{D}}(B_{\mathbb{D}}(\mathbf{x},R))}}{\lambda_{\mathbb{D}}(B_{\mathbb{D}}(\mathbf{x},R))} \lambda_{\mathbb{D}}(d\mathbf{x}),
\end{align*}
where $B\subseteq \mathbb{D}$.
\end{customprop}
\begin{proof}
Let us first start with $X_1$. Define $Y_1\sim PPP(\rho,\mathbb{D})$ to be the homogeneous Poisson process which is thinned to give $X_1$. Then $\forall{B}\subseteq \mathbb{D}$ we can rewrite the counting measure for $X_1$ as,
\begin{align*}
N_{X_1}(B)=\sum_{\mathbf{x}\in Y_1\cap B}\mathbbm{1}[N_{Y_1\setminus\{\mathbf{x}\}}(B_{\mathbb{D}}(\mathbf{x},R))=0],
\end{align*}
where $N_{Y_1\setminus\{\mathbf{x}\}}$, is the random counting measure for the process $Y_1$ without the point $\mathbf{x}$. Then taking expectations and using the Slivnyak-Mecke Theorem,
\begin{align*}
\mathbb{E}[N_{X_1}(B)] &= \mathbb{E}\left[\sum_{\mathbf{x}\in Y_1\cap B}\mathbbm{1}[N_{Y_1\setminus\{\mathbf{x}\}}(B_{\mathbb{D}}(\mathbf{x},R))=0]\right]\\
&= \int_B \mathbb{E}\left[\mathbbm{1}[N_{Y_1}(B_{\mathbb{D}}(\mathbf{x},R))=0]\right]\rho \lambda_{\mathbb{D}}(d\mathbf{x})\\
& = \rho\int_B \mathbb{P}(N_{Y_1}(B_{\mathbb{D}}(\mathbf{x},R))=0)\lambda_{\mathbb{D}}(d\mathbf{x})\\
&= \rho\int_B \mathrm{e}^{-\rho\lambda_{\mathbb{D}}(B_{\mathbb{D}}(\mathbf{x},R))}\lambda_{\mathbb{D}}(d\mathbf{x})
\end{align*}
For $N_{X_2}(B)$ a few more steps are required in order to take into account the mark associated with each point. Similarly to the counting measure for $X_1$, we can rewrite the counting measure for $X_2$ as,
\begin{equation*}
N_{X_2}(B)=\sum_{\mathbf{x}\in Y_2\cap B} \mathbbm{1}[M_{\mathbf{x}}\leq M_{\mathbf{y}}, \forall \mathbf{y}\in (Y_2\setminus\{\mathbf{x}\})\cap B_{\mathbb{D}}(\mathbf{x},R)].
\end{equation*}
By again taking expectations and using the Slivnyak-Mecke Theorem,
\begin{align*}
\mathbb{E}[N_{X_2}(B)]&=\mathbb{E}\left[\sum_{\mathbf{x}\in Y_2\cap B} \mathbbm{1}[M_{\mathbf{x}}\leq M_{\mathbf{y}}, \forall \mathbf{y}\in (Y_2\setminus\{\mathbf{x}\})\cap B_{\mathbb{D}}(\mathbf{x},R)]\right]\\
&=\int_B\mathbb{E}[\mathbbm{1}[M_{\mathbf{x}}\leq M_{\mathbf{y}}, \forall \mathbf{y}\in Y_2\cap B_{\mathbb{D}}(\mathbf{x},R)]]\rho \lambda_{\mathbb{D}}(d\mathbf{x})\\
&=\int_B\mathbb{P}(M_{\mathbf{x}}\leq M_{\mathbf{y}}, \forall \mathbf{y}\in Y_2\cap B_{\mathbb{D}}(\mathbf{x},R))\rho \lambda_{\mathbb{D}}(d\mathbf{x})
\end{align*}
Define $\lambda_{\mathbf{x}}=\rho\lambda_{\mathbb{D}}(B_{\mathbb{D}}(\mathbf{x},R))$, then the probability is calculated as follows,
\begin{align*}
&\mathbb{P}(M_{\mathbf{x}}\leq M_{\mathbf{y}}, \forall \mathbf{y}\in Y_2\cap B_{\mathbb{D}}(\mathbf{x},R))\\
&\phantom{AAAA}= \sum_{n=0}^\infty \mathbb{P}\left(M_{\mathbf{x}}\leq M_{\mathbf{y}}, \forall \mathbf{y}\in Y_2\cap B_{\mathbb{D}}(\mathbf{x},R)| N_{Y_2}(B_{\mathbb{D}}(\mathbf{x},R))=n\right)\mathbb{P}(N_{Y_2}(B_{\mathbb{D}}(\mathbf{x},R))=n)\\
&\phantom{AAAA}=\sum_{n=0}^\infty\frac{\lambda_{\mathbf{x}}^n\mathrm{e}^{-\lambda_{\mathbf{x}}}}{n!}\mathbb{P}\left(M_{\mathbf{x}}\leq M_{\mathbf{y}}, \forall \mathbf{y}\in Y_2\cap B_{\mathbb{D}}(\mathbf{x},R)| N_{Y_2}(B_{\mathbb{D}}(\mathbf{x},R))=n\right)
\end{align*}
Let us label the points $\mathbf{y}_i\in Y_2,$ for $i=1,\dots,n$ to be the $n$ points in $B_{\mathbb{D}}(\mathbf{x},R)$ coming from the process $Y_2$. Then the event $M_{\mathbf{x}}\leq M_{\mathbf{y}}, \forall \mathbf{y}\in Y_2\cap B_{\mathbb{D}}(\mathbf{x},R)$ given that $N_{Y_2}(B_{\mathbb{D}}(\mathbf{x},R))=n$ is identical to the event that the mark associated to the point $\mathbf{x}$ is the smallest of all marks. In other words we are concerned with the event $M_{\mathbf{x}}\leq \mathrm{min}\{M_{\mathbf{y}_1},\dots,M_{\mathbf{y}_n}\}$. Let us define $M_{\mathrm{min}}=\mathrm{min}\{M_{\mathbf{y}_1},\dots,M_{\mathbf{y}_n}\}$ and $f_{M^{\mathrm{min}}}(m)$ to be the density function of $M_{\mathrm{min}}$, then by using order statistics we know that $f_{M^{\mathrm{min}}}(m)=nf_{M_\mathbf{x}}(m)(1-F_{M_{\mathbf{x}}}(m))^{n-1}$, where $f_{M_\mathbf{x}}(m)$  and $F_{M_{\mathbf{x}}}(m)$ are the density and cumulative density functions of an individual markis respectively. Then we have,
\begin{align*}
&\mathbb{P}\left(M_{\mathbf{x}}\leq M_{\mathbf{y}}, \forall \mathbf{y}\in Y_2\cap B_{\mathbb{D}}(\mathbf{x},R)| N_{Y_2}(B_{\mathbb{D}}(\mathbf{x},R))=n\right)\\
&\phantom{AAAA}= \mathbb{P}(M_{\mathbf{x}}\leq M^{\mathrm{min}})\\
&\phantom{AAAA}=\int_{\mathbb{M}^{\mathrm{min}}}\mathbb{P}(M_{\mathbf{x}}\leq m) f_{M^{\mathrm{min}}}(m)dm\\
&\phantom{AAAA}=\int_{\mathbb{M}^{\mathrm{min}}} F_{M_{\mathbf{x}}}(m)nf_{M_\mathbf{x}}(m)(1-F_{M_{\mathbf{x}}}(m))^{n-1} dm\\
&\phantom{AAAA}= \left[F_{M_{\mathbf{x}}}(m)(1-F_{M_{\mathbf{x}}}(m))^n\right]^{\mathbb{M}^{\mathrm{min}}}+\int_{\mathbb{M}^\mathrm{min}} f_{M_{\mathbf{x}}}(m)(1-F_{M_{\mathbf{x}}}(m))^n dm\\
&\phantom{AAAA}=\left[F_{M_{\mathbf{x}}}(m)(1-F_{M_{\mathbf{x}}}(m))^n\right]^{\mathbb{M}^{\mathrm{min}}}+\left[\frac{1}{n+1}(1-F_{M_{\mathbf{x}}}(m))^{n+1}\right]^{\mathbb{M}^{\mathrm{min}}}\\
&\phantom{AAAA}=\frac{1}{n+1},
\end{align*} 
where the last line follows since the range of $\mathbb{M}^{\mathrm{min}}$ is identical to the range of $M_{\mathbf{x}}$. Returning to $\mathbb{P}(\mathbf{x}\in X_2)$,
\begin{align*}
\mathbb{P}(M_{\mathbf{x}}\leq M_{\mathbf{y}}, \forall \mathbf{y}\in Y_2\cap B_{\mathbb{D}}(\mathbf{x},R)) &=\sum_{n=0}^\infty\frac{\lambda_{\mathbf{x}}^n\mathrm{e}^{-\lambda_{\mathbf{x}}}}{n!}\frac{1}{n+1}\\
&=\frac{\mathrm{e}^{-\lambda_{\mathbf{x}}}}{\lambda_{\mathbf{x}}}\sum_{n=0}^\infty \frac{\lambda_{\mathbf{x}}^{n+1}}{(n+1)!}\\
&=\frac{\mathrm{e}^{-\lambda_{\mathbf{x}}}}{\lambda_{\mathbf{x}}}\left[\sum_{n=0}^\infty \frac{\lambda_{\mathbf{x}}^{n}}{n!}-1\right]\\
&=\frac{\mathrm{e}^{-\lambda_{\mathbf{x}}}}{\lambda_{\mathbf{x}}}(\mathrm{e}^{\lambda_{\mathbf{x}}}-1)\\
&=\frac{1-\mathrm{e}^{-\lambda_{\mathbf{x}}}}{\lambda_{\mathbf{x}}}\\
&= \frac{1-\mathrm{e}^{-\rho\lambda_{\mathbb{D}}(B_{\mathbb{D}}(\mathbf{x},R))}}{\rho\lambda_{\mathbb{D}}(B_{\mathbb{D}}(\mathbf{x},R))}.
\end{align*}
Thus returning to the expectation of $N_{X_2}(B)$,
\begin{align*}
\mathbb{E}[N_{X_2}(B)]&=\int_B \frac{1-\mathrm{e}^{-\rho\lambda_{\mathbb{D}}(B_{\mathbb{D}}(\mathbf{x},R))}}{\lambda_{\mathbb{D}}(B_{\mathbb{D}}(\mathbf{x},R))} \lambda_{\mathbb{D}}(d\mathbf{x}).
\end{align*}
\end{proof}
As we are running simulations based on Mat\'{e}rn \rom{2} inhibition processes and given that it is a regular process there is a finite maximum number of points that can arise on the surface $\mathbb{D}$. The following corollary to the previous proposition gives the maximum expected value of $N_{X_2}(\mathbb{D})$ for a fixed hard-core distance.
\begin{customcorollary}{S2}
Let $X$ be a Mat\'{e}rn \rom{2} process over $\mathbb{D}$ with hard-core distance $R$ and defined by a Poisson process with constant intensity function $\rho$. Then,
\begin{equation*}
\mathrm{sup}_{\rho\in\mathbb{R}^+}\mathbb{E}[N_X(\mathbb{D})] = \int_{\mathbb{D}} \frac{1}{\lambda_{\mathbb{D}}(B_{\mathbb{D}}(\mathbf{x},R))}\lambda_{\mathbb{D}}(d\mathbf{x}).
\end{equation*}
\end{customcorollary}
\begin{proof}
Notice that since $\rho\lambda_{\mathbb{D}}(B_{\mathbb{D}}(\mathbf{x},R))$ is always positive this means that $\frac{1-\mathrm{e}^{-\rho\lambda_{\mathbb{D}}(B_{\mathbb{D}}(\mathbf{x},R))}}{\lambda_{\mathbb{D}}(B_{\mathbb{D}}(\mathbf{x},R))}\geq 0,\;\forall \rho\in \mathbb{R}^+$ and $\mathbf{x}\in\mathbb{D}$. Then $\frac{1-\mathrm{e}^{-\rho_1\lambda_{\mathbb{D}}(B_{\mathbb{D}}(\mathbf{x},R))}}{\lambda_{\mathbb{D}}(B_{\mathbb{D}}(\mathbf{x},R))}<\frac{1-\mathrm{e}^{-\rho_2\lambda_{\mathbb{D}}(B_{\mathbb{D}}(\mathbf{x},R))}}{\lambda_{\mathbb{D}}(B_{\mathbb{D}}(\mathbf{x},R))},\;\forall \rho_1<\rho_2$ and so the supremum is when $\rho$ is taken to infinity giving the final result. 
\end{proof}

The following proposition gives the expected number of points for a Thomas-type process on $\mathbb{D}$.

\begin{customprop}{S2}
Let $X$ be a Thomas-type process on $\mathbb{D}$ with concentration parameter $\kappa$ and a constant expected number of offspring, $\lambda$, for each parent point. Then,
\begin{equation*}
\mathbb{E}[N_X(\mathbb{D})]=\lambda_{\mathbb{D}}(\mathbb{D})\rho\lambda.
\end{equation*}
\end{customprop}
\begin{proof}
Let $Y=f(X),$ where $f(\mathbf{x})=\mathbf{x}/||\mathbf{x}||$ then,
\begin{align*}
\mathbb{E}[N_X(\mathbb{D})]&=\mathbb{E}\left[\sum_{\mathbf{c}\in X_p} N_{X_{\mathbf{c}}}(\mathbb{S}^2)\right]\\
&=\int_{\mathbb{D}} \mathbb{E}[N_{X_{\mathbf{c}}}(\mathbb{S}^2)]\rho d\mathbf{c}\\
&=\lambda\rho\int_{\mathbb{D}}d\mathbf{c}\\
&=\lambda_{\mathbb{D}}(\mathbb{D})\rho\lambda.
\end{align*}
\end{proof}

\subsection{Simulation of regular and cluster processes}
In order to examine the power of our hypothesis testing procedure, we need to be able to simulate regular and cluster processes outlined here. Much of it is based on simulation of homogeneous Poisson processes on convex shapes. To simulate homogeneous Poisson processes we can first simulate the number of points in the pattern as $\rho\lambda_{\mathbb{D}}(\mathbb{D})$ and then distribute them uniformly across $\mathbb{D}$. In order to distribute the uniformly across $\mathbb{D}$ we can use the rejection sampler outlined by \cite{Kopytov2012}. 

Simulation of Mat\'{e}rn \rom{1} and \rom{2} then depend upon removing events in an underlying homogeneous Poisson process based on their distance between the events. This depends on being able to calculate the geodesic distance on a given surface. Assuming this can be achieved, then it is simple to simulate Mat\'{e}rn \rom{1} and \rom{2} process using their definitions. For ellipsoidal Mat\'{e}rn \rom{1} and \rom{2} processes we use the \texttt{geographiclib} \cite{Karney2018} available as a \texttt{MATLAB} toolbox.

To simulate a Thomas process this can easily be achieved by rejection sampling again. Simulation of the parents is a homogeneous Poisson process. Then, for each parent, we simulate its random number of offspring and then construct a rejection sampler to sample from
\begin{equation*}
k(\mathbf{x}_p,\mathbf{x})=\frac{1}{\chi(\sigma^2)}\exp\left(-\frac{d^2(\mathbf{x}_p,\mathbf{x})}{2\sigma^2}\right),
\end{equation*}
where $\mathbf{x}_p$ is in the parent process, $\sigma$ is a bandwidth parameter and $\chi(\sigma^2)=\int_{\mathbb{D}}\exp(-d^2(\mathbf{x}_p,\mathbf{x})$ $/2\sigma^2)\lambda_{\mathbb{D}}(d\mathbf{x})$. Again for an ellipsoidal Thomas process we use \texttt{geographiclib} \cite{Karney2018} to calculate geodesic distances on the surface of an ellipsoid.

\subsection{Decreasing power of test statistic}
\begin{figure}[t]
\centering
\includegraphics[width=\textwidth,height=!]{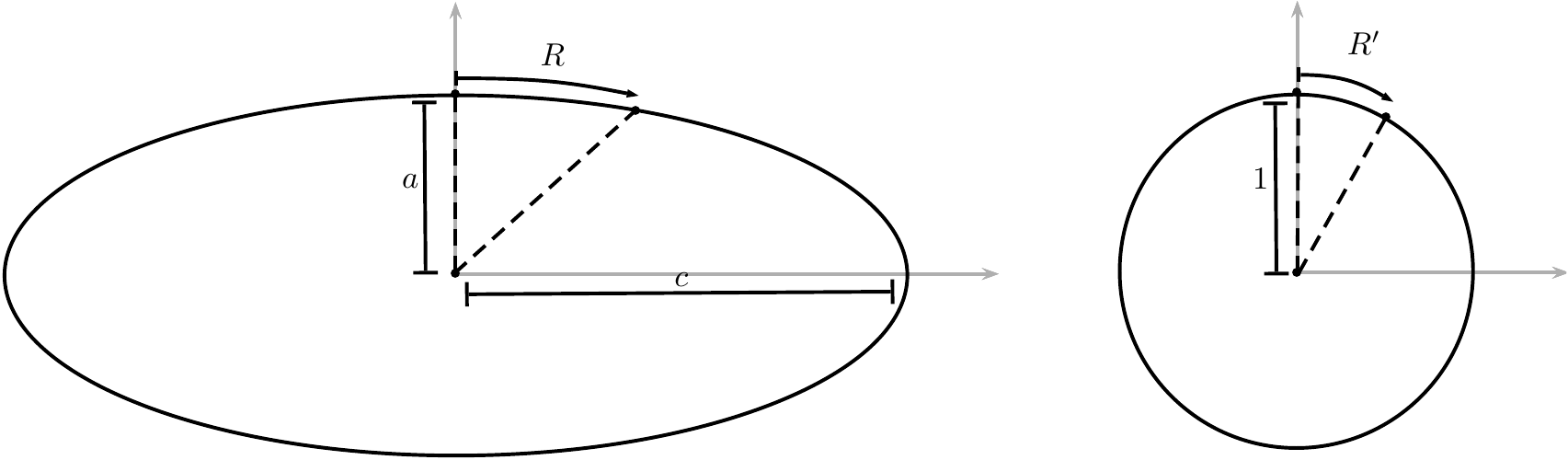}
\caption{Example of hardcore distance reduction due to mapping to a sphere.}
\label{fig:hardcore:red}
\hrulefill
\end{figure}

\begin{figure}[p]
\begin{subfigure}{\textwidth}
\centering
\includegraphics[trim={0 3em 0 3em},clip=TRUE,width=\textwidth,height=11em]{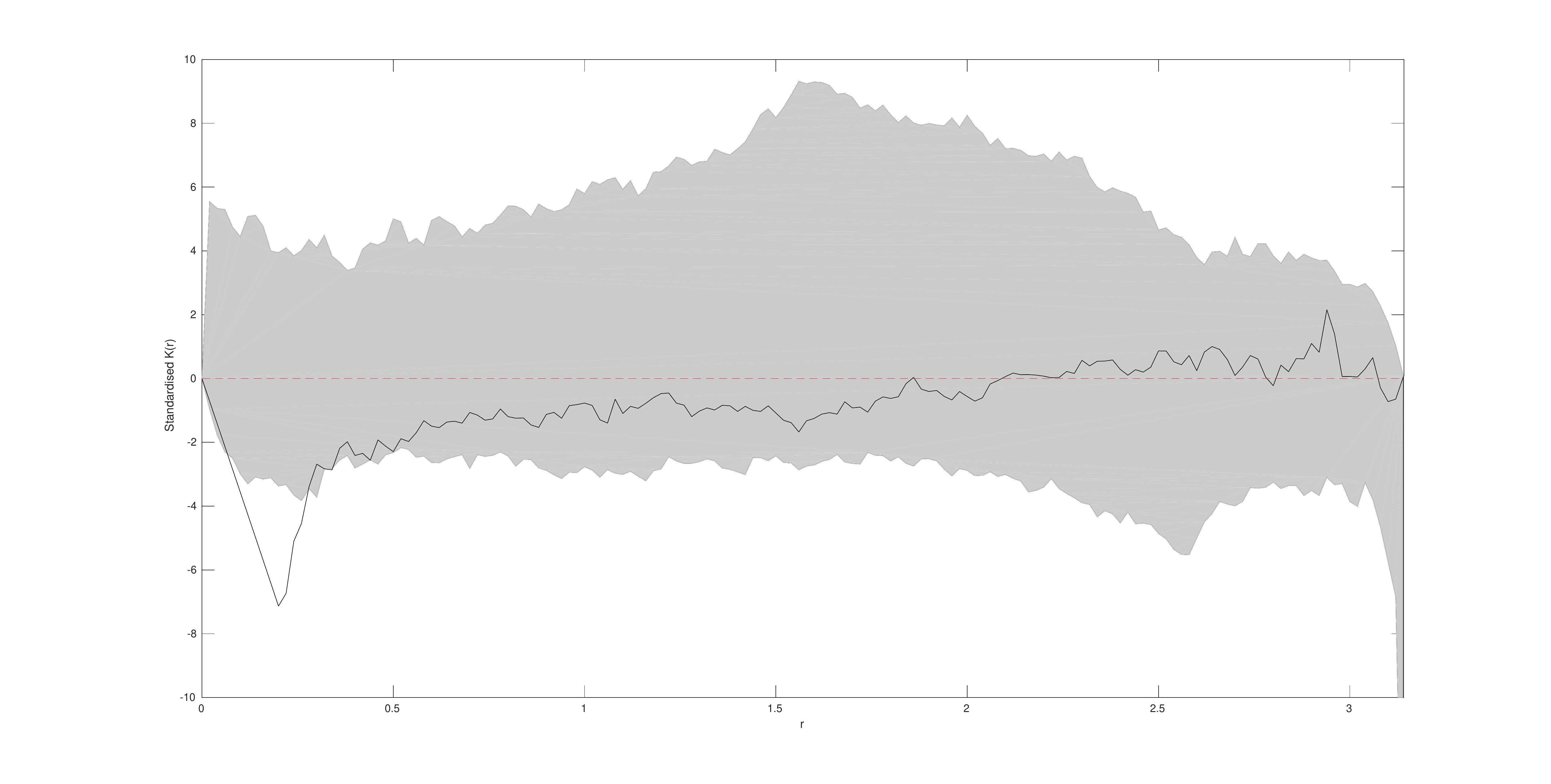}
\end{subfigure}
\begin{subfigure}{\textwidth}
\centering
\includegraphics[trim={0 3em 0 3em},clip=TRUE,width=\textwidth,height=11em]{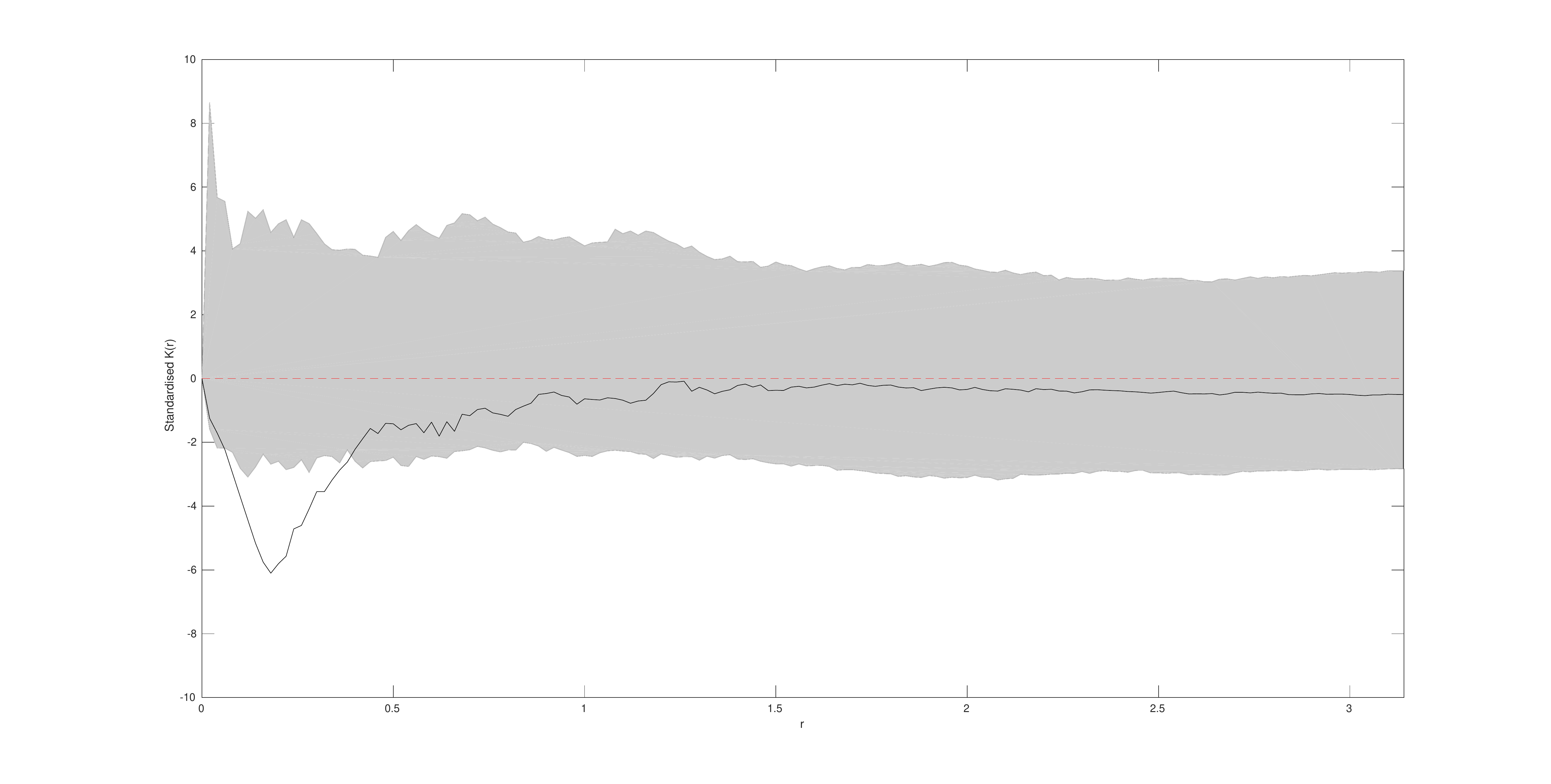}
\end{subfigure}
\begin{subfigure}{\textwidth}
\centering
\includegraphics[trim={0 3em 0 3em},clip=TRUE,width=\textwidth,height=11em]{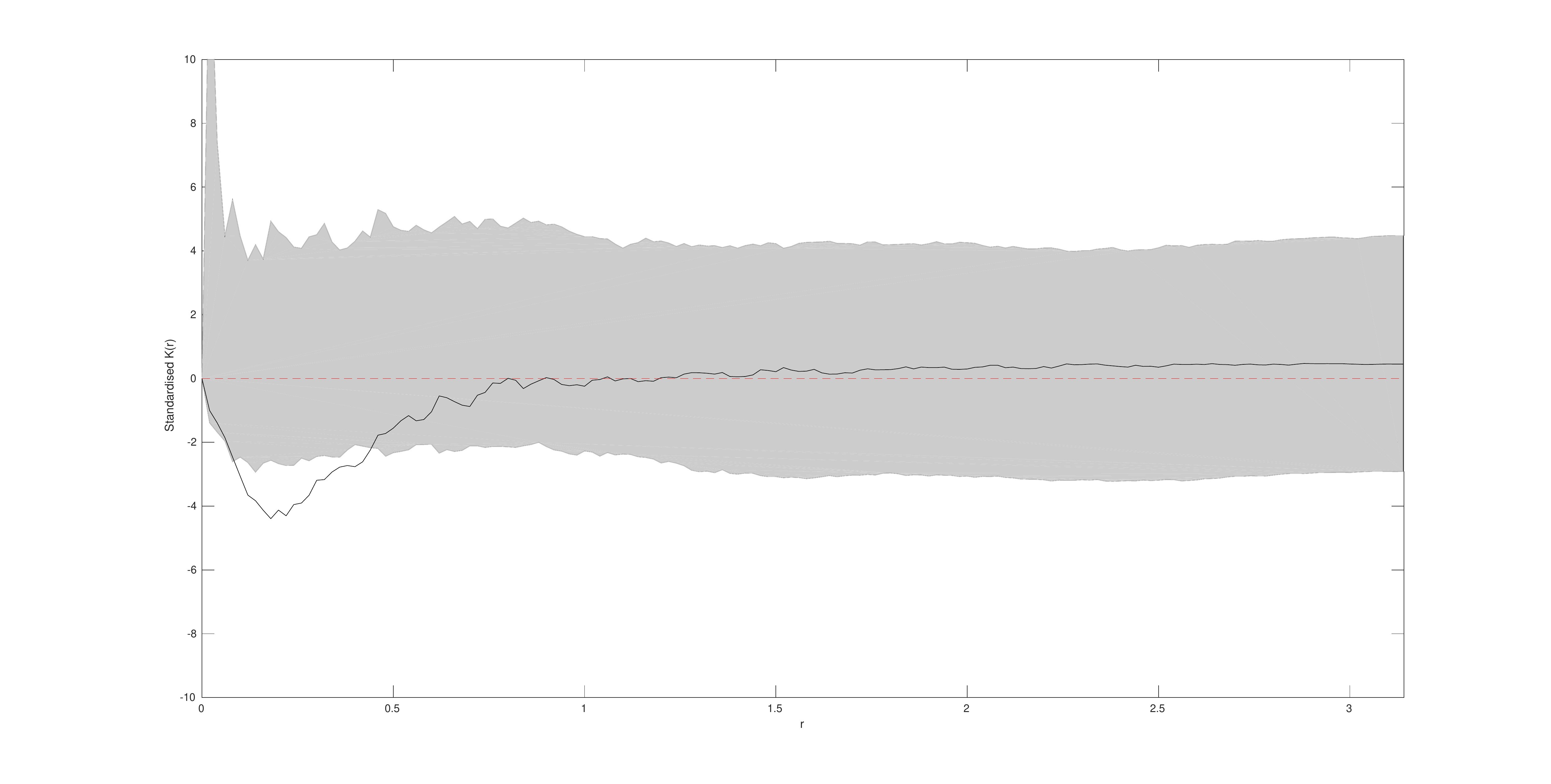}
\end{subfigure}
\begin{subfigure}{\textwidth}
\centering
\includegraphics[trim={0 3em 0 3em},clip=TRUE,width=\textwidth,height=11em]{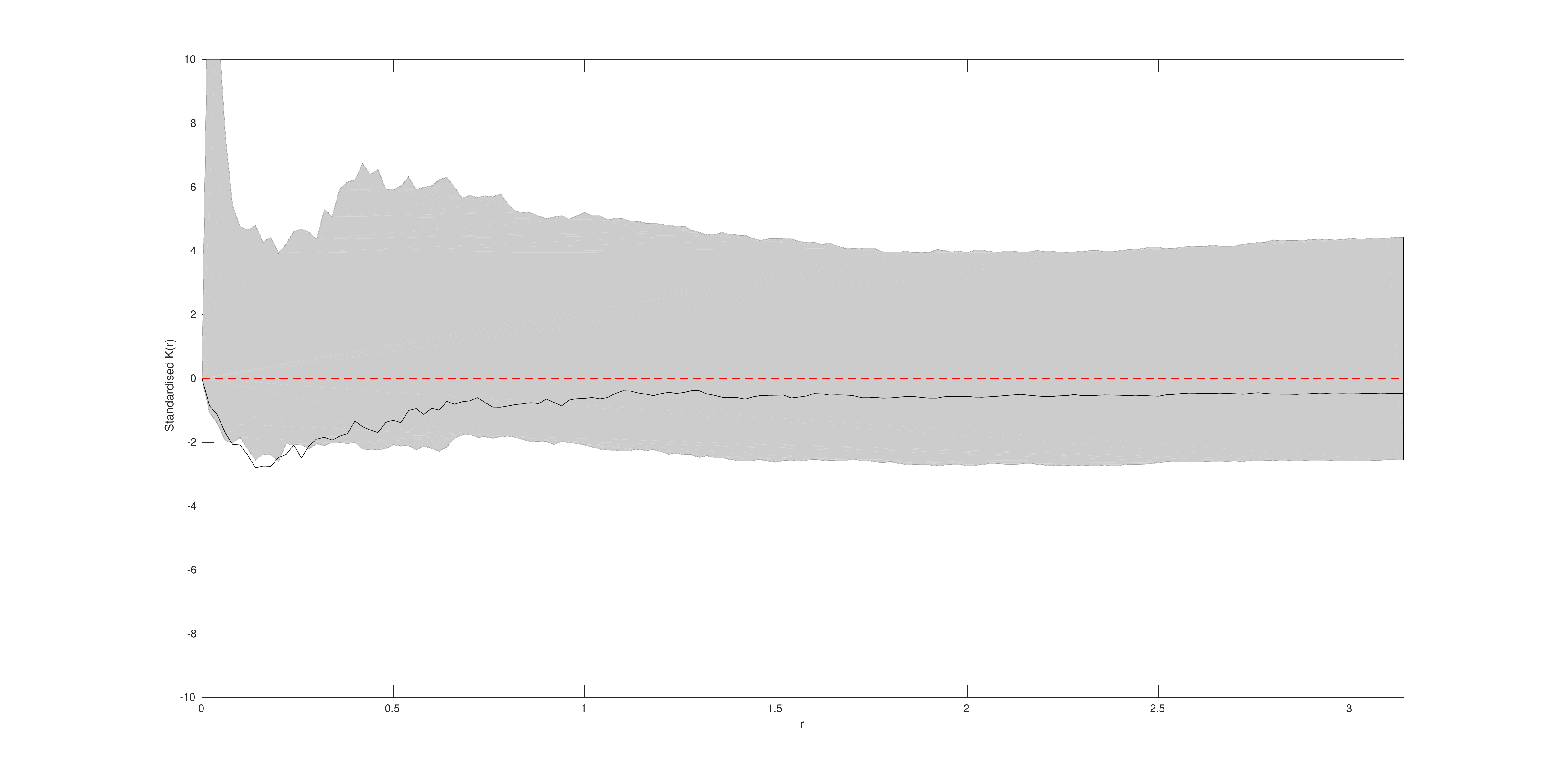}
\end{subfigure}
\caption{Plots of the standardised inhomogeneous $K$-function for a Mat\'{e}rn \rom{2} process for a hardcore distance of $R=0.2$. From top to bottom $a=1$ (sphere), $a=0.8,a=0.6$ and $a=0.4$. Black line is the estimated standardised $\hat{K}_{\text{inhom}}(r)$ for our observed data, dashed red line is the theoretical functional summary statistic for a Poisson process, and the grey shaded area represents the simulation envelope from 999 Monte Carlo simulations of Poisson processes fitted to the observed data.}
\label{fig:kfunction:compare}
\end{figure}

We also see that as $a$ becomes smaller, and therefore $c$ becomes larger a reducing empirical power of our test for both regular and cluster processes, for the same $R$ and $\kappa$ respectively. This effect could be due to a combination of mapping from the ellipsoid to the sphere and an artefact of the test statistic being proposed over a finite grid of points rather than being consider over the entire range of $[0,\pi]$. To see this consider just the Mat\'{e}rn \rom{2} process with a fixed hardcore distance $R$. Further, let us only consider a cross section of the ellipsoid, more specifically the ellipse such that it's major and minor axis lengths are $c$ and $a$ respectively, see Figure \ref{fig:hardcore:red}. Let us consider the point $(0,a)$ lying on the ellipse and take the point to the right of it which is precisely the hardcore distance $R$ away from it, let us label this point $\mathbf{x}$. To find $\mathbf{x}$ use the parametrisation $(x,y)=(c\cos t,a\sin t)$ for $t\in [0,2\pi)$. Then we solve the following equation for $t$ to find $\mathbf{x}$,
\begin{equation*}
R=\int_{t}^{\pi/2} (c^2\sin^2s +a^2\cos^2s)^{1/2}ds.
\end{equation*} 
Let us label $\tilde{t}$ as the solution to this equation and we then find $\mathbf{x}$ as $(c\cos \tilde{t},a\sin \tilde{t})$. Then apply the map onto the unit circle which gives us $(0,a)\mapsto (0,1)$ and $(c\cos \tilde{t},a\sin \tilde{t})\mapsto (\cos \tilde{t},\sin \tilde{t})$ and so our new hardcore distance on the circle is $R'=\cos^{-1}(\sin \tilde{t})$. It should be noted that this calculation is dependent upon where on an ellipsoid the event of interest, $\mathbf{x}$ is as $B_{\mathbb{D}}(\mathbf{x},R)$ does not map to $B_{\mathbb{D}}(f(\mathbf{x}),r)$ for some $r>0$ and where $f$ is our mapping from the ellipsoid to the sphere. Using the example when $a=0.4000$ we have that $c=3.1602$ and thus for a hardcore distance of $R=0.2$, on the cross sectional ellipse we have an effective hardcore distance of $R'=0.0633$, which is an extremely small hardcore distance and since our finite grid of points is too coarse (points are only separated a distance of 0.02 apart) it results in a loss in power of our test. Furthermore, examining the standardised inhomogeneous $K$-function in Figure \ref{fig:kfunction:compare} we can further see the effect of our mapping and taking only a finite grid of points along $[0,\pi]$ as the negative deviation reduces as ellipsoid \emph{deforms} further away from the unit sphere. It should be noted though that even though the power of our test reduces as we move further away from the sphere Figure \ref{fig:kfunction:compare} still indicates evidence that of regularity as for small $r$ for all ellipsoids the observed inhomogeneous $K$-function falls below the simulation envelope. This highlights the importance of a proper examination of graphical representations of functional summary statistics as opposed to the use of formal hypothesis testing \cite{Diggle2003}. Another consideration would be to potential use a two sided hypothesis testing procedure which may provide greater power when the true underlying process is not CSR.

A similar effect also occurs when examining cluster processes. By mapping the pattern from the ellipsoid to the sphere we will distort the isotropic nature of our offspring density relative to its parent. In particular this will cause the cluster size to contract and so if our finite grid of points is too coarse we will struggle to detect deviations away from CSR.

\section{Properties of the estimator of the $K_{\text{inhom}}$-function when $\rho$ is unknown}\label{Supplementary:K:unknown:rho}

In this section we shall rederive the expectation and variance of our estimators for $K_{\text{inhom}}(r)$ in the scenario when $\rho$ is unknown. We restate our estimator as,
\begin{equation}
\tilde{K}_{\text{inhom}}(r)=
\begin{cases}
\frac{\lambda_{\mathbb{D}}^2(\mathbb{D})}{4\pi N_Y(\mathbb{S}^2)(N_Y(\mathbb{S}^2)-1)}\sum_{\mathbf{x}\in Y}\sum_{\mathbf{y}\in Y\setminus \{x\}}\frac{\mathbbm{1}[d(\mathbf{x},\mathbf{y})\leq r]}{\tilde{\rho}(\mathbf{x})\tilde{\rho}(\mathbf{y})}, &\text{if } N_Y(\mathbb{S}^2)>1\\
0, &\text{otherwise}
\end{cases}
\end{equation}
where $Y=f(X)$ and
\begin{equation}\label{rho:tilde:app}
\tilde{\rho}(\mathbf{x})=
\left\{\begin{alignedat}{2}
 l_1(f^{-1}(\mathbf{x})) &J_{(1,f^*)}(\mathbf{x})\sqrt{1-\mathbf{x}_1^2-\mathbf{x}_2^2}, \quad \mathbf{x}\in f(\mathbb{D}_1)\\
&\vdots\\
 l_n(f^{-1}(\mathbf{x})) &J_{(n,f^*)}(\mathbf{x})\sqrt{1-\mathbf{x}_1^2-\mathbf{x}_2^2}, \quad \mathbf{x}\in f(\mathbb{D}_n)
  \end{alignedat}\right.
\end{equation}
Then the following theorem gives the properties of our estimator.

\begin{theorem}\label{thm:variance:K:unknown:rho:app}
The bias and variance of $\tilde{K}_{\text{inhom}}(r)$ are,
\begin{align*}
\text{Bias}(\tilde{K}_{\text{inhom}}(r))=-P(N_Y(\mathbb{S}^2)\leq 1) 2\pi(1-\cos r),
\end{align*}
and,
\begin{equation}\label{variance:K:unknown:rho:2:app}
\begin{split}
&\text{Var}(\tilde{K}_{\text{inhom}}(r))=4\pi^2(1-\cos r)^2(1-P(N_Y(\mathbb{S}^2)\leq 1))P(N_Y(\mathbb{S}^2)\leq 1)\\
&\phantom{AAAA}+\rho^3\lambda_{\mathbb{D}}^4(\mathbb{D})(1-\cos r)^2\left(\int_{\mathbb{S}^2}\frac{1}{\tilde{\rho}(\mathbf{x})}\lambda_{\mathbb{S}^2}(d\mathbf{x})-\frac{16\pi^2}{\lambda_{\mathbb{D}}(\mathbb{D})}\right)\mathbb{E}\left[\frac{1}{(N_Y(\mathbb{S}^2)+3)^2(N_Y(\mathbb{S}^2)+2)^2}\right]\\
&\phantom{AAAA}+\frac{\rho^2\lambda_{\mathbb{D}}^4(\mathbb{D})}{8\pi^2}\left(\int_{\mathbb{S}^2}\int_{\mathbb{S}^2} \frac{\mathbbm{1}[d(\mathbf{x}_1,\mathbf{x}_2)\leq r]}{\tilde{\rho}(\mathbf{x}_1)\tilde{\rho}(\mathbf{x}_2)} \lambda_{\mathbb{S}^2}(d\mathbf{x}_1) \lambda_{\mathbb{S}^2}(d\mathbf{x}_2) - \frac{64\pi^4(1-\cos r)^2}{\lambda_{\mathbb{D}}^2(\mathbb{D})} \right)\\
&\phantom{AAAAAAAA}\times\mathbb{E}\left[\frac{1}{(N_Y(\mathbb{S}^2)+2)^2(N_Y(\mathbb{S}^2)+1)^2}\right],
\end{split}
\end{equation}
where $\tilde{\rho}(\mathbf{x})$ is given by Equation (\ref{rho:tilde:app}).
\end{theorem}

\begin{proof}
We can rewrite our estimator as,
\begin{equation}
\tilde{K}_{\text{inhom}}(r)=\frac{\mathbbm{1}[N_Y(\mathbb{S}^2)>1]\lambda_{\mathbb{D}}^2(\mathbb{D})}{4\pi N_Y(\mathbb{S}^2)(N_Y(\mathbb{S}^2)-1)}\sum_{\mathbf{x}\in Y}\sum_{\mathbf{y}\in Y\setminus \{x\}}\frac{\mathbbm{1}[d(\mathbf{x},\mathbf{y})\leq r]}{\tilde{\rho}(\mathbf{x})\tilde{\rho}(\mathbf{y})}
\end{equation}
Further, note that $N_Y(\mathbb{S}^2)=N_{X}(\mathbb{D}^2),$ where $f(\mathbf{x})=\frac{\mathbf{x}}{|\mathbf{x}|}$. We then take expectations (using iterated expectations conditioning on $N_Y(\mathbb{S}^2)$) of $\tilde{K}_{\text{inhom}}(r)$ to get the bias,
\begin{align*}
\mathbb{E}(\tilde{K}_{\text{inhom}}(r))&=\mathbb{E}\left[\frac{\mathbbm{1}[N_Y(\mathbb{S}^2)>1]\lambda_{\mathbb{D}}^2(\mathbb{D})}{4\pi N_Y(\mathbb{S}^2)(N_Y(\mathbb{S}^2)-1)}\sum_{\mathbf{x}\in Y}\sum_{\mathbf{y}\in Y\setminus \{x\}}\frac{\mathbbm{1}[d(\mathbf{x},\mathbf{y})\leq r]}{\tilde{\rho}(\mathbf{x})\tilde{\rho}(\mathbf{y})}\right]\\
&=\frac{\lambda_{\mathbb{D}}^2(\mathbb{D})}{4\pi}\mathbb{E}\left(\frac{\mathbbm{1}[N_Y(\mathbb{S}^2)>1]}{N_Y(\mathbb{S}^2)(N_Y(\mathbb{S}^2)-1)}\mathbb{E}\left[\left.\sum^{\neq}_{\mathbf{x},\mathbf{y}\in Y}\frac{\mathbbm{1}[d(\mathbf{x},\mathbf{y})\leq r]}{\tilde{\rho}(\mathbf{x})\tilde{\rho}(\mathbf{y})}\right|N_Y(\mathbb{S}^2)=n\right]\right)\\
&=\frac{\lambda_{\mathbb{D}}^2(\mathbb{D})}{4\pi}\mathbb{E}\left(\frac{\mathbbm{1}[N_Y(\mathbb{S}^2)>1]}{N_Y(\mathbb{S}^2)(N_Y(\mathbb{S}^2)-1)}\mathbb{E}\left[\sum^{\neq}_{\mathbf{i},\mathbf{j}\in \{1,\dots,n\}}\frac{\mathbbm{1}[d(\mathbf{Y}_i,\mathbf{Y}_j)\leq r]}{\tilde{\rho}(\mathbf{Y}_i)\tilde{\rho}(\mathbf{Y}_j)}\right]\right)\\
&=\frac{\lambda_{\mathbb{D}}^2(\mathbb{D})}{4\pi}\mathbb{E}\left(\frac{\mathbbm{1}[N_Y(\mathbb{S}^2)>1]N_Y(\mathbb{S}^2)(N_Y(\mathbb{S}^2)-1)}{N_Y(\mathbb{S}^2)(N_Y(\mathbb{S}^2)-1)}\mathbb{E}\left[\frac{\mathbbm{1}[d(\mathbf{Y},\mathbf{Y})\leq r]}{\tilde{\rho}(\mathbf{Y})\tilde{\rho}(\mathbf{Y})}\right]\right)\\
&=\frac{\lambda_{\mathbb{D}}^2(\mathbb{D})}{4\pi}P(N_Y(\mathbb{S}^2)>1)\mathbb{E}\left[\frac{\mathbbm{1}[d(\mathbf{Y},\mathbf{Y})\leq r]}{\tilde{\rho}(\mathbf{Y})\tilde{\rho}(\mathbf{Y})}\right]\\
&=\frac{\lambda_{\mathbb{D}}^2(\mathbb{D})}{4\pi}P(N_Y(\mathbb{S}^2)>1)\int_{\mathbb{S}^2}\int_{\mathbb{S}^2}\frac{\mathbbm{1}[d(\mathbf{x},\mathbf{y})\leq r]}{\tilde{\rho}(\mathbf{x})\tilde{\rho}(\mathbf{y})}\frac{\rho^*(\mathbf{x})\rho^*(\mathbf{y})}{\rho^2\lambda_{\mathbb{D}}^2(\mathbb{D})}\lambda_{\mathbb{S}^2}(d\mathbf{x})\lambda_{\mathbb{S}^2}(d\mathbf{y})\\
&=P(N_Y(\mathbb{S}^2)>1)2\pi(1-\cos r).
\end{align*}
The bias follows by noting that $P(N_Y(\mathbb{S}^2)>1)=1-P(N_Y(\mathbb{S}^2)\leq 1)$. The variance can be calculated using the law of total variance we have,
\begin{equation}\label{variance:law:total}
\text{Var}(\tilde{K}_{\text{inhom}}(r)) = \underbrace{\mathbb{E}[\text{Var}(\tilde{K}_{\text{inhom}}(r)|N_Y(\mathbb{S}^2))]}_{(1)} + \underbrace{\text{Var}[\mathbb{E}[\tilde{K}_{\text{inhom}}(r)|N_Y(\mathbb{S}^2)]]}_{(2)}.
\end{equation}
Considering term $(2)$ first,
\begin{align*}
\mathbb{E}[\tilde{K}_{\text{inhom}}(r)|N_Y(\mathbb{S}^2)=n]&=\frac{\mathbbm{1}[n>1]\lambda_{\mathbb{D}}^2(\mathbb{D})}{4\pi n(n-1)}\mathbb{E}\left[\sum^{\neq}_{\mathbf{x},\mathbf{y}\in Y}\frac{\mathbbm{1}[d(\mathbf{x},\mathbf{y})\leq r]}{\tilde{\rho}(\mathbf{x})\tilde{\rho}(\mathbf{y})}\Bigg|N_Y(\mathbb{S}^2)=n\right]\\
\intertext{Using the fact that given $N_Y(\mathbb{S}^2)=n$, each $\mathbf{x}\in Y$ is independently and identically distributed with density $\rho^*(\mathbf{x})/\mu(\mathbb{D})$, where $\mu(\mathbb{D})=\rho\lambda_{\mathbb{D}}(\mathbb{D})$, the expectation becomes,}
&=\frac{\mathbbm{1}[n>1]\lambda_{\mathbb{D}}^2(\mathbb{D})}{4\pi n(n-1)}\mathbb{E}\left[\sum^{\neq}_{i,j\in \{1,\dots,n\}}\frac{\mathbbm{1}[d(\mathbf{Y}_i,\mathbf{Y}_j)\leq r]}{\tilde{\rho}(\mathbf{Y}_i)\tilde{\rho}(\mathbf{Y}_j)}\Bigg|N_Y(\mathbb{S}^2)=n\right]\\
&=\frac{\mathbbm{1}[n>1]\lambda_{\mathbb{D}}^2(\mathbb{D})}{4\pi n(n-1)}\sum^{\neq}_{i,j\in \{1,\dots,n\}}\mathbb{E}\left[\frac{\mathbbm{1}[d(\mathbf{Y}_i,\mathbf{Y}_j)\leq r]}{\tilde{\rho}(\mathbf{Y}_i)\tilde{\rho}(\mathbf{Y}_j)}\Bigg|N_Y(\mathbb{S}^2)=n\right]\\
&=\frac{\mathbbm{1}[n>1]\lambda_{\mathbb{D}}^2(\mathbb{D})}{4\pi n(n-1)}n(n-1)\mathbb{E}\left[\frac{\mathbbm{1}[d(\mathbf{X},\mathbf{Y})\leq r]}{\tilde{\rho}(\mathbf{X})\tilde{\rho}(\mathbf{Y})}\right]\\
&=\frac{\mathbbm{1}[n>1]\lambda_{\mathbb{D}}^2(\mathbb{D})}{4\pi }\mathbb{E}\left[\frac{\mathbbm{1}[d(\mathbf{X},\mathbf{Y})\leq r]}{\tilde{\rho}(\mathbf{X})\tilde{\rho}(\mathbf{Y})}\right],
\end{align*}
where $\mathbf{X}$ and $\mathbf{Y}$ are independent random vectors distributed in $\mathbb{S}^2$ with density $\rho^*(\mathbf{x})/\mu(\mathbb{D})$. Noting that the joint density of $\mathbf{X}$ and $\mathbf{Y}$ is $\rho^*(\mathbf{x})\rho^*(\mathbf{y})/\mu^2(\mathbb{D})$, the expectation becomes,
\begin{align*}
\mathbb{E}[\tilde{K}_{\text{inhom}}(r)|N_Y(\mathbb{S}^2)=n]&=\frac{\mathbbm{1}[n>1]\lambda_{\mathbb{D}}^2(\mathbb{D})}{4\pi }\int_{\mathbb{S}^2}\int_{\mathbb{S}^2}\frac{\mathbbm{1}[d(\mathbf{x},\mathbf{y})\leq r]}{\tilde{\rho}(\mathbf{y})\tilde{\rho}(\mathbf{y})}\frac{\rho^*(\mathbf{x})\rho^*(\mathbf{y})}{\mu^2(\mathbb{D})}\lambda_{\mathbb{S}^2}(d\mathbf{x}) \lambda_{\mathbb{S}^2}(d\mathbf{y})\\
&=\frac{\mathbbm{1}[n>1]\lambda_{\mathbb{D}}^2(\mathbb{D})}{4\pi }\frac{1}{\lambda_{\mathbb{D}}^2(\mathbb{D})}\int_{\mathbb{S}^2}\int_{\mathbb{S}^2}\mathbbm{1}[d(\mathbf{x},\mathbf{y})\leq r]\lambda_{\mathbb{S}^2}(d\mathbf{x}) \lambda_{\mathbb{S}^2}(d\mathbf{y})\\
&=\frac{\mathbbm{1}[n>1]}{4\pi}\int_{\mathbb{S}^2}\int_{\mathbb{S}^2}\mathbbm{1}[d(\mathbf{x},\mathbf{y})\leq r]\lambda_{\mathbb{S}^2}(d\mathbf{x}) \lambda_{\mathbb{S}^2}(d\mathbf{y})\\
&=\mathbbm{1}[n>1]2\pi(1-\cos r)
\end{align*}
Hence term $(2)$ is,
\begin{align*}
\text{Var}[\mathbb{E}[\tilde{K}_{\text{inhom}}(r)|N_Y(\mathbb{S}^2)]]&=\text{Var}(\mathbbm{1}[N_Y(\mathbb{S}^2)>1]2\pi(1-\cos r))\\
&=4\pi^2(1-\cos r)^2\left[\mathbb{E}(\mathbbm{1}^2[N_Y(\mathbb{S}^2)>1]) - \mathbb{E}^2(\mathbbm{1}[N_Y(\mathbb{S}^2)>1]) \right]\\
&=4\pi^2(1-\cos r)^2\left[P(N_Y(\mathbb{S}^2)>1) - P^2(N_Y(\mathbb{S}^2)>1)\right]\\
&=4\pi^2(1-\cos r)^2P(N_Y(\mathbb{S}^2)>1)\left[1 - P(N_Y(\mathbb{S}^2)>1)\right]\\
&=4\pi^2(1-\cos r)^2(1-P(N_Y(\mathbb{S}^2)\leq 1))P(N_Y(\mathbb{S}^2)\leq 1),
\end{align*}
where, since $N_Y(\mathbb{S}^2)\sim \text{Poisson}(\rho\lambda_{\mathbb{D}}(\mathbb{D}))$,
\begin{equation*}
P(N_Y(\mathbb{S}^2)\leq 1) = 1 - e^{-\rho\lambda_{\mathbb{D}}(\mathbb{D})} - \rho\lambda_{\mathbb{D}}(\mathbb{D})e^{-\rho\lambda_{\mathbb{D}}(\mathbb{D})}.
\end{equation*}
Now consider term $(1)$, we calculate $\text{Var}(\tilde{K}_{\text{inhom}}(r)|N_Y(\mathbb{S}^2))$,
\begin{align*}
\text{Var}(\tilde{K}_{\text{inhom}}(r)|N_Y(\mathbb{S}^2)=n)&=\text{Var}\left(\frac{\mathbbm{1}[N_Y(\mathbb{S}^2)>1]\lambda_{\mathbb{D}}^2(\mathbb{D})}{4\pi N_Y(\mathbb{S}^2)(N_Y(\mathbb{S}^2)-1)}\sum^{\neq}_{\mathbf{x},\mathbf{y}\in Y}\frac{\mathbbm{1}[d(\mathbf{x},\mathbf{y})\leq r]}{\tilde{\rho}(\mathbf{x})\tilde{\rho}(\mathbf{y})}\Bigg|N_Y(\mathbb{S}^2)=n\right)\\
&=\frac{\mathbbm{1}^2[n>1]\lambda_{\mathbb{D}}^4(\mathbb{D})}{16\pi^2 }\text{Var}\left(\frac{1}{n(n-1)}\sum^{\neq}_{\mathbf{x},\mathbf{y}\in Y}\frac{\mathbbm{1}[d(\mathbf{x},\mathbf{y})\leq r]}{\tilde{\rho}(\mathbf{x})\tilde{\rho}(\mathbf{y})}\Bigg|N_Y(\mathbb{S}^2)=n\right)\\
&=\frac{\mathbbm{1}[n>1]\lambda_{\mathbb{D}}^4(\mathbb{D})}{16\pi^2 }\text{Var}\left(\frac{1}{n(n-1)}\sum^{\neq}_{\mathbf{x},\mathbf{y}\in Y}\frac{\mathbbm{1}[d(\mathbf{x},\mathbf{y})\leq r]}{\tilde{\rho}(\mathbf{x})\tilde{\rho}(\mathbf{y})}\Bigg|N_Y(\mathbb{S}^2)=n\right)\\
&=\frac{\mathbbm{1}[n>1]\lambda_{\mathbb{D}}^4(\mathbb{D})}{16\pi^2 }\text{Var}\left(\frac{1}{n(n-1)}\sum^{\neq}_{i,j\in \{1,\dots,n\}}\frac{\mathbbm{1}[d(\mathbf{Y}_i,\mathbf{Y}_j)\leq r]}{\tilde{\rho}(\mathbf{Y}_i)\tilde{\rho}(\mathbf{Y}_j)}\right)
\end{align*}
Here we follow a similar argument to that of \cite{Lang2010} through $U$-statistics. Noting that $\frac{\mathbbm{1}[d(\mathbf{Y}_i,\mathbf{Y}_j)\leq r]}{\tilde{\rho}(\mathbf{Y}_i)\tilde{\rho}(\mathbf{Y}_j)}=\frac{\mathbbm{1}[d(\mathbf{Y}_j,\mathbf{Y}_i)\leq r]}{\tilde{\rho}(\mathbf{Y}_j)\tilde{\rho}(\mathbf{Y}_i)}$, i.e. it is symmetric in its arguments, we rewrite the summation,
\begin{align*}
\frac{1}{n(n-1)}\sum^{\neq}_{i,j\in \{1,\dots,n\}}\frac{\mathbbm{1}[d(\mathbf{Y}_i,\mathbf{Y}_j)\leq r]}{\tilde{\rho}(\mathbf{Y}_i)\tilde{\rho}(\mathbf{Y}_j)}&=\frac{1}{n(n-1)}\sum_{1\leq i < j \leq n}\frac{\mathbbm{1}[d(\mathbf{Y}_i,\mathbf{Y}_j)\leq r]}{\tilde{\rho}(\mathbf{Y}_i)\tilde{\rho}(\mathbf{Y}_j)}+\frac{\mathbbm{1}[d(\mathbf{Y}_j,\mathbf{Y}_i)\leq r]}{\tilde{\rho}(\mathbf{Y}_j)\tilde{\rho}(\mathbf{Y}_i)}\\
&=\frac{2}{n(n-1)}\sum_{1\leq i < j \leq n}\frac{\mathbbm{1}[d(\mathbf{Y}_i,\mathbf{Y}_j)\leq r]}{\tilde{\rho}(\mathbf{Y}_i)\tilde{\rho}(\mathbf{Y}_j)}\\
&={n\choose 2}^{-1}\sum_{1\leq i < j \leq n}\frac{\mathbbm{1}[d(\mathbf{Y}_i,\mathbf{Y}_j)\leq r]}{\tilde{\rho}(\mathbf{Y}_i)\tilde{\rho}(\mathbf{Y}_j)}
\end{align*}
This is form of a $U$-statistic and variances of this class of statistics can be decomposed using the work of \cite{Hoeffding1992}. Using the same notation as \cite{Hoeffding1992}, we define some quantities and derive a number of expectations,
\begin{align*}
U_n&={n\choose 2}^{-1}\sum_{1\leq i < j \leq n}\frac{\mathbbm{1}[d(\mathbf{Y}_i,\mathbf{Y}_j)\leq r]}{\tilde{\rho}(\mathbf{Y}_i)\tilde{\rho}(\mathbf{Y}_j)}\\
\Phi(\mathbf{y}_1,\mathbf{y}_2)&=\frac{\mathbbm{1}[d(\mathbf{y}_1,\mathbf{y}_2)\leq r]}{\tilde{\rho}(\mathbf{y}_1)\tilde{\rho}(\mathbf{y}_2)}\\
\Phi_1(\mathbf{y}_1)&\equiv\Phi_1(\mathbf{y}_1,\mathbf{Y}_2)=\mathbb{E}[\Phi(\mathbf{y}_1,\mathbf{Y}_2)]=\mathbb{E}[\Phi(\mathbf{Y}_1,\mathbf{Y}_2)|\mathbf{Y}_1=\mathbf{y}_1]\\
&=\mathbb{E}\left[\frac{\mathbbm{1}[d(\mathbf{y}_1,\mathbf{Y}_2)\leq r]}{\tilde{\rho}(\mathbf{y}_1)\tilde{\rho}(\mathbf{Y}_2)}\right]\\
&=\int_{\mathbb{S}^2}\frac{\mathbbm{1}[d(\mathbf{y}_1,\mathbf{y}_2)\leq r]}{\tilde{\rho}(\mathbf{y}_1)\tilde{\rho}(\mathbf{y}_2)}\frac{\rho^*(\mathbf{y}_2)}{\rho\lambda_{\mathbb{D}}(\mathbb{D})}\lambda_{\mathbb{S}^2}(d\mathbf{y}_2)\\
&=\frac{1}{\tilde{\rho}(\mathbf{y}_1)\lambda_{\mathbb{D}}(\mathbb{D})}\int_{\mathbb{S}^2}\mathbbm{1}[d(\mathbf{y}_1,\mathbf{y}_2)\leq r]\lambda_{\mathbb{S}^2}(d\mathbf{y}_2)\\
&=\frac{2\pi(1-\cos r)}{\tilde{\rho}(\mathbf{y}_1)\lambda_{\mathbb{D}}(\mathbb{D})}\\
\Phi_2(\mathbf{y}_1,\mathbf{y}_2)&=\mathbb{E}[\Phi(\mathbf{y}_1,\mathbf{y}_2)]=\mathbb{E}[\phi(\mathbf{Y}_1,\mathbf{Y}_2)|\mathbf{Y}_1=\mathbf{y}_1,\mathbf{Y}_2=\mathbf{y}_2]=\Phi(\mathbf{y}_1,\mathbf{y}_2)\\
\mathbb{E}[\Phi_1(\mathbf{Y}_1)]&=\int_{\mathbb{S}^2}  \frac{2\pi(1-\cos r)}{\tilde{\rho}(\mathbf{y}_1)\lambda_{\mathbb{D}}(\mathbb{D})}\frac{\rho^*(\mathbf{y}_1)}{\rho\lambda_{\mathbb{D}}(\mathbb{D})} \lambda_{\mathbb{S}^2}(d\mathbf{y}_1)\\
&=\frac{ 4\pi \cdot 2\pi(1-\cos r)}{\lambda_{\mathbb{D}}^2(\mathbb{D})}\\
\mathbb{E}[\Phi_2(\mathbf{Y}_1,\mathbf{Y}_2)]&=\int_{\mathbb{S}^2}\int_{\mathbb{S}^2} \frac{\mathbbm{1}[d(\mathbf{y}_1,\mathbf{y}_2)\leq r]}{\tilde{\rho}(\mathbf{y}_1)\tilde{\rho}(\mathbf{y}_2)} \frac{\rho^*(\mathbf{y}_1)\rho^*(\mathbf{y}_2)}{\rho^2\lambda_{\mathbb{D}}^2(\mathbb{D})} \lambda_{\mathbb{S}^2}(d\mathbf{y}_1) \lambda_{\mathbb{S}^2}(d\mathbf{y}_2)\\
&=\frac{4\pi\cdot 2\pi(1-\cos r)}{\lambda_{\mathbb{D}}^2(\mathbb{D})}\\
\mathbb{E}[\Phi^2_1(\mathbf{Y}_1)]&=\int_{\mathbb{S}^2}\frac{4\pi^2(1-\cos r)^2}{\tilde{\rho}^2(\mathbf{y}_1)\lambda_{\mathbb{D}}^2(\mathbb{D})}\frac{\rho^*(\mathbf{y}_1)}{\rho\lambda_{\mathbb{D}}(\mathbb{D})} \lambda_{\mathbb{S}^2}(d\mathbf{y}_1)\\
&=\frac{4\pi^2(1-\cos r)^2}{\lambda_{\mathbb{D}}^3(\mathbb{D})}\int_{\mathbb{S}^2}\frac{1}{\tilde{\rho}(\mathbf{y}_1)}\lambda_{\mathbb{S}^2}(d\mathbf{y}_1)\\
\mathbb{E}[\Phi_2^2(\mathbf{Y}_1,\mathbf{Y}_2)]&=\int_{\mathbb{S}^2}\int_{\mathbb{S}^2} \frac{\mathbbm{1}[d(\mathbf{y}_1,\mathbf{y}_2)\leq r]}{\tilde{\rho}^2(\mathbf{y}_1)\tilde{\rho}^2(\mathbf{y}_2)}\frac{\rho^*(\mathbf{y}_1)\rho^*(\mathbf{y}_2)}{\rho^2\lambda_{\mathbb{D}}^2(\mathbb{D})} \lambda_{\mathbb{S}^2}(d\mathbf{y}_1) \lambda_{\mathbb{S}^2}(d\mathbf{y}_2)\\
&=\frac{1}{\lambda_{\mathbb{D}}^2(\mathbb{D})}\int_{\mathbb{S}^2}\int_{\mathbb{S}^2} \frac{\mathbbm{1}[d(\mathbf{y}_1,\mathbf{y}_2)\leq r]}{\tilde{\rho}(\mathbf{y}_1)\tilde{\rho}(\mathbf{y}_2)} \lambda_{\mathbb{S}^2}(d\mathbf{y}_1) \lambda_{\mathbb{S}^2}(d\mathbf{y}_2)\\
\zeta_1 &= \text{Var}(\Phi_1(\mathbf{Y}_1))\\
&=\mathbb{E}[\Phi_1^2(\mathbf{Y}_1)]-\mathbb{E}^2[\Phi_1(\mathbf{Y}_1)]\\
&=\frac{4\pi^2(1-\cos r)^2}{\lambda_{\mathbb{D}}^3(\mathbb{D})}\int_{\mathbb{S}^2}\frac{1}{\tilde{\rho}(\mathbf{y}_1)}\lambda_{\mathbb{S}^2}(d\mathbf{y}_1)-\frac{16\pi^2\cdot 4\pi^2(1-\cos r)^2}{\lambda_{\mathbb{D}}^4(\mathbb{D})}\\
&=\frac{4\pi^2(1-\cos r)^2}{\lambda_{\mathbb{D}}^3(\mathbb{D})}\left(\int_{\mathbb{S}^2}\frac{1}{\tilde{\rho}(\mathbf{y}_1)}\lambda_{\mathbb{S}^2}(d\mathbf{y}_1)-\frac{16\pi^2}{\lambda_{\mathbb{D}}(\mathbb{D})}\right)\\
\zeta_2 &= \text{Var}(\Phi_2(\mathbf{Y}_1,\mathbf{Y}_2))\\
&= \mathbb{E}[\Phi_2^2(\mathbf{Y}_1,\mathbf{Y}_2)]-\mathbb{E}^2[\Phi_2(\mathbf{Y}_1,\mathbf{Y}_2)]\\
&= \frac{1}{\lambda_{\mathbb{D}}^2(\mathbb{D})}\int_{\mathbb{S}^2}\int_{\mathbb{S}^2} \frac{\mathbbm{1}[d(\mathbf{y}_1,\mathbf{y}_2)\leq r]}{\tilde{\rho}(\mathbf{y}_1)\tilde{\rho}(\mathbf{y}_2)} \lambda_{\mathbb{S}^2}(d\mathbf{y}_1) \lambda_{\mathbb{S}^2}(d\mathbf{y}_2)-\frac{16\pi^2\cdot 4\pi^2(1-\cos r)^2}{\lambda_{\mathbb{D}}^4(\mathbb{D})}\\
&=\frac{1}{\lambda_{\mathbb{D}}^2(\mathbb{D})}\left(\int_{\mathbb{S}^2}\int_{\mathbb{S}^2} \frac{\mathbbm{1}[d(\mathbf{y}_1,\mathbf{y}_2)\leq r]}{\tilde{\rho}(\mathbf{y}_1)\tilde{\rho}(\mathbf{y}_2)} \lambda_{\mathbb{S}^2}(d\mathbf{y}_1) \lambda_{\mathbb{S}^2}(d\mathbf{y}_2) - \frac{64\pi^4 (1-\cos r)^2}{\lambda_{\mathbb{D}}^2(\mathbb{D})} \right)
\end{align*}
Then using the variance derived by \cite{Hoeffding1992} for $U$-statistics, the variance of our $U_n$ can be decomposed as,
\begin{align*}
\text{Var}\left(U_n\right)&={n \choose 2}^{-1}\sum_{k=1}^2 {2\choose k} {{n-2}\choose {2-k}}\zeta_k\\
&=\frac{4(n-2)}{n(n-1)}\zeta_1+\frac{2}{n(n-1)}\zeta_2
\end{align*}
Then,
\begin{equation*}
\begin{split}
&\text{Var}(\tilde{K}_{\text{inhom}}(r)|N_Y(\mathbb{S}^2)=n)=\frac{\mathbbm{1}[n>1]\lambda_{\mathbb{D}}^4(\mathbb{D})}{16\pi^2}\left(\frac{4(n-2)}{n(n-1)}\zeta_1+\frac{2}{n(n-1)}\zeta_2 \right)\\
&=\lambda_{\mathbb{D}}(\mathbb{D})(1-\cos r)^2\left(\int_{\mathbb{S}^2}\frac{1}{\tilde{\rho}(\mathbf{y}_1)}\lambda_{\mathbb{S}^2}(d\mathbf{y}_1)-\frac{16\pi^2}{\lambda_{\mathbb{D}}(\mathbb{D})}\right)\cdot\frac{\mathbbm{1}[n>1](n-2)}{n(n-1)}\\
&+\frac{\lambda_{\mathbb{D}}^2(\mathbb{D})}{8\pi^2}\left(\int_{\mathbb{S}^2}\int_{\mathbb{S}^2} \frac{\mathbbm{1}[d(\mathbf{y}_1,\mathbf{y}_2)\leq r]}{\tilde{\rho}(\mathbf{y}_1)\tilde{\rho}(\mathbf{y}_2)} \lambda_{\mathbb{S}^2}(d\mathbf{y}_1) \lambda_{\mathbb{S}^2}(d\mathbf{y}_2) - \frac{64\pi^4(1-\cos r)^2}{\lambda_{\mathbb{D}}^2(\mathbb{D})} \right)\frac{\mathbbm{1}[n>1]}{n(n-1)}
\end{split}
\end{equation*}
Taking expectations gives,
\begin{equation*}
\begin{split}
&\mathbb{E}[\text{Var}(\tilde{K}_{\text{inhom}}(r)|N_Y(\mathbb{S}^2))]\\
&=\lambda_{\mathbb{D}}(\mathbb{D})(1-\cos r)^2\left(\int_{\mathbb{S}^2}\frac{1}{\tilde{\rho}(\mathbf{y}_1)}\lambda_{\mathbb{S}^2}(d\mathbf{y}_1)-\frac{16\pi^2}{\lambda_{\mathbb{D}}(\mathbb{D})}\right)\mathbb{E}\left[\frac{\mathbbm{1}[N_Y(\mathbb{S}^2)>1](N_Y(\mathbb{S}^2)-2)}{N_Y(\mathbb{S}^2)(N_Y(\mathbb{S}^2)-1)}\right]\\
&+\frac{\lambda_{\mathbb{D}}^2(\mathbb{D})}{8\pi^2}\left(\int_{\mathbb{S}^2}\int_{\mathbb{S}^2} \frac{\mathbbm{1}[d(\mathbf{y}_1,\mathbf{y}_2)\leq r]}{\tilde{\rho}(\mathbf{y}_1)\tilde{\rho}(\mathbf{y}_2)} \lambda_{\mathbb{S}^2}(d\mathbf{y}_1) \lambda_{\mathbb{S}^2}(d\mathbf{y}_2) - \frac{64\pi^4(1-\cos r)^2}{\lambda_{\mathbb{D}}^2(\mathbb{D})} \right)\mathbb{E}\left[\frac{\mathbbm{1}[N_Y(\mathbb{S}^2)>1]}{N_Y(\mathbb{S}^2)(N_Y(\mathbb{S}^2)-1)}\right]
\end{split}
\end{equation*}
The expectations can be simplified as follows and defining $\lambda=\rho\lambda_{\mathbb{D}}(\mathbb{D})$,
\begin{align*}
\mathbb{E}\left[\frac{\mathbbm{1}[N_Y(\mathbb{S}^2)>1](N_Y(\mathbb{S}^2)-2)}{N_Y(\mathbb{S}^2)(N_Y(\mathbb{S}^2)-1)}\right]&=\sum_{n=0}^{\infty}\frac{\mathbbm{1}[n>1](n-2)}{n(n-1)}\frac{\lambda^ne^{-\lambda}}{n!}\\
&=\sum_{n=3}^{\infty}\frac{(n-2)}{n(n-1)}\frac{\lambda^ne^{-\lambda}}{n!}\\
&=\sum_{n=3}^{\infty}\frac{1}{n^2(n-1)^2}\frac{\lambda^ne^{-\lambda}}{(n-3)!}\\
&=\sum_{n=0}^{\infty}\frac{1}{(n+3)^2(n+2)^2}\frac{\lambda^{n+3}e^{-\lambda}}{n!}\\
&=\lambda^3\sum_{n=0}^{\infty}\frac{1}{(n+3)^2(n+2)^2}\frac{\lambda^ne^{-\lambda}}{n!}\\
&=\lambda^3\mathbb{E}\left[\frac{1}{(N_Y(\mathbb{S}^2)+3)^2(N_Y(\mathbb{S}^2)+2)^2}\right]
\end{align*}
Similarly the other expectation is,
\begin{equation*}
\mathbb{E}\left[\frac{\mathbbm{1}[N_Y(\mathbb{S}^2)>1]}{N_Y(\mathbb{S}^2)(N_Y(\mathbb{S}^2)-1)}\right]=\lambda^2\mathbb{E}\left[\frac{1}{(N_Y(\mathbb{S}^2)+2)^2(N_Y(\mathbb{S}^2)+1)^2}\right],
\end{equation*}
and so,
\begin{equation*}
\begin{split}
&\mathbb{E}[\text{Var}(\tilde{K}_{\text{inhom}}(r)|N_Y(\mathbb{S}^2))]\\
&\phantom{AAAA}=\rho^3\lambda_{\mathbb{D}}^4(\mathbb{D})(1-\cos r)^2\left(\int_{\mathbb{S}^2}\frac{1}{\tilde{\rho}(\mathbf{y}_1)}\lambda_{\mathbb{S}^2}(d\mathbf{y}_1)-\frac{16\pi^2}{\lambda_{\mathbb{D}}(\mathbb{D})}\right)\mathbb{E}\left[\frac{1}{(N_Y(\mathbb{S}^2)+3)^2(N_Y(\mathbb{S}^2)+2)^2}\right]\\
&\phantom{AAAA}+\frac{\rho^2\lambda_{\mathbb{D}}^4(\mathbb{D})}{8\pi^2}\left(\int_{\mathbb{S}^2}\int_{\mathbb{S}^2} \frac{\mathbbm{1}[d(\mathbf{y}_1,\mathbf{y}_2)\leq r]}{\tilde{\rho}(\mathbf{y}_1)\tilde{\rho}(\mathbf{y}_2)} \lambda_{\mathbb{S}^2}(d\mathbf{y}_1) \lambda_{\mathbb{S}^2}(d\mathbf{y}_2) - \frac{64\pi^4(1-\cos r)^2}{\lambda_{\mathbb{D}}^2(\mathbb{D})} \right)\\
&\phantom{AAAAAAAA}\times\mathbb{E}\left[\frac{1}{(N_Y(\mathbb{S}^2)+2)^2(N_Y(\mathbb{S}^2)+1)^2}\right].
\end{split}
\end{equation*}
Combining everything gives the variance of $\tilde{K}_{\text{inhom}}(r)$,
\begin{equation*}
\begin{split}
&\text{Var}(\tilde{K}_{\text{inhom}}(r))=4\pi^2(1-\cos r)^2(1-P(N_Y(\mathbb{S}^2)\leq 1))P(N_Y(\mathbb{S}^2)\leq 1)\\
&\phantom{AAAA}+\rho^3\lambda_{\mathbb{D}}^4(\mathbb{D})(1-\cos r)^2\left(\int_{\mathbb{S}^2}\frac{1}{\tilde{\rho}(\mathbf{y}_1)}\lambda_{\mathbb{S}^2}(d\mathbf{y}_1)-\frac{16\pi^2}{\lambda_{\mathbb{D}}(\mathbb{D})}\right)\mathbb{E}\left[\frac{1}{(N_Y(\mathbb{S}^2)+3)^2(N_Y(\mathbb{S}^2)+2)^2}\right]\\
&\phantom{AAAA}+\frac{\rho^2\lambda_{\mathbb{D}}^4(\mathbb{D})}{8\pi^2}\left(\int_{\mathbb{S}^2}\int_{\mathbb{S}^2} \frac{\mathbbm{1}[d(\mathbf{y}_1,\mathbf{y}_2)\leq r]}{\tilde{\rho}(\mathbf{y}_1)\tilde{\rho}(\mathbf{y}_2)} \lambda_{\mathbb{S}^2}(d\mathbf{y}_1) \lambda_{\mathbb{S}^2}(d\mathbf{y}_2) - \frac{64\pi^4(1-\cos r)^2}{\lambda_{\mathbb{D}}^2(\mathbb{D})} \right)\\
&\phantom{AAAAAAAA}\times\mathbb{E}\left[\frac{1}{(N_Y(\mathbb{S}^2)+2)^2(N_Y(\mathbb{S}^2)+1)^2}\right]
\end{split}
\end{equation*}
\end{proof}

In order to construct estimates of $\text{Var}(\tilde{K}_{\text{inhom}}(r))$ we need estimators for $P(N_Y(\mathbb{S}^2)\leq 1)$. The following lemma helps us do this.

\begin{lemmaapp}\label{lemma:ratio:unbiased:prob:app}
Let $N\sim \text{Poisson}(\lambda)$ and $k,p\in\mathbb{N}$. Define the following random variable,
\begin{equation*}
R = \frac{N!e^{N-k}}{(N-k)!(e+p)^N}.
\end{equation*}
Then $R$ is ratio-unbiased for $\lambda^ke^{-p\lambda}$.
\end{lemmaapp}
\begin{proof}
Define $S=N!e^{N-k}/(N-k)!$ and $T=(e+p)^N$, then $R=S/T$. Then,
\begin{align*}
\mathbb{E}[S]&=\sum_{n=0}^{\infty} \frac{n!e^{n-k}}{(n-k)!}\frac{\lambda^ne^{-\lambda}}{n!}\\
&=e^{-k}e^{-\lambda}\sum_{n=k}^{\infty}\frac{(e\lambda)^n}{(n-k)!}\\
&=e^{-k}e^{-\lambda}\sum_{l=0}^{\infty}\frac{(e\lambda)^{l+k}}{l!}\\
&=e^{-k}e^{-\lambda}(e\lambda)^{k}\sum_{l=0}^{\infty}\frac{(e\lambda)^l}{l!}\\
&=e^{-k}e^{-\lambda}(e\lambda)^{k}e^{e\lambda}\\
&=\lambda^ke^{\lambda(e-1)}.
\end{align*}
Further,
\begin{align*}
\mathbb{E}[T]&=\sum_{n=0}^{\infty}(e+p)^n\frac{\lambda^ne^{-\lambda}}{n!}\\
&=e^{-\lambda}\sum_{n=0}^{\infty}\frac{(\lambda(e+p))^n}{n!}\\
&=e^{-\lambda}e^{\lambda(e+p)}\\
&=e^{\lambda(e+p-1)}.
\end{align*}
Therefore,
\begin{equation*}
\frac{\mathbb{E}[S]}{\mathbb{E}[T]}=\frac{\lambda^ke^{\lambda(e-1)}}{e^{\lambda(e+p-1)}}=\lambda^ke^{-p\lambda}
\end{equation*}
\end{proof}

\clearpage
\bibliographystyle{unsrt}  
\bibliography{ref3}
